\documentclass[uplatex,12pt]{article}
\usepackage{amsmath}
\usepackage{amssymb}
\usepackage{latexsym}
\usepackage[dvips]{color}
\usepackage{fancybox}
\usepackage{graphics}
\usepackage{amsthm}
\usepackage{tabularx}
\usepackage{amsthm}
\usepackage{subcaption}
\usepackage{mathrsfs}
\usepackage{colortbl}
\usepackage{bm}
\usepackage{url}
\usepackage{mathtools}
\usepackage{amscd}
\usepackage{tikz-cd}
\usepackage{tikz}
\usetikzlibrary{cd}

\usepackage[all]{xy}

\usepackage{hyperref}

\usepackage[a4paper]{geometry}
\usepackage{mathrsfs}
\usepackage{comment}
\usepackage{accsupp}

\usepackage{longtable}
\usepackage{blkarray}
\usepackage{pdfpages}
\usepackage{url}
\usepackage{setspace}

\usepackage{listings}
\usepackage{color}
\def\diam{\mathop{\mathrm{diam}}}
\def\diag{\mathop{\mathrm{diag}}}

\def\argmin{\mathop{\mathrm{arg \ min}}}

\def\Span{\mathop{\mathrm{span}}}

\def\div{\mathop{\mathrm{div}}}

\def\<{\mathop{\textless}}
\def\>{\mathop{\textgreater}}

\def\card{\mathop{\rm{card}}}

\def\conv{\mathop{\mathrm{conv}}}

\def\Tr{\mathop{\rm{Tr}}}
\def\det{\mathop{\rm{det}}}

\def\Int{\mathop{\rm{int}}}

\def\divh{\mathop{\mathrm{div}_h}}

\theoremstyle{definition}
\newtheorem{thr}{Theorem}[section]
\newtheorem{ex}[thr]{Example}
\newtheorem*{ex*}{Example}
\newtheorem{rem}[thr]{Remark}
\newtheorem{defi}[thr]{Definition}
\newtheorem{lem}[thr]{Lemma}

\newtheorem*{pf*}{Proof}
\newtheorem{note}[thr]{Note}
\newtheorem{prop}[thr]{Proposition}
\newtheorem{cor}[thr]{Corollary}
\newtheorem{assume}[thr]{Assumption}
\newtheorem*{ass*}{Assumption}
\newtheorem{Cond}[thr]{Condition}

\newtheorem*{exer*}{Exercise}

\newtheorem*{memo*}{Memo}

\newtheorem{que}[thr]{Question}

\allowdisplaybreaks[3]

\usepackage{makeidx}
\makeindex

\newcounter{sone}
\newcounter{stwo}
\newcounter{sthree}
\newcounter{sfour}
\newcounter{sfive}
\newcounter{ssix}
\newcounter{lone}
\newcounter{ltwo}
\newcounter{lthree}
\newcounter{lfour}
\newcounter{lfive}
\newcounter{lsix}
\newcounter{lseven}
\newcounter{leight}
\makeatletter
    
    \@addtoreset{equation}{section}
  \makeatother

\begin{document}

\setcounter{section}{0}
\title{Interpolation error analysis using a new geometric parameter}
\author{Hiroki ISHIZAKA \thanks{Mail: h.ishizaka005@gmail.com} \thanks{Web: \url{https://teamfem.github.io/hiroki_ishizaka/}} 
}
\date{\today}

\maketitle
\pagestyle{plain}

\begin{abstract}
This article presents novel proof methods for estimating interpolation errors, predicated on the understanding that one has already studied foundational error analysis using the finite element method. This article summarizes References \cite{IshPhD,Ish22,Ish24a,Ish24b,Ish24c,Ish25a,Ish25b,IshKobSuzTsu21,IshKobTsu21a,IshKobTsu21b,IshKobTsu23}. We are also correcting any typos found in each paper as we find them. The purpose is to make an easy-to-understand note of  'Special Topics in Finite Element Methods.' 
\end{abstract}


\tableofcontents

\newpage


\section{Preliminalies}

\subsection{General Convention}
Throughout this article, we denote by $c$ a constant independent of $h$ (defined later) and the angles and aspect ratios of simplices, unless specified otherwise {all constants $c$ are bounded if the maximum angle is bounded}. These values vary across different contexts. 

\subsection{Basic Notation}

\begin{longtable}{l|l}
  \label{BasicNotation}
  \centering
    $d$ & The space dimension, $d \in \{ 2,3 \}$\\
    $\mathbb{R}^d$ & $d$-dimensional real Euclidean space \\
    $\mathbb{N}_0$	& $\mathbb{N}_0 := \mathbb{N} \cup \{ 0 \}$	  \\
     $\mathbb{R}_+$	& The set of positive real numbers \\
    $|\cdot|_d$ & $d$-dimensional Hausdorff measure  \\
     $v|_D$ & Restriction of the function $v$ to the set $D$ \\
     	$\dim(V)$ & Dimension of the vector space $V$ \\
     $\delta_{ij}$	& Kronecker delta: $\delta_{ij} = 1$ if $i=j$ and $0$ otherwise \\
     $(x_1,\ldots,x_d)^T$	 & Cartesian coordinates in $\mathbb{R}^d$  \\
\end{longtable}

\subsection{Vectors and Matrices}

\begin{longtable}{l|l}
  \label{VectorsMatrices}
  \centering
     $(v_1 , \ldots,v_d)^T$ & Cartesian components of the vector $v$ in $\mathbb{R}^d$ \\
      $x \cdot y$	& Euclidean scalar product in $\mathbb{R}^d$: $\displaystyle x \cdot y := \sum_{i=1}^d x_i y_i$  \\
     	$|x|_E$& Euclidean norm in $\mathbb{R}^d$: $|x|_E := (x \cdot x)^{1/2}$  \\
     	$\mathbb{R}^{m \times n}$ & Vector space $m \times n$ matrices with real-valued entries  \\
     	${A}$, ${B}$ & Matrices  \\
     	${A}_{ij}$ or $[{A}]_{ij}$  & Entry of ${A}$ in the $i$th and the $j$th column  \\
     	${A}^{\top}$ & Transpose of the matrix ${A}$ \\
     	$\Tr ({A} )$ & Trace of ${A}$: For ${A} \in \mathbb{R}^{m \times n}$, $\displaystyle \Tr ({A}) := \sum_{i=1}^d {A}_{ii}$  \\
     	 $\det(A)$ & Determinant of $\displaystyle {A}$ \\
     	$\diag({A})$ & Diagonal of ${A}$:   \\
     	& For ${A} \in \mathbb{R}^{m \times n}$, $\diag ({A})_{ij} := \delta_{ij} {A}_{ij}$, $1 \leq i,j \leq d$ \\
     ${A} x$	& Matrix-vector product: \\
     & For ${A} \in \mathbb{R}^{m \times n}$  and $x \in \mathbb{R}^n$, $\displaystyle ({A} x)_i := \sum_{j=1}^d {A}_{ij} x_j$ for $1 \leq i \leq d$  \\
     ${A}:{B}$	& Double contraction:  \\
     	& For ${A} \in \mathbb{R}^{m \times n}$ and ${B} \in \mathbb{R}^{m \times n}$,  $\displaystyle {A}: {B} := \sum_{i=1}^m \sum_{j=1}^n {A}_{ij} {B}_{ij}$  \\
     $\| {A} \|_2$	 & Operator norm of ${A}$: For ${A} \in \mathbb{R}^{d \times d}$, $\displaystyle \| {A} \|_2 := \sup_{0 \neq x \in \mathbb{R}^d} \frac{|{A} x|_E}{|x|_E}$ \\
       $\| {A} \|_{\max}$ & Max norm of ${A}$: For ${A} \in \mathbb{R}^{d \times d}$, $\displaystyle \| {A} \|_{\max} := \max_{1 \leq i,j \leq d} |{A}_{ij}|$ \\
         $O(d)$ & $O(d)$ consists of all orthogonal matrices of determinant $\pm 1$  \\
\end{longtable}

In this article, we use the following facts.

For ${A} \in \mathbb{R}^{m \times n}$, it holds that
\begin{align}
\displaystyle
\| {A} \|_{\max} \leq \| {A} \|_{2} \leq \sqrt{mn} \| {A} \|_{\max},\label{Anorm}
\end{align}
e.g., see \cite[p. 56]{GolLoa96}. For ${A}, {B} \in \mathbb{R}^{m \times m}$, it holds that
\begin{align}
\displaystyle
\| {A} {B} \|_{2} \leq \| {A} \|_{2} \| {B} \|_{2}. \label{Mtriangle}
\end{align}
If ${A}^{\top} {A}$ is a positive definite matrix in $\mathbb{R}^{d \times d}$, the spectral norm of the matrix ${A}^{\top} {A}$ is the largest eigenvalue of ${A}^{\top} {A}$; i.e.,
\begin{align}
\displaystyle
\| {A} \|_2 = \left( \lambda_{\max}({A}^{\top} {A}) \right)^{1/2} = \sigma_{\max}({A}), \label{leigen}
\end{align}
where $\lambda_{\max}({A})$ and $\sigma_{\max}({A})$ are respectively the largest eigenvalues and singular values of ${A}$.

If ${A} \in O(d)$, because ${A}^{\top} = {A}^{-1}$ and
\begin{align*}
\displaystyle
|{A} x|_E^2
&= ({A} x)^{\top} ({A} x) = x^{\top} {A}^{\top} {A} x = x^{\top} {A}^{-1} {A} x = |x|_E^2,
\end{align*}
it holds that
\begin{align*}
\displaystyle
 \| {A} \|_2 &= \sup_{0 \neq x \in \mathbb{R}^d} \frac{|{A} x|_E}{|x|_E}
 = \sup_{0 \neq x \in \mathbb{R}^d} \frac{|x|_E}{|x|_E}  = 1.
\end{align*}

\subsection{Function Spaces}
This article uses standard Sobolev spaces with associated norms (e.g., see \cite{Beretal24,ErnGue04,ErnGue21a}). 

\subsection{Function Space $H(\div;D)$}
Let $D$ be a Lipschitz domain of $\mathbb{R}^d$. We denote the function space by
\begin{align*}
\displaystyle
H(\div;D) :=\left\{ v \in L^2(D)^d ; \ \div v \in L^2(D) \right\},
\end{align*}
which is a Hilbert space with the inner product and norm:
\begin{align*}
\displaystyle
(u,v)_{H(\div;D)} &:= (u,v) + (\div u , \div v), \\
\| v \|_{H(\div;D)} &:= (v,v)_{H(\div;D)}^{1/2} = \left ( \| v \|_{L^2(D)^d}^2 + \| \div v \|_{L^2(D)}^2  \right)^{1/2}.
\end{align*}

\begin{thr}
The space $\mathcal{C}^{\infty}(\overline{D})^d$ is dense in $H(\div;D)$.	
\end{thr}

\begin{pf*}
A proof can be found in \cite[Theorem 2.4]{GirRav86}. The condition of "boundedness" is entered into the assumptions because we use the space $\mathcal{C}^{\infty}(\overline{D})^d$. 
\qed
\end{pf*}

\begin{thr} \label{thrB22}
The trace operator $\gamma^d: \mathcal{C}^{\infty}(\overline{D})^d \to \mathcal{C}^{\infty}(\overline{\partial D})$ which maps $\varphi \mapsto \varphi \cdot n |_{\partial D}$ can be extended to a continuous, linear mapping
\begin{align*}
\displaystyle
\gamma^d: H(\div;D) \to H^{- \frac{1}{2}} (\partial D),
\end{align*}
where $H^{- \frac{1}{2}} (\partial D)$ is the dual space $H^{\frac{1}{2}} (\partial D)$.
\end{thr}

\begin{pf*}
A proof can be found in \cite[Theorem 2.5]{GirRav86}. 
\qed
\end{pf*}

\begin{thr}
The trace theorem is optimal in the sense that
\begin{align*}
\displaystyle
\gamma^d: H(\div;D) \to H^{- \frac{1}{2}} (\partial D),
\end{align*}
is surjective.
\end{thr}

\begin{pf*}
Let $\mu \in  H^{- \frac{1}{2}} (\partial D)$. To show is that there exists $v \in H(\div;D)$ such that
\begin{align*}
\displaystyle
v \cdot n &= \mu \quad \text{on $\partial D$}, \\
\| v \|_{H(\div;D)} &\leq \| v \cdot n \|_{H^{- \frac{1}{2}} (\partial D)}.
\end{align*}
We know that the problem
\begin{align*}
\displaystyle
- \varDelta \varphi + \varphi = 0 \quad \text{in $D$}, \quad \frac{\partial \varphi}{\partial n} = \mu \quad \text{on $\partial D$}
\end{align*}
has a unique solution $\varphi \in H^1(D)$ satisfying
\begin{align}
\displaystyle
\| \varphi \|^2_{H^1(D)} = \langle \mu , \varphi \rangle_{\partial D} \leq \| \mu \|_{H^{- \frac{1}{2}} (\partial D)} \| \varphi \|_{H^1(D)},
\end{align}
see \cite[Section 1.4 and (1.16)]{GirRav86}. Setting $v = \nabla \varphi$, we have $v \in H(\div;D)$, $v \cdot n = \mu$, and
\begin{align*}
\displaystyle
\| v \|_{H(\div;D)}
&= \left( \| v \|_{L^2(D)^d}^2 + \| \div v \|^2_{L^2(D)} \right)^{1/2} = \| \varphi \|_{H^1(D)} \\ 
&\leq \| \mu \|_{H^{- \frac{1}{2}} (\partial D)} 
= \| \mu \|_{H^{- \frac{1}{2}} (\partial D)} =  \| v \cdot n\|_{H^{- \frac{1}{2}} (\partial D)}.
\end{align*}
\qed
\end{pf*}

\begin{thr}
It holds that
\begin{align*}
\displaystyle
H_0(\div;D) := \ker(\gamma^d) = \{ v \in H(\div;D): \ v \cdot n |_{\partial D} = 0 \}.
\end{align*}
\end{thr}

\begin{pf*}
A proof can be found in \cite[Theorem 2.6]{GirRav86}. 
\qed
\end{pf*}

\begin{thr}
Let
\begin{align*}
\displaystyle
H_{\sigma} := \{ v \in H_0(\div;D): \div v = 0 \}.
\end{align*}
It then holds that
\begin{align*}
\displaystyle
L^2(D)^d = H_{\sigma} \oplus H^{\perp},
\end{align*}
where $H^{\perp}$ denotes the orthogonal of $H_{\sigma}$ in $L^2(D)^d$ for the scalar product, that is,
\begin{align*}
\displaystyle
H^{\perp} := \{ v = \nabla q  : \ q \in H^1(D)\}.
\end{align*}
\end{thr}

\begin{pf*}
A proof can be found in \cite[Theorem 2.7]{GirRav86}. Remark that $D$ is open, bounded, connected, and a Lipschitz set, because $D$ is a Lipschitz domain of $\mathbb{R}^d$.
\qed
\end{pf*}

\subsection{Finite-Element-Methods-Related Symbols} \label{MeshSym}
\subsubsection{Symbols}
\begin{longtable}{l|l}
  \label{FunctionSpaces}
  \centering
 $\mathbb{P}^k$ & Vector space of polynomials in the variables $x_1,\ldots,x_d$ of  \\
  & global degree at most $k \in \mathbb{N}_0$ \\
  $N^{(d,k)}$ & $N^{(d,k)} := \dim (\mathbb{P}^k) = \begin{pmatrix}
 d+k \\
 k
\end{pmatrix}
$ \\
$\mathbb{RT}^k$ & The Raviart--Thomas polynomial space of order $k \in \mathbb{N}_0$ as \\
&  $\mathbb{RT}^k := ( \mathbb{P}^k)^d + x \mathbb{P}^k$ for any $x \in \mathbb{R}^d$\\
$N^{(RT)}$ &$N^{(RT)} := \dim RT^k$ \\
$T$, $\widetilde{T}$, $\widehat{T}$, $K$ & Closed simplices in $\mathbb{R}^d$ \\
$ \mathbb{P}^k({T})$, $\mathbb{RT}^k({T})$ & $\mathbb{P}^k({T})$ (or $\mathbb{RT}^k({T})$) is spanned by the restriction to  ${T}$\\
& of polynomials in $\mathbb{P}^k$ (or $\mathbb{RT}^k$) \\
\end{longtable}

\subsubsection{Meshes}
Let $\Omega \subset \mathbb{R}^d$, $d \in \{ 2 , 3 \}$, be a bounded polyhedral domain. Furthermore, we assume that $\Omega$ is convex if necessary. Let $\mathbb{T}_h = \{ T \}$ be a simplicial mesh of $\overline{\Omega}$ made up of closed $d$-simplices, such as
\begin{align*}
\displaystyle
\overline{\Omega} = \bigcup_{T \in \mathbb{T}_h} T,
\end{align*}
with $h := \max_{T \in \mathbb{T}_h} h_{T}$, where $ h_{T} := \diam(T)$. We also use a symbol $\rho_{T}$ which means the radius of the largest ball inscribed in $T$.  We assume that each face of any $d$-simplex $T_1$ in $\mathbb{T}_h$ is either a subset of the boundary $\partial \Omega$ or a face of another $d$-simplex $T_2$ in $\mathbb{T}_h$. That is, $\mathbb{T}_h$ is a simplicial mesh of $\overline{\Omega}$ without hanging nodes. Such mesh $\mathbb{T}_h$ is said to be conformal.  Let $\{ \mathbb{T}_h \}$ be a family of conformal meshes. 

Let $T$ be a simplex of $\mathbb{T}_h$ which is a convex full of $d+1$ vertices, $p_1, \ldots,p_{d+1}$, that do not belong to the same hyperplane. Let $S_i$ be the face of a simplex $T$ opposite to the vertex $p_i$. For $d=3$, angles between faces of a tetrahedron are called \textit{dihedral}, whereas angles between its edges are called \textit{solid}.

\subsubsection{Broken Sobolev Spaces, Mesh faces, Averages and Jumps}
 Let $\mathcal{F}_h^i$ be the set of interior faces, and $\mathcal{F}_h^{\partial}$ be the set of faces on boundary $\partial \Omega$. We set $\mathcal{F}_h := \mathcal{F}_h^i \cup \mathcal{F}_h^{\partial}$. For any $F \in \mathcal{F}_h$, we define the unit normal $n_F$ to $F$ as follows: (\roman{sone}) If  $F \in \mathcal{F}_h^i$ with $F = T_{\natural} \cap T_{\sharp}$, $T_{\natural},T_{\sharp} \in \mathbb{T}_h$, $\natural > \sharp$, let $n_F$ be the unit normal vector from $T_{\natural}$ to  $ T_{\sharp}$.  (\roman{stwo}) If $F \in \mathcal{F}_h^{\partial}$, $n_F$ is {the unit outward normal} $n$ to $\partial \Omega$. We also use the following set. For any $F \in \mathcal{F}_h$,
\begin{align*}
\displaystyle
\mathbb{T}_F := \{ T \in \mathbb{T}_h: \ F \subset T \}.
\end{align*}
Furthermore, for a simplex $T \subset \mathbb{R}^d$, let $\mathcal{F}_{T}$ be the collection of the faces of $T$.

We consider $\mathbb{R}^q$-valued functions for some $q \in \mathbb{N}$. Let $p \in [1,\infty]$ and $s \> 0$ be a positive real number. We define a broken (piecewise) Sobolev space as
\begin{align*}
\displaystyle
W^{s,p}(\mathbb{T}_h;\mathbb{R}^q) := \{ v \in L^p(\Omega;\mathbb{R}^q): \ v|_{T} \in W^{s,p}(T;\mathbb{R}^q) \ \forall T \in \mathbb{T}_h  \}
\end{align*}
with the norms
\begin{align*}
\displaystyle
\| v \|_{W^{s,p}(\mathbb{T}_h;\mathbb{R}^q)} &:= \left( \sum_{T \in \mathbb{T}_h} \| v \|^p_{W^{s,p}(T;\mathbb{R}^q) } \right)^{1/p} \quad \text{if $p \in [1,\infty)$},\\
\| v \|_{W^{s,\infty}(\mathbb{T}_h;\mathbb{R}^q)} &:= \max_{T \in \mathbb{T}_h}  \| v \|_{W^{s,\infty}(T;\mathbb{R}^q)}.
\end{align*}
When $q=1$, we denote $W^{s,p}(\mathbb{T}_h) := W^{s,p}(\mathbb{T}_h;\mathbb{R}) $. When $p=2$, we write $H^s(\mathbb{T}_h)^q := H^s(\mathbb{T}_h;\mathbb{R}^q) := W^{s,2}(\mathbb{T}_h;\mathbb{R}^q)$ and  $H^s(\mathbb{T}_h) := W^{s,2}(\mathbb{T}_h;\mathbb{R})$. We use the norm
\begin{align*}
\displaystyle
| \varphi |_{H^1(\mathbb{T}_h)} &:= \left( \sum_{T \in \mathbb{T}_h}\| \nabla \varphi \|^2_{L^2(T)^d} \right)^{1/2} \quad \varphi \in H^1(\mathbb{T}_h).
\end{align*}
Let $\varphi \in H^1(\mathbb{T}_h)$. Suppose that $F \in \mathcal{F}_h^i$ with $F = T_{\natural} \cap T_{\sharp}$, $T_{\natural},T_{\sharp} \in \mathbb{T}_h$, $\natural > \sharp$. We set $\varphi_{\natural} := \varphi{|_{T_{\natural}}}$ and $\varphi_{\sharp} := \varphi{|_{T_{\sharp}}}$. The jump in $\varphi$ across $F$ is defined as
\begin{align*}
\displaystyle
[\! [ \varphi ]\!] := [\! [ \varphi ]\!]_F := \varphi_{\natural} - \varphi_{\sharp}, \quad \natural > \sharp.
\end{align*}
For a boundary face $F \in \mathcal{F}_h^{\partial}$ with $F = \partial T \cap \partial \Omega$, $[\![\varphi ]\!]_F := \varphi|_{T}$. For any $v \in H^1(\mathbb{T}_h)^d$, the notations
\begin{align*}
\displaystyle
[\![ v \cdot n ]\!] &:= [\![ v \cdot n ]\!]_F := v_{\natural} \cdot n_F - v_{\sharp} \cdot n_F, \quad \natural > \sharp, \\
[\![v]\!] &:=  [\![ v]\!]_F := v_{\natural} - v_{\sharp}, \quad \natural > \sharp,
\end{align*}
denote the jump in the normal component of $v$ and the jump of $v$. Set two nonnegative real numbers $\omega_{T_{\natural},F}$ and $\omega_{T_{\sharp},F}$ such that
\begin{align*}
\displaystyle
\omega_{T_{\natural},F} + \omega_{T_{\sharp},F} = 1.
\end{align*}
The skew-weighted average of $\varphi$ across $F$ is then defined as
\begin{align*}
\displaystyle
\{\! \{ \varphi\} \! \}_{\overline{\omega}} :=  \{\! \{ \varphi\} \! \}_{\overline{\omega},F} := \omega_{T_{\sharp},F} \varphi_{\natural} + \omega_{T_{\natural},F} \varphi_{\sharp}.
\end{align*}
For a boundary face $F \in \mathcal{F}_h^{\partial}$ with $F = \partial T \cap \partial \Omega$, $\{\! \{ \varphi \} \!\}_{\overline{\omega}} := \varphi |_{T}$. Furthermore, 
\begin{align*}
\displaystyle
\{\! \{ v\} \! \}_{\omega} :=  \{\! \{ v \} \! \}_{\omega,F} := \omega_{T_{\natural},F} v_{\natural} + \omega_{T_{\sharp},F} v_{\sharp},
\end{align*}
for the weighted average of $v$. For any $v \in H^{1}(\mathbb{T}_h)^d$ and $\varphi \in H^{1}(\mathbb{T}_h)$,
\begin{align*}
\displaystyle
[\![ (v \varphi) \cdot n ]\!]_F
&=  \{\! \{ v \} \! \}_{\omega,F} \cdot n_F [\! [ \varphi ]\!]_F + [\![ v \cdot n ]\!]_F \{\! \{ \varphi\} \! \}_{\overline{\omega},F}.
\end{align*}

We define a broken gradient operator as follows. Let $p \in [1,\infty]$. For $\varphi \in W^{1,p}(\mathbb{T}_h)$, the broken gradient $\nabla_h: W^{1,p}(\mathbb{T}_h) \to L^p(\Omega)^{d}$ is defined by
\begin{align*}
\displaystyle
(\nabla_h \varphi)|_{T} &:= \nabla (\varphi|_{T}) \quad \forall T \in \mathbb{T}_h,
\end{align*}
and we define the broken $H(\div;T)$ space by
\begin{align*}
\displaystyle
H(\div;\mathbb{T}_h) := \left \{ v \in L^2(\Omega)^d; \ v |_{T} \in H(\div;T) \ \forall T \in \mathbb{T}_h  \right\},
\end{align*}
and the broken divergence operator $\divh : H(\div;\mathbb{T}_h) \to L^2(\Omega)$ such that, for all $v \in H(\div;\mathbb{T}_h)$,
\begin{align*}
\displaystyle
(\divh v)|_{T} := \div (v |_{T}) \quad \forall T \in \mathbb{T}_h.
\end{align*}

\subsubsection{Barycentric Coodinates} \label{Bary}
For a simplex $T \subset \mathbb{R}^d$, let $\{ p_i \}_{i=1}^{d+1}$ be vertices of $T$ and $(x_1^{(i)} , \ldots , x_d^{(i)})^T$ coordinates of $p_i$. We set
\begin{align*}
\displaystyle
\Delta := \det
\begin{pmatrix}
1  & \cdots& 1 \\
 x_1^{(1)} & \cdots &  x_1^{(d+1)}  \\
 \vdots & \vdots & \vdots  \\
 x_d^{(1)} & \cdots & x_d^{(d+1)}  \\
\end{pmatrix}
\> 0.
\end{align*}
The barycentric coordinates $\{ \lambda_i \}_{i=1}^{d+1}: \mathbb{R}^d \to \mathbb{R}$ of the point $p(x_1,\ldots,x_d)$ with respect to  $\{ p_i \}_{i=1}^{d+1}$ are then defined as
\begin{align*}
\displaystyle
\lambda_i (x) := \frac{1}{\Delta} \det \
\begin{blockarray}{cccccc}
 & & i & & \\[-5pt]
 & & \smallsmile & & \\[-5pt]
\begin{block}{(cccccc)}
 1 & \cdots & 1 & \cdots & 1 \\
 x_1^{(1)} & \cdots & x_1 & \cdots & x_1^{(d+1)}  \\
 \vdots &  & \vdots &  & \vdots  \\
x_d^{(1)} & \cdots & x_d  & \cdots  & x_d^{(d+1)}  \\
\end{block}
\end{blockarray}
\ .
\end{align*}
The barycentric coordinates have the following properties:
\begin{align*}
\displaystyle
\lambda_i(p_j) = \delta_{ij}, \quad \sum_{i=1}^{d+1} \lambda_i(x) = 1.
\end{align*}

 \subsection{Useful Tools for Analysis}
 \subsubsection{Jensen-type Inequality}
Let $r,s$ be two nonnegative real numbers and $\{ x_i \}_{i \in I}$ be a finite sequence of nonnegative numbers. It then holds that
\begin{align}
\displaystyle
\begin{cases}
\left( \sum_{i \in I} x_i^s \right)^{\frac{1}{s}} \leq \left( \sum_{i \in I} x_i^r \right)^{\frac{1}{r}} \quad \text{if $r \leq s$},\\
\left( \sum_{i \in I} x_i^s \right)^{\frac{1}{s}} \leq \card(I)^{\frac{r-s}{rs}} \left( \sum_{i \in I} x_i^r \right)^{\frac{1}{r}} \quad \text{if $r \> s$},
\end{cases} \label{jensen}
\end{align}
see \cite[Exercise 12.1]{ErnGue21a}.

\subsubsection{Embedding Theorems}
The following is well known as the Sobolev embedding theorem.

\begin{thr} \label{thr=intro1}
 Let $d \geq 2$, $s \> 0$, and $p \in [1,\infty]$. Let $D \subset \mathbb{R}^d$ be a bounded open subset of $\mathbb{R}^d$. If $D$ is a Lipschitz set, we then have
\begin{align}
\displaystyle
W^{s,p}(D) \hookrightarrow
\begin{cases}
L^q(D) \quad \text{$\forall q \in [p , \frac{pd}{d - sp}]$ if $s p \< d$}, \\
L^q(D) \quad \text{$\forall q \in [p,\infty)$, if $sp=d$},\\
L^{\infty}(D) \cap \mathcal{C}^{0,\xi}(\overline{D}) \quad \text{$\xi = 1 - \frac{d}{s p}$ if $sp \> d$}.
\end{cases} \label{emmed}
\end{align} 
Furthermore, 
\begin{align}
\displaystyle
W^{s,p}(D) \hookrightarrow L^{\infty}(D) \cap \mathcal{C}^{0}(\overline{D}) \quad \text{(case $s=d$ and $p=1$)}. \label{emmed1}
\end{align} 
\end{thr}

\begin{pf*}
See, for example, \cite[Corollary B.43,  Theorem B.40]{ErnGue04} and \cite[Theorem 2.31]{ErnGue21a} and the references therein.
\qed
\end{pf*}

The following is the embedding theorem related to operator from $W^{s,p}(D)$ into $L^q(S_r)$, where $S_r$ is some plane $r$-dimensional piece belonging to $D$ with dimensions $r \< d$.

 \begin{thr} \label{thr=intro2}
  Let $p,q \in [1,+\infty]$ and $s \geq 1$ be an integer. Let $D \subset \mathbb{R}^d$ be a bounded open set having piecewise smooth boundaries. The following embeddings are then continuous:
\begin{align}
\displaystyle
W^{s,p}(D) \hookrightarrow
\begin{cases}
L^q(S_r) \quad \text{if $1 \leq p \< \frac{d}{s}$, $r \> d - sp$ and $q \leq \frac{p r}{d - sp}$}, \\
L^q(S_r) \quad \text{if $p = \frac{d}{s}$ for $q \< + \infty$}.
\end{cases} \label{emmed2}
\end{align} 
\end{thr}
 
\begin{pf*}
See, for example, \cite[Theorem 2.1 (p. 61)]{LadSolUra68} and the references therein.
\qed
\end{pf*}

\subsubsection{Trace Theorem}
The following theorem is well-known, e.g., see \cite{ErnGue21a}.

\begin{thr}[Trace] \label{thr164}
Let $p \in [1,\infty)$. Let $s \> \frac{1}{p}$ if $p \> 1$ or $s \geq 1$ if $p=1$. Let $D$ be a Lipschitz domain (e.g., see \cite[Definition 3.2]{ErnGue21a}) in $\mathbb{R}^d$. There exists a bounded linear operator $\gamma^g:W^{s,p}(D) \to L^p(\partial D)$ such that
\begin{enumerate}
 \item $\gamma^g(\varphi) = \varphi|_{\partial D}$, whenever $\varphi$ is smooth, e.g., $\varphi \in \mathcal{C}(\overline{D})$.
 \item The kernel of $\gamma^g$ is $W^{s,p}_0(D)$.
 \item If $s=1$ and $p=1$, or if $s \in (\frac{1}{2},\frac{3}{2})$ and $p=2$, or if $s \in (\frac{1}{p},1]$ and $p \not\in \{1,2\}$, then $\gamma^g:W^{s,p}(D) \to W^{s - \frac{1}{p},p}(\partial D)$ is bounded and surjective, that is, there exists $C^{\gamma^g}$ such that, for every functions $g \in W^{s-\frac{1}{p},p}(\partial D)$, one can find a function $\varphi_g \in W^{s,p}(D)$, called a lifting of $g$, such that 
\begin{align}
\displaystyle
\gamma^g(\varphi_g) = g, \quad \| \varphi_g \|_{W^{s,p}(D)} \leq C^{\gamma^g} \ell_D^{\frac{1}{p}} \| g \|_{W^{s-\frac{1}{p},p}(\partial D)}, \label{tr65}
\end{align}
where $\ell_D$ is a characteristic length of $D$, e.g., $\ell_D := \diam (D)$.
\end{enumerate}
\end{thr}

\begin{pf*}
See \cite[Theorem 3.10]{ErnGue21a}, and the references therein.
\qed
\end{pf*}

\begin{thr}[Trace on low-dimensional manifolds] \label{low=trace}
Let $p \in [1,\infty)$ and let $D$ be a Lipschitz domain in $\mathbb{R}^d$. Let $M$ be a smooth, or polyhedral, manifold of dimension $r$ in $\overline{D}$, $r \in \{ 0, \ldots ,d\}$. Then, there exists a bounded trace operator from $W^{s,p}(D)$ to $L^p(M)$, provided $sp \> d-r$, or $s \geq d -r$ if $p=1$.
\end{thr}

\begin{pf*}
See \cite[Theorem 3.15]{ErnGue21a}.
\qed
\end{pf*}

\subsubsection{Bramble--Hilbert--type Lemma} \label{sec=BHL}
The Bramble--Hilbert--type lemma (e.g., see \cite{DupSco80,BreSco08}) plays a major role in interpolation error analysis. We use the following estimates on anisotropic meshes proposed in \cite[Lemma 2.1]{Ape99}.

\begin{lem} \label{lem=intro1}
Let $D \subset \mathbb{R}^d$ be a connected open set that is star-shaped concerning balls $B$. Let $\gamma$ be a multi-index with $m := |\gamma|$ and $\varphi \in L^1(D)$ be a function with $\partial^{\gamma} \varphi \in W^{\ell -m,p}(D)$, where $\ell \in \mathbb{N}$, $m \in \mathbb{N}_0$, $0 \leq m \leq \ell$, $p \in [1,\infty]$. It then holds that
\begin{align}
\displaystyle
\| \partial^{\gamma} (\varphi - Q^{(\ell)} \varphi) \|_{W^{\ell -m,p}(D)} \leq C^{BH} |\partial^{\gamma} \varphi|_{W^{\ell-m,p}(D)},  \label{BH=1}
\end{align}
where $C^{BH}$ depends only on $d$, $\ell$, $\diam D$, and $\diam B$, and $ Q^{(\ell)} \varphi$ is defined as
\begin{align}
\displaystyle
(Q^{(\ell)} \varphi)(x) := \sum_{|\delta| \leq \ell -1} \int_B \eta(y) (\partial^{\delta}\varphi)(y) \frac{(x-y)^{\delta}}{\delta !} dy \in \mathbb{P}^{\ell -1},  \label{BH=2}
\end{align}
where $\eta \in \mathcal{C}_0^{\infty}(B)$ is a given function with $\int_B \eta dx = 1$.
\end{lem}

To give local interpolation error estimates on isotropic meshes, we use the inequalities given in  \cite[Theorem 1.1]{DekLev04} and \cite{DupSco80,BreSco08,Ver99} which are variants of the Bramble--Hilbert lemma.

\begin{lem} \label{DekLev=lem}
Let ${D} \subset \mathbb{R}^d$ be a bounded convex domain. Let $\varphi \in W^{m,p}({D})$ with $m \in \mathbb{N}$ and $1 \leq p \leq \infty$. There exists a polynomial $\eta \in \mathbb{P}^{m-1}$ such that
\begin{align}
\displaystyle
 | \varphi - \eta |_{ W^{k,p}({D})} \leq C^{BH}(d,m) \diam({D})^{m-k} | \varphi |_{ W^{m,p}({D})}, \quad k=0,1,\ldots,m.  \label{DekLev}
\end{align}
\end{lem}

\begin{pf*}
The proof is found in \cite[Theorem 1.1]{DekLev04}.
\qed
\end{pf*}

\begin{rem} \label{BA=rem}
In \cite[Lemma 4.3.8]{BreSco08}, the Bramble--Hilbert lemma is given as follows. Let $B$ be a ball in $D \subset \mathbb{R}^d$ such that $D$ is star-shaped with respect to $B$ and its radius $r \> \frac{1}{2} r_{\max}$, where $r_{\max} := \sup \{ r: D$ is star-shaped with respect to a ball of radius $r \}$. Let $\varphi \in W^{m,p}(D)$ with $m \in \mathbb{N}$ and $1 \leq p \leq \infty$. There exists a polynomial $\eta \in \mathbb{P}^{m-1}$ such that
\begin{align}
\displaystyle
 | \varphi - \eta |_{ W^{k,p}({D})} \leq C^{BH}(d,m,\gamma) \diam({D})^{m-k} | \varphi |_{ W^{m,p}({D})}, \quad k=0,1,\ldots,m. \label{BS}
\end{align}
Here, $\gamma$ is called the chunkiness parameter of $D$, which is defined by
\begin{align*}
\displaystyle
\gamma := \frac{\diam({D})}{r_{\max}}.
\end{align*}
The main drawback is that the constant $C^{BH}(d,m,\gamma)$ depends on the chunkiness parameter. Meanwhile, the constant $ C^{BH}(d,m)$ of the estimate \eqref{DekLev} does not depend on the geometric parameter $\gamma$.
\end{rem}

\begin{rem}
For general Sobolev spaces $W^{m,p}(\Omega)$, the upper bounds on the constant $C^{BH}(d,m)$ are not given, as far as we know. However, when $p=2$, the following result has been obtained by Verf$\it{\ddot{u}}$rth \cite{Ver99}.

Let ${D} \subset \mathbb{R}^d$ be a bounded convex domain.  Let $\varphi \in H^{m}({D})$ with $m \in \mathbb{N}$.  There exists a polynomial $\eta \in \mathbb{P}^{m-1}$ such that
\begin{align}
\displaystyle
 | \varphi - \eta |_{ H^{k}({D})} \leq C^{BH}(d,k,m) \diam({D})^{m-k} | \varphi |_{ H^{m}({D})}, \quad k=0,1,\ldots,m-1. \label{Ver}
\end{align}
Verf$\it{\ddot{u}}$rth has given upper bounds on the constants in the estimates such that
\begin{align*}
\displaystyle
 C^{BH}(d,k,m) \leq \pi^{k-m}
\begin{pmatrix}
 d + k -1 \\
 k
\end{pmatrix}
^{1/2} 
\frac{\{ (m-k) ! \}^{1/2}}{\{ \left[ \frac{m-k}{d} \right] !  \}^{d/2}},
\end{align*}
where $[x]$ denotes the largest integer less than or equal to $x$. 

As an example, let us consider the case $d=3$, $k=1$, and $m=2$. We then have
\begin{align*}
\displaystyle
 C^{BH}(3,1,2) \leq \frac{\sqrt{3}}{\pi},
\end{align*}
thus on the standard reference element $\widehat{T}$, we obtain
\begin{align*}
\displaystyle
 | \hat{\varphi} - \hat{\eta} |_{ H^{1}({\widehat{T}})} \leq  \frac{\sqrt{6}}{\pi} | \hat{\varphi} |_{ H^{2}({\widehat{T}})} \quad \forall \hat{\varphi} \in H^{2}({\widehat{T}}),
\end{align*}
because $\diam (\widehat{T}) = \sqrt{2}$.

\end{rem}

\subsubsection{Poincar\'e inequality}

\begin{thr}[Poincar\'e inequality] \label{thr=poincare}
Let $D \subset \mathbb{R}^d$  be a convex domain with diameter $\diam(D)$. It then holds that, for $\varphi \in H^1(D)$ with $\int_D \varphi dx = 0$,  
\begin{align}
\displaystyle
\| \varphi \|_{L^2(D)} \leq \frac{\diam (D)}{\pi} |\varphi|_{H^1(D)}.\label{poincare}
\end{align}
\end{thr}

\begin{pf*}
The proof is found in \cite[Theorem 3.2]{Mar03}, also see \cite{PayWei60}.
\qed
\end{pf*}

\begin{rem}
The coefficient $\frac{1}{\pi}$ of \eqref{poincare} may be improved. 
\end{rem}

\subsection{Abbreviated expression}

\begin{longtable}{l|l}
  \label{abbreviation}
  \centering
    FE & Finite Element\\
    FEMs & Finite Element Methods\\
    CR & Crouzeix--Raviart \\
    RT & Raviart--Thomas\\
\end{longtable}

\section{Isotropic and Anisotropic Mesh Elements} \label{isoani=sec}
In the context of FEMs, mesh elements can be classified based on their geometric properties. \textit{An isotropic mesh element} has equal or nearly equal edge lengths and angles, resulting in a balanced shape. In contrast, \textit{an anisotropic mesh element} features significant variation in edge lengths and angles. 

Consider the following examples: Let $s,\delta \in \mathbb{R}_+$, and $\varepsilon \geq 1$, $\varepsilon \in \mathbb{R}$.

\begin{ex} \label{sec2=ex2=1}
In the case of the simplex $T \subset \mathbb{R}^2$ defined by the vertices $p_1 := (0,0)^{\top}$, $p_2 := (2s,0)^{\top}$, and $p_3 := (s , \delta s)^{\top}$, the triangle is classified as follows:
\begin{itemize}
  \item If $\delta \approx 1$, the triangle $T$ is considered an isotropic mesh element.
  \item Conversely, if $\delta$ is much less than 1, i.e., $\delta \ll 1$, the triangle $T$ becomes an anisotropic mesh element.
\end{itemize}
\end{ex}

\begin{ex} \label{sec2=ex2=2}
In this case, consider the simplex $T \subset \mathbb{R}^2$ defined by the vertices $p_1 := (0,0)^{\top}$, $p_2 := (2s,0)^{\top}$, and $p_3 := (s ,s^{\varepsilon})^{\top}$. Here, the vertex $p_3$ introduces a parameter $\varepsilon$ that can influence the shape of the simplex. The classification of this simplex as isotropic or anisotropic depends on the value of $\varepsilon$:
\begin{itemize}
 \item If $\varepsilon = 1$, the triangle maintains a balanced shape, making it isotropic.
 \item If $\varepsilon \> 1$, the triangle becomes flat when $s \ll 1$, resulting in an anisotropic mesh element.
\end{itemize}
\end{ex}

\begin{ex} \label{sec2=ex2=3}
Consider the simplex $T \subset \mathbb{R}^2$ defined by the vertices $p_1 := (0,0)^{\top}$, $p_2 := (s,0)^{\top}$, and $p_3 := (0 , \delta s)^{\top}$. In this configuration, the classification of the simplex as isotropic or anisotropic depends on the value of $\delta$:
\begin{itemize}
 \item  If $\delta \approx 1$, the triangle is an isotropic mesh element.
 \item  If  $\delta \ll 1$, the triangle becomes an anisotropic mesh element.
\end{itemize}
\end{ex}

\begin{ex} \label{sec2=ex2=4}
Let $T \subset \mathbb{R}^2$ be the simplex with vertices $p_1 := (0,0)^{\top}$, $p_2 := (s,0)^{\top}$, and $p_3 := (0 ,s^{\varepsilon})^{\top}$. In this case, the classification of the simplex as isotropic or anisotropic depends on the value of $\varepsilon$:
\begin{itemize}
 \item  If $\varepsilon = 1$, the triangle is isotropic because the height from $p_3$ is equal to the base length.
 \item If $\varepsilon \> 1$, the triangle will be classified as anisotropic, as the edge lengths will differ significantly when $s \ll 1$.
\end{itemize}
\end{ex}

\begin{ex} \label{sec2=ex2=5}
Let $T \subset \mathbb{R}^2$ be the simplex with vertices $p_1 := (0,0)^{\top}$, $p_2 := (s,0)^{\top}$ and $p_3 := (s^{\delta},s^{\varepsilon})^{\top}$.  The classification of the simplex into two types of anisotropic structures is determined by the values of $\delta$ and $\varepsilon$:
\begin{itemize}
\item If $1 \< \varepsilon \< \delta $, the triangle is flattened so that the point $p_3$ approaches the point $p_1$, i.e. the origin as $s \to 0$.
\item If $1 \<  \delta \< \varepsilon$, the triangle is flattened so that point $p_3$ approaches a point on the straight line $\overline{p_1 p_2}$ that does not include points $p_1$ and $p_2$ as $s \to 0$.
\end{itemize}
\end{ex}

\section{Classical Geometric Conditions}

\subsection{Classical Interpolation Error Estimate}
Let $\widehat{T} \subset \mathbb{R}^d$ and  $T \subset \mathbb{R}^d$ be a reference element and a simplex, respectively. Let these two elements be affine equivalent. Let us consider two finite elements $\{ \widehat{T},\widehat{P}, \widehat{\Sigma} \}$ and $\{ {T},{P}, {\Sigma} \}$ with associated normed vector spaces $V(\widehat{T})$ and $V(T)$. The transformation $\Phi_{T}$ takes the form
\begin{align*}
\displaystyle
&\Phi_{T}: \widehat{T} \ni \hat{x}  \mapsto \Phi_{T}(\hat{x}) := {B}_{T} \hat{x} + b_{T} \in T, 
\end{align*}
where ${B}_{T} \in \mathbb{R}^{d \times d}$ is an invertible matrix and $b_{T} \in \mathbb{R}^{d}$. Let $I_{T}: V(T) := W^{2,p}(T) \to P := \mathbb{P}^1(T)$ with $p \in [1,\infty]$ be an interpolation on $T$ with $I_{T} p = p$ for any $p \in \mathcal{P}^1(T)$. According to the classical theory (e.g., see \cite{Cia02,ErnGue04}), there exists a positive constant $c$, independent of $h_{T}$, such that
\begin{align*}
\displaystyle
|\varphi - I_{T} \varphi|_{W^{1,p}(T)} &\leq c \left( \| {B}_{T} \|_2 \| {B}_{T}^{-1} \|_2 \right ) \| {B}_{T} \|_2 | \varphi |_{W^{2,p}(T)}.
\end{align*}
Here, the quantity $ \| {B}_{T} \|_2 \| {B}_{T}^{-1} \|_2$ is called the \textit{Euclidean condition number} of ${B}_{T}$. By standard estimates (e.g., see \cite[Lemma 1.100]{ErnGue04}), we have
\begin{align*}
\displaystyle
\| {B}_{T} \|_2 \| {B}_{T}^{-1} \|_2 \leq c \frac{h_{T}}{\rho_{T}}, \quad  \| {B}_{T} \|_2 \leq c h_{T}.
\end{align*}
It thus holds that
\begin{align}
\displaystyle
|\varphi - I_{T} \varphi|_{W^{1,p}(T)} &\leq c \frac{h_{T}}{\rho_{T}} h_{T} | \varphi |_{W^{2,p}(T)}. \label{int=sec3=1}
\end{align}
As a geometric condition,  the \textit{shape-regularity condition} is well known to obtain global interpolation error estimates. This condition is stated as follows.

\begin{Cond}[Shape-regularity condition] \label{Cond1}
There exists a constant $\gamma_1 \> 0$ such that
\begin{align}
\displaystyle
\rho_{T} \geq \gamma_1 h_{T} \quad \forall \mathbb{T}_h \in \{ \mathbb{T}_h \}, \quad \forall T \in \mathbb{T}_h. \label{geo1}
\end{align}
\end{Cond}

Under Condition \ref{Cond1}, that is, when the quantity $ \frac{h_{T}}{\rho_{T}}$ is bounded on each $T$, it holds that
\begin{align*}
\displaystyle
|\varphi - I_{h} \varphi|_{W^{1,p}(\Omega)} \leq c h | \varphi |_{W^{2,p}(\Omega)},
\end{align*}
where $I_h \varphi$ is the standard global linear interpolation of $\varphi$ on $\mathbb{T}_h$.

\subsection{Regular Mesh Conditions} \label{regmesh=sec3=2}
Geometric conditions equivalent to the shape-regularity condition are known; that is, the following three conditions are equivalent to the shape-regularity condition \eqref{geo1}. A proof can be found in \cite[Theorem 1]{BraKorKri08}.

\begin{Cond}[Zl\'amal's condition] \label{Cond2}
There exists a constant $\gamma_2 \> 0$ such that for any $\mathbb{T}_h \in \{ \mathbb{T}_h \}$, any simplex $T \in \mathbb{T}_h$ and any dihedral angle $\psi$ and for $d=3$, also any solid angle $\theta$ of $T$, we have
\begin{align}
\displaystyle
\psi \geq \gamma_2, \quad \theta \geq \gamma_2. \label{geo2}
\end{align}
\end{Cond}

\begin{Cond} \label{Cond3}
There exists a constant $\gamma_3 \> 0$ such that for any $\mathbb{T}_h \in \{ \mathbb{T}_h \}$ and any simplex $T \in \mathbb{T}_h$, we have
\begin{align}
\displaystyle
|T|_d \geq \gamma_3 h_{T}^d. \label{geo3}
\end{align}
\end{Cond}

\begin{Cond} \label{Cond4}
There exists a constant $\gamma_4 \> 0$ such that for any $\mathbb{T}_h \in \{ \mathbb{T}_h \}$ and any simplex $T \in \mathbb{T}_h$, we have
\begin{align}
\displaystyle
|T|_d \geq \gamma_4 |B^{T}_d|,  \label{geo4}
\end{align}
where $B^{T} \supset T$ is the circumscribed ball of $T$.
\end{Cond}

\begin{note}
If Condition \ref{Cond1} or \ref{Cond2} or \ref{Cond3} or \ref{Cond4} holds, a family of simplicial partitions is called \textit{regular}.
\end{note}

\begin{note}
Condition \ref{Cond2} was presented by Zl\'amal \cite{Zla68} in 1968. The condition is called the \textit{minimum-angle condition} and guarantees the convergence of finite element methods for linear elliptic problems on $\mathbb{R}^2$. Zl\'amal's condition can be generalised into $\mathbb{R}^n$ for any $n \in \{ 2,3,\ldots \}$. Later, the shape-regularity condition (the inscribed ball condition) was introduced; see \cite{Cia02}. Triangles or tetrahedra cannot be too flat in a shape-regular family of triangulations. 
\end{note}

\begin{note}
Condition \ref{Cond3} seems to be simpler than Condition \ref{Cond1},  Condition \ref{Cond2} and Condition \ref{Cond4}. Therefore, it may be useful to analyse theoretical finite element methods and implement finite element codes to keep nondegenerate mesh partitions.
\end{note}

\subsection{What happens when anisotropic meshes are used?}
Using the equivalence conditions in Section \ref{regmesh=sec3=2}, the error estimate \eqref{int=sec3=1} is rewritten as
\begin{align}
\displaystyle
|\varphi - I_{T} \varphi|_{W^{1,p}(T)} &\leq c \frac{h_{T}^2}{|T|_2} h_{T} | \varphi |_{W^{2,p}(T)}. \label{int=sec3=6}
\end{align}

We considered the following five anisotropic elements as in Section \ref{isoani=sec}: Let $0 \< s, \delta \ll 1$, $s,\delta \in \mathbb{R}$, and $\varepsilon \> 1$, $\varepsilon \in \mathbb{R}$.

\begin{ex} \label{sec3=ex3=8}
Let $T \subset \mathbb{R}^2$ be the simplex with vertices $p_1 := (0,0)^{\top}$, $p_2 := (2s,0)^{\top}$, and $p_3 := (s , \delta s)^{\top}$. Then, we have that $h_T = 2s$, $|T|_2 = \delta s^2$, and
\begin{align*}
\displaystyle
\frac{h_T^2}{|T|_2} = \frac{4}{\delta} \< + \infty.
\end{align*}
Therefore, the shape regularity is satisfied. The estimate \eqref{int=sec3=6} is as follows:
\begin{align*}
\displaystyle
|\varphi - I_{T} \varphi|_{W^{1,p}(T)} &\leq \frac{c}{\delta} h_{T} | \varphi |_{W^{2,p}(T)}.
\end{align*}
When $\delta \ll 1$, the interpolation error \eqref{int=sec3=6} may be large.
\end{ex}

\begin{ex} \label{sec3=ex3=9}
Let $T \subset \mathbb{R}^2$ be the simplex with vertices $p_1 := (0,0)^{\top}$, $p_2 := (2s,0)^{\top}$, and $p_3 := (s ,s^{\varepsilon})^{\top}$. Then, we have that $h_T = 2 s$, $|T|_2 = s^{1+\varepsilon}$ and
\begin{align*}
\displaystyle
\frac{h_T^2}{|T|_2} = 4 s^{1-\varepsilon} \to \infty \ \text{as $s \to 0$}.
\end{align*}
Therefore, the shape-regularity is not satisfied. In this case, when $\varepsilon \> 2$, the estimate \eqref{int=sec3=6} diverges as $s \to 0$.
\end{ex}

\begin{ex} \label{sec3=ex3=10}
Let $T \subset \mathbb{R}^2$ be the simplex with vertices $p_1 := (0,0)^{\top}$, $p_2 := (s,0)^{\top}$, and $p_3 := (0 , \delta s)^{\top}$. Then, we have that $h_T = s \sqrt{1 + \delta^2} \approx s$,  $|T|_2 = \frac{1}{2} \delta s^2$ and
\begin{align*}
\displaystyle
\frac{h_T^2}{|T|_2} =\frac{2(1+\delta^2)}{\delta} \< + \infty.
\end{align*}
Therefore, the shape regularity is satisfied. The estimate \eqref{int=sec3=6} is as follows:
\begin{align}
\displaystyle
|\varphi - I_{T} \varphi|_{W^{1,p}(T)} &\leq \frac{c}{\delta} h_{T} | \varphi |_{W^{2,p}(T)}. \label{int=sec3=7}
\end{align}
It is implied that the interpolation error \eqref{int=sec3=7} may be large when $\delta \ll 1$.
\end{ex}

\begin{ex} \label{sec3=ex3=11}
Let $T \subset \mathbb{R}^2$ be the simplex with vertices $p_1 := (0,0)^{\top}$, $p_2 := (s,0)^{\top}$, and $p_3 := (0 ,s^{\varepsilon})^{\top}$. Subsequently, we obtain $h_T = \sqrt{s^2 + s^{2 \varepsilon}} \approx s$, $|T|_2 = \frac{1}{2} s^{1+\varepsilon}$ and
\begin{align*}
\displaystyle
\frac{h_T^2}{|T|_2} = \frac{2 (s^2 + s^{2 \varepsilon})}{s^{1+\varepsilon}} \to \infty \ \text{as $s \to 0$}.
\end{align*}
Therefore, the shape-regularity condition is not satisfied. In this case, it is implied that the estimate \eqref{int=sec3=6} diverges as $s \to 0$.
\end{ex}

\begin{ex} \label{sec3=ex3=12}
Let $T \subset \mathbb{R}^2$ be the simplex with vertices $p_1 := (0,0)^{\top}$, $p_2 := (s,0)^{\top}$ and $p_3 := (s^{\delta},s^{\varepsilon})^{\top}$. If $1 \< \varepsilon \< \delta $, we have $h_T =  \sqrt{(s - s^{\delta})^2 + s^{2 \varepsilon}}$, $|T|_2 = \frac{1}{2} s^{1 + \varepsilon}$ and 
\begin{align*}
\displaystyle
\frac{h_T^2}{|T|_2} = \frac{2 {(s - s^{\delta})^2 + s^{2 \varepsilon}}}{s^{1 + \varepsilon}} \leq c s^{1 - \varepsilon} \to \infty \ \text{as $s \to 0$}.
\end{align*}
Therefore, the shape-regularity condition is not satisfied. In this case, it is implied that the estimate \eqref{int=sec3=6} diverges as $s \to 0$. If $1 \<  \delta \< \varepsilon$, we have $h_T =  \sqrt{(s - s^{\delta})^2 + s^{2 \varepsilon}}$, $|T|_2 = \frac{1}{2} s^{1 + \varepsilon}$ and 
\begin{align*}
\displaystyle
\frac{h_T^2}{|T|_2} = \frac{2 {(s - s^{\delta})^2 + s^{2 \varepsilon}}}{s^{1 + \varepsilon}} \leq c s^{1 - \varepsilon} \to \infty \ \text{as $s \to 0$}.
\end{align*}
Therefore, the shape-regularity condition is not satisfied. In this case, it is implied that the estimate \eqref{int=sec3=6} diverges as $s \to 0$.
\end{ex}

\begin{rem}
As will be explained later, the factor $\frac{1}{\delta}$ in Example \ref{sec3=ex3=10} is violated. The interpolation error estimate converges in the cases of Example \ref{sec3=ex3=11} and Example \ref{sec3=ex3=12} with $1 \< \varepsilon \< \delta $ using new precise interpolation error estimates under more relaxed geometric conditions.
\end{rem}

\section{Classical Relaxed  Geometric Conditions}

\subsection{Semi-regular Mesh Conditions for $d=2$}
In 1957, Synge \cite[Section 3.8]{Syn57} proposed the following condition.
\begin{Cond}[Synge's condition] \label{Cond5}
There exists $\frac{\pi}{3} \leq \gamma_5 \< \pi$ such that, for any $\mathbb{T}_h \in \{ \mathbb{T}_h \}$ and any simplex $T \in \mathbb{T}_h$, 
\begin{align}
\displaystyle
\theta_{T,\max} \leq \gamma_5,  \label{geo5}
\end{align}
where $\theta_{T,\max}$ is the maximal angle of $T$. 
\end{Cond}
Under Condition \ref{Cond5}, Synge proved an optimal interpolation error estimate as follows.
\begin{align*}
\displaystyle
\| \varphi - I_{h} \varphi \|_{W^{1,p}(\Omega)} \leq c h | \varphi |_{W^{2,p}(\Omega)} \quad \text{for $p = \infty$}.
\end{align*}
The inequality \eqref{geo5} is called \textit{Synge's condition} or the \textit{maximum-angle condition}. In 1976, several author's \cite{BabAzi76,BarGre76,Gre76,Jam76} independently proved the convergence of finite element for $p \< \infty$. It ensures that finite elements converge effectively when the minimum angle approaches zero as the mesh size decreases. If this condition is not met, the accuracy of interpolation for linear triangular elements can suffer, similar to the absence of Zl\'amal's condition, see e.g. \cite[p. 223]{BabAzi76}. This underscores the importance of keeping proper geometric constraints to ensure reliable outcomes in numerical methods. Synge's condition is essential in finite element analysis. 

In \cite{Kri91}, K$\rm{\check{r}}$\'{i}$\rm{\check{z}}$ek proposed the following circumscribed ball condition for $d=2$ which is equivalent to Synge's condition.

\begin{Cond} \label{Cond6}
There exists $\gamma_6\> 0$ such that, for any $\mathbb{T}_h \in \{ \mathbb{T}_h \}$ and any simplex $T \in \mathbb{T}_h$, 
\begin{align}
\displaystyle
\frac{R_2}{h_{T}} \leq \gamma_6,  \label{geo6}
\end{align}
where $R_2$ is the radius of the circumscribed ball of $T \subset \mathbb{R}^2$.
\end{Cond}

\begin{note}
If Condition \ref{Cond5} or \ref{Cond6} holds, the associated families of partitions are called \textit{semi-regular}.
\end{note}

\begin{rem}
Assume that Condition \ref{Cond3} holds, that is, there exists a constant $\gamma_3 \> 0$ such that for any $\mathbb{T}_h \in \{ \mathbb{T}_h \}$ and any simplex $T \in \mathbb{T}_h$, we have
\begin{align*}
\displaystyle
|T|_d \geq \gamma_3 h_{T}^2. 
\end{align*}
Let $T \subset \mathbb{R}^2$ be the triangle with vertices $P_1$, $P_2$ and $P_3$  such that the maximum angle $\theta_{T,\max}$ of $T$ is $\angle P_2 P_1 P_3$. We then have $h_{T} = |P_2 P_3|$ and
\begin{align*}
\displaystyle
\frac{R_2}{h_{T}}
&= \frac{|P_2 P_3|}{2 h_{T} \sin \theta_{T,\max}} = \frac{|P_1 P_2| |P_1 P_3|}{2 |P_1 P_2| |P_1 P_3| \sin \theta_{T,\max}} \leq c \frac{h_{T}^2}{|T|_2} \leq \frac{c}{\gamma_3} =: \gamma_6.
\end{align*}
This implies that each regular family is semi-regular. However, the converse implication does not hold. 	
\end{rem}

\subsection{Semi-regular Mesh Conditions for $d=3$}
Synge's condition \eqref{geo5} is extended to the case of tetrahedra in \cite{Kri92}.

\begin{Cond} \label{Cond7}
There exists a constant $0 \< \gamma_7 \< \pi$ such that for any $\mathbb{T}_h \in \{ \mathbb{T}_h \}$ and any simplex $T \in \mathbb{T}_h$, 
\begin{subequations} \label{geo7}
\begin{align}
\displaystyle
\theta_{T, \max} &\leq \gamma_7, \label{geo7a}\\
\psi_{T, \max} &\leq \gamma_7, \label{geo7b}
\end{align}
\end{subequations}
where $\theta_{T, \max}$ is the maximum angle of all triangular faces of the tetrahedron $T$ and $\psi_{T, \max}$ is the maximum dihedral angle of $T$. 
\end{Cond}

\begin{rem}
The theory of anisotropic interpolation has been advanced through extensive research (\cite{ApeDob92,Ape99,CheShiZha04}).	
\end{rem}

\begin{que}
Is there a semi-regularity condition which equivalent to Synge's condition \eqref{geo7} for $d=3$?	
\end{que}

\begin{rem}
This article introduces a novel geometric condition intended to serve as an alternative to Synge's condition specifically for three-dimensional cases.
\end{rem}

\section{Settings for New Interpolation Theory}

\subsection{Reference Elements} \label{reference}
We first define the reference elements $\widehat{T} \subset \mathbb{R}^d$.

\subsubsection*{Two-dimensional case} \label{reference2d}
Let $\widehat{T} \subset \mathbb{R}^2$ be a reference triangle with vertices $\hat{p}_1 := (0,0){^{\top}}$, $\hat{p}_2 := (1,0){^{\top}}$, and $\hat{p}_3 := (0,1){^{\top}}$. 

\subsubsection*{Three-dimensional case} \label{reference3d}
In the three-dimensional case, we consider the following two cases: (\roman{sone}) and (\roman{stwo}); see Condition \ref{cond2}.

Let $\widehat{T}_1$ and $\widehat{T}_2$ be reference tetrahedra with the following vertices:
\begin{description}
   \item[(\roman{sone})] $\widehat{T}_1$ has vertices $\hat{p}_1 := (0,0,0){^{\top}}$, $\hat{p}_2 := (1,0,0){^{\top}}$, $\hat{p}_3 := (0,1,0){^{\top}}$, and $\hat{p}_4 := (0,0,1)^T$;
 \item[(\roman{stwo})] $\widehat{T}_2$ has vertices $\hat{p}_1 := (0,0,0){^{\top}}$, $\hat{p}_2 := (1,0,0){^{\top}}$, $\hat{p}_3 := (1,1,0){^{\top}}$, and $\hat{p}_4 := (0,0,1)^T$.
\end{description}
Therefore, we set $\widehat{T} \in \{ \widehat{T}_1 , \widehat{T}_2 \}$. 
Note that the case (\roman{sone}) is called \textit{the regular vertex property}, see \cite{AcoDur99}.

\subsection{Two-step Affine Mapping} \label{two=step}
To an affine simplex $T \subset \mathbb{R}^d$, we construct two affine mappings $\Phi_{\widetilde{T}}: \widehat{T} \to \widetilde{T}$ and $\Phi_{T}: \widetilde{T} \to T$. First, we define the affine mapping $\Phi_{\widetilde{T}}: \widehat{T} \to \widetilde{T}$ as
\begin{align}
\displaystyle
\Phi_{\widetilde{T}}: \widehat{T} \ni \hat{x} \mapsto \tilde{x} := \Phi_{\widetilde{T}}(\hat{x}) := {A}_{\widetilde{T}} \hat{x} \in  \widetilde{T}, \label{aff=1}
\end{align}
where ${A}_{\widetilde{T}} \in \mathbb{R}^{d \times d}$ is an invertible matrix. We then define the affine mapping $\Phi_{T}: \widetilde{T} \to T$ as follows:
\begin{align}
\displaystyle
\Phi_{T}: \widetilde{T} \ni \tilde{x} \mapsto x := \Phi_{T}(\tilde{x}) := {A}_{T} \tilde{x} + b_{T} \in T, \label{aff=2}
\end{align}
where $b_{T} \in \mathbb{R}^d$ is a vector and ${A}_{T} \in O(d)$ denotes the rotation and mirror-imaging matrix. We define the affine mapping $\Phi: \widehat{T} \to T$ as
\begin{align*}
\displaystyle
\Phi := {\Phi}_{T} \circ {\Phi}_{\widetilde{T}}: \widehat{T} \ni \hat{x} \mapsto x := \Phi (\hat{x}) =  ({\Phi}_{T} \circ {\Phi}_{\widetilde{T}})(\hat{x}) = {A} \hat{x} + b_{T} \in T, 
\end{align*}
where ${A} := {A}_{T} {A}_{\widetilde{T}} \in \mathbb{R}^{d \times d}$.

\subsubsection*{Construct mapping $\Phi_{\widetilde{T}}: \widehat{T} \to \widetilde{T}$} \label{sec221} 
We consider the affine mapping \eqref{aff=1}. We define the matrix $ {A}_{\widetilde{T}} \in \mathbb{R}^{d \times d}$ as follows. We first define the diagonal matrix as
\begin{align}
\displaystyle
\widehat{A} :=  \diag (h_1,\ldots,h_d), \quad h_i \in \mathbb{R}_+ \quad \forall i,\label{aff=3}
\end{align}
where $\mathbb{R}_+$ denotes the set of positive real numbers.

For $d=2$, we define the regular matrix $\widetilde{A} \in \mathbb{R}^{2 \times 2}$ as
\begin{align}
\displaystyle
\widetilde{A} :=
\begin{pmatrix}
1 & s \\
0 & t \\
\end{pmatrix}, \label{aff=4}
\end{align}
with the parameters
\begin{align*}
\displaystyle
s^2 + t^2 = 1, \quad t \> 0.
\end{align*}
For the reference element $\widehat{T}$, let $\mathfrak{T}^{(2)}$ be a family of triangles.
\begin{align*}
\displaystyle
\widetilde{T} &= \Phi_{\widetilde{T}}(\widehat{T}) = {A}_{\widetilde{T}} (\widehat{T}), \quad {A}_{\widetilde{T}} := \widetilde {A} \widehat{A}
\end{align*}
with the vertices $\tilde{p}_1 := (0,0)^{\top}$, $\tilde{p}_2 := (h_1,0)^{\top}$ and $\tilde{p}_3 :=(h_2 s , h_2 t)^{\top}$. Then, $h_1 = |\tilde{p}_1 - \tilde{p}_2| \> 0$ and $h_2 = |\tilde{p}_1 - \tilde{p}_3| \> 0$. 

For $d=3$, we define the regular matrices $\widetilde{A}_1, \widetilde{A}_2 \in \mathbb{R}^{3 \times 3}$ as follows:
\begin{align}
\displaystyle
\widetilde{A}_1 :=
\begin{pmatrix}
1 & s_1 & s_{21} \\
0 & t_1  & s_{22}\\
0 & 0  & t_2\\
\end{pmatrix}, \
\widetilde{A}_2 :=
\begin{pmatrix}
1 & - s_1 & s_{21} \\
0 & t_1  & s_{22}\\
0 & 0  & t_2\\
\end{pmatrix} \label{aff=5}
\end{align}
with the parameters
\begin{align*}
\displaystyle
\begin{cases}
s_1^2 + t_1^2 = 1, \ s_1 \> 0, \ t_1 \> 0, \ h_2 s_1 \leq h_1 / 2, \\
s_{21}^2 + s_{22}^2 + t_2^2 = 1, \ t_2 \> 0, \ h_3 s_{21} \leq h_1 / 2.
\end{cases}
\end{align*}
Therefore, we set $\widetilde{A} \in \{ \widetilde{A}_1 , \widetilde{A}_2 \}$. For the reference elements $\widehat{T}_i$, $i=1,2$, let $\mathfrak{T}_i^{(3)}$, $i=1,2$, be a family of tetrahedra.
\begin{align*}
\displaystyle
\widetilde{T}_i &= \Phi_{\widetilde{T}_i} (\widehat{T}_i) =  {A}_{\widetilde{T}_i} (\widehat{T}_i), \quad {A}_{\widetilde{T}_i} := \widetilde {A}_i \widehat{A}, \quad i=1,2,
\end{align*}
with the vertices
\begin{align*}
\displaystyle
&\tilde{p}_1 := (0,0,0)^{\top}, \ \tilde{p}_2 := (h_1,0,0)^{\top}, \ \tilde{p}_4 := (h_3 s_{21}, h_3 s_{22}, h_3 t_2)^{\top}, \\
&\begin{cases}
\tilde{p}_3 := (h_2 s_1 , h_2 t_1 , 0)^{\top} \quad \text{for case (\roman{sone})}, \\
\tilde{p}_3 := (h_1 - h_2 s_1, h_2 t_1,0)^{\top} \quad \text{for case (\roman{stwo})}.
\end{cases}
\end{align*}
Subsequently, $h_1 = |\tilde{p}_1 - \tilde{p}_2| \> 0$, $h_3 = |\tilde{p}_1 - \tilde{p}_4| \> 0$, and
\begin{align*}
\displaystyle
h_2 =
\begin{cases}
|\tilde{p}_1 - \tilde{p}_3| \> 0  \quad \text{for case (\roman{sone})}, \\
|\tilde{p}_2 - \tilde{p}_3| \> 0  \quad \text{for case (\roman{stwo})}.
\end{cases}
\end{align*}

\subsubsection*{Construct mapping $\Phi_{T}: \widetilde{T} \to T$}  \label{sec322}
We determine the affine mapping \eqref{aff=2} as follows. Let ${T} \in \mathbb{T}_h$ have vertices ${p}_i$ ($i=1,\ldots,d+1$). Let $b_{T} \in \mathbb{R}^d$ be the vector and ${A}_{T} \in O(d)$ be the rotation and mirror imaging matrix such that
\begin{align*}
\displaystyle
p_{i} = \Phi_T (\tilde{p}_i) = {A}_{T} \tilde{p}_i + b_T, \quad i \in \{1, \ldots,d+1 \},
\end{align*}
where vertices $p_{i}$ ($i=1,\ldots,d+1$) satisfy the following conditions:

\begin{Cond}[Case in which $d=2$] \label{cond1}
Let ${T} \in \mathbb{T}_h$ have vertices ${p}_i$ ($i=1,\ldots,3$). We assume that $\overline{{p}_2 {p}_3}$ is the longest edge of ${T}$, that is, $ h_{{T}} := |{p}_2 - {p}_ 3|$. We set $h_1 = |{p}_1 - {p}_2|$ and $h_2 = |{p}_1 - {p}_3|$. We then assume that $h_2 \leq h_1$. {Because $\frac{1}{2} h_T < h_1 \leq h_T$, ${h_1 \approx h_T}$.} 
\end{Cond}

\begin{Cond}[Case in which $d=3$] \label{cond2}
Let ${T} \in \mathbb{T}_h$ have vertices ${p}_i$ ($i=1,\ldots,4$). Let ${L}_i$ ($1 \leq i \leq 6$) be the edges of ${T}$. We denote by ${L}_{\min}$  the edge of ${T}$ with the minimum length; that is, $|{L}_{\min}| = \min_{1 \leq i \leq 6} |{L}_i|$. We set $h_2 := |{L}_{\min}|$ and assume that 
\begin{align*}
\displaystyle
&\text{the endpoints of ${L}_{\min}$ are either $\{ {p}_1 , {p}_3\}$ or $\{ {p}_2 , {p}_3\}$}.
\end{align*}
Among the four edges sharing an endpoint with ${L}_{\min}$, we consider the longest edge ${L}^{({\min})}_{\max}$. Let ${p}_1$ and ${p}_2$ be the endpoints of edge ${L}^{({\min})}_{\max}$. Thus, we have
\begin{align*}
\displaystyle
h_1 = |{L}^{(\min)}_{\max}| = |{p}_1 - {p}_2|.
\end{align*}
We consider cutting $\mathbb{R}^3$ with a plane that contains the midpoint of the edge ${L}^{(\min)}_{\max}$ and is perpendicular to the vector ${p}_1 - {p}_2$. Thus, there are two cases. 
\begin{description}
  \item[(Type \roman{sone})] ${p}_3$ and ${p}_4$  belong to the same half-space;
  \item[(Type \roman{stwo})] ${p}_3$ and ${p}_4$  belong to different half-spaces.
\end{description}
In each case, we set
\begin{description}
  \item[(Type \roman{sone})] ${p}_1$ and ${p}_3$ as the endpoints of ${L}_{\min}$, that is, $h_2 =  |{p}_1 - {p}_3| $;
  \item[(Type \roman{stwo})] ${p}_2$ and ${p}_3$ as the endpoints of ${L}_{\min}$, that is, $h_2 =  |{p}_2 - {p}_3| $.
\end{description}
Finally, we set $h_3 = |{p}_1 - {p}_4|$. We implicitly assume that ${p}_1$ and ${p}_4$ belong to the same half-space. Additionally, note that ${h_1 \approx h_T}$.
\end{Cond}

\begin{note}
As an example, we {define} the matrices $A_{T}$ as 
\begin{align*}
\displaystyle
A_{T} := 
\begin{pmatrix}
\cos \theta  & - \sin \theta \\
 \sin \theta & \cos \theta
\end{pmatrix}, \quad 
{A}_{T} := 
\begin{pmatrix}
 \cos \theta  & - \sin \theta & 0\\
 \sin \theta & \cos \theta & 0 \\
 0 & 0 & 1 \\
\end{pmatrix},
\end{align*}
where $\theta$ denotes the angle. 
\end{note}

\begin{note}
None of the lengths of the edges of a simplex or the measures of the simplex {are} changed by the transformation, i.e.,
\begin{align}
\displaystyle
h_i \leq  h_{T}, \quad i=1,\ldots,d. \label{ineq331}
\end{align}	
\end{note}

\subsection{Additional Notations and Assumptions} \label{addinot}
\begin{figure}[tbhp]
\vspace{-7cm}
  \includegraphics[bb=0 0 944 702,scale=0.55]{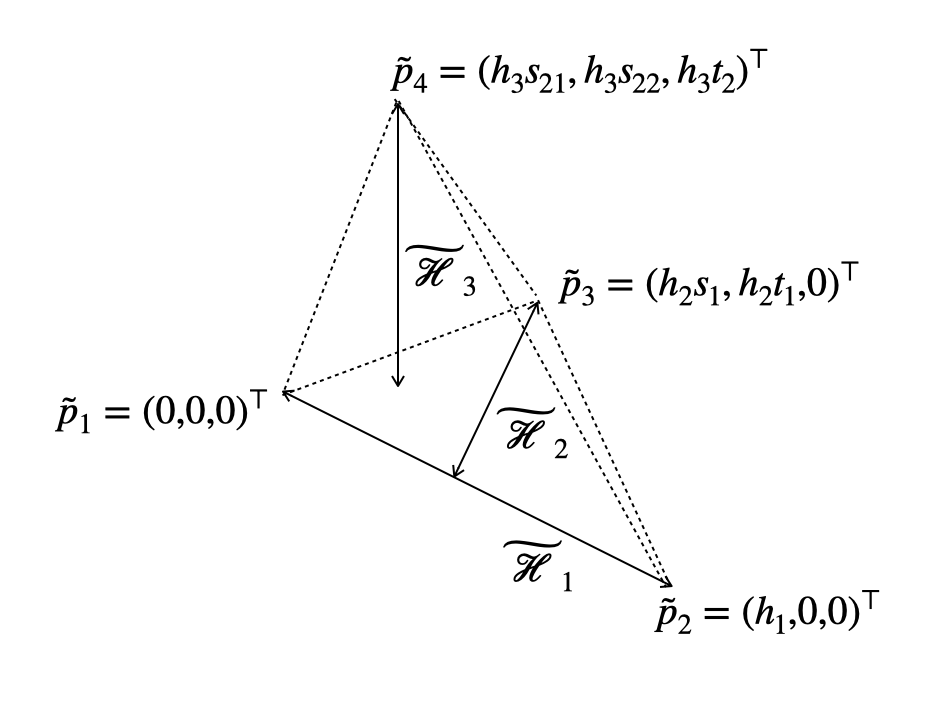}
\caption{New parameters $\widetilde{\mathscr{H}}_i$, $i=1,2,3$}
\label{mathscrH}
\end{figure}

For convenience, we introduce the following additional notation.
We {define} a parameter $\widetilde{\mathscr{H}}_i$, $i=1,\ldots,d$, as
\begin{align}
\displaystyle
\begin{cases} \label{sym=mthscrH}
\widetilde{\mathscr{H}}_1 := h_1, \quad \widetilde{\mathscr{H}}_2 := h_2 t \quad \text{if $d=2$}, \\
\widetilde{\mathscr{H}}_1 := h_1, \quad \widetilde{\mathscr{H}}_2 := h_2 t_1, \quad \widetilde{\mathscr{H}}_3 := h_3 t_2 \quad \text{if $d=3$},
\end{cases}
\end{align}
see Fig. \ref{mathscrH}.

\begin{assume} \label{assume0}
In an anisotropic interpolation error analysis, we impose a geometric condition for the simplex $\widetilde{T}$:
\begin{enumerate}
 \item If $d=2$, there are no additional conditions;
 \item If $d=3$, there exists a positive constant $M$ independent of $h_{\widetilde{T}}$ such that $|s_{22}| \leq M \frac{h_2 t_1}{h_3}$. Note that if $s_{22} \neq 0$, this condition means that the order concerning $h_T$ of $h_3$ coincides with the order of $h_2$, and if $s_{22} = 0$, the order of $h_3$ may be different from that of $h_2$. 
\end{enumerate}
\end{assume}

\begin{figure}[htbp]
  \includegraphics[keepaspectratio, scale=0.45]{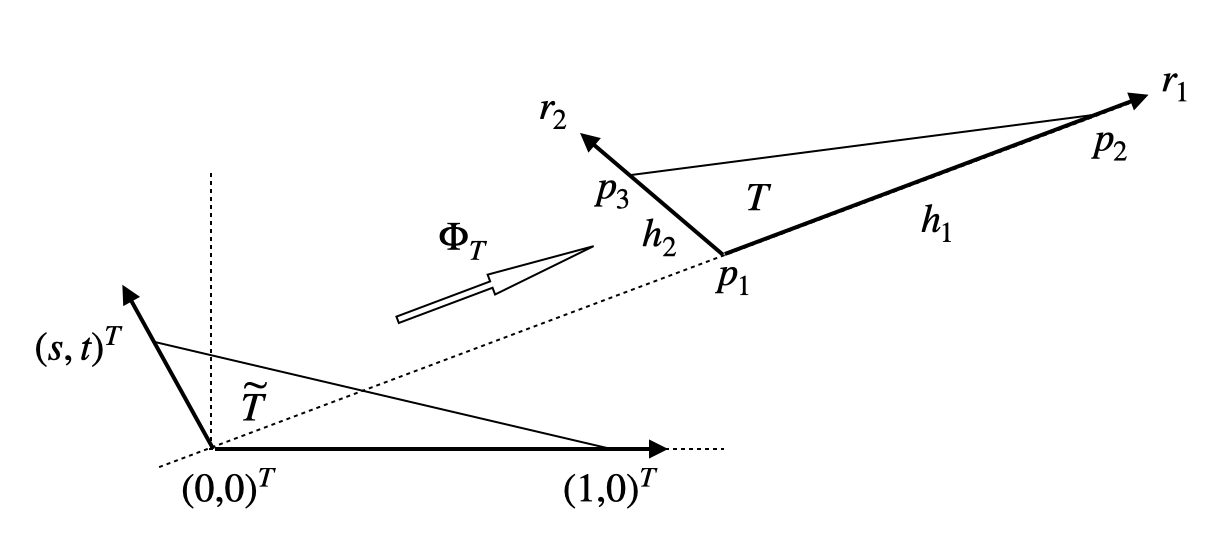}
\caption{Affine mapping $\Phi_{T}$ and vectors $r_i$, $i=1,2$}
\label{affine_2d}
\end{figure}

\begin{figure}[htbp]
  \begin{minipage}[b]{0.4\linewidth}
    \centering
    \includegraphics[keepaspectratio, scale=0.45]{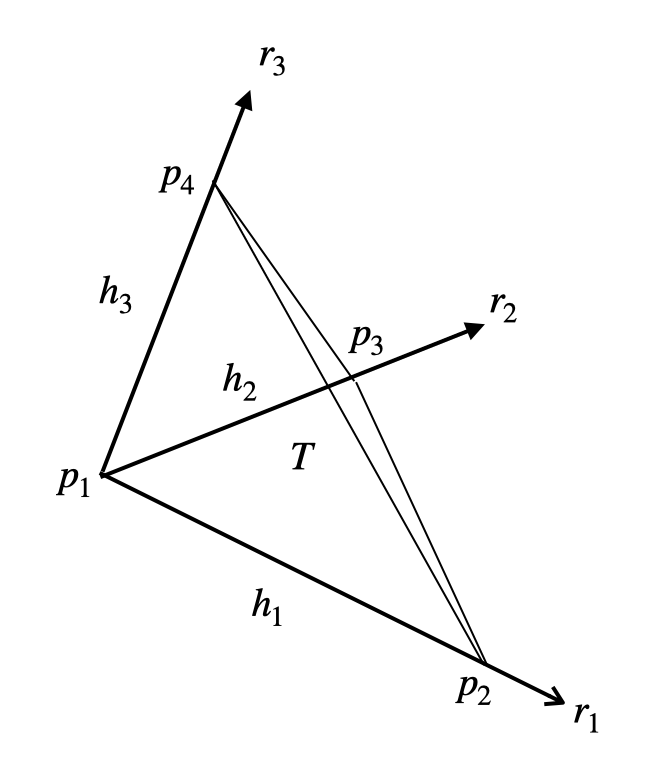}
    \caption{(Type \roman{sone}) Vectors $r_i$, $i=1,2,3$}
     \label{affine_3d_1}
  \end{minipage}
  \begin{minipage}[b]{0.4\linewidth}
    \centering
    \includegraphics[keepaspectratio, scale=0.45]{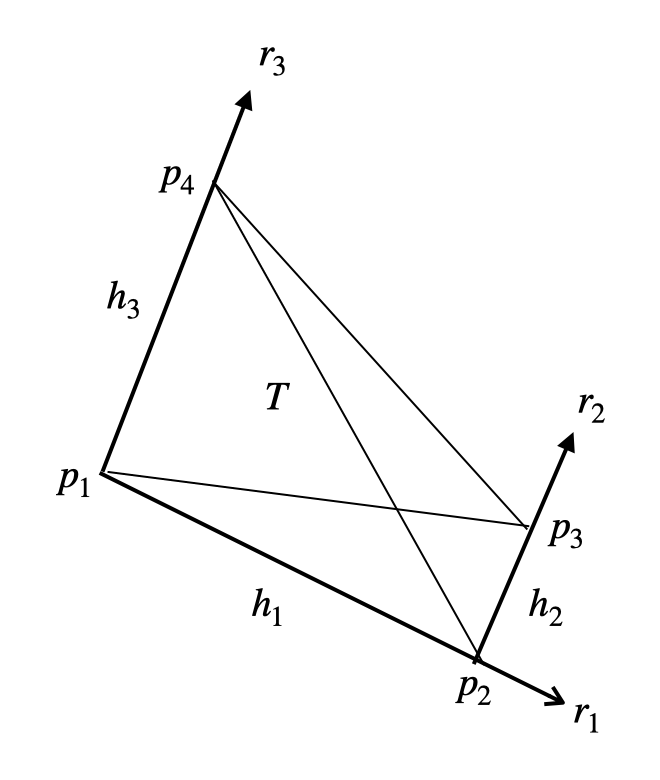}
    \caption{(Type \roman{stwo}) Vectors $r_i$, $i=1,2,3$}
     \label{affine_3d_2}
  \end{minipage}
\end{figure}

We define the vectors ${r}_n \in \mathbb{R}^d$ and $n=1,\ldots,d$ as follows: If $d=2$,
\begin{align*}
\displaystyle
{r}_1 := \frac{p_2 - p_1}{|p_2 - p_1|}, \quad {{r}_2} := \frac{p_3 - p_1}{|p_3 - p_1|},
\end{align*}
see Fig. \ref{affine_2d}, and if $d=3$,
\begin{align*}
\displaystyle
&{r}_1 := \frac{p_2 - p_1}{|p_2 - p_1|}, \quad {r}_3 := \frac{p_4 - p_1}{|p_4 - p_1|}, \quad
\begin{cases}
\displaystyle
{r}_2 := \frac{p_3 - p_1}{|p_3 - p_1|}, \quad \text{for case (\roman{sone})}, \\
\displaystyle
{r}_2 := \frac{p_3 - p_2}{|p_3 - p_2|} \quad \text{for case (\roman{stwo})},
\end{cases}
\end{align*}
see Fig \ref{affine_3d_1} for (Type \roman{sone}) and Fig \ref{affine_3d_2} for (Type \roman{stwo}). Furthermore, we define the vectors $\tilde{r}_n \in \mathbb{R}^d$ and $n=1,\ldots,d$ as follows. If $d=2$,
\begin{align*}
\displaystyle
\tilde{r}_1 := (1 , 0){^{\top}}, \quad \tilde{r}_2 := (s,t){^{\top}},
\end{align*}
and if $d=3$,
\begin{align*}
\displaystyle
&\tilde{r}_1 := (1 , 0,0){^{\top}}, \quad \tilde{r}_3 := ( s_{21}, s_{22} , t_2){^{\top}}, \quad
\begin{cases}
\tilde{r}_2 := ( s_1 ,  t_1 , 0){^{\top}} \quad \text{for case (\roman{sone})}, \\
\tilde{r}_2 := (- s_1,  t_1,0){^{\top}} \quad \text{for case (\roman{stwo})}.
\end{cases}
\end{align*}

\begin{rem}
The vectors $\tilde{r}_i$, $i \in \{ 1,\ldots,d \}$ are unit vectors. Indeed, if $d=2$,
\begin{align*}
\displaystyle
|\tilde{r}_1|_E = 1, \quad |\tilde{r}_2|_E = \sqrt{s^2+t^2} = 1,
\end{align*}
if $d=3$,
\begin{align*}
\displaystyle
|\tilde{r}_1|_E = 1, \quad |\tilde{r}_2|_E = \sqrt{s_1^2+t_1^2} = 1, \quad |\tilde{r}_3|_E = \sqrt{s_{21}^2 + s_{22}^2 + t_2^2} = 1.
\end{align*}
\end{rem}

\begin{rem}
Let ${A}_{T} \in O(d)$ be the orthogonal matrix defined in \eqref{aff=2}. Then,
\begin{align}
\displaystyle
r_i = A_T \tilde{r}_i, \quad i = 1, \ldots, d. \label{vec=ATvec}
\end{align}
\end{rem}

For a sufficiently smooth function $\varphi$ and vector function $v := (v_{1},\ldots,v_{d}){^{\top}}$, we define the directional derivative of $i \in \{ 1, \ldots, d \}$ as:
\begin{align*}
\displaystyle
\frac{\partial \varphi}{\partial {r_i}} &:= ( {r}_i \cdot  \nabla_{x} ) \varphi = \sum_{i_0=1}^d ({r}_i)_{i_0} \frac{\partial \varphi}{\partial x_{i_0}^{}}, \\
\frac{\partial v}{\partial r_i} &:= \left(\frac{\partial v_{1}}{\partial r_i}, \ldots, \frac{\partial v_{d}}{\partial r_i} \right){^{\top}} 
= ( ({r}_i  \cdot \nabla_{x}) v_{1}, \ldots, ({r}_i  \cdot \nabla_{x} ) v_{d} ){^{\top}},
\end{align*}
and for a sufficiently smooth function $\tilde{\varphi}$ and vector function $\tilde{v} := (\tilde{v}_{1},\ldots,\tilde{v}_{d}){^{\top}}$, 
\begin{align*}
\displaystyle
\frac{\partial \tilde{\varphi}}{\partial {\tilde{r}_i}} &:= ( \tilde{r}_i \cdot  \nabla_{\tilde{x}} ) \tilde{\varphi} = \sum_{i_0=1}^d (\tilde{r}_i)_{i_0} \frac{\partial \tilde{\varphi}}{\partial \tilde{x}_{i_0}^{}}, \\
\frac{\partial \tilde{v}}{\partial \tilde{r}_i} &:= \left(\frac{\partial \tilde{v}_{1}}{\partial \tilde{r}_i}, \ldots, \frac{\partial \tilde{v}_{d}}{\partial \tilde{r}_i} \right){^{\top}} 
= ( (\tilde{r}_i  \cdot \nabla_{\tilde{x}}) \tilde{v}_{1}, \ldots, (\tilde{r}_i  \cdot \nabla_{\tilde{x}} ) \tilde{v}_{d} ){^{\top}}.
\end{align*}

For a multi-index $\beta = (\beta_1,\ldots,\beta_d) \in \mathbb{N}_0^d$, we use the following notation.
\begin{align*}
\displaystyle
\partial^{\beta} \varphi := \frac{\partial^{|\beta|} \varphi}{\partial x_1^{\beta_1} \ldots \partial x_d^{\beta_d}}, \quad \partial^{\beta}_{r} \varphi := \frac{\partial^{|\beta|} \varphi}{\partial r_1^{\beta_1} \ldots \partial r_d^{\beta_d}}.
\end{align*}
Note that $\partial^{\beta} \varphi \neq  \partial^{\beta}_{r} \varphi$.

\section{New Semi-regularity Condition}
\subsection{New Geometric Parameter and Condition} \label{new=geeocon=sec}
We proposed a new geometric parameter $H_{T}$ in \cite{IshKobTsu21a}.
 \begin{defi} \label{defi1}
 Parameter $H_T$ is defined as follows:
\begin{align*}
\displaystyle
H_T := \frac{\prod_{i=1}^d h_i}{|T|_d} h_T.
\end{align*}
\end{defi}

We introduce geometric conditions to obtain the optimal convergence rate of the anisotropic error estimates.

\begin{Cond} \label{Cond1=sec6}
A family of meshes $\{ \mathbb{T}_h\}$ is semi-regular if there exists $\gamma_0 \> 0$ such that
\begin{align}
\displaystyle
\frac{H_{T}}{h_{T}} \leq \gamma_0 \quad \forall \mathbb{T}_h \in \{ \mathbb{T}_h \}, \quad \forall T \in \mathbb{T}_h. \label{NewGeo}
\end{align}
\end{Cond}

\begin{rem}
The geometric condition in \eqref{NewGeo}	is equivalent to the maximum angle condition (Section \ref{geo=sec7}).
\end{rem}

\begin{rem} \label{goodel}
We consider the good elements on the meshes in Section \ref{goodbad}. On anisotropic meshes, good elements may satisfy the following conditions:
\begin{description}
  \item[($d=2$)] $h_2 \approx h_2 t$;
  \item[($d=3$)] $h_2 \approx h_2 t_1$ and $h_3 \approx h_3 t_2$.
\end{description}
\end{rem}

\subsection{Properties of the New Geometric Parameter}
We first show the relation between $h_{T}$ and $H_{T}$.

\begin{lem}
It holds that 
\begin{align}
\displaystyle
h_{T} \leq \frac{1}{2} H_{T} \quad \text{if $d=2$}, \label{ineq421=a} \\
h_{T} \< \frac{1}{6} H_{T} \quad \text{if $d=3$}. \label{ineq421=b}
\end{align}
\end{lem}

\begin{pf*}
We consider for each dimension, $d=2,3$.

\textbf{Two-dimensional case.} By constructing the standard element in the two-dimensional case, the angle $\theta_{\max}:= \angle p_2 p_1 p_3$ is the maximum angle of $T$. We then have $\frac{\pi}{3} \< \theta_{\max} \< \pi$, that is, $0 \< \sin \theta_{\max} \leq 1$. Therefore, it holds that
\begin{align*}
\displaystyle
H_{T} = \frac{h_1 h_2}{|T|_2} h_{T} = \frac{2}{ \sin \theta_{\max}} h_{T} \geq 2 h_{T}.
\end{align*}
We here used the fact that $|T|_2 = \frac{1}{2} h_1 h_2 \sin \theta_{\max}$.

\textbf{Three-dimensional case.} We denote by $\phi_{T}$ the angle between the base $\triangle p_1 p_2 p_3$ of $T$ and the segment $\overline{p_1 p_4}$. Recall that there are two types of standard elements, (Type \roman{sone}) or (Type \roman{stwo}). We denote by $\theta_{T}$
\begin{description}
  \item[(Type \roman{sone})]the angle between the segments $\overline{p_1 p_2}$ and $\overline{p_1 p_3}$, that is, $\theta_{T} := \angle p_2 p_1 p_3$, or
  \item[(Type \roman{stwo})] the angle between the segments $\overline{p_2 p_1}$ and $\overline{p_2 p_3}$, that is, $\theta_{T} := \angle p_1 p_2 p_3$.
\end{description}
We set $t_1 := \sin \theta_{T}$ and $t_2 := \sin \phi_{T}$. By constructing the standard element in the three-dimensional case, the angle $ \angle p_1 p_3 p_2$ is the maximum angle of the base $\triangle p_1 p_2 p_3$ of $T$. Therefore, we have $0 \< \theta_{T} \< \frac{\pi}{2}$. Because $0 \< \phi_{T} \< \pi$, it holds that
\begin{align*}
\displaystyle
H_{T} = \frac{h_1 h_2 h_3}{|T|_3} h_{T} = \frac{6}{ \sin \theta_{T}  \sin \phi_{T}} h_{T} \> 6 h_{T}.
\end{align*}
We here used the fact that $|T|_3 = \frac{1}{6} h_1 h_2 h_3 \sin \theta_{T}  \sin \phi_{T}$.
\qed
\end{pf*}

We introduce another geometric parameter regarding Definition \ref{defi1}.

\begin{defi}[Another parameter $H^*_{T}$] \label{defi411}
For $T \in \mathbb{T}_h$,  we denote by $L_i$ edges of the simplex $T$. We define the new parameter $H^*_{T}$ as
 \begin{align}
\displaystyle
H^*_{T} := \frac{h_{T}^2}{|T|_2} \min_{1 \leq i \leq 3} |L_i|  \quad \text{if $d=2$}, \quad H^*_{T} := \frac{h_{T}^2}{|T|_3} \min_{1 \leq i , j \leq 6, i \neq j} |L_i| |L_j| \quad \text{if $d=3$}. \label{NewPar1}
\end{align}
\end{defi}

The parameters $H^*_{T}$ and $H_{T}$ are equivalent.

\begin{lem} \label{lem422}
It holds that
\begin{align}
\displaystyle
\frac{1}{2} H^*_{T} \< H_{T}  \< 2 H^*_{T}. \label{ineq422}
\end{align}
Furthermore, $H^*_{T}$ is equivalent to the circumradius $R_2$ of $T$ in the two-dimensional case. 
\end{lem}

\begin{pf*}
We consider for each dimension, $d=2,3$.

\textbf{Two-dimensional case.} Let  $L_i$ $(i=1,2,3)$ denote edges of the triangle $T$ with $|L_1| \leq |L_2| \leq |L_3|$.  It obviously holds that $h_2 = |L_1|$ and $h_{T} = |L_3| = h_{T}$. Because $h_2 \leq h_1 \< 2 h_{T}$ and $h_{T} \< h_1 + h_2 \leq 2 h_1$ for the triangle $\triangle p_1 p_2 p_3$, it holds that
\begin{align*}
\displaystyle
 \frac{1}{2} h_{T} \< h_1 = |L_2| \< 2 h_{T} = 2 h_{T}. 
\end{align*}
We thus have
\begin{align*}
\displaystyle
\frac{1}{2} H^*_{T} = \frac{1}{2} \frac{|L_1|}{|T|_2} h_{T}^2 \< H_{T} = \frac{h_1 h_2}{|T|_2} h_{T} \< 2 \frac{|L_1|}{|T|_2} h_{T}^2 = 2 H^*_{T}.
\end{align*}
Furthermore, it holds that
\begin{align*}
\displaystyle
2 R_2 = 2 \frac{|L_1| |L_2| |L_3|}{4 |T|_2} \< H^*_{T} = \frac{|L_1|}{|T|_2} h_{T}^2 \< 8 \frac{|L_1| |L_2| |L_3|}{4 |T|_2} = 8 R_2.
\end{align*}

\textbf{Three-dimensional case.} Let  $L_i$ $(i=1,\ldots,6)$ denote edges of the triangle $T$ with $|L_1| \leq |L_2| \leq \cdots \leq |L_6|$.  It obviously holds that $h_2 = |L_1|$ and $h_{T} = |L_6|$. Recall that there are two types of standard elements, (Type \roman{sone}) or (Type \roman{stwo}). 
\begin{description}
  \item[(Type \roman{sone})] We set $h_4 := |p_3 - p_4|$, $h_5 := |p_2 - p_4|$ and $h_6 := |p_2- p_3|$.  Because $h_1 = |E^{(\min)}_{\max}| = |p_1 - p_2|$ is the longest edge among the four edges that share an endpoint with $L_{1}$, it holds that
\begin{align}
\displaystyle
h_2 \leq \min \{ h_3, h_4, h_6 \} \leq  \max \{ h_3, h_4, h_6 \} \leq h_1. \label{alpha346}
\end{align}
Because $p_1$ and $p_4$ belong to the same half-space for the triangle $\triangle p_1 p_2 p_4$,  it holds that
\begin{align*}
\displaystyle
\begin{cases}
h_3 \leq h_5 \leq h_1 = h_{T}  \quad \text{or}\\
h_3 \leq h_1 \leq h_5 = h_{T}.
\end{cases}
\end{align*} 
We thus have
\begin{align*}
\displaystyle
\begin{cases}
h_3 \leq h_5 \leq h_1 = h_{T} \quad \text{or}\\
h_3 \leq h_1 \leq h_{T}  \<  2 h_1, \quad \frac{1}{2} h_{T} \< h_1 \leq h_{T}.
\end{cases}
\end{align*}
Because $h_3 \leq h_5$, the length of the edge $L_2$ is equal to the one of $h_3$, $h_4$ or $h_6$. 

Assume that $|L_2| = h_3$. We then have 
\begin{align*}
\displaystyle
\frac{1}{2}  H^*_{T} = \frac{1}{2} \frac{|L_1| |L_2| }{|T|_3} h_{T}^2  \< H_{T} = \frac{h_1 h_2 h_3}{|T|_3} h_{T} \leq \frac{|L_1| |L_2| }{|T|_3} h_{T}^2 = H^*_{T} (\< 2 H^*_{T}).
\end{align*}

Assume that $|L_2| = h_4$. We consider the triangle $\triangle p_1 p_3 p_4$. From the assumption, we have $h_2 \leq h_4 \leq h_3$ and $\frac{1}{2} h_3 \< h_4 \leq h_3$. We then obtain
\begin{align*}
\displaystyle
\frac{1}{2}  H^*_{T} = \frac{1}{2} \frac{|L_1| |L_2 |}{|T|_3} h_{T}^2  \< H_{T} = \frac{h_1 h_2 h_3}{|T|_3} h_{T^s} \< 2 \frac{|L_1| |L_2| }{|T|_3} h_{T}^2 = 2 H^*_{T}.
\end{align*}

Assume that $|L_2| = h_6$. We consider the triangle $\triangle p_1 p_2 p_3$. Because $p_1$ and $p_3$ belong to the same half-space for the triangle $\triangle p_1 p_2 p_3$, it holds that $h_2 \leq h_6 \leq h_1$ and $\frac{1}{2} h_1 \< h_6 \leq h_1$. From \eqref{alpha346}, we have
\begin{align*}
\displaystyle
\frac{1}{2} h_3 \leq \frac{1}{2} h_1 \< h_6 \leq h_1.
\end{align*}
Because $h_6 \leq h_3$, we then obtain
\begin{align*}
\displaystyle
\frac{1}{2}  H^*_{T} = \frac{1}{2} \frac{|L_1| |L_2 |}{|T|_3} h_{T}^2  \< H_{T} = \frac{h_1 h_2 h_3}{|T|_3} h_{T} \< 2 \frac{|L_1| |L_2|}{|T|_3} h_{T}^2 = 2 H^*_{T}.
\end{align*}

  \item[(Type \roman{stwo})] We set $h_4 := |p_3 - p_4|$, $h_5 := |p_2 - p_4|$, and $h_6 := |p_1- p_3|$.  Because $h_1 = |E^{(\min)}_{\max}| = |p_1 - p_2|$ is the longest edge among the four edges that share an endpoint with $L_{1}$, it holds that
\begin{align}
\displaystyle
h_2 \leq \min \{ h_4, h_5, h_6 \} \leq  \max \{ h_4, h_5, h_6 \} \leq h_1. \label{alpha456}
\end{align}
Because $p_1$ and $p_4$ belong to the same half-space for the triangle $\triangle p_1 p_2 p_4$ and \eqref{alpha456}, it holds that
\begin{align*}
\displaystyle
h_3 \leq h_5 \leq h_1.
\end{align*}
This implies that $h_1 = h_{T}$. Therefore, the length of the edge $L_2$ is equal to the one of $h_3$, $h_4$, or $h_6$. 

Assume that $|L_2| = h_3$. We then have 
\begin{align*}
\displaystyle
\left(\frac{1}{2} H^*_{T} \< \right) H^*_{T} = \frac{|L_1| |L_2 |}{|T|_3} h_{T}^2 = H_{T} &= \frac{h_1 h_2 h_3}{|T|_3} h_{T} \\
&= \frac{|L_1| |L_2| }{|T|_3} h_{T}^2 = H^*_{T} (\< 2 H_{T}).
\end{align*}

Assume that $|L_2| = h_4$. For the triangle $\triangle p_2 p_3 p_4$, we have
\begin{align*}
\displaystyle
h_2 \leq h_4 \leq h_5 \< 2 h_4.
\end{align*}
Because $h_3 \leq h_5$ and $h_4 \leq h_3$, it holds that
\begin{align*}
\displaystyle
\left(\frac{1}{2} H^*_{T} \< \right) H^*_{T} = \frac{|L_1| |L_2 |}{|T|_3} h_{T}^2  \leq H_{T} &= \frac{h_1 h_2 h_3}{|T|_3} h_{T} \\
&\< 2 \frac{|L_1| |L_2| }{|T|_3} h_{T}^2 = 2 H^*_{T}.
\end{align*}

Assume that $|L_2| = h_6$. We have $h_1 \< h_2 + h_6 \< 2 h_6$ for the triangle $\triangle p_1 p_2 p_3$.  Therefore, since $h_6 \leq h_3 \leq h_1$, we obtain
\begin{align*}
\displaystyle
\left(\frac{1}{2} H^*_{T} \< \right) H^*_{T} = \frac{|L_1| |L_2 |}{|T|_3} h_{T}^2  \leq H_{T} &= \frac{h_1 h_2 h_3}{|T|_3} h_{T} \\
&\< 2 \frac{|L_1| |L_2| }{|T|_3} h_{T}^2 = 2 H^*_{T}.
\end{align*}
\end{description}
\qed
\end{pf*}

\subsection{Euclidean Condition Number}
Examining the Euclidean condition number is useful for deriving appropriate interpolation error estimates.

\begin{lem} \label{lem351}
It holds that
\begin{subequations} \label{CN331}
\begin{align}
\displaystyle
\| \widehat{{A}} \|_2 &\leq  h_{T}, \quad \| \widehat{{A}} \|_2 \| \widehat{{A}}^{-1} \|_2 = \frac{\max \{h_1 , \cdots, h_d \}}{\min \{h_1 , \cdots, h_d \}}, \label{CN331a} \\
\| \widetilde{{A}} \|_2 &\leq 
\begin{cases}
\sqrt{2} \quad \text{if $d=2$}, \\
2  \quad \text{if $d=3$},
\end{cases}
\quad \| \widetilde{{A}} \|_2 \| \widetilde{{A}}^{-1} \|_2 \leq
\begin{cases}
\frac{h_1 h_2}{|T|_2} = \frac{H_{T}}{h_{T}} \quad \text{if $d=2$}, \\
\frac{2}{3} \frac{h_1 h_2 h_3}{|T|_3} = \frac{2}{3} \frac{H_{T}}{h_{T}} \quad \text{if $d=3$},
\end{cases}
\label{CN331b} \\
\| {A}_{T} \|_2 &= 1, \quad \| {A}_{T}^{-1} \|_2 = 1. \label{CN331c}
\end{align}
\end{subequations}
where a parameter $H_{T}$ is defined in Definition \ref{defi1}. Furthermore, we have
\begin{align}
\displaystyle
| \det ({A}_{\widetilde{T}}) | = | \det(\widetilde{{A}}) | | \det (\widehat{{A}}) | = \frac{|T|_d}{|\widetilde{T}|_d} \frac{|\widetilde{T}|_d}{|\widehat{T}|_d} = d ! |T|_d, \quad | \det ({A}_{T}) | = 1. \label{CN332}
\end{align}
\end{lem}

\begin{pf*}
We first show the equality \eqref{CN332}. Because
\begin{align*}
\displaystyle
\int_{T} dx = \int_{\widetilde{T}} |\det({A}_T)| d \tilde{x}, \quad  \int_{\widetilde{T}}  d \tilde{x} = \int_{\widehat{T}} | \det ({A}_{\widetilde{T}})  | d \hat{x},
\end{align*}
and $|T|_d = |\widetilde{T}|_d$, we conclude \eqref{CN332}.

We show the equality \eqref{CN331a}. From
\begin{align*}
\displaystyle
(\widehat{{A}})^{\top} \widehat{{A}} = \diag (h_1^2,\ldots,h_d^2), \quad \widehat{{A}}^{-1} \widehat{\mathcal{A}}^{- \top} = \diag (h_1^{-2},\ldots,h_d^{-2}),
\end{align*}
we have
\begin{align*}
\displaystyle
\| \widehat{{A}} \|_2
&= \lambda_{\max} ( \widehat{{A}}^{\top} \widehat{{A}})^{\frac{1}{2}} = \max \{h_1 , \cdots, h_d \} \leq h_{T},
\end{align*}
and
\begin{align*}
\displaystyle
\| \widehat{{A}} \|_2 \| \widehat{{A}}^{-1} \|_2
&=  \lambda_{\max} ( \widehat{{A}}^{\top} \widehat{{A}} )^{\frac{1}{2}}  \lambda_{\max} ( \widehat{{A}}^{-1} \widehat{{A}}^{- \top})^{\frac{1}{2}}
= \frac{\max \{h_1 , \cdots, h_d \}}{\min \{h_1 , \cdots, h_d \}},
\end{align*}
which leads to \eqref{CN331a}.

We next show the equality \eqref{CN331b}. We consider for each dimension, $d=2,3$.

\textbf{Two-dimensional case.} Because
\begin{align*}
\displaystyle
\widetilde{{A}}^{\top} \widetilde{{A}} =
\begin{pmatrix}
1 & s \\
s & 1 \\
\end{pmatrix}, \quad
\widetilde{{A}}^{-1} \widetilde{{A}}^{-\top} = \frac{1}{t^2}
\begin{pmatrix}
1 & -s \\
-s & 1 \\
\end{pmatrix}, \quad
|s| \leq 1,
\end{align*}
we have
\begin{align*}
\displaystyle
 \| \widetilde{{A}} \|_2 
 &= \lambda_{\max}(\widetilde{{A}}^{\top} \widetilde{{A}})^{\frac{1}{2}} \leq (1+|s|)^{\frac{1}{2}} \leq \sqrt{2},
\end{align*}
and
\begin{align*}
\displaystyle
 \| \widetilde{{A}} \|_2 \| \widetilde{{A}}^{-1} \|_2
 &=  \lambda_{\max}(\widetilde{{A}}^{\top} \widetilde{{A}})^{\frac{1}{2}} \lambda_{\max} (\widetilde{{A}}^{-1} \widetilde{{A}}^{- \top} )^{\frac{1}{2}}
 \leq \frac{2}{t}  = \frac{h_1 h_2}{|T|_d},
\end{align*}
which leads to \eqref{CN331b} for $d=2$. Here, we used the fact that $|\widetilde{T}|_d = \frac{1}{2} h_1 h_2 t$ and  $|T|_d = |\widetilde{T}|_d$.

\textbf{Three-dimensional case.} The matrices $\widetilde{{A}}_1$ and $\widetilde{{A}}_2$ introduced in \eqref{aff=5} can be decomposed as $\widetilde{{A}}_1 = \widetilde{{M}}_0 \widetilde{{M}}_1$ and $\widetilde{{A}}_2 = \widetilde{{M}}_0 \widetilde{{M}}_2$ with
\begin{align*}
\displaystyle
\widetilde{{M}}_0 :=
\begin{pmatrix}
1 & 0 & s_{21} \\
0 & 1  & s_{22}\\
0 & 0  & t_2\\
\end{pmatrix}, \
\widetilde{{M}}_1 :=
\begin{pmatrix}
1 &  s_1 & 0 \\
0 & t_1  & 0\\
0 & 0  & 1\\
\end{pmatrix}, \
\widetilde{{M}}_2 :=
\begin{pmatrix}
1 & -s_1 & 0 \\
0 & t_1  & 0\\
0 & 0  & 1\\
\end{pmatrix}.
\end{align*}
\noindent
The eigenvalues of $\widetilde{{M}}_2^{\top} \widetilde{{M}}_2$ coincide with those of $\widetilde{{M}}_1^{\top} \widetilde{{M}}_1$, and only Case  (\roman{sone}) is shown.

We have the inequalities
\begin{align*}
\displaystyle
 \| \widetilde{{A}}_1 \|_2
&= \lambda_{\max}(\widetilde{{A}}_1^{\top} \widetilde{{A}}_1)^{\frac{1}{2}}
\leq \lambda_{\max}(\widetilde{{M}}_0^{\top} \widetilde{{M}}_0)^{\frac{1}{2}} \lambda_{\max}(\widetilde{{M}}_1^{\top} \widetilde{{M}}_1)^{\frac{1}{2}} \notag\\
&\leq \left(1 +  \sqrt{s_{21}^2 + s_{22}^2} \right)^{\frac{1}{2}} (1 + |s_1|)^{\frac{1}{2}} \leq {2},
\end{align*}
and
\begin{align*}
\displaystyle
 \| \widetilde{{A}}_1 \|_2 \| \widetilde{{A}}_1^{-1} \|_2
&= \lambda_{\max}(\widetilde{{A}}_1^{\top} \widetilde{{A}}_1)^{\frac{1}{2}} \lambda_{\max} (\widetilde{{A}}_1^{-1} \widetilde{{A}}_1^{- \top} )^{\frac{1}{2}} \notag\\
&\leq \frac{\left(1 +  \sqrt{s_{21}^2 + s_{22}^2} \right) (1 + |s_1| )}{t_1 t_2} 
\leq \frac{4}{t_1 t_2} = \frac{2}{3} \frac{h_1 h_2 h_3}{|T|_d}, 
\end{align*}
where we used the fact that $|\widetilde{T}|_d = \frac{1}{6} h_1 h_2 h_3 t_1 t_2$ and  $|T|_d = |\widetilde{T}|_d$.

Because the length of all edges of a simplex and measure of the simplex is not changed by a rotation and mirror imaging matrix and ${A}_{T}, {A}_{T}^{-1} \in O(d)$, 
\begin{align*}
\displaystyle
\| {A}_{T} \|_2 = 1, \quad  \| {A}_{T}^{-1} \|_2 =1,
\end{align*}
which is \eqref{CN331c}.
\qed
\end{pf*}

\section{New Geometric Mesh Condition and the Maximum-angle Condition} \label{geo=sec7}

\subsection{Statements}
We state the following theorems concerning the new condition.

\begin{thr} \label{thr=2d}
Condition \ref{Cond1=sec6} holds if and only if Condition \ref{Cond5} holds when $d=2$.
\end{thr}

\begin{pf*}
In the case of $d=2$, we use the previous result presented in \cite{Kri91}; i.e., there exists a constant $\gamma_6  \> 0$ such that
\begin{align*}
\displaystyle
\frac{R_2}{h_{T}} \leq \gamma_6 \quad \forall \mathbb{T}_h \in \{ \mathbb{T}_h \}, \quad \forall T \in \mathbb{T}_h,
\end{align*}
if and only if Condition \ref{Cond5} is satisfied. Combining this result with $H_{T}$ being equivalent to the circumradius $R_2$ of $T$ (Lemma \ref{lem422}), we have the desired conclusion.
\qed	
\end{pf*}

\begin{thr} \label{thr=3d}
Condition \ref{Cond1=sec6} holds if and only if  Condition \ref{Cond7} holds when $d=3$.
\end{thr}

The proof can be found in  \cite{IshKobSuzTsu21}. Preparation is needed to prove the three-dimensional case. \textcolor{red}{The following subsection shows the symbols used only in this section.}

\subsection{Notation}
\begin{figure}[htbp]
\includegraphics[keepaspectratio, scale=0.85]{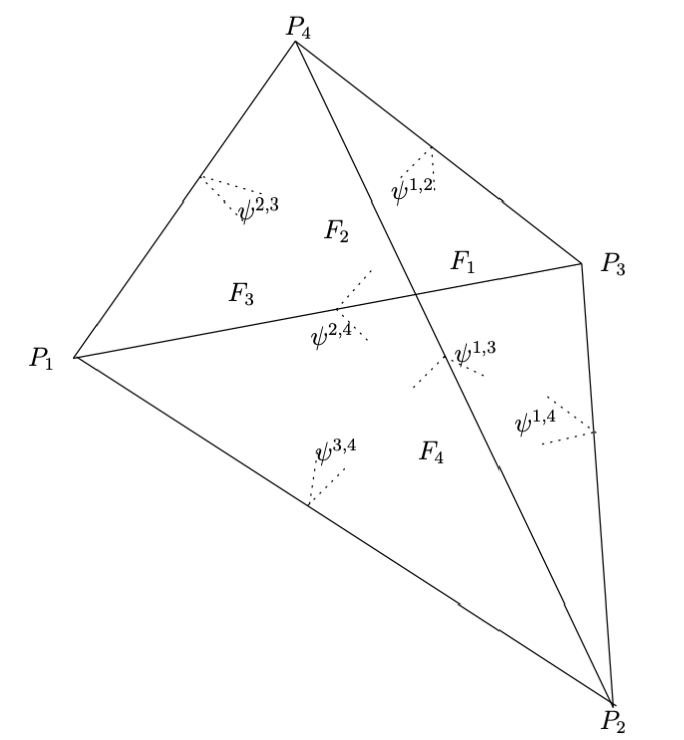}
\caption{Tetrahedra}
\label{tetra1}
\end{figure}

Let $T \in \mathbb{T}_h$ be the standard element in $\mathbb{R}^3$ with vertices, $P_1$, $P_{2}$, $P_{3}$ and $P_{4}$. Let $F_i$ be the face of a simplex $T$ opposite to the vertex $P_i$. We denote by $\psi^{i,j}$ (Table \ref{table=psi}) the angle between the face $F_i$ and the face $F_j$, see Figure \ref{tetra1}. Note that $\psi^{i,j} = \psi^{j,i}$. Furthermore, we denote by $\theta^i_j$ (Table \ref{table=theta}) the internal angle at the vertex $P_j$ on the face $F_i$ and by $\phi^i_j$ (Table \ref{table=phi}) the angle between the face $F_i$ and the segment $\overline{P_j P_i}$. 

\begin{table}[htbp]
\begin{tabular}{c}

\begin{minipage}{0.55\hsize}
\begin{center}
\caption{$\psi^{i,j}$}
\begin{tabular}{c||c|c|c|c} \hline
    & $F_1$ & $F_2$ & $F_3$ & $F_4$ \\ \hline \hline
    $F_1$ &      -      &   $\psi^{1,2}$          &   $\psi^{1,3}$    &  $\psi^{1,4}$     \\ \hline
    $F_2$ &      $\psi^{2,1}$      &      -       &   $\psi^{2,3}$    &  $\psi^{2,4}$    \\ \hline
    $F_3$ &    $\psi^{3,1}$        &    $\psi^{3,2}$         &    -   &   $\psi^{3,4}$    \\ \hline
    $F_4$ &    $\psi^{4,1}$        &      $\psi^{4,2}$       &   $\psi^{4,3}$    &  -     \\ \hline
  \end{tabular} \label{table=psi}
\end{center}
\end{minipage}

\begin{minipage}{0.55\hsize}
\begin{center}
\caption{$\theta^{i}_j$}
\begin{tabular}{c||c|c|c|c} \hline
    & $F_1$ & $F_2$ & $F_3$ & $F_4$ \\ \hline \hline
    $P_1$ &   -        &        $\theta^2_1$      &   $\theta^3_1$     &   $\theta^4_1$    \\ \hline
    $P_2$ &      $\theta^1_2$       &         -    &  $\theta^3_2$      &    $\theta^4_2$  \\ \hline
    $P_3$ &     $\theta^1_3$       &   $\theta^2_3$          &    -   &   $\theta^4_3$   \\ \hline
    $P_4$ &      $\theta^1_4$      &      $\theta^2_4$       &     $\theta^3_4$  &    -   \\ \hline
  \end{tabular} \label{table=theta}
\end{center}
\end{minipage}
\end{tabular}
\end{table}

\begin{table}[htbp]
\begin{center}
\caption{$\phi^{i}_j$}
  \begin{tabular}{c||c|c|c|c} \hline
    & $F_1$ & $F_2$ & $F_3$ & $F_4$ \\ \hline \hline
    $P_1$ &    -        &      $\phi_1^2$       &    $\phi_1^3$   &   $\phi_1^4$   \\ \hline
    $P_2$ &    $\phi_2^1$        &      -       &   $\phi_2^3$    &  $\phi_2^4$    \\ \hline
    $P_3$ &      $\phi_3^1$      &     $\phi_3^2$        &   -    &  $\phi_3^4$    \\ \hline
    $P_4$ &    $\phi_4^1$        &       $\phi_4^2$      &   $\phi_4^3$    &   -   \\ \hline
  \end{tabular}  \label{table=phi}
 \end{center}
\end{table}

\begin{figure}[htbp]
\begin{minipage}{0.5\hsize}
\begin{center}
\includegraphics[keepaspectratio, scale=0.85]{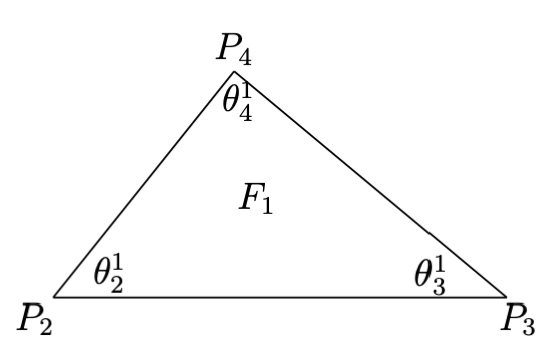}
\end{center}
\caption{Face 1}
\end{minipage}
\begin{minipage}{0.5\hsize}
\begin{center}
\includegraphics[keepaspectratio, scale=0.85]{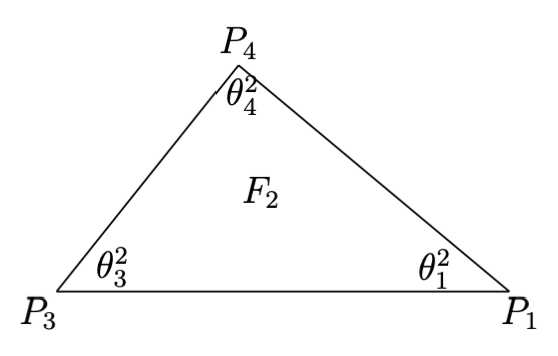}
\end{center}
\caption{Face 2}
\end{minipage} 
\end{figure}

\begin{figure}[htbp]
\begin{minipage}{0.5\hsize}
\begin{center}
\includegraphics[keepaspectratio, scale=0.85]{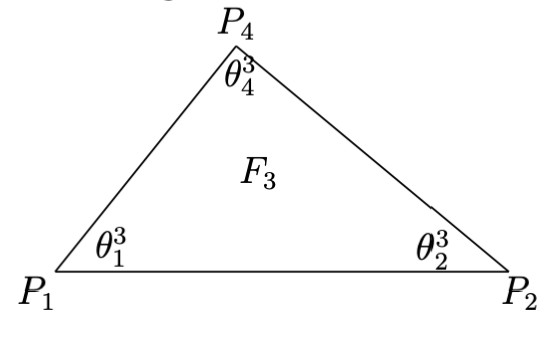}
\end{center}
\caption{Face 3}
\end{minipage}
\begin{minipage}{0.5\hsize}
\begin{center}
\includegraphics[keepaspectratio, scale=0.85]{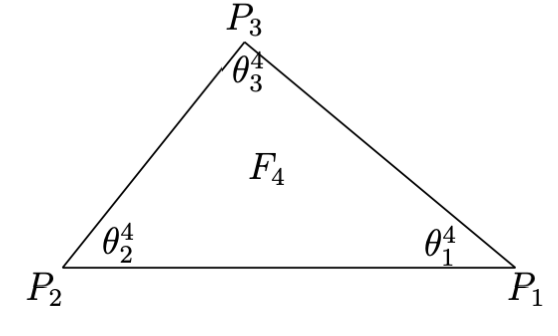}
\end{center}
\caption{Face 4}
\end{minipage}
\end{figure}

\subsection{Preliminaries: Part 1}
We introduce three lemmata.

\begin{lem} \label{lem53}
Let $K \subset \mathbb{R}^2$ be a simplex and let $\theta_1$, $\theta_2$ and $\theta_3$ be internal angles of $K$ with $\theta_1 \leq \theta_2 \leq \theta_3$. If there exists $0 \< \theta_0 \< \pi$, $\theta_0 \in \mathbb{R}$, such that $\theta_{3} \leq \theta_0$, we then have
\begin{align*}
\displaystyle
\sin \theta_{2}, \  \sin \theta_{3} \geq \min \left \{ \sin \frac{\pi - \theta_0}{2}, \sin \theta_0 \right \}.
\end{align*}
\end{lem}

\begin{pf*}
Because $\theta_{1} + \theta_{2} + \theta_{3} = \pi$ and $\theta_{1} \leq \theta_{2} \leq \theta_{3}$, we have
\begin{align*}
\displaystyle
\theta_0 \geq \theta_{3} \geq \theta_{2}  \geq \frac{\theta_{1} + \theta_{2}}{2} \geq \frac{\pi - \theta_{3}}{2} \geq \frac{\pi - \theta_0}{2},
\end{align*}
which leads to the target inequality.
\qed
\end{pf*}

\begin{lem} \label{lem54}
Let $K \subset \mathbb{R}^2$ be a simplex with internal angles $\theta_1$, $\theta_2$ and $\theta_3$. For any fixed $\gamma \in \mathbb{R}$ with $0 \< \gamma \< \pi$, we assume that $\pi - \gamma \leq \theta_i$, $i \in \{1,2,3 \}$. We then have $\theta_{i+1}, \ \theta_{i+2} \leq \gamma$, where the indices $i$, $i + 1$ and $i + 2$ have to be understood ”mod 3".
\end{lem}

\begin{pf*}
Because $\theta_{1} + \theta_{2} + \theta_{3} = \pi$, we have
\begin{align*}
\displaystyle
\theta_{i+1} = \pi - \theta_i - \theta_{i+2} \< \pi - \theta_i \leq \pi - (\pi - \gamma) = \gamma.
\end{align*}
\qed
\end{pf*}

\begin{lem} \label{lem55}
Let $\gamma \in \mathbb{R}$ with $\frac{\pi}{3} \leq \gamma \< \pi$. It then holds that
\begin{align*}
\displaystyle
0 \< \frac{\cos \gamma + 1}{\sin \frac{\gamma}{2} + 1} \leq 1 .
\end{align*}
\end{lem}

\begin{pf*}
Because $\cos \gamma = 1 - 2 \sin^2 \frac{\gamma}{2}$, we have
\begin{align*}
\displaystyle
\frac{\cos \gamma + 1}{\sin \frac{\gamma}{2} + 1}
&= \frac{2 - 2 \sin^2 \frac{\gamma}{2} }{\sin \frac{\gamma}{2} + 1} = 2 \left( 1 - \sin \frac{\gamma}{2} \right).
\end{align*}
Therefore, for $\frac{\pi}{3} \leq \gamma \< \pi$, the target inequality holds.
\qed
\end{pf*}

\subsection{Preliminaries: Part 2}

\begin{lem}[Cosine rules for the sides and for the angles]
It holds that
 \begin{subequations} \label{semi17}
\begin{align}
\displaystyle
\cos \theta^{j+3}_j &= \cos \theta^{j+1}_j \cos \theta^{j+2}_j + \sin \theta^{j+1}_j \sin \theta^{j+2}_j \cos \psi^{{j+1},{j+2}}, \label{semi17a}\\
\cos \theta^{j+1}_j &= \cos \theta^{j+2}_j \cos \theta^{j+3}_j + \sin \theta^{j+2}_j \sin \theta^{j+3}_j \cos \psi^{{j+2},{j+3}},\label{semi17b} \\
\cos \theta^{j+2}_j &= \cos \theta^{j+3}_j \cos \theta^{j+1}_j + \sin \theta^{j+3}_j \sin \theta^{j+1}_j \cos \psi^{{j+3},{j+1}},\label{semi17c} \\
\cos \psi^{{j+1},{j+2}} &= \sin \psi^{{j+2},{j+3}} \sin \psi^{{j+3},{j+1}} \cos \theta^{j+3}_j - \cos \psi^{{j+2},{j+3}} \cos \psi^{{j+3},{j+1}}, \label{semi17d}\\
\cos \psi^{{j+2},{j+3}} &= \sin \psi^{{j+3},{j+1}} \sin \psi^{{j+1},{j+2}} \cos \theta^{j+1}_j - \cos \psi^{{j+3},{j+1}} \cos \psi^{{j+1},{j+2}}, \label{semi17e}\\
\cos \psi^{{j+3},{j+1}} &= \sin \psi^{{j+1},{j+2}} \sin \psi^{{j+2},{j+3}} \cos \theta^{j+2}_j - \cos \psi^{{j+1},{j+2}} \cos \psi^{{j+2},{j+3}}, \label{semi17f}
\end{align}
\end{subequations}
where the indices $j$, $j + 1$, $j + 2$ and $j+3$ have to be understood ”mod 4".
 \end{lem}
 
\begin{pf*}
A proof can be found in \cite{GelCot75,Tod1886}.
\qed
\end{pf*}

\begin{lem} \label{lem56}
Let $\gamma_{\max} \in \mathbb{R}$ with $\frac{\pi}{3} \leq \gamma_{\max} \< \pi$ satisfy Condition \ref{Cond7} for the maximum solid $\theta_{T,\max}$ and the maximum dihedral $\psi_{T,\max}$ of $T$. Assume that for each $j=1,2$, $\theta^4_j$ is not the minimum angle of $\triangle P_1 P_2 P_3$ and $\theta^4_j \< \frac{\pi}{2}$. Then, setting $\delta := \delta(\gamma_{\max})$, $0 \< \delta \leq \frac{\pi}{2}$ such that
\begin{align*}
\displaystyle
\sin \delta = \left( \frac{\cos \gamma_{\max} + 1}{\sin \frac{\gamma_{\max}}{2} + 1} \right)^{1/2},
\end{align*}
it holds that
\begin{align}
\displaystyle
\psi^{j+1,4} \geq \delta, \quad \text{or} \quad \psi^{3,4} \geq \delta, \label{
semi18}
\end{align}
where the indices $j$ and $j + 1$ have to be understood ”mod 2".
\end{lem}

\begin{pf*}
From Lemma \ref{lem55}, we have
\begin{align*}
\displaystyle
0 \< \frac{\cos \gamma_{\max} + 1}{\sin \frac{\gamma_{\max}}{2} + 1} \leq 1,
\end{align*}
because $\frac{\pi}{3} \leq \gamma_{\max} \< \pi$. Therefore, $\delta$ is well-defined.

We use proof by contradiction. Suppose that 
\begin{align*}
\displaystyle
0 \< \psi^{j+1,4} \< \delta, \quad 0 \< \psi^{3,4} \< \delta,
\end{align*}
that is,
\begin{align*}
\displaystyle
0 \< \sin  \psi^{j+1,4} \sin \psi^{3,4} \< \sin^2 \delta, \ \text{and} \ 1 > \cos  \psi^{j+1,4} \cos \psi^{3,4} \> \cos^2 \delta \geq 0. 
\end{align*}
From Lemma \ref{lem53} and assumption, we have
\begin{align*}
\displaystyle
\frac{\pi - \gamma_{\max}}{2} \leq \theta^4_j \< \frac{\pi}{2},
\end{align*}
which implies
\begin{align*}
\displaystyle
0 \< \cos  \theta^4_j \leq \cos \left( \frac{\pi - \gamma_{\max}}{2} \right) = \sin \frac{\gamma_{\max}}{2}.
\end{align*}
We thus obtain
\begin{align*}
\displaystyle
\sin  \psi^{j+1,4} \sin \psi^{3,4} \cos  \theta^4_j \< \sin^2 \delta \sin \frac{\gamma_{\max}}{2}.
\end{align*}
Using the cosine rule \eqref{semi17d} with $j=1$ and the above inequalities yield
\begin{align*}
\displaystyle
\cos \psi_{{2},{3}} &= \sin \psi^{{3},{4}} \sin \psi^{{4},{2}} \cos \theta^{4}_1 - \cos \psi^{{3},{4}} \cos \psi^{{4},{2}} \\
&\< \sin^2 \delta \sin \frac{\gamma_{\max}}{2} - (1 - \sin^2 \delta ) \\
&= \frac{\cos \gamma_{\max} + 1}{\sin \frac{\gamma_{\max}}{2} + 1} \left( \sin \frac{\gamma_{\max}}{2} + 1 \right) - 1 =  \cos \gamma_{\max}.
\end{align*}
This is contradiction for the maximum-angle condition $0 \< \psi^{{2},{3}} \leq \gamma_{\max} < \pi $, that is, $\cos \psi^{{2},{3}} \geq \cos \gamma_{\max}$.

Analogously, using the cosine rule \eqref{semi17f} with $j=2$ and the above inequalities yield
\begin{align*}
\displaystyle
\cos \psi^{{1},{3}} &= \sin \psi^{{3},{4}} \sin \psi^{{4},{1}} \cos \theta^{4}_2 - \cos \psi^{{3},{4}} \cos \psi^{{4},{1}} \\
&<  \cos \gamma_{\max}.
\end{align*}
This is contradiction for the maximum-angle condition $0 \< \psi^{{1},{3}} \leq \gamma_{\max} < \pi $, that is, $\cos \psi^{{1},{3}} \geq \cos \gamma_{\max}$.
\qed
\end{pf*}

\begin{cor} \label{coro57}
For each $j=1,2$, under assumptions in Lemma \ref{lem56}, it holds that setting $C_0 := \min \{ \delta , \gamma_{\max} \}$,
\begin{align*}
\displaystyle
\sin \psi^{j+1,4} \geq C_0, \quad \text{or} \quad \sin \psi^{3,4} \geq C_0 
\end{align*}
where the indices $j$ and $j + 1$ have to be understood ”mod 2".
\end{cor}

\begin{lem} \label{lem58}
For any $i,j \in \{ 1,2,3,4 \}$, $i \neq j$ and $k \in \{ 1,2,3,4 \}$, $k \neq i,j$, it holds that
\begin{align*}
\displaystyle
\sin {\phi_j^i} = \sin \theta_j^k \sin \psi^{k,i}.
\end{align*}
\end{lem}

\begin{pf*}
We only show the case $i=4$, $j=1$ and $k=2$. We then have
\begin{align*}
\displaystyle
\sin {\phi_1^4} = |\overline{P_1 P_4}|  \sin \theta_1^2 \times \frac{1}{|\overline{P_1 P_4}|}  \sin \psi^{2,4} =  \sin \theta_1^2  \sin \psi^{2,4}.
\end{align*}
\qed
\end{pf*}

\begin{lem} \label{lem59}
Assume that there exists a positive constant $M_j$ ($j=1,2$) with $0 \< M_j \< 1$ such that
\begin{align*}
\displaystyle
\sin \theta^4_j \sin \phi^4_1 &\> M_j, \quad j=1,2.
\end{align*}
Setting $\gamma (M_j) := \pi - \sin^{-1} M_j$ ($j=1,2$), we have $\frac{\pi}{2} \< \gamma(M_j) \< \pi$ and it holds that for each $j=1,2$,
\begin{align*}
\displaystyle
&\theta^4_1, \ \theta^4_2, \ \theta^4_3 \< \gamma(M_j),\\
&\theta_{2}^3, \ \theta_4^3, \ \theta_3^2, \ \theta^2_4, \ \theta^2_1, \ \theta^3_1, \ \psi^{2,4}, \ \psi^{3,4} \< \gamma(M_j).
\end{align*}
\end{lem}

\begin{pf*}
From assumption, we have, for each $j=1,2$,
\begin{align*}
\displaystyle
&\sin \theta^4_j \geq \sin \theta^4_j \sin \phi^4_1 \> M_j, \\
&\sin \phi^4_1 \> M_j.
\end{align*}
The definition of $\gamma(M_j)$ and Lemma \ref{lem54} yield, for each $j=1,2$,
\begin{align*}
\displaystyle
&\pi - \gamma \< \theta^4_j \< \gamma(M_j), \quad \theta^4_{j+1} \< \gamma(M_j), \quad \theta^4_{j+2} \< \gamma(M_j), \\
&\pi - \gamma \< \phi^4_1 \< \gamma(M_j),
\end{align*}
where the indices $j$, $j + 1$ and $j + 2$ have to be understood ”mod 3".

We obtain, from Lemma \ref{lem58},
\begin{align*}
\displaystyle
\sin {\phi_1^4} = \sin \theta_1^2 \sin \psi^{2,4} = \sin \theta_1^3 \sin \psi^{3,4} \> M_j, \quad j=1,2.
\end{align*}
We then have, for each $j=1,2$,
\begin{align*}
\displaystyle
\sin \theta_1^2, \ \sin \psi^{2,4}, \ \sin \theta_1^3, \ \sin \psi^{3,4}  \> M_j,
\end{align*}
that is,
\begin{align*}
\displaystyle
\pi - \gamma(M_j) \< \theta^2_1, \ \theta^3_1, \ \psi^{2,4}, \ \psi^{3,4} \< \gamma(M_j).
\end{align*}
On $\triangle P_1 P_2 P_4$ and $\triangle P_1 P_3 P_4$, using Lemma \ref{lem54} yields
\begin{align*}
\displaystyle
\theta_{2}^3, \ \theta_4^3, \ \theta_3^2, \ \theta^2_4 \< \gamma(M_j), \quad j=1,2.
\end{align*}
\qed
\end{pf*}

By analogous argument with Lemma \ref{lem59}, we get the following two lemmata.

\begin{lem} \label{lem510}
Assume that there exists $M_3$ with $0 \< M_3 \< 1$ such that
\begin{align*}
\displaystyle
\sin \theta^1_3 \sin \phi^1_3 &\> M_3.
\end{align*}
Setting $\gamma (M_3) := \pi - \sin^{-1} M_3$, we have $\frac{\pi}{2} \< \gamma(M_3) \< \pi$ and it holds that
\begin{align*}
\displaystyle
\theta^2_3, \ \theta^4_3, \ \theta^1_2, \ \theta^1_4, \ \theta^1_3, \ \psi^{2,1}, \ \psi^{4,1} \< \gamma(M_3).
\end{align*}
\end{lem}

\begin{pf*}
From assumption, we have
\begin{align*}
\displaystyle
&\sin \theta^1_3 \geq \sin \theta^1_3 \sin \phi^1_3 \> M_3, \quad \sin \phi^1_3 \> M_3.
\end{align*}
Using the definition of $\gamma(M_3)$ yields
\begin{align*}
\displaystyle
&\pi - \gamma \< \theta^3_1 \< \gamma(M_3), \quad \pi - \gamma \< \phi^3_1 \< \gamma(M_3).
\end{align*}
We obtain, from Lemma \ref{lem58},
\begin{align*}
\displaystyle
\sin {\phi_3^1} = \sin \theta_3^2 \sin \psi^{2,1} = \sin \theta_3^4 \sin \psi^{4,1} \> M_3.
\end{align*}
We then have
\begin{align*}
\displaystyle
\sin \theta_3^2, \ \sin \psi^{2,1}, \ \sin \theta_3^4, \ \sin \psi^{4,1}  \> M_3,
\end{align*}
that is,
\begin{align*}
\displaystyle
\pi - \gamma(M_3) \< \theta^2_3, \ \theta^4_3, \ \psi^{2,1}, \ \psi^{4,1} \< \gamma(M_3).
\end{align*}

Meanwhile, on $\triangle P_2 P_3 P_4$, using Lemma \ref{lem54}, we have
\begin{align*}
\displaystyle
\theta^1_2, \ \theta^1_4  \< \gamma(M_3).
\end{align*}
\qed
\end{pf*}

\begin{lem} \label{lem511}
Assume that there exists $M_4$ with $0 \< M_4 \< 1$ such that
\begin{align*}
\displaystyle
\sin \theta^1_2 \sin \phi^1_4 &\> M_4.
\end{align*}
Setting $\gamma (M_4) := \pi - \sin^{-1} M_4$, we have $\frac{\pi}{2} \< \gamma(M_4) \< \pi$ and it holds that
\begin{align*}
\displaystyle
\theta^2_4, \ \theta^3_4, \ \theta^1_2, \ \theta^1_3, \ \theta^1_4, \ \psi^{1,2}, \ \psi^{1,3} \< \gamma(M_4).
\end{align*}
\end{lem}

\begin{pf*}
The proof is obtained by using an analogous argument with Lemma \ref{lem510}.
\qed
\end{pf*}

\subsection{Proof of Theorem \ref{thr=3d} in (Type \roman{sone})}

\subsubsection{Condition \ref{Cond7} $\Rightarrow$ Condition \ref{Cond1=sec6}}
We set $t_1 := \sin \theta_1^4$ and $t_2 := \sin \phi_1^4$. We then have
 \begin{align*}
\displaystyle
\frac{H_{T}}{h_{T}} = \frac{h_1 h_2 h_3}{|T|_3} = \frac{6}{\sin \theta_1^4  \sin \phi_1^4}.
\end{align*}
We here used the fact that $|T|_3 = \frac{1}{6} h_1 h_2 h_3 \sin \theta_1^4  \sin \phi_1^4$. By construct of the standard element (Type \roman{sone}), the angle $\theta_3^4$ and $\theta_2^4$ are respectively the maximum angle and the minimum angle of the base $\triangle P_1 P_2 P_3$ of $T$. We hence have $\theta_1^4 \< \frac{\pi}{2}$. From Lemma \ref{lem53}, we have
\begin{align*}
\displaystyle
\frac{\pi - \gamma_{11}}{2} \leq \theta_1^4 \leq \gamma_{11}, \quad  \sin \theta_1^4 \geq \min \left \{ \sin \frac{\pi - \gamma_{11}}{2}, \sin \gamma_{11} \right \} =: C_1.
\end{align*}
Due to Lemma \ref{lem56}, setting $\delta := \delta(\gamma_{11})$, $0 \< \delta \leq \frac{\pi}{2}$ such that
\begin{align*}
\displaystyle
\sin \delta = \left( \frac{\cos \gamma_{11} + 1}{\sin \frac{\gamma_{11}}{2} + 1} \right)^{1/2},
\end{align*}
it holds that
\begin{align*}
\displaystyle
\psi^{2,4} \geq \delta, \quad \text{or} \quad \psi^{3,4} \geq \delta.
\end{align*}
Suppose that $\psi^{2,4} \geq \delta$. By Corollary \ref{coro57} and Lemma \ref{lem58}, we have 
\begin{align*}
\displaystyle
\sin \phi_1^4 = \sin \theta^2_1 \sin \psi^{2,4} \geq C_0 \sin \theta^2_1.
\end{align*}
By construct of the standard element (Type \roman{sone}), the angle $\theta_1^2$ is not the minimum angle of $\triangle P_1 P_3 P_4$. From Lemma \ref{lem53}, we have
\begin{align*}
\displaystyle
\frac{\pi - \gamma_{11}}{2} \leq \theta_1^2 \leq \gamma_{11}, \quad  \sin \theta_1^2 \geq  C_1.
\end{align*}
We thus obtain
\begin{align*}
\displaystyle
\sin \phi_1^4 \geq C_0 C_1.
\end{align*}
Suppose that $\psi^{3,4} \geq \delta$. By Corollary \ref{coro57} and Lemma \ref{lem58}, we have 
\begin{align*}
\displaystyle
\sin \phi_1^4 = \sin \theta^3_1 \sin \psi^{3,4} \geq C_0 \sin \theta^3_1.
\end{align*}
By construct of the standard element (Type \roman{sone}), the angle $\theta_1^3$ is not the minimum angle of $\triangle P_1 P_2 P_4$. From Lemma \ref{lem53}, we have
\begin{align*}
\displaystyle
\frac{\pi - \gamma_{11}}{2} \leq \theta_1^3 \leq \gamma_{11}, \quad  \sin \theta_1^3 \geq  C_1.
\end{align*}
We thus obtain
\begin{align*}
\displaystyle
\sin \phi_1^4 \geq C_0 C_1.
\end{align*}
In both cases
\begin{align*}
\displaystyle
\psi^{2,4} \geq \delta, \quad \text{or} \quad \psi^{3,4} \geq \delta,
\end{align*}
gathering the above results yield
\begin{align*}
\displaystyle
\frac{H_{T}}{h_{T}} = \frac{6}{\sin \theta_1^4  \sin \phi_1^4} \leq \frac{6}{C_0 C_1^2} =: D_1 \> 0,
\end{align*}
that is, Condition \ref{Cond1=sec6} holds.
\qed

\subsubsection{Condition \ref{Cond1=sec6} $\Rightarrow$ Condition \ref{Cond7}}
From assumption, it holds that
\begin{align*}
\displaystyle
\frac{H_{T}}{h_{T}} = \frac{h_1 h_2 h_3}{|T|_3} = \frac{6}{\sin \theta_1^4  \sin \phi_1^4}  \leq \gamma_9.
\end{align*}
Remark that $\frac{6}{\gamma_9} \< 1$ because $\theta^4_1 \< \frac{\pi}{2}$ and $\sin \theta_1^4  \sin \phi_1^4 \< 1$. Therefore, we have
\begin{align*}
\displaystyle
\sin \theta_1^4  \sin \phi_1^4 \geq \frac{6}{\gamma_9} =: C_2.
\end{align*}
From Lemma \ref{lem59} with $j=1$, setting $\gamma (C_2) := \pi - \sin^{-1} C_2$, we have $\frac{\pi}{2} \< \gamma(C_2) \< \pi$ and it holds that
\begin{align*}
\displaystyle
&\theta^4_1, \ \theta^4_2, \ \theta^4_3 \< \gamma(C_2),\\
&\theta_{2}^3, \ \theta_4^3, \ \theta^3_1, \ \theta_3^2, \ \theta^2_4, \ \theta^2_1, \ \psi^{2,4}, \ \psi^{3,4} \< \gamma(C_2).
\end{align*}
Furthermore, we write
\begin{align*}
\displaystyle
|T|_3 &= \frac{1}{3} \times \frac{1}{2} |\overline{P_2 P_3}| |\overline{P_3 P_4}| \sin \theta^1_3 \times h_2 \sin \phi_3^1
= \frac{1}{6} h_2  |\overline{P_2 P_3}| |\overline{P_3 P_4}|  \sin \theta^1_3 \sin \phi_3^1 \\
&\< \frac{1}{3} h_1 h_2 h_3 \sin \theta^1_3 \sin \phi_3^1,
\end{align*}
where we used the fact that $|\overline{P_3 P_4}| \< |\overline{P_1 P_4}| + |\overline{P_1 P_3}| \leq 2 h_3$ on $\triangle P_1 P_3 P_4$ and $|\overline{P_2 P_3}| \leq h_1$. We thus have
\begin{align*}
\displaystyle
\gamma_9 \geq \frac{H_{T}}{h_{T}} \> \frac{3}{\sin \theta^1_3 \sin \phi_3^1},
\end{align*}
that is,
\begin{align*}
\displaystyle
\sin \theta^1_3 \sin \phi_3^1 \> \frac{3}{\gamma_9} =: C_3.
\end{align*}
From Lemma \ref{lem510}, setting $\gamma (C_3) := \pi - \sin^{-1} C_3$, we have $\frac{\pi}{2} \< \gamma(C_3) \< \pi$ and it holds that
\begin{align*}
\displaystyle
\theta^1_2, \ \theta^1_4, \ \theta^1_3, \ \psi^{2,1}, \ \psi^{4,1} \< \gamma(C_3).
\end{align*}
Due to the cosine rule \eqref{semi17f} with $j=2$, we get
\begin{align*}
\displaystyle
\cos \psi^{{1},{3}} &= \sin \psi^{{3},{4}} \sin \psi^{{4},{1}} \cos \theta^{4}_2 - \cos \psi^{{3},{4}} \cos \psi^{{4},{1}}.
\end{align*}
By constructing the standard element (Type \roman{sone}), the angle $\theta_2^4$ is the minimum angle of $\triangle P_1 P_2 P_3$. Therefore, we have
\begin{align*}
\displaystyle
&\cos \theta_2^4 \geq \frac{1}{2} \quad \text{because $\theta_2^4 \leq \frac{\pi}{3}$},\\
&\sin \psi^{{3},{4}} \sin \psi^{{4},{1}} \cos \theta^{4}_2 \> 0, \quad \text{because $\sin \psi^{{3},{4}} \sin \psi^{{4},{1}}  \> 0$},
\end{align*}
and thus
\begin{align*}
\displaystyle
\cos \psi^{{1},{3}} \> - \cos \psi^{{3},{4}} \cos \psi^{{4},{1}}.
\end{align*}
Using $\sin \psi^{{3},{4}} \> C_2$ and $\sin \psi^{{4},{1}} \> C_3$ yields
\begin{align*}
\displaystyle
\cos \psi^{{1},{3}} &\> - \cos \psi^{{3},{4}} \cos \psi^{{4},{1}} \\
&\geq - | \cos \psi^{{3},{4}} | | \cos \psi^{{4},{1}} | = - \sqrt{1 - \sin^2 \psi^{{3},{4}}} \sqrt{1 - \sin^2 \psi^{{4},{1}}} \\
&\> - \sqrt{1 - C_2^2} \sqrt{1 - C_3^2}=: C_4 \> -1.
\end{align*}
Setting $\gamma(C_4) := \cos^{-1} C_4$, it holds that
\begin{align*}
\displaystyle
\psi^{{1},{3}} \< \gamma(C_4) \< \pi.
\end{align*}
Due to the cosine rule \eqref{semi17d} with $j=1$, we get
\begin{align*}
\displaystyle
\cos \psi^{{2},{3}} &= \sin \psi^{{3},{4}} \sin \psi^{{4},{2}} \cos \theta^{4}_1 - \cos \psi^{{3},{4}} \cos \psi^{{4},{2}}.
\end{align*}
By constructing the standard element (Type \roman{sone}), the angle $\theta_3^4$ and $\theta_2^4$ are respectively the maximum angle and the minimum angle of the base $\triangle P_1 P_2 P_3$ of $T^s$. We hence have $\theta_1^4 \< \frac{\pi}{2}$. Therefore, we have
\begin{align*}
\displaystyle
&\cos \theta_1^4 \> 0 \quad \text{because $\theta_1^4 \leq \frac{\pi}{2}$},\\
&\sin \psi^{{3},{4}} \sin \psi^{{4},{2}} \cos \theta^{4}_1 \> 0, \quad \text{because $\sin \psi^{{3},{4}} \sin \psi^{{4},{2}}  \> 0$},
\end{align*}
and thus
\begin{align*}
\displaystyle
\cos \psi^{{2},{3}} \> - \cos \psi^{{3},{4}} \cos \psi^{{4},{2}}.
\end{align*}
Using $\sin \psi^{{3},{4}} \> C_2$ and $\sin \psi^{{4},{2}} \> C_2$ yield
\begin{align*}
\displaystyle
\cos \psi^{{2},{3}} &\> - \cos \psi^{{3},{4}} \cos \psi^{{4},{2}} \\
&\geq - | \cos \psi^{{3},{4}} | | \cos \psi^{{4},{2}} | = - \sqrt{1 - \sin^2 \psi^{{3},{4}}} \sqrt{1 - \sin^2 \psi^{{4},{2}}} \\
&\> - (1 - C_2^2)=: C_5 \> -1.
\end{align*}
Setting $\gamma(C_5) := \cos^{-1} C_5$, it holds that
\begin{align*}
\displaystyle
\psi^{{2},{3}} \< \gamma(C_5) \< \pi.
\end{align*}
We set $\gamma_{\max} := \max \{ \gamma(C_3), \gamma(C_4) , \gamma(C_5) \}$. We then have $0 \< \gamma_{\max} \< \pi$, that is, Condition \ref{Cond7} holds.
\qed

\subsection{Proof of Theorem \ref{thr=3d} in (Type \roman{stwo})}

\subsubsection{Condition \ref{Cond7} $\Rightarrow$ Condition \ref{Cond1=sec6}}
We set $t_1 := \sin \theta_2^4$ and $t_2 := \sin \phi_1^4$. We then have
 \begin{align*}
\displaystyle
\frac{H_{T}}{h_{T}} = \frac{h_1 h_2 h_3}{|T|_3} = \frac{6}{\sin \theta_2^4  \sin \phi_1^4}.
\end{align*}
We here used the fact that $|T|_3 = \frac{1}{6} h_1 h_2 h_3 \sin \theta_2^4  \sin \phi_1^4$. By construct of the standard element (Type \roman{stwo}), the angle $\theta_3^4$ and $\theta_1^4$ are respectively the maximum angle and the minimum angle of the base $\triangle P_1 P_2 P_3$ of $T^s$. We hence have $\theta_2^4 \< \frac{\pi}{2}$. From Lemma \ref{lem53}, we have
\begin{align*}
\displaystyle
\frac{\pi - \gamma_{11}}{2} \leq \theta_2^4 \leq \gamma_{11}, \quad  \sin \theta_2^4 \geq C_1.
\end{align*}
Due to Lemma \ref{lem56}, it holds that
\begin{align*}
\displaystyle
\psi^{1,4} \geq \delta, \quad \text{or} \quad \psi^{3,4} \geq \delta.
\end{align*}
Suppose that $\psi^{1,4} \geq \delta$. By Corollary \ref{coro57} and Lemma \ref{lem58}, we have 
\begin{align*}
\displaystyle
\sin \phi_2^4 = \sin \theta^1_2 \sin \psi^{1,4} \geq C_0 \sin \theta^1_2.
\end{align*}
Furthermore, it holds that
\begin{align*}
\displaystyle
\sin \phi_1^4 = \frac{|\overline{P_2 P_4}| \sin \phi_2^4}{h_3}.
\end{align*}
By construct of the standard element (Type \roman{stwo}), the angle $\theta_2^1$ is not the minimum angle of $\triangle P_2 P_3 P_4$. From Lemma \ref{lem53}, we have
\begin{align*}
\displaystyle
\frac{\pi - \gamma_{11}}{2} \leq \theta_2^1 \leq \gamma_{11}, \quad  \sin \theta_2^1 \geq  C_1.
\end{align*}
Because $h_3 = |\overline{P_1 P_4}| \< |\overline{P_2 P_4}|$ on $\triangle P_1 P_2 P_4$, we thus obtain
\begin{align*}
\displaystyle
\sin \phi_1^4 &= \frac{|\overline{P_2 P_4}|}{h_3} \sin \phi_2^4
\> C_0 C_1.
\end{align*}
Suppose that $\psi^{3,4} \geq \delta$. By Corollary \ref{coro57} and Lemma \ref{lem58}, we have 
\begin{align*}
\displaystyle
\sin \phi_1^4 = \sin \theta^3_1 \sin \psi^{3,4} \geq C_0 \sin \theta^3_1.
\end{align*}
By constructing the standard element (Type \roman{stwo}), the angle $\theta_1^3$ is not the minimum angle of $\triangle P_1 P_2 P_4$. From Lemma \ref{lem53}, we have
\begin{align*}
\displaystyle
\frac{\pi - \gamma_{11}}{2} \leq \theta_1^3 \leq \gamma_{11}, \quad  \sin \theta_1^3 \geq  C_1.
\end{align*}
We thus obtain
\begin{align*}
\displaystyle
\sin \phi_1^4\>  C_0 C_1.
\end{align*}
In both cases
\begin{align*}
\displaystyle
\psi^{1,4} \geq \delta, \quad \text{or} \quad \psi^{3,4} \geq \delta,
\end{align*}
gathering the above results yields
\begin{align*}
\displaystyle
\frac{H_{T}}{h_{T}} = \frac{6}{\sin \theta_2^4  \sin \phi_1^4} \leq \frac{6}{C_0 C_1^2} = D_1 \> 0,
\end{align*}
that is, Condition \ref{Cond1=sec6} holds.
\qed

\subsubsection{Condition \ref{Cond1=sec6} $\Rightarrow$ Condition \ref{Cond7}}
From assumption, it holds that
\begin{align*}
\displaystyle
\frac{H_{T}}{h_{T}} = \frac{h_1 h_2 h_3}{|T|_3} = \frac{6}{\sin \theta_2^4  \sin \phi_1^4}  \leq \gamma_9.
\end{align*}
Remark that $\frac{6}{\gamma_9} \< 1$ because $\theta^4_2 \< \frac{\pi}{2}$ and $\sin \theta_2^4  \sin \phi_1^4 \< 1$. Therefore, we have
\begin{align*}
\displaystyle
\sin \theta_2^4  \sin \phi_1^4 \geq \frac{6}{\gamma_9} = C_2.
\end{align*}
From Lemma \ref{lem59} with $j=2$, it holds that
\begin{align*}
\displaystyle
&\theta^4_1, \ \theta^4_2, \ \theta^4_3 \< \gamma(C_2),\\
&\theta_{2}^3, \ \theta_4^3, \ \theta^3_1, \ \theta_3^2, \ \theta^2_4, \ \theta^2_1, \ \psi^{2,4}, \ \psi^{3,4} \< \gamma(C_2).
\end{align*}
Furthermore, we write
\begin{align*}
\displaystyle
|T|_3 &= \frac{1}{3} \times \frac{1}{2} |\overline{P_2 P_4}| |\overline{P_2 P_3}| \sin \theta^1_2 \times h_3 \sin \phi_4^1 \\
&\< \frac{1}{6} h_1 h_2 h_3 \sin \theta^1_2 \sin \phi_4^1,
\end{align*}
where we used the fact that $|\overline{P_3 P_2}| = h_2$ and $|\overline{P_2 P_4}| \leq h_1$. We thus have
\begin{align*}
\displaystyle
\gamma_9 \geq \frac{H_{T^s}}{h_{T^s}} \> \frac{6}{\sin \theta^1_2 \sin \phi_4^1},
\end{align*}
that is,
\begin{align*}
\displaystyle
\sin \theta^1_2 \sin \phi_4^1 \> \frac{6}{\gamma_9} = C_2.
\end{align*}
From Lemma \ref{lem511}, it holds that
\begin{align*}
\displaystyle
\theta^1_2, \ \theta^1_4, \ \theta^1_3, \ \psi^{1,2}, \ \psi^{1,3} \< \gamma(C_2).
\end{align*}
Due to the cosine rule \eqref{semi17e} with $j=2$, we get
\begin{align*}
\displaystyle
\cos \psi^{{4},{1}} &= \sin \psi^{{1},{3}} \sin \psi^{{3},{4}} \cos \theta^{3}_2 - \cos \psi^{{1},{3}} \cos \psi^{{3},{4}}.
\end{align*}
By constructing the standard element (Type \roman{stwo}), the angle $\theta_2^3$ is the minimum angle of $\triangle P_1 P_2 P_4$. Therefore, we have
\begin{align*}
\displaystyle
&\cos \theta_2^3 \geq \frac{1}{2} \quad \text{because $\theta_2^3 \leq \frac{\pi}{3}$},\\
&\sin \psi^{{1},{3}} \sin \psi^{{3},{4}} \cos \theta^{3}_2 \> 0, \quad \text{because $\sin \psi^{{1},{3}} \sin \psi^{{3},{4}}  \> 0$},
\end{align*}
and thus
\begin{align*}
\displaystyle
\cos \psi^{{4},{1}} \> - \cos \psi^{{1},{3}} \cos \psi^{{3},{4}}.
\end{align*}
Using $\sin \psi^{{1},{3}} \> C_2$ and $\sin \psi^{{3},{4}} \> C_2$ yield
\begin{align*}
\displaystyle
\cos \psi^{{4},{1}} &\> - \cos \psi^{{1},{3}} \cos \psi^{{3},{4}} \\
&\geq - \sqrt{1 - \sin^2 \psi^{{1},{3}}} \sqrt{1 - \sin^2 \psi^{{3},{4}}} \\
&\> - (1 - C_2^2)= C_5 \> -1.
\end{align*}
It then holds that
\begin{align*}
\displaystyle
\psi^{{4},{1}} \< \gamma(C_5) \< \pi.
\end{align*}
Due to the cosine rule \eqref{semi17d} with $j=1$, we get
\begin{align*}
\displaystyle
\cos \psi^{{2},{3}} &= \sin \psi^{{3},{4}} \sin \psi^{{4},{2}} \cos \theta^{4}_1 - \cos \psi^{{3},{4}} \cos \psi^{{4},{2}}.
\end{align*}
By constructing the standard element (Type \roman{stwo}), the angle $\theta_1^4$ is the minimum angle of $\triangle P_1 P_2 P_3$. We hence have $\theta_1^4 \< \frac{\pi}{3}$. Therefore, we have
\begin{align*}
\displaystyle
&\cos \theta_1^4 \geq \frac{1}{2} \quad \text{because $\theta_1^4 \leq \frac{\pi}{3}$},\\
&\sin \psi^{{3},{4}} \sin \psi^{{4},{2}} \cos \theta^{4}_1 \> 0, \quad \text{because $\sin \psi^{{3},{4}} \sin \psi^{{4},{2}}  \> 0$},
\end{align*}
and thus
\begin{align*}
\displaystyle
\cos \psi^{{2},{3}} \> - \cos \psi^{{3},{4}} \cos \psi^{{4},{2}}.
\end{align*}
Using $\sin \psi^{{3},{4}} \> C_2$ and $\sin \psi^{{4},{2}} \> C_2$ yield
\begin{align*}
\displaystyle
\cos \psi^{{2},{3}} &\> - \cos \psi^{{3},{4}} \cos \psi^{{4},{2}} \\
&\geq - \sqrt{1 - \sin^2 \psi^{{3},{4}}} \sqrt{1 - \sin^2 \psi^{{4},{2}}} \\
&\> - (1 - C_2^2)= C_5 \> -1.
\end{align*}
It then holds that
\begin{align*}
\displaystyle
\psi^{{2},{3}} \< \gamma(C_5) \< \pi.
\end{align*}
We set $\gamma_{\max} := \max \{ \gamma(C_2), \gamma(C_5) \}$. We then have $0 \< \gamma_{\max} \< \pi$, that is, Condition \ref{Cond7} holds.
\qed

\section{Good Elements or not for $d=2,3$?} \label{goodbad}
In this subsection, we consider good elements on meshes. In this paper, we define 'good elements' on meshes as the existence of a positive constant $\gamma_0 > 0$ satisfying \eqref{NewGeo}. We treat a "Right-angled triangle", "Blade" and "Dagger" for $d=2$, and "Spire", "Spear", "Spindle", "Spike", "Splinter" and "Sliver" for $d=3$, which are introduced in \cite{Cheetal00}. We give the quantities $h_{\max} / h_{\min}$ and $H_{T}/h_{T}$ for those elements. The parameters $h_{\max}$ and $h_{\min}$ are defined as
\begin{align}
\displaystyle
h_{\max} := \max \{ h_1 , \ldots, h_d \}, \quad h_{\min} := \min \{ h_1 , \ldots, h_d \}. \label{hmin=hmax}
\end{align}

\subsection{Isotropic Mesh Elements} \label{iso=mesh=ex2d}
Recall that an isotropic mesh element has equal or nearly equal edge lengths and angles, resulting in a balanced shape. Then, the geometric condition \eqref{geo3} is satisfied. Therefore, it holds that
\begin{align*}
\displaystyle
\frac{H_{T}}{h_{T}} \leq \frac{h_{T}^d}{|T|_d} \leq \frac{1}{\gamma_3}, \quad \frac{h_{\max}}{h_{\min}} \leq  c \frac{h_{T}^d}{|T|_d} \leq \frac{c}{\gamma_3}.
\end{align*}
In this case, elements satisfying the geometric condition \eqref{geo3} are "good."


\subsection{Anisotropic mesh: two-dimensional case} \label{ani=mesh=ex2d}
Let $S \subset \mathbb{R}^2$ be a triangle. Let $0 \< s \ll 1$, $s \in \mathbb{R}$ and $\varepsilon,\delta,\gamma \in \mathbb{R}$.

\begin{ex}[Right-angled triangle] \label{Note2}
Let $S \subset \mathbb{R}^2$ be the simplex with vertices $p_1 := (0,0)^{\top}$, $p_2 := (s,0)^{\top}$ and $p_3 := (0,s^{\varepsilon})^{\top}$ with $1 \< \varepsilon$. We then have $h_1 = s$ and $h_2 =  s^{\varepsilon}$; i.e.,
\begin{align*}
\displaystyle
\frac{h_{\max}}{h_{\min}} &\leq s^{1 - \varepsilon} \to \infty \quad \text{as $s \to 0$}, \quad \frac{H_{S}}{h_{S}} =2.
\end{align*}
In this case, the element $S$ is "good."
\end{ex}

\begin{ex}[Dagger] \label{Note1}
Let $S \subset \mathbb{R}^2$ be the simplex with vertices $p_1 := (0,0)^{\top}$, $p_2 := (s,0)^{\top}$ and $p_3 := (s^{\delta},s^{\varepsilon})^{\top}$ with $1 \< \varepsilon \< \delta $. We then have $h_1 = \sqrt{(s - s^{\delta})^2 + s^{2 \varepsilon}}$ and $h_2 = \sqrt{s^{2 \delta} + s^{2 \varepsilon}}$; i.e.,
\begin{align*}
\displaystyle
\frac{h_{\max}}{h_{\min}} &= \frac{ \sqrt{(s - s^{\delta})^2 + s^{2 \varepsilon}}}{\sqrt{s^{2 \delta} + s^{2 \varepsilon}}} \leq c s^{1 - \varepsilon} \to \infty \quad \text{as $s \to 0$},\\
\frac{H_{S}}{h_{S}} &= \frac{ \sqrt{(s- s^{\delta})^2 + s^{2 \varepsilon}} \sqrt{s^{2 \delta} + s^{2 \varepsilon}}}{\frac{1}{2} s^{1 + \varepsilon}} \leq c.
\end{align*}
In this case, the element $S$ is "good."
\end{ex}

\begin{rem}
In the above examples, $h_2 \approx \widetilde{\mathscr{H}}_2$ holds. That is, the good element $S \subset \mathbb{R}^2$ may satisfy conditions such as $h_2 \approx \widetilde{\mathscr{H}}_2$.
\end{rem}

\begin{ex}[Blade] \label{Ex4}
Let $S \subset \mathbb{R}^2$ be the simplex with vertices $p_1 := (0,0)^{\top}$, $p_2 := (2s,0)^{\top}$ and $p_3 := (s ,s^{\varepsilon})^{\top}$ with $1 \< \varepsilon $. We then have $h_1 = h_2 = \sqrt{s^{2} + s^{2 \varepsilon}}$; i.e.,
\begin{align*}
\displaystyle
\frac{h_{\max}}{h_{\min}} = 1, \quad \frac{H_{S}}{h_{S}} = \frac{s^{2} + s^{2 \varepsilon}}{s^{1 + \varepsilon}} \to \infty \quad \text{as $s \to 0$}.
\end{align*}
In this case, the element $S$ is "not good."
\end{ex}

\begin{ex}[Dagger] \label{Ex5}
Let $S \subset \mathbb{R}^2$ be the simplex with vertices $p_1 := (0,0)^{\top}$, $p_2 := (s,0)^{\top}$ and $p_3 := (s^{\delta},s^{\varepsilon})^{\top}$ with $1 \< \delta \< \varepsilon $. We then have $h_1 = \sqrt{(s - s^{\delta})^2 + s^{2 \varepsilon}}$ and $h_2 = \sqrt{s^{2 \delta} + s^{2 \varepsilon}}$; i.e.,
\begin{align*}
\displaystyle
\frac{h_{\max}}{h_{\min}} &= \frac{ \sqrt{(s - s^{\delta})^2 + s^{2 \varepsilon}}}{\sqrt{s^{2 \delta} + s^{2 \varepsilon}}} \leq c s^{1 - \delta} \to \infty \quad \text{as $s \to 0$},\\
\frac{H_{S}}{h_{S}} &= \frac{ \sqrt{(s- s^{\delta})^2 + s^{2 \varepsilon}} \sqrt{s^{2 \delta} + s^{2 \varepsilon}}}{\frac{1}{2} s^{1 + \varepsilon}} \leq c s^{\delta - \varepsilon} \to \infty \quad \text{as $s \to 0$}.
\end{align*}
In this case, the element $S$ is "not good."
\end{ex}

Anisotropic elements in the next two examples are also "good." However, these examples differ slightly from Examples \ref{Note2} and \ref{Ex4}.

\begin{ex}[Right-angled triangle] \label{Ad=Ex6}
Let $S \subset \mathbb{R}^2$ be the simplex with vertices $p_1 := (0,0)^{\top}$, $p_2 := (s,0)^{\top}$ and $p_3 := (0, \delta s)^{\top}$ with $\delta \ll 1$. We then have $h_1 = s$ and $h_2 = \delta s$; i.e.,
\begin{align*}
\displaystyle
\frac{h_{\max}}{h_{\min}} &= \frac{1}{\delta} , \quad \frac{H_{S}}{h_{S}} =2.
\end{align*}
In this case, the element $S$ is "good." However, the factor $\frac{1}{\delta}$ is very large.
\end{ex}

\begin{ex}[Blade] \label{Ad=Ex7}
Let $S \subset \mathbb{R}^2$ be the simplex with vertices $p_1 := (0,0)^{\top}$, $p_2 := (2s,0)^{\top}$ and $p_3 := (s , \delta s)^{\top}$ with $\delta \ll 1$. We then have $h_1 = h_2 = s \sqrt{1 + \delta^2}$; i.e.,
\begin{align*}
\displaystyle
\frac{h_{\max}}{h_{\min}} = 1, \quad \frac{H_{S}}{h_{S}} = \frac{s^2 (1 + \delta^2)}{\delta s^2} \leq \frac{c}{\delta},
\end{align*}
In this case, the element $S$ is "good." However, the factor $\frac{1}{\delta}$ is very large.
\end{ex}

\subsection{Anisotropic mesh: three-dimensional case}

\begin{ex} \label{3d=good}
Let $T \subset \mathbb{R}^3$ be a tetrahedron. Let $S$ be the base of $T$; i.e., $S = \triangle p_1 p_2 p_3$. Recall that
\begin{align}
\displaystyle
\frac{H_{T}}{h_{T}} = \frac{h_1 h_2 h_3}{|T|_3} = \frac{h_1 h_2}{\frac{1}{2} h_1 h_2 t_1} \frac{h_3}{\frac{1}{3} h_3 t_2} \leq  \frac{H_{S}}{h_{S}}  \frac{h_3}{\frac{1}{3} \widetilde{\mathscr{H}}_3}. \label{bad56}
\end{align}
If the triangle $S$ is "not good" such as in Examples \ref{Ex4} and \ref{Ex5}, the quantity \eqref{bad56} may diverge. In the following, we consider the case that the triangle $S$ is "good".

Assume that there exists a positive constant $M$ such that $\frac{H_{S}}{h_{S}} \leq M$. For simplicity, we set $p_1 := (0,0,0)^{\top}$, $p_2 := (2s,0,0)^{\top}$, and $p_3 := (2s - \sqrt{4 s^2 - s^{2 \gamma}}, s^{\gamma},0)^{\top}$ with $1 \< \gamma$. Then,
\begin{align*}
\displaystyle
h_1 = 2s, \quad h_2 = \sqrt{\frac{4 s^{2 \gamma}}{2 + \sqrt{4  - s^{2 \gamma -2} }}},
\end{align*}
and because $h_{\max} \approx c s$,
\begin{align*}
\displaystyle
\frac{h_{\max}}{h_{\min}} &\leq \frac{c s}{h_2} \leq c s^{1 - \gamma} \to \infty \quad \text{as $s \to 0$}.
\end{align*}

If we set $p_4 := (s,0,s^{\varepsilon})^{\top}$ with $1 \< \varepsilon$, the triangle $\triangle p_1 p_2 p_4$ is the blade (Example \ref{Ex4}). Then,
\begin{align*}
\displaystyle
h_3 = \sqrt{s^2 + s^{2 \varepsilon}}.
\end{align*}
We thus have
\begin{align*}
\displaystyle
\frac{H_{T}}{h_{T}} \leq c \frac{s^{2 + \gamma}}{s^{1+\gamma+\varepsilon}} \leq c s^{1-\varepsilon} \to \infty \quad \text{as $s \to 0$}.
\end{align*}
In this case, the element $T$ is "not good."

If we set $p_4 := (s^{\delta},0,s^{\varepsilon})^{\top}$ with $1 \< \delta \< \varepsilon \< \gamma$, the triangle $\triangle p_1 p_2 p_4$ is the dagger (Example \ref{Ex5}, Fig. \ref{pic=3d=good}). Then,
\begin{align*}
\displaystyle
h_3 = \sqrt{s^{2 \delta} + s^{2 \varepsilon}}.
\end{align*}
We thus have
\begin{align*}
\displaystyle
\frac{H_{T}}{h_{T}} \leq c \frac{s^{1+\gamma + \delta}}{s^{1+\gamma+\varepsilon}} \leq c s^{\delta - \varepsilon} \to \infty \quad \text{as $s \to 0$}.
\end{align*}
In this case, the element $T$ is "not good."

If we set $p_4 := (s^{\delta},0,s^{\varepsilon})^{\top}$ with $1  \< \varepsilon \< \delta \< \gamma$, the triangle $\triangle p_1 p_2 p_4$ is the dagger (Example \ref{Note1}). Then,
\begin{align*}
\displaystyle
h_3 = \sqrt{s^{2 \delta} + s^{2 \varepsilon}}.
\end{align*}
We thus have
\begin{align*}
\displaystyle
\frac{H_{T}}{h_{T}} \leq c \frac{s^{1+\gamma + \varepsilon} }{ s^{1+\gamma+\varepsilon}} \leq c.
\end{align*}
In this case, the element $T$ is "good" and  $h_3 \approx h_3 t_2 = \widetilde{\mathscr{H}}_3$ holds.

\begin{figure}[htbp]
\vspace{-5cm}
  \includegraphics[bb=0 0 774 530,scale=0.55]{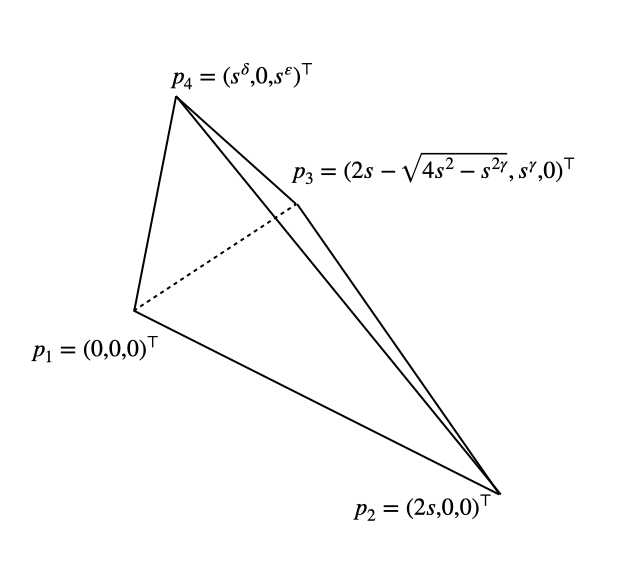}
\caption{Example \ref{3d=good}}
\label{pic=3d=good}
\end{figure}

\end{ex}

\begin{ex}
In \cite{Cheetal00}, the spire has a cycle of three daggers among its four triangles; see Figure \ref{FIG5}. The splinter has four daggers; see Figure \ref{FIG9}. The spear and spike have two daggers and two blades as triangles; see Figures \ref{FIG6}, \ref{FIG8}. The spindle has four blades as triangles; see Figure \ref{FIG7}.

\begin{figure}[htbp]
\vspace{-2cm}
  \begin{minipage}[b]{0.3\linewidth}
    \centering
    \includegraphics[bb=0 0 520 320,scale=0.35]{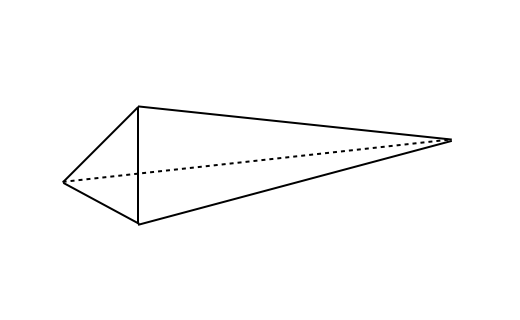}
    \caption{Spire}
    \label{FIG5}
  \end{minipage}
  \begin{minipage}[b]{0.3\linewidth}
    \centering
     \includegraphics[bb=0 0 520 320,scale=0.35]{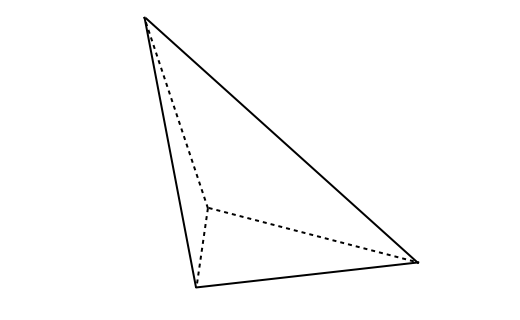}
    \caption{Spear}
    \label{FIG6}
  \end{minipage}
  \begin{minipage}[b]{0.3\linewidth}
    \centering
     \includegraphics[bb=0 0 520 320,scale=0.35]{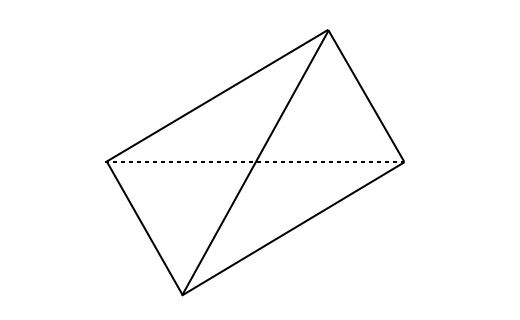}
    \caption{Spindle}
    \label{FIG7}
  \end{minipage}
\end{figure}

\begin{figure}[htbp]
\vspace{-2cm}
  \begin{minipage}[b]{0.3\linewidth}
    \centering
    \includegraphics[bb=0 0 520 320,scale=0.35]{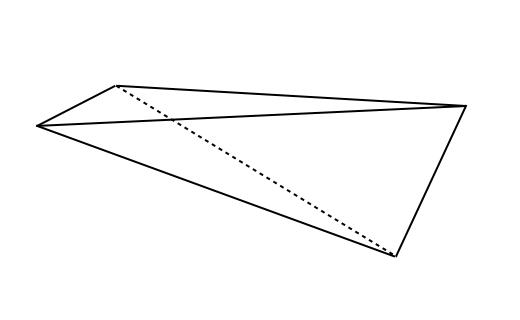}
    \caption{Spike}
    \label{FIG8}
  \end{minipage}
  \begin{minipage}[b]{0.3\linewidth}
    \centering
     \includegraphics[bb=0 0 520 320,scale=0.35]{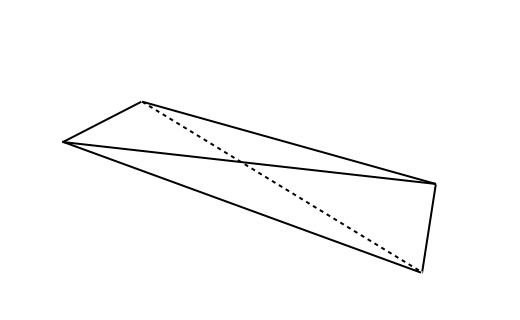}
    \caption{Splinter}
    \label{FIG9}
  \end{minipage}
\end{figure}

\end{ex}

\begin{rem}
The above examples reveal that the good element $T \subset \mathbb{R}^3$ may satisfy conditions such as $h_2 \approx \tilde{\mathscr{H}}_2$ and $h_3 \approx  \tilde{\mathscr{H}}_3$.
\end{rem}

\begin{figure}[htbp]
 \vspace{-3.5cm}
   \includegraphics[bb=0 0 522 454,scale=0.6]{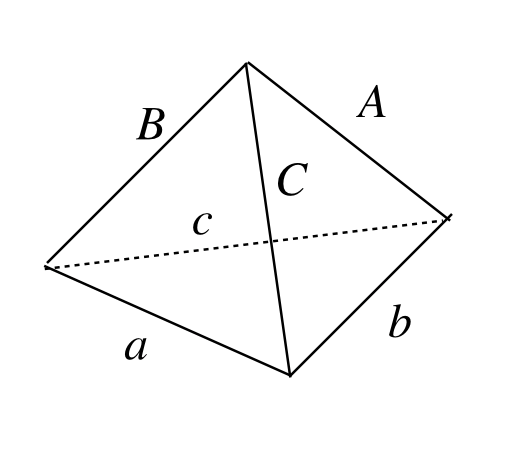}
    \caption{$R_3$}
  \label{fig:one}
\end{figure}

\begin{figure}[htbp]
\vspace{-10cm}
   \includegraphics[bb=0 0 690 642,scale=0.85]{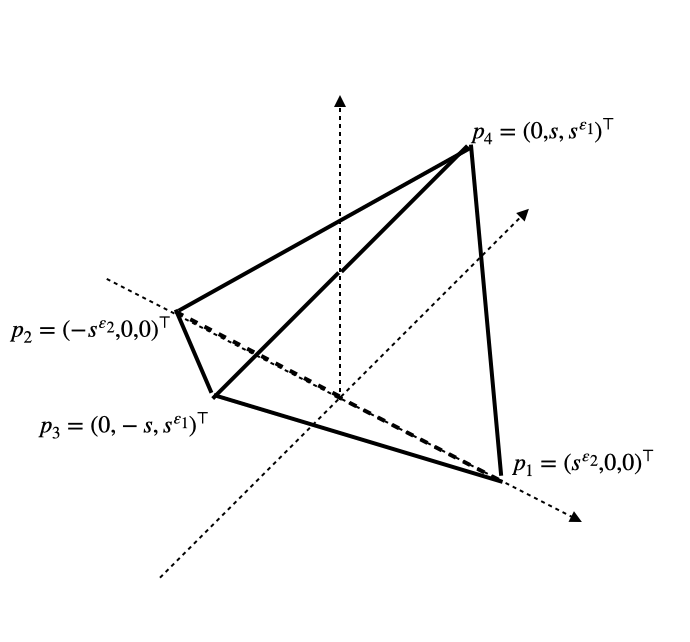}
  \caption{Sliver}
  \label{Sliver}
\end{figure}

\begin{ex} \label{Ex_sliver}
Using an element $T$ called \textit{Sliver}, we compare the three quantities $\frac{h_{T}^3}{|{T}|_3}$, $\frac{H_{T}}{h_{T}}$, and $\frac{R_3}{h_{T}}$, where the formulation of the circumradius $R_3$ of a tetrahedron $T$ is as follows, e.g., see \cite{HitKur14}. Let $a$, $b$ and $c$ be the lengths of the three edges of $T$ and $A$, $B$, $C$ the length of the opposite edges of $a$, $b$, $c$, respectively. Then, 
\begin{align*}
\displaystyle
R_3 = \frac{\sqrt{(a A+b B+c C)(a A+b B-c C)(a A-b B+c C)(-a A+b B+c C)}}{24 |T|_3},
\end{align*}
see Fig. \ref{fig:one}.

Let $T \subset \mathbb{R}^3$ be the simplex with vertices $p_1 := (s^{\varepsilon_2},0,0)^{\top}$, $p_2 := (-s^{\varepsilon_2},0,0)^{\top}$, $p_3 := (0,-s,s^{\varepsilon_1})^{\top}$, and $p_4 := (0,s,s^{\varepsilon_1})^{\top}$ ($\varepsilon_1, \varepsilon_2 \> 1$), where $s := \frac{1}{N}$, $N \in \mathbb{N}$, see Fig. \ref{Sliver}. Let $L_i$ ($1 \leq i \leq 6$) be the edges of $T$ with $h_{\min} = L_1 \leq L_2 \leq \cdots \leq L_6 = h_{T}$ . Recall that $h_{\max} \approx h_{T}$ and
\begin{align*}
\displaystyle
\frac{h_{\max}}{h_{\min}} \leq c \frac{L_6}{L_1}, \quad \frac{H_{T}}{h_{T}} = \frac{L_1 L_2}{|T|_3} h_{T}.
\end{align*}

\begin{table}[htbp]
\caption{$h_{T}^3/{|T|_3}$, $H_{T}/h_{T}$ and $R_3/h_{T}$ ($\varepsilon_1 = 1.5$, $\varepsilon_2 = 1.0$)}
\centering
\begin{tabular}{l | l | l | l | l | l } \hline
$N$ &  $s$ & $L_6 / L_1$ & $h_{T}^3/{|T|_3}$  & $H_{T}/h_{T}$ & $R_3/h_{T}$\\ \hline \hline
32 & 3.1250e-02 & 1.4033 & 6.7882e+01   & 3.4471e+01  &  5.0195e-01     \\
64 & 1.5625e-02  & 1.4087 & 9.6000e+01  & 4.8375e+01  &  5.0098e-01  \\
128 &7.8125e-03 &  1.4115 & 1.3576e+02  &6.8147e+01 &  5.0049e-01   \\
\hline
\end{tabular}
\label{tab=sliver1}
\end{table}

\begin{table}[ht]
\caption{$h_{T}^3/{|T|_3}$, $H_{T}/h_{T}$ and $R_3/h_{T}$ ($\varepsilon_1 = 1.0$, $\varepsilon_2 = 1.5$)}
\centering
\begin{tabular}{l | l | l | l | l | l } \hline
$N$ &  $s$ & $L_6 / L_1$ & $h_{T}^3/{|T|_3}$  & $H_{T}/h_{T}$ & $R_3/h_{T}$\\ \hline \hline
32 & 3.1250e-02 & 5.6569 & 6.7882e+01  & 8.5513  & 5.0006e-01  \\
64 & 1.5625e-02 & 8.0000  & 9.6000e+01 &8.5184   & 5.0002e-01 \\
128 &7.8125e-03 & 1.1314e+01  & 1.3576e+02  & 8.5018 & 5.0000e-01     \\
\hline
\end{tabular}
\label{tab=sliver2}
\end{table}

\begin{table}[ht]
\caption{$h_{T}^3/{|T|_3}$, $H_{T}/h_{T}$ and $R_3/h_{T}$ ($\varepsilon_1 = 1.5$, $\varepsilon_2 = 1.5$)}
\centering
\begin{tabular}{l | l | l | l | l | l } \hline
$N$ &  $s$ & $L_6 / L_1$ & $h_{T}^3/{|T|_3}$  & $H_{T}/h_{T}$ & $R_3/h_{T}$\\ \hline \hline
32 & 3.1250e-02 & 5.6569 & 3.8400e+02 & 3.4986e+01 &  1.4170  \\
64 & 1.5625e-02  & 8.0000  & 7.6800e+02 & 4.8744e+01 & 2.0010 \\
128 &7.8125e-03& 1.1314e+01 & 1.5360e+03 & 6.8411e+01& 2.8288 \\
\hline
\end{tabular}
\label{tab=sliver3}
\end{table}

In Table \ref{tab=sliver1}, the angle between $\triangle p_1 p_2 p_3$ and $\triangle p_1 p_2 p_4$ tends to $\pi$ as $s \to 0$, and the simplex $T$ is "not good." In Table \ref{tab=sliver2}, the angle between $\triangle p_1 p_3 p_4$ and $\triangle p_2 p_3 p_4$ tends to $0$ as $s \to 0$, the simplex $T$ is "good." In Table \ref{tab=sliver3}, from the numerical result, the simplex $T$ is "not good."
\end{ex}

\section{Examples of Anisotropic Mesh Partitions used in Numerical Calculations} 

\subsection{Examples of Good Elements}
\begin{figure}[htbp]
  \begin{minipage}[b]{0.45\linewidth}
    \centering
    \includegraphics[keepaspectratio, scale=0.15]{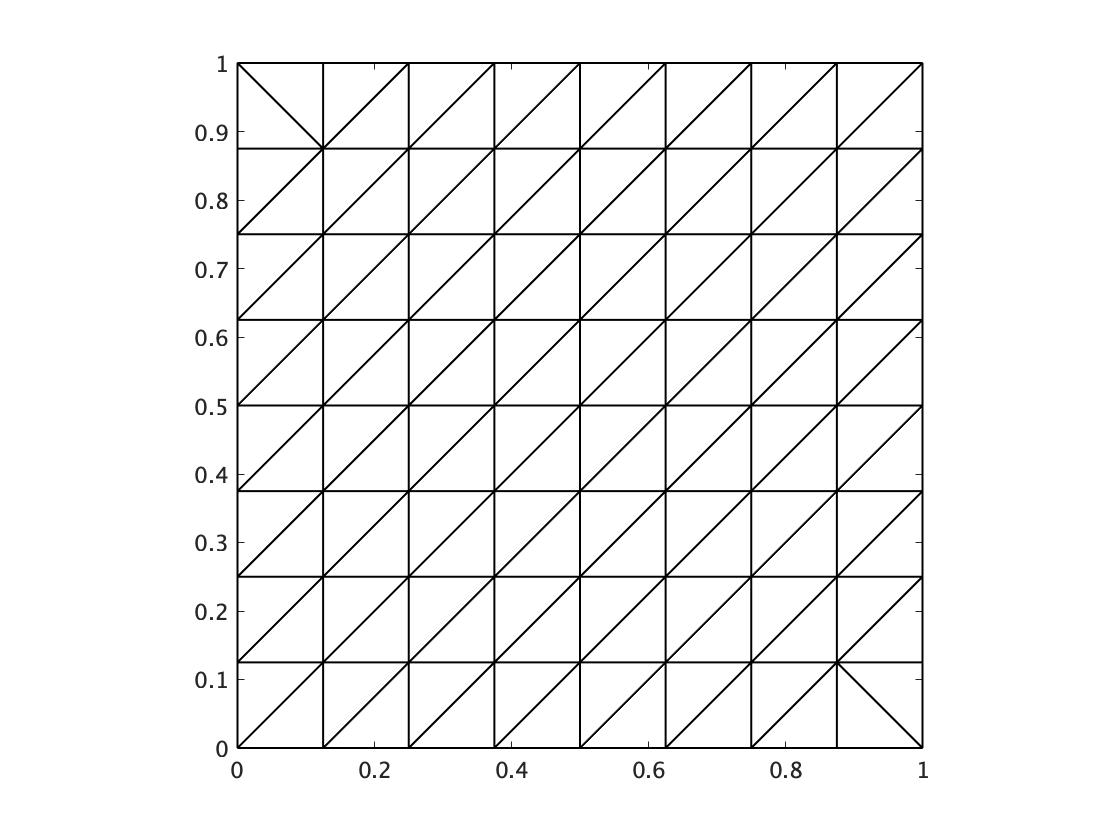}
    \caption{(\Roman{lone}) Standard mesh}
     \label{fig1}
  \end{minipage}
  \begin{minipage}[b]{0.45\linewidth}
    \centering
   \includegraphics[keepaspectratio, scale=0.15]{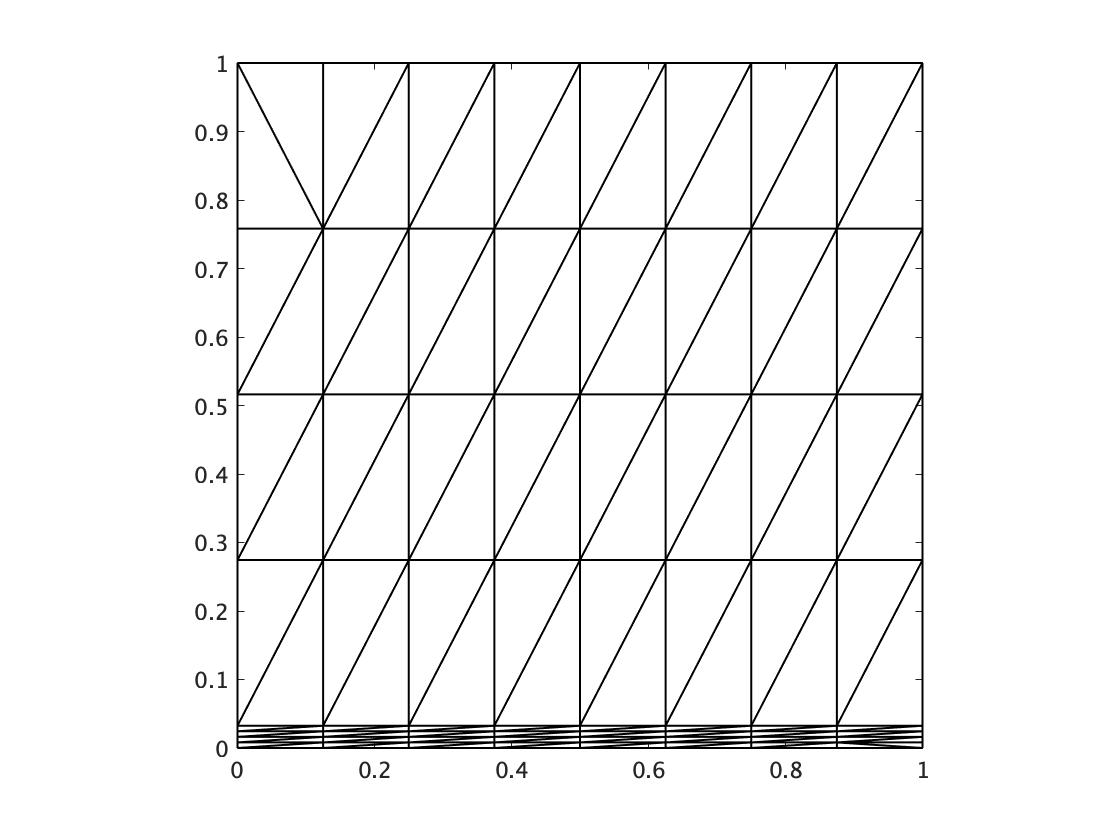}
    \caption{(\Roman{ltwo}) Shishkin mesh, $\delta = \frac{1}{128}$}
     \label{fig2}
  \end{minipage}
\end{figure}

\begin{figure}[htbp]
  \begin{minipage}[b]{0.45\linewidth}
    \centering
   \includegraphics[keepaspectratio, scale=0.15]{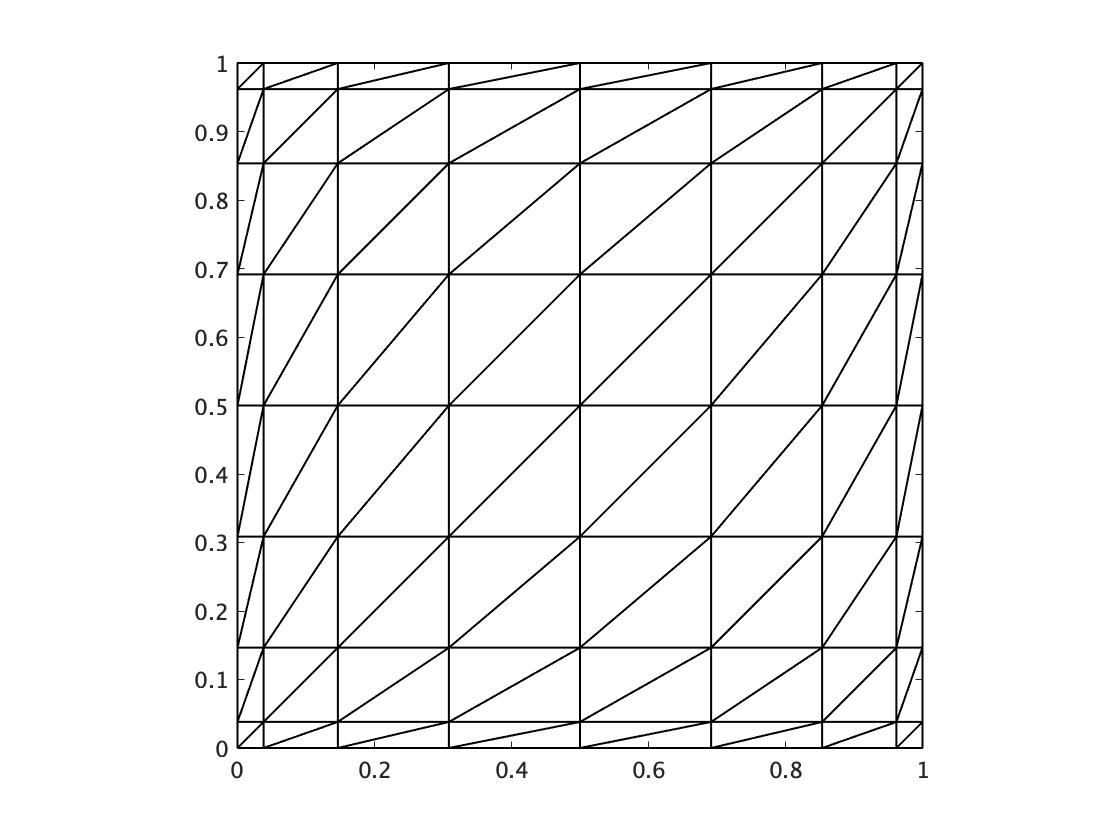}
    \caption{(\Roman{lthree}) Anisotropic mesh}
     \label{fig33}
  \end{minipage}
  \begin{minipage}[b]{0.45\linewidth}
    \centering
  \includegraphics[keepaspectratio, scale=0.15]{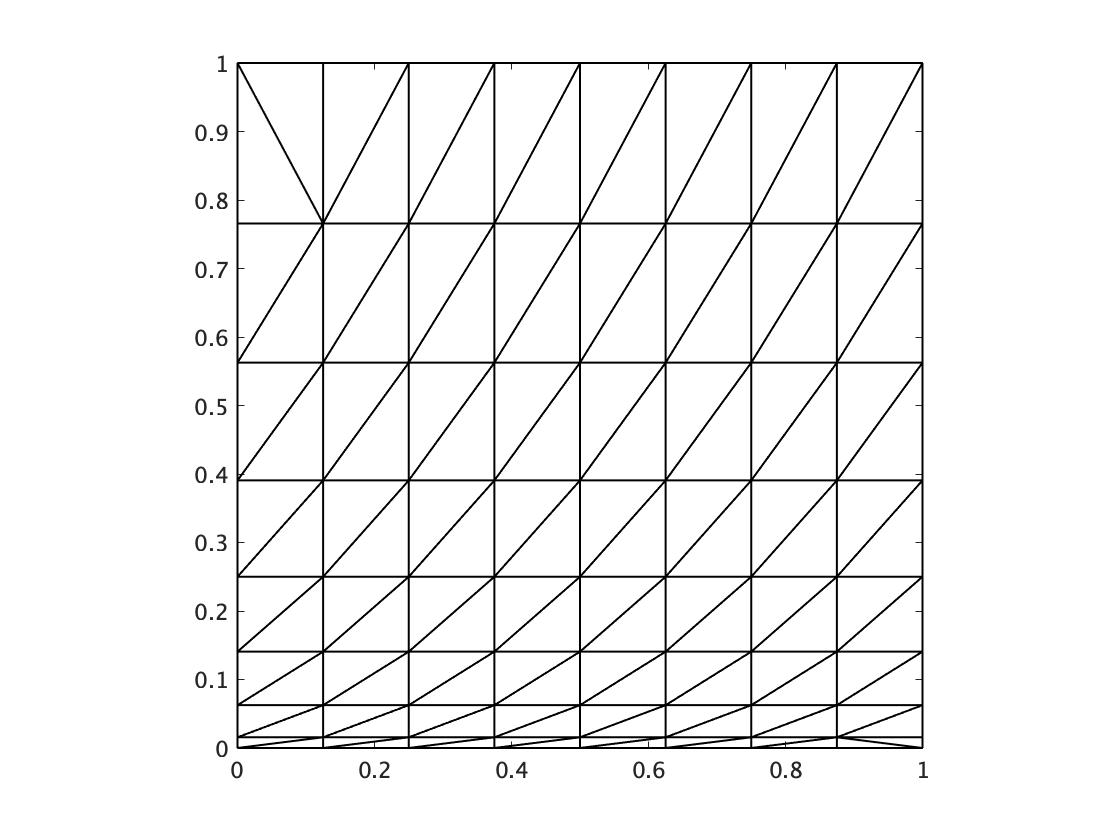}
    \caption{(\Roman{lfour}) Graded mesh}
     \label{fig44}
  \end{minipage}
\end{figure}

Let $N$ be the division number of each side of the bottom and the height edges of $\Omega$. We consider four types of mesh partitions. Let $(x_1^i, x_2^i)^{\top}$ be grip points of triangulations $\mathbb{T}_h$ defined as follows. Let $i \in \mathbb{N}$.
\begin{description}
  \item[(\Roman{lone}) Standard mesh (Fig. \ref{fig1})] 
\begin{align*}
\displaystyle
x_1^i := \frac{i}{N}, \quad x_2^i := \frac{i}{N}, \quad  i \in \{0, \ldots , N \}.
\end{align*}
  \item[(\Roman{ltwo}) Shishkin mesh (Fig. \ref{fig2})]
\begin{align*}
\displaystyle
x_1^i &:= \frac{i}{N}, \quad  i \in \{0 , \ldots , N \}, \\
 x_2^i &:=
\begin{cases}
\tau \frac{2}{N} i, \quad  i \in \left\{0, \ldots , \frac{N}{2} \right \}, \\
\tau + (1 - \tau) \frac{2}{N} \left( i - \frac{N}{2} \right), \quad i \in \left\{ \frac{N}{2}+1 , \ldots , N \right\},
\end{cases}
\end{align*}
where $\tau := 2 \delta | \ln(N) |$ with $\delta = \frac{1}{128}$, see \cite[Section 2.1.2]{Lin10}.
\item[(\Roman{lthree}) Anisotropic mesh {from \cite{CheLiuQia10}} (Fig. \ref{fig33})]
\begin{align*}
\displaystyle
x_1^i := \frac{1}{2}\left( 1 - \cos \left( \frac{i \pi}{N} \right) \right), \quad x_2^i := \frac{1}{2}\left( 1 - \cos \left( \frac{i \pi}{N} \right) \right), \quad  i \in \{0, \ldots, N \}.
\end{align*}
 \item[(\Roman{lfour}) Graded mesh (Fig. \ref{fig44})]
 \begin{align*}
\displaystyle
x_1^i := \frac{i}{N}, \quad x_2^i := \left ( \frac{i}{N} \right)^{2}, \quad  i \in \{0, \ldots, N \}.
\end{align*}
\end{description}

\subsection{Are geometric mesh conditions fulfilled?}
As described in Section \ref{regmesh=sec3=2}, the shape-regularity condition is equivalent to the following condition. There exists a constant $\gamma_2 \> 0$ such that for any $\mathbb{T}_h \in \{ \mathbb{T}_h \}$ and simplex $T \in \mathbb{T}_h$, we have
\begin{align*}
\displaystyle
|T|_2 \geq \gamma_3 h_{T}^2. 
\end{align*}
We presented the new geometric mesh condition in Section \ref{new=geeocon=sec}, i.e., there exists $\gamma_0 \> 0$ such that
\begin{align*}
\displaystyle
\frac{H_{T}}{h_{T}} \leq \gamma_0 \quad \forall \mathbb{T}_h \in \{ \mathbb{T}_h \}, \quad \forall T \in \mathbb{T}_h. 
\end{align*}
The following parameters are computed.
\begin{align*}
\displaystyle
\textit{MinAngle} := \max_{T \in \mathbb{T}_h} \frac{|L_3|^2}{|T|_2}, \quad \textit{MaxAngle} := \max_{T \in \mathbb{T}_h} \frac{|L_1| |L_2|}{|T|_2},
\end{align*}
where $L_i$, $i=1,2,3,$ denote the edges of the simplex $T \in \mathbb{T}_h$ with $|L_1| \leq |L_2| \leq |L_3|$. 

\begin{table}[htbp]
\caption{Mesh conditions}
\centering
\begin{tabular}{l|l|l|l} \hline
 Mesh &$N$ &  $\textit{MinAngle} $ & $\textit{MaxAngle}$   \\ \hline \hline
 \Roman{lone} & 32  & 4.00000 & 2.00000   \\
 & 64  &  4.00000  & 2.00000     \\
 & 128  & 4.00000 & 2.00000   \\
\hline
 \Roman{ltwo} & 32  & 1.86831e+01   & 2.00000    \\
 & 64  &  1.56487e+01 & 2.00000     \\
& 128  &  1.34936e+01   & 2.00000    \\
\hline
 \Roman{lthree}  & 32  &4.08092e+01 & 2.00000   \\
& 64  &  8.15201e+01   &  2.00000    \\
& 128  & 1.62991e+02   & 2.00000   \\
\hline
 \Roman{lfour} & 32  & 6.40625e+01 & 2.00000   \\
 & 64  & 1.28031e+02   & 2.00000    \\
& 128  &2.56016e+02   & 2.00000   \\
\hline
\end{tabular}
\label{table=mesh}
\end{table}

Notably, a sequence with meshes (\Roman{lone}) or (\Roman{ltwo}) satisfies the shape-regularity condition, but a sequence with meshes (\Roman{lthree}) or (\Roman{lfour}) does not fulfil the shape-regularity condition. See Table \ref{table=mesh}.

\subsection{Bad Elements}
\begin{figure}[htbp]
  \begin{minipage}[b]{0.45\linewidth}
    \centering
   \includegraphics[keepaspectratio, scale=0.15]{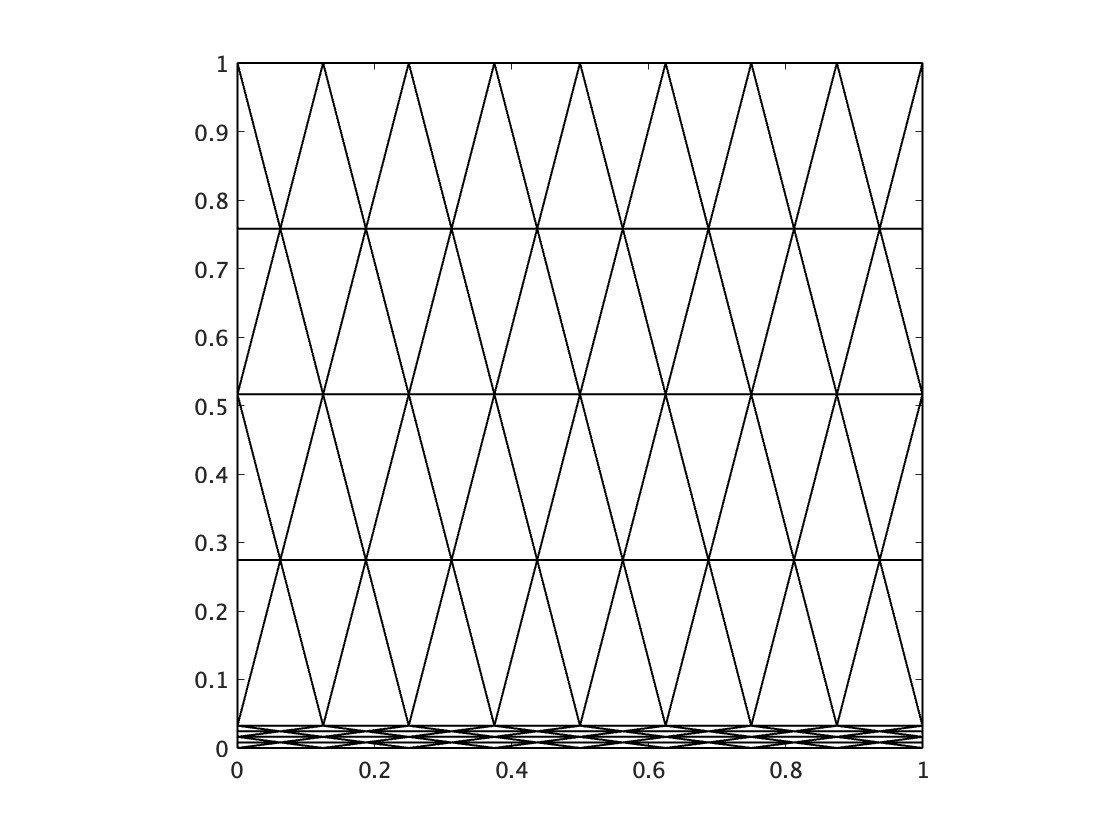}
    \caption{(\Roman{lfive}) Shishkin   mesh}
     \label{fig112}
  \end{minipage}
  \begin{minipage}[b]{0.45\linewidth}
    \centering
  \includegraphics[keepaspectratio, scale=0.15]{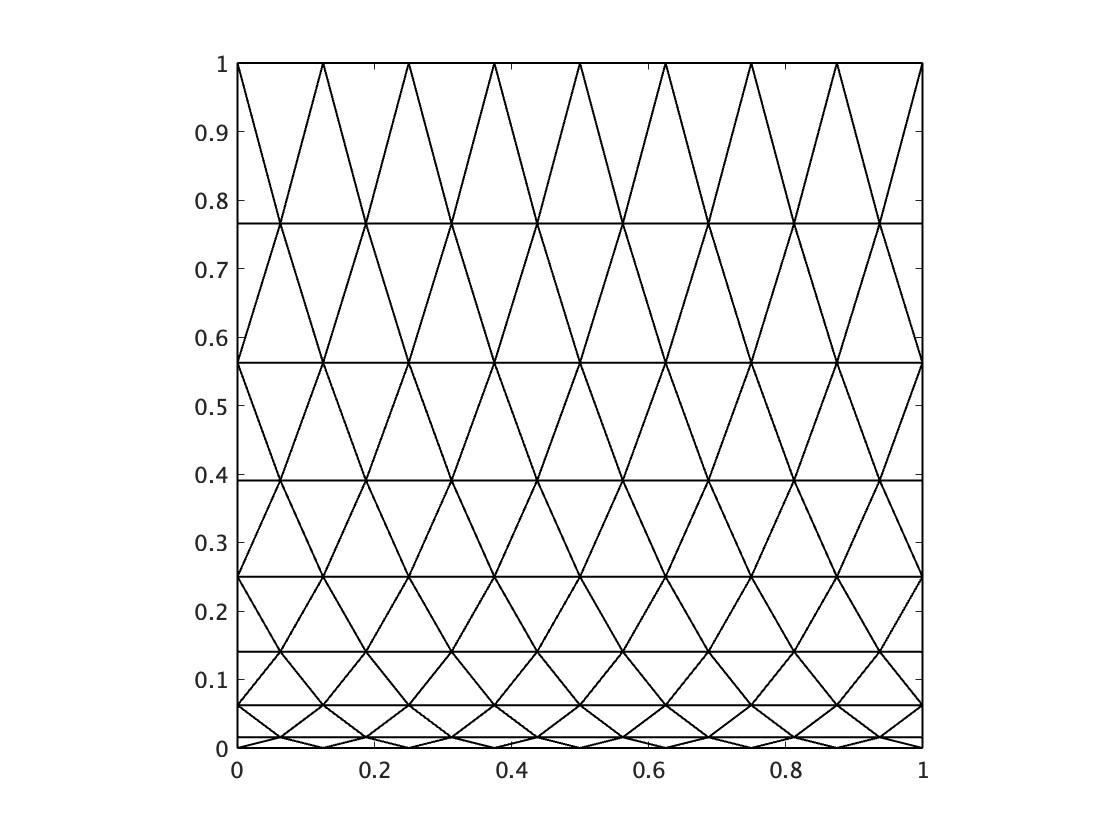}
    \caption{(\Roman{lsix}) Graded  mesh}
     \label{fig113}
  \end{minipage}
\end{figure}

We consider the following mesh partitions:
\begin{description}
\item[(\Roman{lfive}) Shishkin mesh (Fig. \ref{fig112})],
\item[(\Roman{lsix}) Graded mesh (Fig. \ref{fig113})].
\end{description}

\begin{table}[htbp]
\caption{Mesh conditions}
\centering
\begin{tabular}{l|l|l|l} \hline
 Mesh &$N$ &  $\textit{MinAngle} $ & $\textit{MaxAngle}$   \\ \hline \hline
 \Roman{lfive} & 32  &1.84665e+01  &4.83323\\
 & 64  &  1.53887e+01 &  4.10712   \\
& 128  & 1.31904e+01  & 3.60084    \\
\hline
 \Roman{lsix} & 32  & 6.40000e+01  & 1.60625e+01    \\
 & 64  & 1.28000e+02   &3.20312e+01    \\
& 128  &2.56000e+02   &6.40156e+01     \\
\hline
\end{tabular}
\label{table=mesh2}
\end{table}

Notably, a sequence with the mesh (\Roman{lfive}) satisfies the shape-regularity condition, but a sequence with the mesh (\Roman{lsix}) does not satisfy either the shape-regularity condition or the semi-regular mesh condition. See Table \ref{table=mesh2}. When Example (\Roman{lfive}) is used, interpolation error estimates approach zero as the mesh division becomes finer, see Example \ref{ex14=8}. It is 'good' in the sense described in Section \ref{goodbad}. As discussed in Example \ref{ex14=8}, the error may become larger.

\section{FE Generation} \label{FEMG}
We follow the procedure described in \cite[Chapter 9]{ErnGue21a} and \cite[Section 1.4.1 and 1.2.1]{ErnGue04}; also see \cite[Section 3.5]{IshKobTsu21a}. The definition of a \textit{finite element} can be found in \cite[p. 78]{Cia02} and \cite[Definition 5.2]{ErnGue21a}.

For the reference element $\widehat{T}$ defined in Sections \ref{reference}, let $\{ \widehat{T} , \widehat{{P}} , \widehat{\Sigma} \}$ be a fixed reference finite element, where $\widehat{{P}} $ is a vector space of functions $\hat{q}: \widehat{T} \to \mathbb{R}^n$ for some positive integer $n$ (typically $n=1$ or $n=d$) and $\widehat{\Sigma}$ is a set of $n_{0}$ linear forms $\{ \hat{\chi}_1 , \ldots , \hat{\chi}_{n_0} \}$ such that
\begin{align*}
\displaystyle
\widehat{{P}} \ni \hat{q} \mapsto (\hat{\chi}_1(\hat{q}) , \ldots , \hat{\chi}_{n_0}(\hat{q}))^{\top} \in \mathbb{R}^{n_0}
\end{align*}
is bijective; i.e., $\widehat{\Sigma}$ is a basis for $\mathcal{L}(\widehat{P};\mathbb{R})$. Further, we denote by $\{ \hat{\theta}_1 , \ldots, \hat{\theta}_{n_0} \}$ in $\widehat{{P}}$ the local ($\mathbb{R}^n$-valued) shape functions such that
\begin{align*}
\displaystyle
\hat{\chi}_i(\hat{\theta}_j) = \delta_{ij}, \quad 1 \leq i,j \leq n_0.
\end{align*}

Let $V(\widehat{T})$ be a normed vector space of functions $\hat{\varphi}: \widehat{T} \to \mathbb{R}^n$ such that $\widehat{P} \subset V(\widehat{T})$ and the linear forms $\{ \hat{\chi}_1 , \ldots , \hat{\chi}_{n_0} \}$ can be extended to $V(\widehat{T})^{\prime}$, i.e., there exist $\{ \bar{\chi}_1 , \ldots , \bar{\chi}_{n_0} \}$ and $c_{\chi}$ such that $\hat{\chi}_i(\hat{q}) = \bar{\chi}_i (\hat{q})$ for any $\hat{q} \in \widehat{{P}}$, and $|  \bar{\chi}_i  (\hat{v})| \leq c_{\chi} \| \hat{v} \|_{V(\widehat{T})}$ and for $i \in \{ 1,\ldots,n_0 \}$. We use the same symbol $\hat{\chi}_i$ instead of $\bar{\chi}_i$. The local interpolation operator ${I}_{\widehat{T}}$ is then defined by
\begin{align}
\displaystyle
{I}_{\widehat{T}} : V(\widehat{T}) \ni \hat{\varphi} \mapsto \sum_{i=1}^{n_0} \hat{\chi}_i (\hat{\varphi}) \hat{\theta}_i \in \widehat{{P}}. \label{int1}
\end{align}
It obviously holds that, for any $\hat{\varphi} \in V(\widehat{T})$,
\begin{align}
\displaystyle
 \hat{\chi}_i ({I}_{\widehat{T}} \hat{\varphi}) =  \hat{\chi}_i (\hat{\varphi}) \quad i=1,\ldots,n_0. \label{int512}
\end{align}

\begin{prop}
$\widehat{P}$ is invariant under $I_{\widehat{T}}$, that is, 
\begin{align}
\displaystyle
I_{\widehat{T}} \hat{q} = \hat{q} \quad \forall \hat{q} \in \widehat{P}. \label{int513}
\end{align}
\end{prop}

\begin{pf*}
Let $\hat{q} = \sum_{j=1}^{n_0} \alpha_j  \hat{\theta}_j$ for $\alpha_j \in \mathbb{R}$, $1 \leq j \leq n_0$. Then,
\begin{align*}
\displaystyle
I_{\widehat{T}} \hat{q} = \sum_{i,j=1}^{n_0} \alpha_j \hat{\chi}_i(\hat{\theta}_j) \hat{\theta}_i = \hat{q}.
\end{align*}
\qed
\end{pf*}
Let $\Phi_{\widetilde{T}}: \widehat{T} \to \widetilde{T}$ and $\Phi_{T}: \widetilde{T} \to T$ be the two affine mappings defined in Section \ref{two=step}. For any ${T} \in \mathbb{T}_h$ with $T = \Phi(\widehat{T}) = ({\Phi}_{T} \circ {\Phi}_{\widetilde{T}})(\widehat{T})$, we define a Banach space $V(T)$ of $\mathbb{R}^n$-valued functions that is the counterpart of $V(\widehat{T})$ and define a linear bijection mapping by
\begin{align*}
\displaystyle
\psi := \psi_{\widehat{T}} \circ \psi_{\widetilde{T}}: V(T) \ni \varphi \mapsto \hat{\varphi} := \psi(\varphi) := \varphi \circ \Phi \in   V(\widehat{T}),
\end{align*}
with two linear bijection mappings:
\begin{align*}
\displaystyle
&\psi_{\widetilde{T}} : V({T}) \ni \varphi \mapsto \tilde{\varphi} := \psi_{\widetilde{T}}(\varphi) := \varphi \circ {\Phi_T} \in   V(\widetilde{T}), \\
&\psi_{\widehat{T}} : V(\widetilde{T}) \ni \tilde{\varphi} \mapsto \hat{\varphi} :=  \psi_{\widehat{T}} (\tilde{\varphi}) :=  \tilde{\varphi} \circ {\Phi_{\widetilde{T}}} \in V(\widehat{T}).
\end{align*}

Triples $\{ \widetilde{T} , \widetilde{P} , \widetilde{\Sigma}\}$ and $\{ T , {P}, \Sigma \}$ are defined as follows:
\begin{align*}
\displaystyle
\begin{cases}
\displaystyle
\widetilde{T} = \Phi_{\widetilde{T}}(\widehat{T}); \\
\displaystyle
 \widetilde{P} = \{ \psi_{\widehat{T}}^{-1}(\hat{q}) ; \ \hat{q} \in \widehat{{P}}\}; \\
\displaystyle
\widetilde{\Sigma} = \{ \{ \tilde{\chi}_{i} \}_{1 \leq i \leq n_0}; \ \tilde{\chi}_{i} = \hat{\chi}_i(\psi_{\widehat{T}}(\tilde{q})), \forall \tilde{q} \in \widetilde{P}, \hat{\chi}_i \in \widehat{\Sigma} \},
\end{cases}
\end{align*}
and
\begin{align*}
\displaystyle
\begin{cases}
\displaystyle
T = {\Phi}_T(\widetilde{T}); \\
\displaystyle
{P} = \{ \psi_{\widetilde{T}}^{-1}(\tilde{q}) ; \ \tilde{q} \in \widetilde{{P}}\}; \\
\displaystyle
\Sigma = \{ \{ \chi_{i} \}_{1 \leq i \leq n_0}; \ \chi_{i} = \tilde{\chi}_i(\psi_{\widetilde{T}}(q)), \forall q \in {P}, \tilde{\chi}_i \in \widetilde{\Sigma} \}.
\end{cases}
\end{align*}

\begin{prop}
The triples $\{ \widetilde{T} , \widetilde{P} , \widetilde{\Sigma}\}$ and $\{ T , {P}, \Sigma \}$ are finite elements.
\end{prop}

\begin{pf*}
A proof can be obtained similarly for \cite[Proposition 9.2]{ErnGue21a}.
\qed
\end{pf*}

The local shape functions are $\tilde{\theta}_{i} = \psi_{\widehat{T}}^{-1}(\hat{\theta}_i)$ and $\theta_{i} = \psi_{\widetilde{T}}^{-1}(\tilde{\theta}_i)$, $1 \leq i \leq n_0$, and the associated local interpolation operators are respectively defined by
\begin{align}
\displaystyle
{I}_{\widetilde{T}} : V(\widetilde{T}) \ni \tilde{\varphi} \mapsto {I}_{\widetilde{T}} \tilde{\varphi} &:= \sum_{i=1}^{n_0} \tilde{\chi}_{i}(\tilde{\varphi}) \tilde{\theta}_{i} \in \widetilde{P}, \label{int2} \\
{I}_{T} : V(T) \ni \varphi \mapsto {I}_{T} \varphi &:= \sum_{i=1}^{n_0} \chi_{i}(\varphi) \theta_{i} \in {P}. \label{int3}
\end{align}

The following diagrams play an important role in analysing the interpolation error.

\begin{prop}[Commuting diagrams] \label{prop512}
The diagrams
\begin{align*}
\displaystyle
\xymatrix{
V(T)\ar[r]^-{\psi_{\widetilde{T}}}\ar[d]_-{I_{T}}\ar@{}|{}&V(\widetilde{T})\ar[r]^-{\psi_{\widehat{T}}}\ar[d]_-{{I}_{\widetilde{T}}}\ar@{}|{}&V(\widehat{T})\ar[d]^-{{I}_{\widehat{T}}} \\
P\ar[r]_-{\psi_{\widetilde{T}}}&\widetilde{{P}}\ar[r]_-{\psi_{\widehat{T}}}&\widehat{{P}}
}
\end{align*}
commute. Furthermore, $\widetilde{P}$ and $P$ are respectively invariant under ${I}_{\widetilde{T}}$ and ${I}_{T}$.
\end{prop}

\begin{pf*}
A proof can be obtained similarly for \cite[Proposition 9.3]{ErnGue21a}.

Let $\tilde{\varphi} \in  V(\widetilde{T})$. The definition of $\{ \widetilde{T} , \widetilde{P} , \widetilde{\Sigma}\}$ implies that
\begin{align*}
\displaystyle
I_{\widehat{T}}(\psi_{\widehat{T}}(\tilde{\varphi}))
=  \sum_{i=1}^{n_0} \hat{\chi}_i (\psi_{\widehat{T}}(\tilde{\varphi})) \hat{\theta}_i
=  \sum_{i=1}^{n_0} \tilde{\chi}_i (\tilde{\varphi})   \psi_{\widehat{T}} (\tilde{\theta}_{i})
= \psi_{\widehat{T}} ({I}_{\widetilde{T}} \tilde{\varphi}).
\end{align*}
Here, we used the linearity of $\psi_{\widehat{T}}$. Therefore, the right diagram commutes. 

Let $\tilde{q} \in \widetilde{P}$. Because $\psi_{\widehat{T}}(\tilde{q}) \in \widehat{P}$ and $\widehat{P}$ is invariant under ${I}_{\widehat{T}}$, 
\begin{align*}
\displaystyle
{I}_{\widetilde{T}}(\tilde{q})
= \psi_{\widehat{T}}^{-1} (I_{\widehat{T}}(\psi_{\widehat{T}}(\tilde{q})) )  
= \psi_{\widehat{T}}^{-1} (\psi_{\widehat{T}}(\tilde{q}))  = \tilde{q}.
\end{align*}

Another diagram can be proved in the same way.
\qed
\end{pf*}

\begin{ex}
Let  $\{ \widehat{T} , \widehat{{P}} , \widehat{\Sigma} \}$ be a finite element. 
\begin{enumerate}
 \item For the Lagrange finite element of degree $k$, we set $V(\widehat{T}) := \mathcal{C}^0(\widehat{T})$.
 \item For the Hermite finite element, we set $V(\widehat{T}) := \mathcal{C}^1(\widehat{T})$.
 \item For the Crouzeix--Raviart finite element with $k=1$, we set $V(\widehat{T}) := W^{1,1}(\widehat{T})$. 
\end{enumerate}
\end{ex}

\section{New Scaling Argument: Part 1} \label{sec54}
This section gives estimates related to a scaling argument corresponding to \cite[Lemma 1.101]{ErnGue04}.

\subsection{Preliminalies}
\subsubsection{Additional New Condition} \label{add=new=cond=10}
The following condition is used for obtaining optimal interpolation error estimates.

\begin{Cond} \label{Cond333}
In anisotropic interpolation error analysis, we impose the following geometric condition for the simplex $T$:
\begin{enumerate}
 \item If $d=2$, there are no additional conditions;
 \item If $d=3$, there must exist a positive constant $M$ independent of $h_{T}$ such that $|s_{22}| \leq M \frac{h_2 t_1}{h_3}$. Note that if $s_{22} \neq 0$, this condition means that the order of $h_3$ with respect to $h_{T}$ coincides with the order of $h_2$, and if $s_{22} = 0$, the order of $h_3$ may be different from that of $h_2$. 
\end{enumerate}
\end{Cond}

Recall that
\begin{align*}
\displaystyle
&|s| \leq 1, \quad h_2 \leq h_1 \quad \text{if $d=2$},\\
&|s_1|\leq 1, \ |s_{21}| \leq 1, \quad  h_2 \leq h_3\leq h_1 \quad \text{if $d=3$}.
\end{align*}
When $d=3$, if Condition \ref{Cond333} is imposed, there exists a positive constant $M$ independent of $h_T$ such that $|s_{22}| \leq M \frac{h_2 t_1}{h_3}$. 
We thus have, if $d=2$,
\begin{align*}
\displaystyle
&h_1 | [\widetilde{{A}}]_{j1} | \leq \widetilde{\mathscr{H}}_j, \quad h_2 | [\widetilde{{A}}]_{j2} | \leq \widetilde{\mathscr{H}}_j, \quad j=1,2,
\end{align*}
and, if $d=3$, for $\widetilde{{A}} \in \{ \widetilde{{A}}_1 , \widetilde{{A}}_2  \}$ and $j=1,2,3$,
\begin{align*}
\displaystyle
&h_1 | [\widetilde{{A}}]_{j1} | \leq \widetilde{\mathscr{H}}_j, \quad h_2 | [\widetilde{{A}}]_{j2} | \leq \widetilde{\mathscr{H}}_j, \quad h_3 | [\widetilde{{A}}]_{j3} | \leq \max \{ 1,M\} \widetilde{\mathscr{H}}_j,  \quad j=1,2,3.
\end{align*}

\subsubsection{Calculations 1}
We use the following calculations in \eqref{scaling2}. Recall that $ \tilde{x} = {A}_{\widetilde{T}} \hat{x}$ with ${A}_{\widetilde{T}} = \widetilde {A} \widehat{A}$ and $x = {A}_{T} \tilde{x} + b_{T}$. For any multi-indices $\beta$ and $\gamma$, we have
\begin{align*}
\displaystyle
\partial^{\beta + \gamma}_{\hat{x}} &= \frac{\partial^{|\beta| + |\gamma|}}{\partial \hat{x}_1^{\beta_1} \cdots \partial \hat{x}_d^{\beta_d} \partial \hat{x}_1^{\gamma_1} \cdots \partial \hat{x}_d^{\gamma_d}} \notag\\
&\hspace{-1.5cm}  =  \underbrace{\sum_{i_1^{(1)},i_1^{(0,1)} = 1}^d h_1 [\widetilde{{A}}]_{i_1^{(1)} 1} [{A}_{T}]_{i_1^{(0,1)} i_1^{(1)}} \cdots   \sum_{i_{\beta_1}^{(1)},i_{\beta_1}^{(0,1)} = 1}^d h_1 [\widetilde{{A}}]_{i_{\beta_1}^{(1)} 1} [{A}_{T}]_{i_{\beta_1}^{(0,1)} i_{\beta_1}^{(1)}}   }_{\beta_1 \text{times}} \cdots \notag \\
&\hspace{-1.2cm} \underbrace{ \sum_{i_1^{(d)} , i_1^{(0,d)} = 1}^d h_d [\widetilde{{A}}]_{i_1^{(d)} d} [{A}_{T}]_{i_1^{(0,d)} i_1^{(d)}}  \cdots \sum_{i_{\beta_d}^{(d)} , i_{\beta_d}^{(0,d)}= 1}^d h_d [\widetilde{{A}}]_{i_{\beta_d}^{(d)} d} [{A}_{T}]_{i_{\beta_d}^{(0,d)} i_{\beta_d}^{(d)}}  }_{\beta_d \text{times}}  \notag \\
&\hspace{-1.2cm}  \underbrace{\sum_{j_1^{(1)} , j_1^{(0,1)} = 1}^d h_1 [\widetilde{{A}}]_{j_1^{(1)} 1} [{A}_{T}]_{j_1^{(0,1)} j_1^{(1)}}  \cdots  \sum_{j_{\gamma_1}^{(1)} , j_{\gamma_1}^{(0,1)}= 1}^d h_1 [\widetilde{{A}}]_{j_{\gamma_1}^{(1)} 1} [{A}_{T}]_{j_{\gamma_1}^{(0,1)} j_{\gamma_1}^{(1)}} }_{\gamma_1 \text{times}}  \cdots   \notag \\
&\hspace{-1.2cm}  \underbrace{\sum_{j_1^{(d)} , j_1^{(0,d)} = 1}^d h_d [\widetilde{{A}}]_{j_1^{(d)} d} [{A}_{T}]_{j_1^{(0,d)} j_1^{(d)}}  \cdots \sum_{j_{\gamma_d}^{(d)} , j_{\gamma_d}^{(0,d)}= 1}^d h_d [\widetilde{{A}}]_{j_{\gamma_d}^{(d)} d} [{A}_{T}]_{j_{\gamma_d}^{(0,d)} j_{\gamma_d}^{(d)}}  }_{\gamma_d \text{times}} \notag \\
&\hspace{-1.2cm}  \underbrace{\frac{\partial^{\beta_1}}{\partial {x}_{i_1^{(0,1)}}^{} \cdots \partial {x}_{i_{\beta_1}^{(0,1)}}^{}}}_{\beta_1 \text{times}} \cdots \underbrace{\frac{\partial^{\beta_d}}{\partial {x}_{i_1^{(0,d)}}^{} \cdots \partial {x}_{i_{\beta_d}^{(0,d)}}^{}}}_{\beta_d \text{times}}  \underbrace{ \frac{\partial^{\gamma_1}}{ \partial {x}_{j_1^{(0,1)}}^{} \cdots \partial {x}_{j_{\gamma_1}^{(0,1)}}^{} } }_{\gamma_1 \text{times}} \cdots \underbrace{\frac{\partial^{\gamma_d}}{ \partial {x}_{j_1^{(0,d)}}^{} \cdots \partial {x}_{j_{\gamma_d}^{(0,d)}}^{} }}_{\gamma_d \text{times}}.
\end{align*}
Let $\hat{\varphi} \in \mathcal{C}^{2}(\widehat{T})$ with $\tilde{\varphi} = \hat{\varphi} \circ {\Phi}_{\widetilde{T}}^{-1}$ and ${\varphi} = \tilde{\varphi} \circ {\Phi}_{T}^{-1}$. Then, for $1 \leq i \leq d$, 
\begin{align*}
\displaystyle
\left| \frac{\partial \hat{\varphi}}{\partial \hat{x}_i  } \right|
&= \left|  \sum_{i_1^{(1)}=1}^d  \sum_{i_1^{(0,1)} =1}^d h_i  [\widetilde{{A}}]_{i_1^{(1)} i} [{A}_{T}]_{i_1^{(0,1)} i_1^{(1)}}  \frac{\partial \varphi}{\partial x_{i_1^{(0,1)}}^{}}  \right| \\
&= h_i   \left |  \sum_{i_1^{(1)}=1}^d  \sum_{i_1^{(0,1)}=1}^d  [{A}_{T}]_{i_1^{(0,1)} i_1^{(1)}} (\tilde{r}_i)_{i_1^{(1)}}  \frac{\partial \varphi}{\partial x_{i_1^{(0,1)}}^{}} \right| = h_i \left |  \frac{\partial \varphi}{\partial r_{i}^{}} \right| \\
&\leq h_i \| \widetilde{{A}} \|_{\max} \| {A}_{T} \|_{\max}  \sum_{i_1^{(1)}=1}^d  \sum_{i_1^{(0,1)}=1}^d \left | \frac{\partial \varphi}{\partial x_{i_1^{(0,1)}}^{}} \right|,
\end{align*}
and for $1 \leq i,j \leq d$, 
\begin{align*}
\displaystyle
\left| \frac{\partial^2 \hat{\varphi}}{\partial \hat{x}_i \partial \hat{x}_j } \right|
&= \Biggl | \sum_{i_1^{(1)}, j_1^{(1)}=1}^d  \sum_{i_1^{(0,1)}, j_1^{(0,1)}=1}^d h_i  h_j  [\widetilde{{A}}]_{i_1^{(1)} i} [\widetilde{{A}}]_{j_1^{(1)} j}  \\
&\quad \quad [{A}_{T}]_{i_1^{(0,1)} i_1^{(1)}} [{A}_{T}]_{j_1^{(0,1)} j_1^{(1)}}  \frac{\partial^2 \varphi}{\partial x_{i_1^{(0,1)}}^{} \partial x_{j_1^{(0,1)}}^{}} \Biggr | = h_i  h_j \left |  \frac{\partial^2 \varphi}{\partial r_{i}^{} \partial r_{j}^{}} \right| \\
&\leq h_i  h_j \sum_{j_1^{(1)} =1}^d | [\widetilde{{A}}]_{j_1^{(1)} j} |    \Biggl | \sum_{j_1^{(0,1)}=1}^d [{A}_{T}]_{j_1^{(0,1)} j_1^{(1)}}  \frac{\partial^2 \varphi}{\partial r_{i}^{} \partial x_{j_1^{(0,1)}}^{}} \Biggr | \\
&\leq h_i  h_j  \| \widetilde{{A}} \|_{\max}   \| {{A}}_{T} \|_{\max} \sum_{j_1^{(0,1)}=1}^d \Biggl | \frac{\partial^2 \varphi}{\partial r_{i}^{} \partial x_{j_1^{(0,1)}}^{}} \Biggr | \\
&\leq h_i  h_j  \| \widetilde{{A}} \|_{\max}^2   \| {{A}}_{T} \|_{\max}^2  \sum_{i_1^{(0,1)}, j_1^{(0,1)}=1}^d \left| \frac{\partial^2 \varphi}{\partial x_{i_1^{(0,1)}}^{} \partial x_{j_1^{(0,1)}}^{}} \right|.
\end{align*}

\subsubsection{Calculations 2}
We use the following calculations in \eqref{scaling3}. Recall that $ \tilde{x} = {A}_{\widetilde{T}} \hat{x}$ with ${A}_{\widetilde{T}} = \widetilde {A} \widehat{A}$. For any multi-indices $\beta$ and $\gamma$, we have
\begin{align*}
\displaystyle
\partial^{\beta + \gamma}_{\hat{x}} &= \frac{\partial^{|\beta| + |\gamma|}}{\partial \hat{x}_1^{\beta_1} \cdots \partial \hat{x}_d^{\beta_d} \partial \hat{x}_1^{\gamma_1} \cdots \partial \hat{x}_d^{\gamma_d}} \notag\\
&\hspace{-1.5cm}  =  \underbrace{\sum_{i_1^{(1)} = 1}^d h_1 [\widetilde{{A}}]_{i_1^{(1)} 1} \cdots   \sum_{i_{\beta_1}^{(1)} = 1}^d h_1 [\widetilde{{A}}]_{i_{\beta_1}^{(1)} 1}   }_{\beta_1 \text{times}} \cdots \underbrace{ \sum_{i_1^{(d)}  = 1}^d h_d [\widetilde{{A}}]_{i_1^{(d)} d}  \cdots \sum_{i_{\beta_d}^{(d)} = 1}^d h_d [\widetilde{{A}}]_{i_{\beta_d}^{(d)} d}  }_{\beta_d \text{times}}  \notag \\
&\hspace{-1.2cm}  \underbrace{\sum_{j_1^{(1)}  = 1}^d h_1 [\widetilde{{A}}]_{j_1^{(1)} 1} \cdots  \sum_{j_{\gamma_1}^{(1)}= 1}^d h_1 [\widetilde{{A}}]_{j_{\gamma_1}^{(1)} 1} }_{\gamma_1 \text{times}}  \cdots  \underbrace{\sum_{j_1^{(d)} = 1}^d h_d [\widetilde{{A}}]_{j_1^{(d)} d}  \cdots \sum_{j_{\gamma_d}^{(d)} = 1}^d h_d [\widetilde{{A}}]_{j_{\gamma_d}^{(d)} d}  }_{\gamma_d \text{times}} \notag \\
&\hspace{-1.2cm}  \underbrace{\frac{\partial^{\beta_1}}{\partial \tilde{x}_{i_1^{1)}} \cdots \partial \tilde{x}_{i_{\beta_1}^{(1)}}}}_{\beta_1 \text{times}} \cdots \underbrace{\frac{\partial^{\beta_d}}{\partial \tilde{x}_{i_1^{(d)}} \cdots \partial \tilde{x}_{i_{\beta_d}^{(d)}}}}_{\beta_d \text{times}} \underbrace{ \frac{\partial^{\gamma_1}}{ \partial \tilde{x}_{j_1^{(1)}} \cdots \partial \tilde{x}_{j_{\gamma_1}^{(1)}} } }_{\gamma_1 \text{times}} \cdots \underbrace{\frac{\partial^{\gamma_d}}{ \partial \tilde{x}_{j_1^{(d)}} \cdots \partial \tilde{x}_{j_{\gamma_d}^{(d)}} }}_{\gamma_d \text{times}}.
\end{align*}
Let $\hat{\varphi} \in \mathcal{C}^{2}(\widehat{T})$ with $\tilde{\varphi} = \hat{\varphi} \circ {\Phi}_{\widetilde{T}}^{-1}$. Then, for $1 \leq i \leq d$, 
\begin{align*}
\displaystyle
\left| \frac{\partial \hat{\varphi}}{\partial \hat{x}_i  } \right|
&\leq  \sum_{i_1^{(1)}=1}^d h_i  \left| [\widetilde{{A}}]_{i_1^{(1)} i}  \right|  \left|\frac{\partial \tilde{\varphi}}{\partial \tilde{x}_{i_1^{(1)}}}  \right| 
\leq 
\begin{cases}	
 h_i  \| \widetilde{{A}}\|_{\max}  \sum_{i_1^{(1)}=1}^d \left| \frac{\partial \tilde{\varphi}}{\partial \tilde{x}_{i_1^{(1)}}} \right|  \quad \text{or},\\
 c \sum_{i_1^{(1)}=1}^d  \widetilde{\mathscr{H}}_{i_1^{(1)}}  \left| \frac{\partial \tilde{\varphi}}{\partial \tilde{x}_{i_1^{(1)}}} \right|,
\end{cases}
\end{align*}
and for $1 \leq i,j \leq d$, 
\begin{align*}
\displaystyle
\left| \frac{\partial^2 \hat{\varphi}}{\partial \hat{x}_i \partial \hat{x}_j } \right|
&= \left | \sum_{i_1^{(1)}, j_1^{(1)}=1}^d h_i h_j  [\widetilde{{A}}]_{i_1^{(1)} i} [\widetilde{{A}}]_{j_1^{(1)} j}  \frac{\partial^2 \tilde{\varphi}}{\partial \tilde{x}_{i_1^{(1)}} \partial \tilde{x}_{j_1^{(1)}}} \right | \\
&\leq
\begin{cases}	
h_i h_j \| \widetilde{{A}}\|_{\max}^2  \sum_{i_1^{(1)}, j_1^{(1)}=1}^d \left |  \frac{\partial^2 \tilde{\varphi}}{\partial \tilde{x}_{i_1^{(1)}} \partial \tilde{x}_{j_1^{(1)}}} \right | \quad \text{or}, \\
h_j  \sum_{j_1^{(1)}=1}^d  |[\widetilde{{A}}]_{j_1^{(1)} j}|   \left|  \sum_{i_1^{(1)} = 1}^d  h_i [\widetilde{{A}}]_{i_1^{(1)} i} \frac{\partial^2 \tilde{\varphi}}{\partial \tilde{x}_{i_1^{(1)}} \partial \tilde{x}_{j_1^{(1)}}}   \right|  \\
\quad \leq c h_j \| \widetilde{{A}} \|_{\max} \sum_{j_1^{(1)}=1}^d \sum_{i_1^{(1)} = 1}^d  \widetilde{\mathscr{H}}_{i_1^{(1)}} \left | \frac{\partial^2 \tilde{\varphi}}{\partial \tilde{x}_{i_1^{(1)}} \partial \tilde{x}_{j_1^{(1)}}}   \right| \quad \text{or}, \\
c \sum_{i_1^{(1)} = 1}^d \sum_{j_1^{(1)}=1}^d  \widetilde{\mathscr{H}}_{i_1^{(1)}}  \widetilde{\mathscr{H}}_{j_1^{(1)}} \left | \frac{\partial^2 \tilde{\varphi}}{\partial \tilde{x}_{i_1^{(1)}} \partial \tilde{x}_{j_1^{(1)}}}   \right|.
\end{cases}
\end{align*}

\subsubsection{Calculations 3} \label{sec=cal3}
We use the following calculations in \eqref{scaling1}.  Recall that $\hat{x} = {A}_{\widetilde{T}}^{-1} \tilde{x}$ with ${A}_{\widetilde{T}} = \widetilde {A} \widehat{A}$. For any multi-indices $\beta$, we have
\begin{align*}
\displaystyle
\partial^{\beta}_{\tilde{x}} &= \frac{\partial^{|\beta|} }{\partial \tilde{x}_1^{\beta_1} \cdots \partial \tilde{x}_d^{\beta_d} } \notag\\
&\hspace{-1.5cm}  =  \underbrace{\sum_{i_1^{(1)} = 1}^d h_{i_1^{(1)}}^{-1} [\widetilde{{A}}^{-1}]_{i_1^{(1)} 1} \cdots   \sum_{i_{\beta_1}^{(1)} = 1}^d h_{i_{\beta_1}^{(1)}}^{-1} [\widetilde{{A}}^{-1}]_{i_{\beta_1}^{(1)} 1}   }_{\beta_1 \text{times}} \cdots \underbrace{ \sum_{i_1^{(d)}  = 1}^d h_{i_1^{(d)}}^{-1} [\widetilde{{A}}^{-1}]_{i_1^{(d)} d}  \cdots \sum_{i_{\beta_d}^{(d)} = 1}^d h_{i_{\beta_d}^{(d)}}^{-1} [\widetilde{{A}}^{-1}]_{i_{\beta_d}^{(d)} d}  }_{\beta_d \text{times}}  \notag \\
&\hspace{-1.2cm}  \underbrace{\frac{\partial^{\beta_1}}{\partial \hat{x}_{i_1^{1)}} \cdots \partial \hat{x}_{i_{\beta_1}^{(1)}}}}_{\beta_1 \text{times}} \cdots \underbrace{\frac{\partial^{\beta_d}}{\partial \hat{x}_{i_1^{(d)}} \cdots \partial \hat{x}_{i_{\beta_d}^{(d)}}}}_{\beta_d \text{times}}.
\end{align*}
Let $\tilde{\varphi} \in \mathcal{C}^{2}(\widetilde{T})$ with $\hat{\varphi} = \tilde{\varphi} \circ {\Phi}_{\widetilde{T}}$. Then, for $1 \leq i \leq d$, 
\begin{align*}
\displaystyle
\left| \frac{\partial \tilde{\varphi}}{\partial \tilde{x}_i  } \right|
&\leq \sum_{i_1^{(1)} = 1}^d h_{i_1^{(1)}}^{-1} \left| [\widetilde{{A}}^{-1}]_{i_1^{(1)} i} \right|  \left|\frac{\partial \hat{\varphi}}{\partial \hat{x}_{i_1^{(1)}}}\right| 
\leq  \| \widetilde{{A}}^{-1}\|_{\max}  \sum_{i_1^{(1)} = 1}^d h_{i_1^{(1)}}^{-1} \left|\frac{\partial \hat{\varphi}}{\partial \hat{x}_{i_1^{(1)}}}\right|,
\end{align*}
and for $1 \leq i,j \leq d$, 
\begin{align*}
\displaystyle
\left| \frac{\partial^2 \tilde{\varphi}}{\partial \tilde{x}_i \partial \tilde{x}_j } \right|
&= \left | \sum_{i_1^{(1)}, j_1^{(1)}=1}^d h^{-1}_{i_1^{(1)}} h^{-1}_{j_1^{(1)}}  [\widetilde{{A}}^{-1}]_{i_1^{(1)} i} [\widetilde{{A}}^{-1}]_{j_1^{(1)} j}  \frac{\partial^2 \hat{\varphi}}{\partial \hat{x}_{i_1^{(1)}} \partial \hat{x}_{j_1^{(1)}}} \right | \\
&\leq  \| \widetilde{{A}}^{-1}\|_{\max}^2 \sum_{i_1^{(1)}, j_1^{(1)}=1}^d h^{-1}_{i_1^{(1)}} h^{-1}_{j_1^{(1)}} \left| \frac{\partial^2 \hat{\varphi}}{\partial \hat{x}_{i_1^{(1)}} \partial \hat{x}_{j_1^{(1)}}} \right|.
\end{align*}

\subsection{Main Results}
\begin{lem} \label{lem541}
Let $m,\ell \in \mathbb{N}_0$ with $\ell \geq m$. Let  $\beta := (\beta_1,\ldots,\beta_d) \in \mathbb{N}_0^d$ and $\gamma := (\gamma_1,\ldots,\gamma_d) \in \mathbb{N}_0^d$ be multi-indices with $|\beta| = m$ and $|\gamma| = \ell - m$. Then, for any $\hat{\varphi} \in W^{m,p}(\widehat{T})$ with $\tilde{\varphi} = \hat{\varphi} \circ {\Phi}_{\widetilde{T}}^{-1}$, it holds that
\begin{subequations} \label{scaling1}
\begin{align}
\displaystyle
|\tilde{\varphi}|_{W^{m,p}(\widetilde{T})}
&\leq c |\det({A}_{\widetilde{T}})|^{\frac{1}{p}} \| \widetilde{{A}}^{-1} \|^m_2 \left( \sum_{|\beta| = m}  ( h^{- \beta} )^p \|  \partial_{\hat{x}}^{\beta} \hat{\varphi} \|^p_{L^p(\widehat{T})} \right)^{\frac{1}{p}} \quad \text{if $p \in [1,\infty)$}, \label{scaling1a} \\
|\tilde{\varphi}|_{W^{m,\infty}(\widetilde{T})} 
&\leq c \| \widetilde{{A}}^{-1} \|^m_2 \max_{|\beta| = m} \left(  h^{- \beta}  \|  \partial^{\beta}_{\hat{x}} \hat{\varphi} \|_{L^{\infty}(\widehat{T})} \right) \quad \text{if $p = \infty$}. \label{scaling1b}
\end{align}
\end{subequations}

Let $p \in [0,\infty]$. Furthermore, for any $\hat{\varphi} \in W^{\ell,p}(\widehat{T})$ with $\tilde{\varphi} = \hat{\varphi} \circ {\Phi}_{\widetilde{T}}^{-1}$ and ${\varphi} = \tilde{\varphi} \circ {\Phi}_{T}^{-1}$, it holds that
\begin{align}
\displaystyle
\| \partial_{\hat{x}}^{\beta} \partial_{\hat{x}}^{\gamma}  \hat{\varphi}\|_{L^p(\widehat{T})}
&\leq c |\det({A}_{\widetilde{T}})|^{-\frac{1}{p}} \| \widetilde{{A}}  \|_2^{m} h^{\beta} \sum_{|\epsilon| =  |\gamma|} {h}^{\varepsilon} | \partial_r^{\epsilon} \varphi |_{W^{m,p}(T)}. \label{scaling2}
\end{align}
In particular, if Condition \ref{Cond333} is imposed, then for any $\hat{\varphi} \in W^{\ell,p}(\widehat{T})$ with $\tilde{\varphi} = \hat{\varphi} \circ {\Phi}_{\widetilde{T}}^{-1}$, it holds that
\begin{align}
\displaystyle
\| \partial_{\hat{x}}^{\beta} \partial_{\hat{x}}^{\gamma} \hat{\varphi} \|_{L^p(\widehat{T})}
&\leq c |\det({A}_{\widetilde{T}})|^{-\frac{1}{p}} \| \widetilde{{A}} \|^m_2 h^{\beta} \sum_{|\epsilon| =  |\gamma|} \widetilde{\mathscr{H}}^{\varepsilon} | \partial_{\tilde{x}}^{\varepsilon} \tilde{\varphi} |_{W^{m,p}(\widetilde{T})}.  \label{scaling3}
\end{align}
Here, for $p = \infty$ and any positive real $x$, $x^{- \frac{1}{p}} = 1$.
\end{lem}

\begin{pf*}
We divide the proof into three parts.

\textbf{Proof of \eqref{scaling1}.} 
Let $p \in [1,\infty)$. Because the space $\mathcal{C}^{m}(\widehat{T})$ is dense in the space ${W}^{m,p}(\widehat{T})$, we show \eqref{scaling1} for $\hat{\varphi} \in \mathcal{C}^{m}(\widehat{T})$ with $\tilde{\varphi} = \hat{\varphi} \circ {\Phi}_{\widetilde{T}}^{-1}$. Through the calculation (Section \ref{sec=cal3}) and \eqref{Anorm}, we have for any multi-index $\gamma$ with $|\gamma| = m$,
\begin{align*}
\displaystyle
|\partial^{\gamma}_{\tilde{x}} \tilde{\varphi}|
&\leq c  \| \widetilde{{A}}^{-1}\|_{2}^m \sum_{|\beta| = m} h^{- \beta } |\partial_{\hat{x}}^{\beta} \hat{\varphi}|.
\end{align*}
Through a change in a variable, we obtain
\begin{align*}
\displaystyle
|\tilde{\varphi}|_{W^{m,p}(\widetilde{T})}^p
&= \sum_{|\gamma| = m}  \| \partial_{\tilde{x}}^{\gamma} \tilde{\varphi} \|^p_{L^p(\widetilde{T})} 
\leq c |\det ({A}_{\widetilde{T}})|   \| \widetilde{{A}}^{-1}\|_{2}^{mp}  \sum_{|\beta| = m}   (h^{- \beta } )^p \|  \partial_{\hat{x}}^{\beta} \hat{\varphi} \|^p_{L^p(\widehat{T})},
\end{align*}
which leads to the inequality \eqref{scaling1a}. We consider the case that $p = \infty$. A function $\hat{\varphi} \in W^{m,\infty}(\widehat{T})$ belongs to the space $W^{m,p}(\widehat{T})$ for any $p \in [1,\infty)$.  It therefore holds that $\tilde{\varphi} \in W^{m,p}(\widetilde{T})$ for any $p \in [1,\infty)$ and, from \eqref{jensen},
\begin{align}
\displaystyle
  \| \partial_{\tilde{x}}^{\gamma} \tilde{\varphi} \|_{L^p(\widetilde{T})} 
  &\leq |  \tilde{\varphi} |_{W^{|\gamma|,p}(\widetilde{T})} \notag \\
  &\leq c |\det({A}_{\widetilde{T}})|^{\frac{1}{p}} \| \widetilde{{A}}^{-1} \|^m_2 \left( \sum_{|\beta| = |\gamma|}  ( h^{- \beta} )^p \|  \partial_{\hat{x}}^{\beta} \hat{\varphi} \|^p_{L^p(\widehat{T})} \right)^{\frac{1}{p}} \notag \\
  &\leq c \left( \sup_{1 \leq p} |\det({A}_{\widetilde{T}})|^{\frac{1}{p}} \right)  \| \widetilde{{A}}^{-1} \|^m_2 \sum_{|\beta| = |\gamma|} h^{- \beta}  \|  \partial_{\hat{x}}^{\beta} \hat{\varphi} \|_{L^p(\widehat{T})} \notag\\
    &\leq c \left( \sup_{1 \leq p} |\det({A}_{\widetilde{T}})|^{\frac{1}{p}} \right)  \| \widetilde{{A}}^{-1} \|^m_2 \sum_{|\beta| = |\gamma|} h^{- \beta}  \|  \partial_{\hat{x}}^{\beta} \hat{\varphi} \|_{L^{\infty}(\widehat{T})} \< + \infty,  \label{infty31}
\end{align}
for multi-index $\gamma \in \mathbb{N}_0^d$ with $|\gamma| \leq m$. This implies that the function $\partial_{\tilde{x}}^{\gamma} \tilde{\varphi}$ is in the space $L^{\infty}(\widetilde{T})$ for each $|\gamma| \leq m$. We therefore have $\tilde{\varphi} \in W^{m,\infty}(\widetilde{T})$. By passing to the limit $p \to \infty$ in \eqref{infty31} and because $\lim_{p \to \infty} \| \cdot \|_{L^p(\widetilde{T})} = \| \cdot\|_{L^{\infty}(\widetilde{T})}$, we have
\begin{align*}
\displaystyle
| \tilde{\varphi} |_{W^{m,\infty}(\widetilde{T})} \leq c \| \widetilde{{A}}^{-1} \|^m_2  \max_{|\beta| = m} \left( h^{- \beta}  \|  \partial_{\hat{x}}^{\beta} \hat{\varphi} \|_{L^{\infty}(\widehat{T})} \right), 
\end{align*}
which is \eqref{scaling1b}.

\textbf{Proof of \eqref{scaling3}.}
Let $\varepsilon = (\varepsilon_1,\ldots,\varepsilon_d) \in \mathbb{N}_0^d$ and $\delta = (\delta_1,\ldots,\delta_d) \in \mathbb{N}_0^d$ be multi-indies with $|\varepsilon| = |\gamma|$ and $|\delta| = |\beta|$. Let $p \in [1,\infty)$. Because the space $\mathcal{C}^{\ell}(\widehat{T})$ is dense in the space ${W}^{\ell,p}(\widehat{T})$, we show \eqref{scaling3} for $\hat{\varphi} \in \mathcal{C}^{\ell}(\widehat{T})$ with $\tilde{\varphi} = \hat{\varphi} \circ {\Phi}_{\widetilde{T}}^{-1}$. Through a simple calculation, we have
\begin{align}
\displaystyle
|\partial_{\hat{x}}^{\beta + \gamma} \hat{\varphi} | &= \left| \frac{\partial^{\ell} \hat{ \varphi}}{\partial \hat{x}_1^{\beta_1} \cdots \partial \hat{x}_d^{\beta_d} \partial \hat{x}_1^{\gamma_1} \cdots \partial \hat{x}_d^{\gamma_d}} \right| \notag\\
&\hspace{-1.5cm}  \leq c h^{\beta}  \| \widetilde{{A}} \|_{\max}^{|\beta|}  \underbrace{\sum_{i_1^{(1)} = 1}^d  \cdots   \sum_{i_{\beta_1}^{(1)} = 1}^d }_{\beta_1 \text{times}} \cdots \underbrace{ \sum_{i_1^{(d)} = 1}^d  \cdots \sum_{i_{\beta_d}^{(d)} = 1}^d  }_{\beta_d \text{times}}  \underbrace{\sum_{j_1^{(1)} = 1}^d  \cdots  \sum_{j_{\gamma_1}^{(1)} = 1}^d  }_{\gamma_1 \text{times}}  \cdots   \underbrace{\sum_{j_1^{(d)} = 1}^d \cdots \sum_{j_{\gamma_d}^{(d)} = 1}^d }_{\gamma_d \text{times}} \notag \\
&\hspace{-1.2cm}  \underbrace{\widetilde{\mathscr{H}}_{j_1^{(1)}} \cdots \widetilde{\mathscr{H}}_{j_{\varepsilon_1}^{(1)}}}_{\gamma_1 \text{times}} \cdots  \underbrace{ \widetilde{\mathscr{H}}_{j_1^{(d)}} \cdots \widetilde{\mathscr{H}}_{j_{\varepsilon_d}^{(d)}}}_{\gamma_d \text{times}} \notag \\
&\hspace{-1.2cm} \Biggl | \underbrace{\frac{\partial^{\beta_1}}{\partial \tilde{x}_{i_1^{(1)}} \cdots \partial \tilde{x}_{i_{\beta_1}^{(1)}}}}_{\beta_1 \text{times}} \cdots \underbrace{\frac{\partial^{\beta_d}}{\partial \tilde{x}_{i_1^{(d)}} \cdots \partial \tilde{x}_{i_{\beta_d}^{(d)}}}}_{\beta_d \text{times}} \underbrace{ \frac{\partial^{\gamma_1}}{ \partial \tilde{x}_{j_1^{(1)}} \cdots \partial \tilde{x}_{j_{\gamma_1}^{(1)}} } }_{\gamma_1 \text{times}} \cdots \underbrace{\frac{\partial^{\gamma_d}}{ \partial \tilde{x}_{j_1^{(d)}} \cdots \partial \tilde{x}_{j_{\gamma_d}^{(d)}} }}_{\gamma_d \text{times}} \tilde{\varphi} \Biggr |  \notag \\
&\hspace{-1.5cm}  \leq c h^{\beta}  \| \widetilde{{A}} \|_{\max}^{|\beta|}   \sum_{|\delta| = |\beta|} \sum_{|\varepsilon| = |\gamma|} \widetilde{\mathscr{H}}^{\varepsilon}  | { \partial_{\tilde{x}}^{\delta} \partial_{\tilde{x}}^{\varepsilon} \tilde{\varphi}}  |. \notag
\end{align}
We then have, using \eqref{Anorm},
\begin{align*}
\displaystyle
\int_{\widehat{T}} | \partial_{\hat{x}}^{\beta} \partial_{\hat{x}}^{ \gamma} \hat{\varphi}|^p d \hat{x}
&\leq c  \| \widetilde{{A}}  \|_2^{mp} h^{\beta p}  \sum_{|\delta| = |\beta|} \sum_{|\epsilon| = |\gamma|} \widetilde{\mathscr{H}}^{\varepsilon p} \int_{\widehat{T}} | { \partial_{\tilde{x}}^{\delta} \partial_{\tilde{x}}^{\varepsilon} \tilde{\varphi}}  |^p d \hat{x} \\
&= c |\det({A}_{\widetilde{T}})|^{-1} \| \widetilde{{A}}  \|_2^{mp} h^{\beta p}  \sum_{|\delta| = |\beta|} \sum_{|\epsilon| = |\gamma|} \widetilde{\mathscr{H}}^{\varepsilon p} \int_{\widetilde{T}} |{\partial_{\tilde{x}}^{\delta} \partial_{\tilde{x}}^{\varepsilon} \tilde{\varphi}}  |^p d \tilde{x}.
\end{align*}
Therefore, using \eqref{jensen}, we have
\begin{align*}
\displaystyle
\|  \partial_{\hat{x}}^{\beta} \partial_{\hat{x}}^{\gamma} \hat{\varphi}\|_{L^p(\widehat{T})}
&\leq c |\det({A}_{\widetilde{T}}))|^{-\frac{1}{p}} \| \widetilde{{A}}  \|_2^{m} h^{\beta} \sum_{|\epsilon| =  |\gamma|} \widetilde{\mathscr{H}}^{\varepsilon} | \partial_{\tilde{x}}^{\epsilon} \tilde{\varphi} |_{W^{m,p}(\widetilde{T})},
\end{align*}
which concludes \eqref{scaling3}. We consider the case that $p = \infty$. A function $\varphi \in W^{\ell,\infty}(T)$ belongs to the space $W^{\ell,p}(T)$ for any $p \in [1,\infty)$. It therefore holds that $\hat{\varphi} \in W^{\ell,p}(\widehat{T})$ for any $p \in [1,\infty)$ and thus
\begin{align}
\displaystyle
\|  \partial_{\hat{x}}^{\beta} \partial_{\hat{x}}^{\gamma} \hat{\varphi}\|_{L^p(\widehat{T})}
&\leq c |\det({A}_{\widetilde{T}})|^{-\frac{1}{p}} \| \widetilde{{A}}  \|_2^{m} h^{\beta} \sum_{|\epsilon| =  |\gamma|} \widetilde{\mathscr{H}}^{\varepsilon} | \partial_{\tilde{x}}^{\epsilon} \tilde{\varphi} |_{W^{m,p}(\widetilde{T})} \notag \\
&\leq c \| \widetilde{{A}}  \|_2^{m} h^{\beta} \sum_{|\epsilon| =  |\gamma|} \widetilde{\mathscr{H}}^{\varepsilon} | \partial_{\tilde{x}}^{\epsilon} \tilde{\varphi} |_{W^{m,\infty}(\widetilde{T})} \< \infty. \label{infty32}
\end{align}
This implies that the function $\partial_{\hat{x}}^{\beta} \partial_{\hat{x}}^{\gamma} \hat{\varphi}$ is in the space $L^{\infty}(\widehat{T})$.
Inequality  \eqref{scaling3} for $p=\infty$ is obtained by passing to the limit $p \to \infty$ in \eqref{infty32} on the basis that $\lim_{p \to \infty} \| \cdot \|_{L^p(\widehat{T})} = \| \cdot\|_{L^{\infty}(\widehat{T})}$.

\textbf{Proof of \eqref{scaling2}.} 
We follow the proof of \eqref{scaling3}. Let $p \in [1,\infty)$. Because the space $\mathcal{C}^{\ell}(\widehat{T})$ is dense in the space ${W}^{\ell,p}(\widehat{T})$, we show \eqref{scaling2} for $\hat{\varphi} \in \mathcal{C}^{\ell}(\widehat{T})$ with $\tilde{\varphi} = \hat{\varphi} \circ {\Phi}_{\widetilde{T}}^{-1}$ and ${\varphi} = \tilde{\varphi} \circ {\Phi}_{T}^{-1}$,  it holds that, for $1 \leq i,k \leq d$, 
\begin{align*}
\displaystyle
\left| \partial_{\hat{x}}^{\beta + \gamma} \hat{\varphi} \right| 
&\leq c h^{\beta}  \| \widetilde{{A}} \|_{\max}^{|\beta|}  \| {A}_{T} \|_{\max}^{|\beta|} \sum_{|\delta| = |\beta|}  \sum_{|\varepsilon| = |\gamma|} h^{\varepsilon}  \left|  \partial_x^{\delta}  \partial_{r^{}}^{\varepsilon} \varphi  \right|.
\end{align*}
Using \eqref{CN331c} and \eqref{Anorm}, we obtain \eqref{scaling2} for $p \in [1,\infty]$ by an argument analogous to the proof of \eqref{scaling3}.
\qed
\end{pf*}

\begin{rem}
In inequality \eqref{scaling3}, it is possible to obtain the estimates in $T$ by specifically determining the matrix $\mathcal{A}_{T}$.

Let $\ell=2$, $m=1$ and $p=q=2$. Recall that
\begin{align*}
\displaystyle
\Phi_{T}: \widetilde{T} \ni \tilde{x} \mapsto x = \Phi_{T}(\tilde{x}) = {A}_{T} \tilde{x} + b_{T} \in T.
\end{align*}
For ${\tilde\varphi} \in \mathcal{C}^{2}(\widetilde{T})$ with ${\varphi} = \tilde{\varphi} \circ {\Phi}_{T}^{-1}$ and $1 \leq i,j \leq d$, 
we have
\begin{align*}
\displaystyle
\left| \frac{\partial^2 \tilde{\varphi}}{\partial \tilde{x}_i  {\partial \tilde{x}_j}} ({x}) \right|
&= \left| \sum_{i_1^{(1)} , j_1^{(1)}=1}^2  [{A}_{T}]_{i_1^{(1)} i} [{A}_{T}]_{j_1^{(1)} j} \frac{\partial^2 \varphi}{\partial x_{i_1^{(1)}} \partial x_{j_1^{(1)}} } (x)  \right|.
\end{align*}

Let $d=2$. We define the matrix ${A}_{T}$ as 
\begin{align*}
\displaystyle
{A}_{T} := 
\begin{pmatrix}
\cos \frac{\pi}{2}  & - \sin \frac{\pi}{2} \\
 \sin \frac{\pi}{2} & \cos \frac{\pi}{2}
\end{pmatrix}.
\end{align*}
Because $ \| {A}_{T} \|_{\max} = 1$, we have
\begin{align*}
\displaystyle
\left| \frac{\partial^2 \tilde{\varphi}}{\partial \tilde{x}_i  {\partial \tilde{x}_j}} ({x}) \right|
&\leq \left|  \frac{\partial^2 \varphi}{\partial x_{i+1}^{} \partial x_{j+1}^{}} (x)  \right|,
\end{align*}
where the indices $i$, $i+1$ and $j$, $j+1$ have to be understood mod 2. Because $|\det ({A}_{T})| = 1$, it holds that
\begin{align*}
\displaystyle
\left \| \frac{\partial^2 \tilde{\varphi}}{\partial \tilde{x}_i  {\partial \tilde{x}_j}} \right \|_{L^2(\widetilde{T})} \leq  \left \|  \frac{\partial^2 \varphi}{\partial x_{i_{i+1}}^{} \partial x_{j+1}^{}}  \right \|_{L^2(T)}.
\end{align*}
We then have
\begin{align*}
\displaystyle
\sum_{j = 1}^2 \widetilde{\mathscr{H}}_j \left | \frac{\partial {\tilde{\varphi}}}{ {\partial \tilde{x}_j}} \right |_{H^1(\widetilde{T})}
&\leq \sum_{j = 1}^2 \widetilde{\mathscr{H}}_j \left |  \frac{\partial \varphi}{ \partial x_{j+1}^{}}  \right |_{H^1(T)},
\end{align*}
where the indices $j$, $j+1$ have to be understood mod 2. 

We define the matrix ${A}_{T}$ as 
\begin{align*}
\displaystyle
{A}_{T} := 
\begin{pmatrix}
\cos \frac{\pi}{4}  & - \sin \frac{\pi}{4} \\
 \sin \frac{\pi}{4} & \cos \frac{\pi}{4}
\end{pmatrix}.
\end{align*}
We then have
\begin{align*}
\displaystyle
\left| \frac{\partial^2 \tilde{\varphi}}{\partial \tilde{x}_i  {\partial \tilde{x}_j}} ({x}) \right|
&\leq \frac{1}{\sqrt{2}}  \sum_{i_1^{(1)} , j_1^{(1)}=1}^2  \left|  \frac{\partial^2 \varphi}{\partial x_{i_1^{(1)}}^{} \partial x_{j_1^{(1)}}^{}} (x)  \right|,
\end{align*}
which leads to
\begin{align*}
\displaystyle
\left \| \frac{\partial^2 \tilde{\varphi}}{\partial \tilde{x}_i  {\partial \tilde{x}_j}} \right \|_{L^2(\widetilde{T})}^2 \leq  c \sum_{i_1^{(1)} , j_1^{(1)}=1}^2 \left \|  \frac{\partial^2 \varphi}{\partial x_{i_1^{(1)}}^{} \partial x_{j_1^{(1)}}^{}}  \right \|_{L^2(T)}^2 \leq c |\varphi|^2_{H^2(T)}.
\end{align*}
We then have, using \eqref{jensen},
\begin{align*}
\displaystyle
\sum_{j = 1}^2 \widetilde{\mathscr{H}}_j \left | \frac{\partial {\varphi}^s}{ {\partial {x}^s_j}} \right |_{H^1(T^s)}
&\leq \sum_{j = 1}^2 \widetilde{\mathscr{H}}_j |\varphi|_{H^2(T)} \leq c h_{T} |\varphi|_{H^2(T)}.
\end{align*}
In this case, anisotropic interpolation error estimates cannot be obtained.
\end{rem}

\begin{rem}
We consider a general case. Let $p=q=2$. The space $\mathcal{C}^{1}(\widetilde{T})$ is dense in the space ${H}^{1}(\widetilde{T})$. For $\tilde{\varphi} \in \mathcal{C}^{1}(\widetilde{T})$ with $\varphi = \tilde{\varphi} \circ \Phi_{T}^{-1}$ and $1 \leq i \leq d$, we have
\begin{align*}
\displaystyle
\left| \frac{\partial \tilde{\varphi}}{\partial \tilde{x}_i} (\tilde{x}) \right|
&= \left| \sum_{i_1^{(1)}=1}^d  [A_{T}]_{i_1^{(1)} i} \frac{\partial \varphi}{\partial x_{i_1^{(1)}}} (x)  \right|.
\end{align*}
Let $d=2$.  We define a rotation matrix $A_{T}$ as
\begin{align*}
\displaystyle
A_{T} := 
\begin{pmatrix}
\cos \theta  & - \sin \theta \\
 \sin \theta & \cos \theta
\end{pmatrix},
\end{align*}
where $\theta$ denotes the angle. We then have
\begin{align*}
\displaystyle
\left| \frac{\partial \tilde{\varphi}}{\partial \tilde{x}_1} (\tilde{x}) \right|
&=  \left| \cos \theta  \frac{\partial \varphi}{\partial x_{1}} (x)  +  \sin \theta  \frac{\partial \varphi}{\partial x_{2}} (x) \right|, \\
\left| \frac{\partial \tilde{\varphi}}{\partial \tilde{x}_2} (\tilde{x}) \right|
&=  \left| - \sin \theta  \frac{\partial \varphi}{\partial x_{1}} (x)  +  \cos \theta  \frac{\partial \varphi}{\partial x_{2}} (x) \right|.
\end{align*}
If $|\sin \theta| \leq c \frac{\widetilde{\mathscr{H}}_2}{\widetilde{\mathscr{H}}_1}$ and $\widetilde{\mathscr{H}}_2 \leq c \widetilde{\mathscr{H}}_1$, we can deduce
\begin{align*}
\displaystyle
\left| \frac{\partial \tilde{\varphi}}{\partial \tilde{x}_1} (\tilde{x}) \right|
&\leq \left| \frac{\partial \varphi}{\partial x_{1}} (x) \right| + c \frac{\widetilde{\mathscr{H}}_2}{\widetilde{\mathscr{H}}_1} \left| \frac{\partial \varphi}{\partial x_{2}} (x) \right| , \\
\left| \frac{\partial \tilde{\varphi}}{\partial \tilde{x}_2} (\tilde{x}) \right|
&\leq c \left| \frac{\partial \varphi}{\partial x_{1}} (x) \right| +  \left| \frac{\partial \varphi}{\partial x_{2}} (x) \right|.
\end{align*}
As $|\det (A_{T})| = 1$, it holds that for $i=1,2$,
\begin{align*}
\displaystyle
\widetilde{\mathscr{H}}_i \left \| \frac{\partial \tilde{\varphi}}{\partial \tilde{x}_i} \right \|_{L^2(\widetilde{T})}
&\leq c \sum_{j=1}^2 \widetilde{\mathscr{H}}_j  \left \|  \frac{\partial \varphi}{\partial x_{{j}}}  \right \|_{L^2(T)}.
\end{align*}
\end{rem}

\begin{lem} \label{lem542}
Let $\Phi_{T}$ be the affine mapping defined in \eqref{aff=2}. Let $s \geq 0$ and $1 \leq p \leq \infty$. There exists positive constants $c_1$ and $c_2$ such that, for all $T \in \mathbb{T}_h$ and $\varphi \in W^{s,p}(T)$,
\begin{align}
\displaystyle
c_1  |\varphi |_{W^{s,p}(T)} \leq |\tilde{\varphi}|_{W^{s,p}(\widetilde{T})} &\leq c_ 2 | \varphi |_{W^{s,p}(T)}, \label{gint562}
\end{align}
with $\tilde{\varphi} = \varphi \circ {\Phi}_{T}$.

\end{lem}

\begin{pf*}
The following inequalities are found in \cite[Lemma 1.101]{ErnGue04}. There exists a positive constant $c$ such that, for all $T \in \mathbb{T}_h$ and $\varphi \in W^{s,p}(T)$,
\begin{align}
\displaystyle
|\tilde{\varphi}|_{W^{s,p}(\widetilde{T})} &\leq c \| {A}_{T} \|_2^s |\det ({A}_{T})|^{ - \frac{1}{p}} |\varphi|_{W^{s,p}(T)}, \label{gint563}\\
|\varphi|_{W^{s,p}(T)} &\leq c \| {A}_{T}^{-1} \|_2^s |\det ({A}_{T})|^{ \frac{1}{p}} |\tilde{\varphi}|_{W^{s,p}(\widetilde{T})}. \label{gint564}
\end{align}

Because the length of all edges of a simplex and measure of the simplex is not changed by a rotation and mirror imaging matrix and ${A}_{T}, {A}_{T}^{-1} \in O(d)$, 
\begin{align}
\displaystyle
&|\det ({A}_{T})| = \frac{|T|_d}{|\widetilde{T}|_d} = 1, \quad \| {A}_{T} \|_2 = 1, \quad  \| {A}_{T}^{-1} \|_2 =1.\label{gint565}
\end{align}
From \eqref{gint563}, \eqref{gint564}, and \eqref{gint565}, we obtain the desired inequality \eqref{gint562}.
\qed
\end{pf*}

\section{Classical Interpolation Error Estimates} 
\subsection{Local Interpolation Error Estimates}
The following theorem is another representation of the standard interpolation error estimates, e.g., see \cite[Theorem 1.103]{ErnGue04}.

\begin{thr} \label{thr521}
Let $1 \leq p \leq \infty$ and assume that there exists a nonnegative integer $k$ such that
\begin{align*}
\displaystyle
\mathbb{P}^{k} \subset \widehat{{P}} \subset W^{k+1,p}(\widehat{T}) \subset V(\widehat{T}).
\end{align*}
Let $\ell$ ($0 \leq \ell \leq k$) be such that $W^{\ell+1,p}(\widehat{T}) \subset V(\widehat{T})$ with continuous embedding. Furthermore, assume that $\ell, m \in \mathbb{N} \cup \{ 0 \}$ and $p , q \in [1,\infty]$ such that $0 \leq m \leq \ell + 1$ and 
\begin{align}
\displaystyle
W^{\ell +1,p}(\widehat{T}) \hookrightarrow W^{m,q} (\widehat{T}). \label{cla40}
\end{align}
It holds that, for any $m \in \{ 0 , \ldots, \ell+1 \}$ and any $\varphi \in W^{\ell+1,p}(T)$,
\begin{align}
\displaystyle
 |\varphi - {I}_{T} \varphi|_{W^{m,q}(T)}
&\leq C_*^I |T|_d^{\frac{1}{q} - \frac{1}{p}} \left( \frac{h_{\max}}{h_{\min}} \right)^m \left( \frac{H_{T}}{h_{T}} \right)^m h_{T}^{\ell+1-m} | \varphi |_{W^{\ell+1,p}(T)}, \label{ineq523}
\end{align}
where $C_*^I$ is a positive constant independent of $h_{T}$ and $H_{T}$, and the parameters $h_{\max}$ and $h_{\min}$ are defined by \eqref{hmin=hmax}, that is,
\begin{align*}
\displaystyle
h_{\max} = \max \{ h_1 , \ldots, h_d \}, \quad h_{\min} = \min \{ h_1 , \ldots, h_d \}.
\end{align*}
\end{thr}

\begin{pf*}
Let $\hat{\varphi} \in W^{\ell+1,p}(\widehat{T})$. Because $0 \leq \ell \leq k$, $\mathbb{P}^{\ell} \subset \mathbb{P}^k \subset \widehat{{P}}$. Therefore, for any $\hat{\eta} \in \mathbb{P}^{\ell}$, we have $ I_{\widehat{T}} \hat{\eta} = \hat{\eta}$. Using \eqref{int513} and \eqref{cla40}, we obtain
\begin{align*}
\displaystyle
| \tilde{\varphi} - I_{\widehat{T}} \hat{\varphi} |_{W^{m,q}(\widehat{T})}
&\leq | \tilde{\varphi} -\hat{\eta} |_{W^{m,q}(\widehat{T})} + |  I_{\widehat{T}} ( \hat{\eta} -  \hat{\varphi} ) |_{W^{m,q}(\widehat{T})} \\
&\leq c  \| \tilde{\varphi} -\hat{\eta} \|_{W^{\ell+1,p}(\widehat{T})},
\end{align*}
where we used the stability of the interpolation operator $I_{\widehat{T}}$, that is,
\begin{align*}
\displaystyle
|  I_{\widehat{T}} ( \hat{\eta} -  \hat{\varphi} ) |_{W^{m,q}(\widehat{T})} 
&\leq \sum_{i=1}^{n_0} |\hat{\chi}_i (\hat{\eta} -  \hat{\varphi})| |\hat{\theta}_i|_{W^{m,q}(\widehat{T})}  
\leq c \| \hat{\eta} - \tilde{\varphi} \|_{W^{\ell+1,p}(\widehat{T})}.
\end{align*}
Using the classic Bramble--Hilbert--type lemma (e.g., \cite[Lemma 4.3.8]{BreSco08}), we obtain
\begin{align}
\displaystyle
| \tilde{\varphi} - I_{\widehat{T}} \hat{\varphi} |_{W^{m,q}(\widehat{T})}
&\leq c \inf_{\hat{\eta} \in \mathbb{P}^{\ell}} \| \hat{\eta} - \tilde{\varphi} \|_{W^{\ell+1,p}(\widehat{T})} \leq c |\hat{\varphi}|_{W^{\ell+1,p}(\widehat{T})}. \label{cla43}
\end{align}

The inequalities \eqref{gint562}, \eqref{scaling1}, \eqref{jensen}, and \eqref{cla43} yield
\begin{align}
\displaystyle
&| \varphi - {I}_{T} \varphi |_{W^{m,q}({T})} \leq c | \tilde{\varphi} - {I}_{\widetilde{T}} \tilde{\varphi} |_{W^{m,q}(\widetilde{T})} \notag \\
&\ \leq c  |\det({A}_{\widetilde{T}})|^{\frac{1}{q}} \| \widetilde{{A}}^{-1} \|^m_2  \left( \sum_{|\beta| = m}  ( h^{- \beta} )^q \|  \partial^{\beta}  (\hat{\varphi} - I_{\widehat{T}} \hat{\varphi} )  \|^q_{L^q(\widehat{T})} \right)^{\frac{1}{q}} \notag \\
&\ \leq c |\det({A}_{\widetilde{T}})|^{\frac{1}{q}}  \| \widetilde{{A}}^{-1} \|_2^m \max \{ h_1^{-1} , \ldots, h_d^{-1} \}^{|\beta|} | \tilde{\varphi} - I_{\widehat{T}} \hat{\varphi} |_{W^{m,q}(\widehat{T})} \notag \\
&\ \leq c |\det({A}_{\widetilde{T}})|^{\frac{1}{q}} \| \widetilde{{A}}^{-1} \|_2^m  h_{\min}^{- |\beta|}  |\hat{\varphi}|_{W^{\ell+1,p}(\widehat{T})}. \label{cla44}
\end{align}
Using the inequalities \eqref{jensen}, \eqref{gint562} and \eqref{scaling2}, we have
\begin{align}
\displaystyle
|\hat{\varphi}|_{W^{\ell+1,p}(\widehat{T})}
&\leq  \sum_{|\gamma| = \ell + 1 -m} \sum_{|\beta| = m} \| \partial^{\beta} \partial^{\gamma} \hat{\varphi} \|_{L^{p}(\widehat{T})} \notag \\
&\leq c  |\det({A}_{\widetilde{T}})|^{-\frac{1}{p}} \| \widetilde{{A}}  \|_2^{m}  \sum_{|\gamma| = \ell + 1 -m} \sum_{|\beta| = m} h^{\beta} \sum_{|\epsilon| =  |\gamma|} {h}^{\varepsilon} | \partial_r^{\epsilon} \varphi |_{W^{m,p}(T)} \notag \\
&\leq c  |\det({A}_{\widetilde{T}})|^{-\frac{1}{p}} \| \widetilde{{{A}}}  \|_2^{m}  \max \{ h_1 , \ldots, h_d \}^{|\beta|} h_{T}^{\ell+1-m} |\varphi|_{W^{\ell+1,p} (T)} \notag \\
&\leq c  |\det({A}_{\widetilde{T}})|^{-\frac{1}{p}} \| \widetilde{{{A}}}  \|_2^{m} h_{\max}^{|\beta|}  h_{T}^{\ell+1-m} |\varphi|_{W^{\ell+1,p} (T)}.  \label{cla45}
\end{align}
From \eqref{cla44} and  \eqref{cla45} together with \eqref{CN331} and \eqref{CN332}, we have the desired estimate \eqref{ineq523}.
\qed	
\end{pf*}

\begin{rem}
We introduced the estimate \eqref{DekLev}, a variant of the Bramble--Hilbert lemma. However, because we prove estimate \eqref{cla43} with $p=q$ using the reference element,  it is sufficient to use the standard estimate  (e.g., \cite{DupSco80,BreSco08}) to achieve our goal.
\end{rem}

\begin{ex}
As the examples in \cite[Example 1.106]{ErnGue04}, we get local interpolation error estimates for a Lagrange finite element of degree $k$, a more general finite element, and the Crouzeix--Raviart finite element with $k=1$.
\begin{enumerate}
 \item For a Lagrange finite element of degree $k$, we set $V(\widehat{T}) := \mathcal{C}^0(\widehat{T})$. The condition on $\ell$ in Theorem \ref{thr521} is $\frac{d}{p} - 1 \< \ell \leq k$ because $W^{\ell+1,p}(\widehat{T}) \subset \mathcal{C}^0(\widehat{T})$ if $\ell+1 \> \frac{d}{p}$ according to the Sobolev imbedding theorem. 
 \item For a general finite element with $V(\widehat{T}) := \mathcal{C}^t(\widehat{T})$ and $t \in \mathbb{N}$. The condition on $\ell$ in Theorem \ref{thr521} is $\frac{d}{p} - 1 + t \< \ell \leq k$. When $t=1$, there is a Hermite finite element.
 \item For the Crouzeix--Raviart finite element with $k=1$, we set $V(\widehat{T}) := W^{1,1}(\widehat{T})$. The condition on $\ell$ in Theorem \ref{thr521} is $0 \leq \ell \leq 1$.
\end{enumerate}
\end{ex}

\subsection{Examples of Anisotropic Elements}
When $m=\ell = 1$ and $q=p$ in \eqref{ineq523} of Theorem \ref{thr521}, the estimate is written as
\begin{align}
\displaystyle
 |\varphi - {I}_{T} \varphi|_{W^{1,p}(T)}
&\leq C_*^I  \frac{h_{\max}}{h_{\min}} \frac{H_{T}}{h_{T}}  h_{T} | \varphi |_{W^{2,p}(T)}. \label{ineq523=ex}
\end{align}
Let $T \subset \mathbb{R}^2$ be a triangle. As described in Section \ref{iso=mesh=ex2d}, an isotropic mesh element has equal or nearly equal edge lengths and angles, resulting in a balanced shape. Then, the estimate \eqref{ineq523=ex} becomes
\begin{align}
\displaystyle
 |\varphi - {I}_{T} \varphi|_{W^{1,p}(T)}
&\leq c h_{T} | \varphi |_{W^{2,p}(T)}. \label{ineq523=ex2}
\end{align}

We introduce typical examples of the quantities $\frac{h_{\max} }{ h_{\min}}$ and $\frac{H_{T}}{h_{T}}$ in anisotropic elements. We considered the following five anisotropic elements as in Section \ref{ani=mesh=ex2d}:  Let $0 \< s \ll 1$, $s \in \mathbb{R}$ and $\varepsilon,\delta,\gamma \in \mathbb{R}$.

\begin{ex}[Right-angled triangle] \label{ex11=5}
Let $T \subset \mathbb{R}^2$ be the simplex with vertices $p_1 := (0,0)^{\top}$, $p_2 := (s,0)^{\top}$ and $p_3 := (0,s^{\varepsilon})^{\top}$ with $1 \< \varepsilon$. We then have $h_1 = s$, $h_2 =  s^{\varepsilon}$ and $h_T = \sqrt{s^2 + s^{2 \varepsilon}}$; i.e.,
\begin{align*}
\displaystyle
\frac{h_{\max}}{h_{\min}} &\leq s^{1 - \varepsilon} \to \infty \quad \text{as $s \to 0$}, \quad \frac{H_{T}}{h_{T}} =2.
\end{align*}
In this case, the estimate \eqref{ineq523=ex} becomes
\begin{align*}
\displaystyle
 |\varphi - {I}_{T} \varphi|_{W^{1,p}(T)}
&\leq 2 C_*^I  s^{2 - \varepsilon} | \varphi |_{W^{2,p}(T)}. 
\end{align*}
When $\varepsilon > 2$, this implies that the estimate diverges as $s \to 0$. However, new interpolation error estimates will be shown to converge, see Example \ref{ex14=3}.
\end{ex}

\begin{ex}[Dagger] \label{ex11=6}
Let $T \subset \mathbb{R}^2$ be the simplex with vertices $p_1 := (0,0)^{\top}$, $p_2 := (s,0)^{\top}$ and $p_3 := (s^{\delta},s^{\varepsilon})^{\top}$ with $1 \< \varepsilon \< \delta $. We then have $h_1 = \sqrt{(s - s^{\delta})^2 + s^{2 \varepsilon}}$, $h_2 = \sqrt{s^{2 \delta} + s^{2 \varepsilon}}$ and $h_T=s$; i.e.,
\begin{align*}
\displaystyle
\frac{h_{\max}}{h_{\min}} &= \frac{ \sqrt{(s - s^{\delta})^2 + s^{2 \varepsilon}}}{\sqrt{s^{2 \delta} + s^{2 \varepsilon}}} \leq c s^{1 - \varepsilon} \to \infty \quad \text{as $s \to 0$},\\
\frac{H_{T}}{h_{T}} &= \frac{ \sqrt{(s- s^{\delta})^2 + s^{2 \varepsilon}} \sqrt{s^{2 \delta} + s^{2 \varepsilon}}}{\frac{1}{2} s^{1 + \varepsilon}} \leq c.
\end{align*}
In this case, the estimate \eqref{ineq523=ex} becomes
\begin{align*}
\displaystyle
 |\varphi - {I}_{T} \varphi|_{W^{1,p}(T)}
&\leq c s^{2 - \varepsilon} | \varphi |_{W^{2,p}(T)}. 
\end{align*}
When $\varepsilon > 2$, this implies that the estimate diverges as $s \to 0$. However, new interpolation error estimates will be shown to converge, see Example \ref{ex14=4}.
\end{ex}

\begin{ex}[Blade] \label{ex11=7}
Let $T \subset \mathbb{R}^2$ be the simplex with vertices $p_1 := (0,0)^{\top}$, $p_2 := (2s,0)^{\top}$ and $p_3 := (s ,s^{\varepsilon})^{\top}$ with $1 \< \varepsilon $. We then have $h_1 = h_2 = \sqrt{s^{2} + s^{2 \varepsilon}}$ and $h_T = 2s$; i.e.,
\begin{align*}
\displaystyle
\frac{h_{\max}}{h_{\min}} = 1, \quad \frac{H_{T}}{h_{T}} = \frac{s^{2} + s^{2 \varepsilon}}{s^{1 + \varepsilon}} \to \infty \quad \text{as $s \to 0$}.
\end{align*}
In this case, the estimate \eqref{ineq523=ex} becomes
\begin{align*}
\displaystyle
 |\varphi - {I}_{T} \varphi|_{W^{1,p}(T)}
&\leq c s^{2 - \varepsilon} | \varphi |_{W^{2,p}(T)}. 
\end{align*}
When $\varepsilon > 2$, this implies that the estimate diverges as $s \to 0$. In this case, the interpolation error estimate can not be improved, , see Example \ref{ex14=5}.
\end{ex}

\begin{ex}[Dagger] \label{ex11=8}
Let $T \subset \mathbb{R}^2$ be the simplex with vertices $p_1 := (0,0)^{\top}$, $p_2 := (s,0)^{\top}$ and $p_3 := (s^{\delta},s^{\varepsilon})^{\top}$ with $1 \< \delta \< \varepsilon $. We then have $h_1 = \sqrt{(s - s^{\delta})^2 + s^{2 \varepsilon}}$, $h_2 = \sqrt{s^{2 \delta} + s^{2 \varepsilon}}$ and $h_T=s$; i.e.,
\begin{align*}
\displaystyle
\frac{h_{\max}}{h_{\min}} &= \frac{ \sqrt{(s - s^{\delta})^2 + s^{2 \varepsilon}}}{\sqrt{s^{2 \delta} + s^{2 \varepsilon}}} \leq c s^{1 - \delta} \to \infty \quad \text{as $s \to 0$},\\
\frac{H_{T}}{h_{T}} &= \frac{ \sqrt{(s- s^{\delta})^2 + s^{2 \varepsilon}} \sqrt{s^{2 \delta} + s^{2 \varepsilon}}}{\frac{1}{2} s^{1 + \varepsilon}} \leq c s^{\delta - \varepsilon} \to \infty \quad \text{as $s \to 0$}.
\end{align*}
In this case, the estimate \eqref{ineq523=ex} becomes
\begin{align*}
\displaystyle
 |\varphi - {I}_{T} \varphi|_{W^{1,p}(T)}
&\leq c s^{2 - \varepsilon} | \varphi |_{W^{2,p}(T)}. 
\end{align*}
When $\varepsilon > 2$, this implies that the estimate diverges as $s \to 0$. In this case, the interpolation error estimate can not be improved, see Example \ref{ex14=6}.
\end{ex}

\begin{ex}[Right-angled triangle] \label{Ad=Ex6}
Let $T \subset \mathbb{R}^2$ be the simplex with vertices $p_1 := (0,0)^{\top}$, $p_2 := (s,0)^{\top}$ and $p_3 := (0, \delta s)^{\top}$ with $\delta \ll 1$. We then have $h_1 = s$, $h_2 = \delta s$ and $h_T = s \sqrt{1 + \delta^2}$; i.e.,
\begin{align*}
\displaystyle
\frac{h_{\max}}{h_{\min}} &= \frac{1}{\delta} , \quad \frac{H_{T}}{h_{T}} =2.
\end{align*}
In this case, the estimate \eqref{ineq523=ex} becomes
\begin{align*}
\displaystyle
 |\varphi - {I}_{T} \varphi|_{W^{1,p}(T)}
&\leq \frac{c}{\delta} s | \varphi |_{W^{2,p}(T)}. 
\end{align*}
This implies that the estimate converges as $s \to 0$ and the error may be large. However, new interpolation error estimates remove the factor $\frac{1}{\delta}$, see Example \ref{ex14=7}.
\end{ex}

\begin{ex}[Blade] \label{Ad=Ex7}
Let $T \subset \mathbb{R}^2$ be the simplex with vertices $p_1 := (0,0)^{\top}$, $p_2 := (2s,0)^{\top}$ and $p_3 := (s , \delta s)^{\top}$ with $\delta \ll 1$. We then have $h_1 = h_2 = s \sqrt{1 + \delta^2}$ and $h_T = 2s$; i.e.,
\begin{align*}
\displaystyle
\frac{h_{\max}}{h_{\min}} = 1, \quad \frac{H_{T}}{h_{T}} = \frac{s^2 (1 + \delta^2)}{\delta s^2} \leq \frac{c}{\delta},
\end{align*}
In this case, the estimate \eqref{ineq523=ex} becomes
\begin{align*}
\displaystyle
 |\varphi - {I}_{T} \varphi|_{W^{1,p}(T)}
&\leq \frac{c}{\delta} s | \varphi |_{W^{2,p}(T)}. 
\end{align*}
This implies that the estimate converges as $s \to 0$ and the error may be large. Unfortunately, new interpolation error estimates do not remove the factor $\frac{1}{\delta}$, see Example \ref{ex14=8}.
\end{ex}

\begin{ex}[$\mathbb{P}^1$ + bubble finite element in $\mathbb{R}^2$] \label{ex526}
We give a numerical example which is not optimal in the usual sense. Let $T \subset \mathbb{R}^2$ be the triangle with vertices $p_1 := (0,0)^{\top}$, $p_2 := (s,0)^{\top}$, $p_3 := (0, s^{\varepsilon})^{\top}$ (Example \ref{ex11=5}), where $s := \frac{1}{N}$, $N \in \mathbb{N}$ and $\varepsilon \in \mathbb{R}$, $1 \< \varepsilon \leq 2$.  Let $p_{4}$ be the barycentre of  $T$. 

Using the barycentric coordinates $ \lambda_i: \mathbb{R}^2 \to \mathbb{R}$, $i=1,\ldots,3$, we define the local basis functions as
\begin{align*}
\displaystyle
\theta_4(x) := 27 \lambda_1(x) \lambda_2(x) \lambda_3(x), \quad \theta_i(x) := \lambda_i(x) - \frac{1}{3} \theta_4(x), \quad i=1,2,3.
\end{align*}
The interpolation operator $I_T^b$ defined by
\begin{align*}
\displaystyle
I_{T}^{b}:H^2(T) \ni \varphi  \mapsto I_{T}^{b} \varphi := \sum_{i=1}^{4} \varphi (x_{i}) \theta_i \in \Span \{ \theta_1, \theta_2, \theta_3 , \theta_4 \}.
\end{align*}
From Theorem \ref{thr521}, we have
\begin{align*}
\displaystyle
 |\varphi - {I}_{T}^b \varphi|_{H^{1}(T)}
&\leq c h_{T}^{2 - \varepsilon} | \varphi |_{H^{2}(T^s)} \quad \forall \varphi \in H^{2}(T).
\end{align*}
Let $\varphi$ be a function such that
\begin{align*}
\displaystyle
\varphi(x,y) := 2x^2 - xy + 3y^2.
\end{align*}
We compute the convergence order concerning the $H^1$ norm defined by
\begin{align*}
\displaystyle
&Err_s^{b}(H^1) := \frac{| \varphi - I_{T}^{b} \varphi |_{H^1(T)}}{| \varphi |_{H^2(T)}}, 
\end{align*}
for the cases: $\varepsilon = 1.5$ (Table \ref{modify=table1}) and $\varepsilon = 2.0$ (Table \ref{modify=table2}). The convergence indicator $r$ is defined by
\begin{align*}
\displaystyle
r = \frac{1}{\log(2)} \log \left( \frac{Err_t^{b}(H^1)}{Err_{t/2}^{b}(H^1)} \right).
\end{align*}

\begin{table}[htbp]
\caption{Error of the local interpolation operator ($\varepsilon = 1.5$)}
\centering
\begin{tabular}{l | l | l | l } \hline
$N$ &  $s$  & $Err_s^{b}(H^1)$ & $r$  \\ \hline \hline
128 & 7.8125e-03 & 2.9951e-02  &       \\
256 & 3.9062e-03  & 2.1101e-02 & 5.0529e-01   \\
512 & 1.9531e-03  & 1.4874e-02 & 5.0452e-01     \\
1024 & 9.7656e-04  & 1.0491e-02 &  5.0364e-01   \\
\hline
\end{tabular}
\label{modify=table1}
\end{table}

\begin{table}[htbp]
\caption{Error of the local interpolation operator ($\varepsilon = 2.0$)}
\centering
\begin{tabular}{l | l | l | l } \hline
$N$ &  $s$  & $Err_s^{b}(H^1)$ & $r$  \\ \hline \hline
128 & 7.8125e-03 & 3.3397e-01 &       \\
256 & 3.9062e-03  & 3.3366e-01& 1.3398e-03   \\
512 & 1.9531e-03  & 3.3350e-01 & 6.9198e-04     \\
1024 & 9.7656e-04  & 3.3341e-01 & 3.8939e-04    \\
\hline
\end{tabular}
\label{modify=table2}
\end{table}
\end{ex}

\begin{rem}
If we are concerned with anisotropic elements, it would be desirable to remove the quantity $h_{\max} / h_{\min}$  from estimate \eqref{ineq523}.	
\end{rem}

\section{Anisotropic Interpolation on the Reference Element} \label{sec53}
We introduce estimates on the reference element due to \cite{ApeDob92,Ape99} to obtain anisotropic interpolation error estimates.

For the reference element $\widehat{T}$ defined in Sections \ref{reference2d} and \ref{reference3d}, let the triple $\{ \widehat{T}, \widehat{P} , \widehat{\Sigma} \}$ be the reference finite element with associated normed vector space $V(\widehat{T})$.

\begin{thr} \label{thrApel}
Let ${I}_{\widehat{T}} : \mathcal{C}({\widehat{T}}) \to  \mathbb{P}^{k}(\widehat{T})$ be a linear operator. Fix $m,\ell \in \mathbb{N}$ and $p,q \in [1,\infty]$ such that $0 \leq m \leq \ell  \leq k+1$ and
\begin{align}
\displaystyle
W^{\ell - m,p}(\widehat{T}) \hookrightarrow L^q(\widehat{T}). \label{Sobolev511}
\end{align}
Let $\beta$ be a multi-index with $|\beta| = m$. We set $j := \dim (\partial_{\hat{x}}^{\beta} \mathcal{P}^{k})$.  Assume that there exist linear functionals $\mathscr{F}_i$, $i=1,\ldots,j$, such that
\begin{subequations} \label{functionals}
\begin{align}
\displaystyle
&\mathscr{F}_i \in W^{\ell  - m ,p} (\widehat{T})^{\prime},\quad \forall i = 1,\ldots,j, \label{functionals=a}\\
&\mathscr{F}_i ( \partial_{\hat{x}}^{\beta} (\hat{\varphi} - I_{\widehat{T}} \hat{\varphi}) ) = 0 \quad  \forall i = 1,\ldots,j, \quad \forall \hat{\varphi} \in \mathcal{C}({\widehat{T}}): \ \partial_{\hat{x}}^{\beta} \hat{\varphi} \in W^{\ell - m , p} (\widehat{T}), \label{functionals=b} \\
&\hat{\eta} \in \mathbb{P}^{k}, \quad \mathscr{F}_i(\partial_{\hat{x}}^{\beta} \hat{\eta}) = 0 \quad \forall i = 1,\ldots,j \quad \Rightarrow \quad \partial^{\beta} \hat{\eta} = 0. \label{functionals=c}
\end{align}
\end{subequations}
It holds that for all $\hat{\varphi} \in \mathcal{C}({\widehat{T}})$ with $\partial^{\beta} \hat{\varphi} \in W^{\ell  - m , p} (\widehat{T})$,
\begin{align}
\displaystyle
\| \partial_{\hat{x}}^{\beta} (\hat{\varphi} - I_{\widehat{T}} \hat{\varphi}) \|_{L^q(\widehat{T})} \leq C^{F} | \partial_{\hat{x}}^{\beta} \hat{\varphi} |_{W^{\ell  - m , p} (\widehat{T})}. \label{ineq512}
\end{align}
\end{thr}

\begin{pf*}
We follow \cite[Lemma 2.2]{Ape99}.

For all $\hat{\eta} \in \mathbb{P}^{\ell -1}$, we have
\begin{align}
\displaystyle
\| \partial_{\hat{x}}^{\beta} (\hat{\varphi} - I_{\widehat{T}} \hat{\varphi}) \|_{L^q(\widehat{T})}
\leq \| \partial_{\hat{x}}^{\beta} (\hat{\varphi} - \hat{\eta}) \|_{L^q(\widehat{T})} + \| \partial_{\hat{x}}^{\beta} (\hat{\eta} - I_{\widehat{T}} \hat{\varphi}) \|_{L^q(\widehat{T})}. \label{ineq513}
\end{align}
Note that $\hat{\eta} - I_{\widehat{T}} \hat{\varphi} \in  \mathbb{P}^{k}$, because $\ell  \leq k + 1$. That is, $\partial_{\hat{x}}^{\beta} ( \hat{\eta} - I_{\widehat{T}} \hat{\varphi} ) \in \partial_{\hat{x}}^{\beta} \mathbb{P}^{k}$. Because the polynomial spaces are finite-dimensional all norms are equivalent, that is, by the fact $\sum_{i=1}^j |\mathscr{F}_i(\hat{\eta})|$ is a norm on $\partial_{\hat{x}}^{\beta} \mathbb{P}^{k}$, together with \eqref{functionals=a}, \eqref{functionals=b} and \eqref{functionals=c}, we have for any $\hat{\eta} \in \mathbb{P}^{\ell -1}$,
\begin{align*}
\displaystyle
\| \partial_{\hat{x}}^{\beta} (\hat{\eta} - I_{\widehat{T}} \hat{\varphi}) \|_{L^q(\widehat{T})}
&\leq c \sum_{i=1}^j | \mathscr{F}_i (\partial_{\hat{x}}^{\beta} ( \hat{\eta} - I_{\widehat{T}} \hat{\varphi} ) )| 
= c \sum_{i=1}^j |\mathscr{F}_i (\partial_{\hat{x}}^{\beta} ( \hat{\eta} - \hat{\varphi} ) )| \\
&\leq c \| \partial_{\hat{x}}^{\beta} ( \hat{\eta} - \hat{\varphi} ) \|_{W^{\ell  - m ,p}(\widehat{T})}.
\end{align*}
Using  \eqref{ineq513} and \eqref{Sobolev511}, it holds that for any $\hat{\eta} \in \mathbb{P}^{\ell -1}$,
\begin{align*}
\displaystyle
\| \partial_{\hat{x}}^{\beta} (\hat{\varphi} - I_{\widehat{T}} \hat{\varphi}) \|_{L^q(\widehat{T})}
&\leq \| \partial_{\hat{x}}^{\beta} (\hat{\varphi} - \hat{\eta}) \|_{L^q(\widehat{T})} + \| \partial_{\hat{x}}^{\beta} (\hat{\eta} - I_{\widehat{T}} \hat{\varphi}) \|_{L^q(\widehat{T})} \\
&\leq c \| \partial_{\hat{x}}^{\beta} ( \hat{\eta} - \hat{\varphi} ) \|_{W^{\ell  - m ,p}(\widehat{T})}.
\end{align*}
By Lemma \ref{BH=1}, we have
\begin{align*}
\displaystyle
\| \partial_{\hat{x}}^{\beta} (\hat{\varphi} - I_{\widehat{T}} \hat{\varphi}) \|_{L^q(\widehat{T})}
&\leq c \inf_{\hat{\eta} \in \mathbb{P}^{\ell-1}} \| \partial_{\hat{x}}^{\beta} ( \hat{\eta} - \hat{\varphi} ) \|_{W^{\ell  - m ,p}(\widehat{T})} \\
&\leq c | \partial_{\hat{x}}^{\beta} \hat{\varphi} |_{W^{\ell  - m ,p}(\widehat{T})}.
\end{align*}
\qed
\end{pf*}

\begin{rem}
Note that it is not required $I_{\widehat{T}} \hat{\eta} = \hat{\eta}$ for any $\hat{\eta} \in \mathbb{P}^{\ell -1}$.
\end{rem}

\section{Remarks on Anisotropic Interpolation Analysis} \label{sec_mis}
Let $\widehat{T} \subset \mathbb{R}^2$ be the reference element defined in Section \ref{reference2d}. We set $k = m  = 1$, $\ell = 2$, and $p=2$. For $\hat{\varphi} \in H^2(\widehat{T})$, we set $\tilde{\varphi} = \hat{\varphi} \circ {\Phi}_{\widetilde{T}}^{-1}$ and ${\varphi} = \tilde{\varphi} \circ {\Phi}_{T}^{-1}$. Inequalities \eqref{scaling1} and \eqref{gint562} yield
\begin{align}
\displaystyle
|{\varphi} - I_T \varphi|_{H^{1}({T})}
&\leq  c |\det({A}_{\widetilde{T}})|^{\frac{1}{2}} \| \widetilde{{A}}^{-1} \|_2 \left( \sum_{i=1}^2 h^{- 2}_i \|  \partial_{\hat{x}_i} ( \hat{\varphi} - I_{\widehat{T}} \hat{\varphi} ) \|^2_{L^2(\widehat{T})} \right)^{\frac{1}{2}}. \label{rem36}
\end{align}
The coefficient $h_i^{-2}$ appears on the right-hand side of Eq.~\eqref{rem36}. A further assumption is required for this. Using  Eq.~\eqref{BH=2} and the triangle inequality, we have
\begin{align*}
\displaystyle
 \|  \partial_{\hat{x}_i} ( \hat{\varphi} - I_{\widehat{T}} \hat{\varphi} ) \|^2_{L^2(\widehat{T})}
 &\leq 2 \|  \partial_{\hat{x}_i} ( \hat{\varphi} - Q^{(2)} \hat{\varphi} ) \|^2_{L^2(\widehat{T})} + 2 \|  \partial_{\hat{x}_i} ( Q^{(2)} \hat{\varphi}  -   I_{\widehat{T}} \hat{\varphi} ) \|^2_{L^2(\widehat{T})}.
\end{align*}
We use inequality \eqref{BH=1} to remove the coefficient $h_i^{-2}$. To this end, we have to show that
\begin{align}
\displaystyle
 \|  \partial_{\hat{x}_i} ( Q^{(2)} \hat{\varphi}  -   I_{\widehat{T}} \hat{\varphi} ) \|_{L^2(\widehat{T})}
 &\leq c \|  \partial_{\hat{x}_i} ( \hat{\varphi} - Q^{(2)} \hat{\varphi} ) \|_{H^1(\widehat{T})}. \label{rem37}
\end{align}
However, this is unlikely to hold because Eqs.~\eqref{int1} and \eqref{int513} yield
\begin{align*}
\displaystyle
 \|  \partial_{\hat{x}_i} ( Q^{(2)} \hat{\varphi}  -   I_{\widehat{T}} \hat{\varphi} ) \|_{L^2(\widehat{T})}
 &=  \|  \partial_{\hat{x}_i} (    I_{\widehat{T}} ( Q^{(2)} \hat{\varphi} )  -   I_{\widehat{T}} \hat{\varphi} ) \|_{L^2(\widehat{T})} \\
 &\leq c \| Q^{(2)} \hat{\varphi} -   \hat{\varphi} \|_{H^2(\widehat{T})} \leq c |\hat{\varphi}|_{H^2(\widehat{T})}.
\end{align*}
Using the classical scaling argument (see \cite[Lemma 1.101]{ErnGue04}), we have
\begin{align*}
\displaystyle
|\hat{\varphi}|_{H^2(\widehat{T})} \leq c |\det(A)|^{- \frac{1}{2}} \| {A} \|_2 |\varphi|_{H^2(T)},
\end{align*}
which does not include the quantity $h_i$. Therefore, the quantity $h_i^{-1}$ in Eq.~\eqref{rem36} remains. 

To overcome this problem, we use Theorem \ref{thrApel}. That is, we assume that there exists a linear functional $\mathscr{F}_1$ such that
\begin{align*}
\displaystyle
&\mathscr{F}_1 \in H^{1} (\widehat{T})^{\prime}, \\
&\mathscr{F}_1 ( \partial_{\hat{x}_i} (\hat{\varphi} - I_{\widehat{T}} \hat{\varphi}) ) = 0 \quad i = 1,2, \quad \forall \hat{\varphi} \in \mathcal{C}({\widehat{T}}): \ \partial_{\hat{x}_i} \hat{\varphi} \in H^{1} (\widehat{T}), \\
&\hat{\eta} \in \mathbb{P}^{1}, \quad \mathscr{F}_1(\partial_{\hat{x}_i} \hat{\eta}) = 0 \quad i = 1,2, \quad \Rightarrow \quad \partial_{\hat{x}_i} \hat{\eta} = 0.
\end{align*}
Because the polynomial spaces are finite-dimensional, all norms are equivalent; i.e., because $ | \mathscr{F}_1( \partial_{\hat{x}_i} ( \hat{\eta}  -   I_{\widehat{T}} \hat{\varphi} ) )|$ ($i=1,2$) is a norm on $\mathbb{P}^{0}$, we have that, for $i=1,2$,
\begin{align*}
\displaystyle
\| \partial_{\hat{x}_i} (\hat{\eta} - I_{\widehat{T}} \hat{\varphi}) \|_{L^2(\widehat{T})}
&\leq c  | \mathscr{F}_1 (\partial_{\hat{x}_i} ( \hat{\eta} - I_{\widehat{T}} \hat{\varphi} ) )| 
= c |\mathscr{F}_1 (\partial_{\hat{x}_i} ( \hat{\eta} - \hat{\varphi} ) )| \\
&\leq c \| \partial_{\hat{x}_i} ( \hat{\eta} - \hat{\varphi} ) \|_{H^1(\widehat{T})}.
\end{align*}
Setting $\hat{\eta} := Q^{(2)} \hat{\varphi}$, we obtain Eq.~\eqref{rem37}. Using inequality \eqref{BH=1} yields
\begin{align*}
\displaystyle
 \|  \partial_{\hat{x}_i} ( \hat{\varphi} - I_{\widehat{T}} \hat{\varphi} ) \|^2_{L^2(\widehat{T})}
 &\leq c |  \partial_{\hat{x}_i} \hat{\varphi} |^2_{H^1(\widehat{T})},
\end{align*}
and so inequality \eqref{rem36} together with Eq.~\eqref{jensen} can be written as
\begin{align}
\displaystyle
|{\varphi} - I_T \varphi|_{H^{1}({T})}
&\leq c |\det({A}_{\widetilde{T}})|^{\frac{1}{2}} \| \widetilde{{A}}^{-1} \|_2 \sum_{i,j=1}^2  h^{- 1}_i \|  \partial_{\hat{x}_i}  \partial_{\hat{x}_j}  \hat{\varphi} \|_{L^2(\widehat{T})}. \label{rem38}
\end{align}
Inequality \eqref{scaling2} yields
\begin{align}
\displaystyle
\|  \partial_{\hat{x}_i}  \partial_{\hat{x}_j} \hat{\varphi} \|_{L^2(\widehat{T})}
&\leq c|\det({A}_{\widetilde{T}})|^{-\frac{1}{2}} \| \widetilde{{A}}  \|_2 h_i \sum_{n=1}^2 h_n \left|  \frac{ \partial \varphi}{\partial r_n} \right|_{H^{1}(T)}. \label{rem39}
\end{align}
Therefore, the quantity $h_i^{-1}$ in Eq.~\eqref{rem38} and the quantity $h_i$ in Eq.~\eqref{rem39} cancel out.

\section{New Interpolation Error Estimates}
\subsection{Local Interpolation Error Estimates}
The new scaling arguments in Section \ref{sec54} are the heart of the following local interpolation error estimates.

\begin{thr}[Local interpolation] \label{thr551}
Let $\{ \widehat{T} , \widehat{{P}} , \widehat{\Sigma} \}$ be a finite element with the normed vector space $V(\widehat{T}) := \mathcal{C}({\widehat{T}})$ and $\widehat{{P}} := \mathcal{P}^{k}(\widehat{T})$ with $k \geq 1$.  Let $I_{\widehat{T}} : V({\widehat{T}}) \to   \widehat{{P}}$ be a linear operator. Fix $\ell \in \mathbb{N}$, $m \in \mathbb{N}_0$, and $p,q \in [1,\infty]$ such that $0 \leq m \leq \ell  \leq k+1$, $\ell - m \geq 1$, and the embeddings \eqref{emmed} and \eqref{emmed1} with $s := \ell - m$ hold. Let $\beta$ be a multi-index with $|\beta| = m$. We set $j := \dim (\partial^{\beta} \mathcal{P}^{k})$. Assume that there exist linear functionals $\mathscr{F}_i$, $i=1,\ldots,j$, satisfying the conditions \eqref{functionals}. It then holds that, for all $\hat{\varphi} \in W^{\ell, p} (\widehat{T}) \cap \mathcal{C}({\widehat{T}})$ with ${\varphi} := \hat{\varphi} \circ {\Phi}^{-1}$,
\begin{align}
\displaystyle
| {\varphi} - I_{{T}} {\varphi}|_{W^{m,q}({T})}
\leq  C_1^{I} |T|_d^{\frac{1}{q} - \frac{1}{p}} \left( \frac{H_{T}}{h_{T}} \right)^m \sum_{|\varepsilon| =  \ell-m} {h}^{\varepsilon} | \partial_r^{\varepsilon} \varphi |_{W^{m,p}(T)}, \label{int552}
\end{align}
where $C_1^{I}$ is a positive constant independent of $h_{T}$ and $H_{T}$. In particular, if Condition \ref{Cond333} is imposed, it holds that, for all $\hat{\varphi} \in W^{\ell   , p} (\widehat{T}) \cap \mathcal{C}({\widehat{T}})$ with ${\varphi} = \hat{\varphi} \circ {\Phi}^{-1}$,
\begin{align}
\displaystyle
| {\varphi} - I_{{T}} {\varphi}|_{W^{m,q}({T})}
\leq  C_2^{I} |T|_d^{\frac{1}{q} - \frac{1}{p}} \left( \frac{H_{T}}{h_{T}} \right)^m \sum_{|\varepsilon| = \ell-m} \widetilde{\mathscr{H}}^{\varepsilon}  | \partial_{\tilde{x}}^{\varepsilon} (  {\varphi} \circ \Phi_{T}) |_{W^{ m ,p}({\Phi_{T}^{-1}(T)})}, \label{int553}
\end{align}
where $C_2^{I}$ is a positive constant independent of $h_{T^s}$ and $H_{T^s}$. 

\end{thr}

\begin{pf*}
The introduction of the functionals $\mathscr{F}_i$ follows from \cite{ApeDob92,Ape99}, also see Theorem \ref{thrApel}. Actually, under the same assumptions as in Theorem \ref{thr551}, we have
\begin{align}
\displaystyle
\| \partial_{\hat{x}}^{\beta} (\hat{\varphi} - I_{\widehat{T}} \hat{\varphi}) \|_{L^q(\widehat{T})}
&\leq C^B | \partial_{\hat{x}}^{\beta}  \hat{\varphi} |_{W^{\ell  - m ,p}(\widehat{T})}, \label{int554}
\end{align}
where $|\beta| = m$, $\hat{\varphi} \in \mathcal{C}({\widehat{T}})$, and $\partial_{\hat{x}}^{\beta} \hat{\varphi} \in W^{\ell - m , p} (\widehat{T})$.

The inequalities in \eqref{gint562}, \eqref{jensen}, \eqref{scaling1}, and \eqref{int554} yield
\begin{align}
\displaystyle
&| \varphi - {I}_{T} \varphi |_{W^{m,q}({T})} \leq c| \varphi - {I}_{T} \varphi |_{W^{m,q}({T})} \notag \\
&\ \leq c |\det({A}_{\widetilde{T}})|^{\frac{1}{q}} \| \widetilde{{A}}^{-1} \|_2^m \left( \sum_{|\beta| = m}  ( h^{- \beta} )^q \|  \partial_{\hat{x}}^{\beta}  (\hat{\varphi} - I_{\widehat{T}} \hat{\varphi} )  \|^q_{L^q(\widehat{T})} \right)^{1/q} \notag \\
&\ \leq c |\det({A}_{\widetilde{T}})|^{\frac{1}{q}}  \| \widetilde{{A}}^{-1} \|_2^m  \sum_{|\beta| = m}  ( h^{- \beta} ) \|  \partial_{\hat{x}}^{\beta}  (\hat{\varphi} - I_{\widehat{T}} \hat{\varphi} )  \|_{L^q(\widehat{T})} \notag \\
&\ \leq c |\det({A}_{\widetilde{T}})|^{\frac{1}{q}} \| \widetilde{{A}}^{-1} \|_2^m  \sum_{|\beta| = m}  ( h^{- \beta} )  | \partial_{\hat{x}}^{\beta}  \hat{\varphi} |_{W^{\ell  - m ,p}(\widehat{T})}. \label{int555}
\end{align}
Inequalities \eqref{jensen} and \eqref{scaling2} yield
\begin{align}
\displaystyle
& \sum_{|\beta| = m}  ( h^{- \beta} )  | \partial_{\hat{x}}^{\beta}  \hat{\varphi} |_{W^{\ell  - m ,p}(\widehat{T})} \notag \\
 &\quad \leq  \sum_{|\gamma| = \ell-m} \sum_{|\beta| = m}  ( h^{- \beta} )  \| \partial_{\hat{x}}^{\beta}  \partial_{\hat{x}}^{\gamma}  \hat{\varphi} \|_{L^{p}(\widehat{T})}\notag \\
 &\quad \leq c  |\det({A}_{\widetilde{T}})|^{-\frac{1}{p}} \| \widetilde{{A}}  \|_2^{m}  \sum_{|\gamma| = \ell-m} \sum_{|\beta| = m}  ( h^{- \beta} ) h^{\beta} \sum_{|\epsilon| =  |\gamma|} {h}^{\varepsilon} | \partial_r^{\epsilon} \varphi |_{W^{m,p}(T)} \notag \\
 &\quad \leq c  |\det({A}_{\widetilde{T}})|^{-\frac{1}{p}} \| \widetilde{{A}}  \|_2^{m}  \sum_{|\epsilon| =  \ell-m} {h}^{\varepsilon} | \partial_r^{\varepsilon} \varphi |_{W^{m,p}(T)}. \label{int556}
\end{align}
From \eqref{CN331}, \eqref{CN332}, \eqref{int555}, and  \eqref{int556},
we have
\begin{align*}
\displaystyle
| \varphi - {I}_{T} \varphi |_{W^{m,q}({T})}
\leq c |T|_d^{\frac{1}{q} - \frac{1}{p}} \left( \frac{H_{T}}{h_{T}} \right)^m  \sum_{|\varepsilon| =  \ell-m} {h}^{\varepsilon} | \partial_r^{\varepsilon} \varphi |_{W^{m,p}(T)},
\end{align*}
which is the inequality \eqref{int552}.

Assume that Condition \ref{Cond333} is imposed. Inequality \eqref{scaling3} yields
\begin{align}
\displaystyle
& \sum_{|\beta| = m}  ( h^{- \beta} )  | \partial_{\hat{x}}^{\beta}  \hat{\varphi} |_{W^{\ell  - m ,p}(\widehat{T})} \notag \\
 &\quad \leq  \sum_{|\gamma| = \ell-m} \sum_{|\beta| = m}  ( h^{- \beta} )  \| \partial_{\hat{x}}^{\beta}  \partial_{\hat{x}}^{\gamma}  \hat{\varphi} \|_{L^{p}(\widehat{T})}\notag \\
 &\quad \leq c  |\det({A}_{\widetilde{T}})|^{-\frac{1}{p}} \| \widetilde{{A}}  \|_2^{m}  \sum_{|\gamma| = \ell-m} \sum_{|\beta| = m}  ( h^{- \beta} ) h^{\beta} \sum_{|\epsilon| =  |\gamma|} \widetilde{\mathscr{H}}^{\varepsilon} | \partial_{\tilde{x}}^{\epsilon} \tilde{\varphi} |_{W^{m,p}(\widetilde{T})} \notag \\
 &\quad \leq c  |\det({A}_{\widetilde{T}})|^{-\frac{1}{p}} \| \widetilde{{A}}  \|_2^{m}  \sum_{|\epsilon| =  \ell-m} \widetilde{\mathscr{H}}^{\varepsilon} | \partial_{\tilde{x}}^{\epsilon} \tilde{\varphi} |_{W^{m,p}(\widetilde{T})}. \label{int557}
\end{align}
From \eqref{CN331}, \eqref{CN332}, \eqref{int555}, and  \eqref{int557},
we have
\begin{align*}
\displaystyle
| \varphi - {I}_{T} \varphi |_{W^{m,q}({T})}
\leq c |T|_d^{\frac{1}{q} - \frac{1}{p}} \left( \frac{H_{T}}{h_{T}} \right)^m  \sum_{|\varepsilon| =  \ell-m} \widetilde{\mathscr{H}}^{\varepsilon} | \partial_{\tilde{x}}^{\epsilon} \tilde{\varphi} |_{W^{m,p}(\widetilde{T})},
\end{align*}
which is the inequality \eqref{int553} using $\widetilde{T} = \Phi^{-1}_{T}(T)$ and $\tilde{\varphi} = \varphi \circ \Phi_{T}$.
\qed
\end{pf*}

\subsection{Global Interpolation Error Estimates}  \label{GIEE}
A global interpolation operator ${I}_h$ is constructed as follows (e.g., see \cite[Section 1.4.2]{ErnGue04}). Its domain is defined by
\begin{align*}
\displaystyle
{D}({I}_h) := \{ \varphi \in L^1(\Omega); \  \varphi|_{T} \in V(T), \ \forall T \in \mathbb{T}_h \}.
\end{align*}
 For $T \in \mathbb{T}_h$ and $\varphi \in {D}({I}_h)$, the quantities $\chi_{i}(\varphi |_{T}) $ are meaningful on all the mesh elements and $1 \leq i \leq n_0$. The global interpolation ${I}_h \varphi$ can be specified elementwise using the local interpolation operators, that is,
\begin{align*}
\displaystyle
({I}_h \varphi )|_{T} := {I}_T(\varphi|_{T}) = \sum_{i=1}^{n_0} \chi_{i}(\varphi |_{T}) \theta_{i}  \quad \forall T \in \mathbb{T}_h, \quad \forall \varphi \in {D}({I}_h).
\end{align*}
The global interpolation operator ${I}_h: D(I_h) \to V_h$ is defined as
\begin{align*}
\displaystyle
I_h: D(I_h) \ni \varphi \mapsto I_h \varphi := \sum_{T \in \mathbb{T}_h} \sum_{i=1}^{n_0}  \chi_{i}(\varphi |_{T}) \theta_{i} \in V_h^n,
\end{align*}
where $V_h^n$ is defined as
\begin{align}
\displaystyle
V_h^n := \{ \varphi_h \in L^1(\Omega)^n; \ \varphi_h|_{T} \in {P}, \ \forall T \in \mathbb{T}_h\}. \label{space=Vh}
\end{align}

\begin{cor} \label{newglobal=cor}
Suppose that the assumptions of Theorem \ref{thr551} are satisfied. We impose  Condition \ref{Cond1=sec6}. Let $I_h$ be the corresponding global interpolation operator. It then holds that, for any $ \varphi \in W^{\ell,p}(\Omega)$;
\begin{description}	
  \item[(\Roman{lone})] if Condition \ref{Cond333} is not imposed, 
\begin{align}
\displaystyle
| {\varphi} - I_{{h}} {\varphi}|_{W^{m,q}({\Omega})}
&\leq c  \sum_{T \in \mathbb{T}_h} |T|_d^{\frac{1}{q} - \frac{1}{p}} \sum_{|\varepsilon| =  \ell-m} {h}^{\varepsilon} | \partial_r^{\varepsilon} \varphi |_{W^{m,p}(T)}. \label{gint566}
\end{align}
   \item[(\Roman{ltwo})] if Condition \ref{Cond333} is imposed, 
\begin{align}
\displaystyle
| {\varphi} - I_{{h}} {\varphi}|_{W^{m,q}({\Omega})}
&\leq c  \sum_{T \in \mathbb{T}_h} |T|_d^{\frac{1}{q} - \frac{1}{p}} \sum_{|\varepsilon| =  \ell-m} \widetilde{\mathscr{H}}^{\varepsilon}  | \partial_{\tilde{x}}^{\varepsilon}   ( \varphi \circ \Phi_{T}) |_{W^{ m ,p}(\Phi_{T}^{-1}(T))}. \label{gint567}
\end{align}
\end{description}
\end{cor}

\begin{pf*}
 If Condition \ref{Cond333} is not imposed, then using the local interpolation error \eqref{int552},
\begin{align*}
\displaystyle
| {\varphi} - I_{{h}} {\varphi}|_{W^{m,q}({\Omega})}^q
&= \sum_{T \in \mathbb{T}} | {\varphi} - I_{{T}} {\varphi}|_{W^{m,q}({T})}^q \\
&\leq c  \sum_{T \in \mathbb{T}} |T|_d^{q \left( \frac{1}{q} - \frac{1}{p} \right)} \left( \frac{H_{T}}{h_{T}} \right)^{q m} \left( \sum_{|\varepsilon| =  \ell-m} {h}^{\varepsilon} | \partial_r^{\varepsilon} \varphi |_{W^{m,p}(T)} \right)^q,
\end{align*}
which leads to the desired result together with \eqref{jensen} and Condition \ref{Cond1=sec6}.

 If Condition \ref{Cond333} is not imposed, then using the local interpolation error \eqref{int553},
 \begin{align*}
\displaystyle
| {\varphi} - I_{{h}} {\varphi}|_{W^{m,q}({\Omega})}^q
&= \sum_{T \in \mathbb{T}} | {\varphi} - I_{{T}} {\varphi}|_{W^{m,q}({T})}^q \\
&\leq c  \sum_{T \in \mathbb{T}} |T|_d^{q \left( \frac{1}{q} - \frac{1}{p} \right)} \left( \frac{H_{T}}{h_{T}} \right)^{q m} \left( \sum_{|\varepsilon| =  \ell-m}  \widetilde{\mathscr{H}}^{\varepsilon}  | \partial_{\tilde{x}}^{\varepsilon} (  {\varphi} \circ \Phi_{T}) |_{W^{ m ,p}({\Phi_{T}^{-1}(T)})} \right)^q,
\end{align*}
 which leads to the desired result together with \eqref{jensen} and Condition \ref{Cond1=sec6}.
\qed
\end{pf*}

\subsection{Examples of Anisotropic Elements}
When $k=1$, $ell = 2$, $m=1$ and $q=p$ in \eqref{int552} of Theorem \ref{thr551}, the estimate is written as
\begin{align}
\displaystyle
| {\varphi} - I_{{T}} {\varphi}|_{W^{1,p}({T})}
\leq  C_1^{I} \frac{H_{T}}{h_{T}} \sum_{i=1}^d  {h}_i \left |  \frac{\partial \varphi}{\partial r_i} \right|_{W^{1,p}(T)}. \label{int552=ex}
\end{align}
Let $T \subset \mathbb{R}^2$ be a triangle. As described in Section \ref{iso=mesh=ex2d}, an isotropic mesh element has equal or nearly equal edge lengths and angles, resulting in a balanced shape. Then, the estimate \eqref{int552=ex} becomes
\begin{align}
\displaystyle
 |\varphi - {I}_{T} \varphi|_{W^{1,p}(T)}
&\leq c h_{T} | \varphi |_{W^{2,p}(T)}. \label{int552=ex2}
\end{align}

We considered the following five anisotropic elements as in Section \ref{ani=mesh=ex2d}:  Let $0 \< s \ll 1$, $s \in \mathbb{R}$ and $\varepsilon,\delta,\gamma \in \mathbb{R}$.

\begin{ex}[Right-angled triangle] \label{ex14=3}
Let $T \subset \mathbb{R}^2$ be the simplex with vertices $p_1 := (0,0)^{\top}$, $p_2 := (s,0)^{\top}$ and $p_3 := (0,s^{\varepsilon})^{\top}$ with $1 \< \varepsilon$. We then have $h_1 = s$, $h_2 =  s^{\varepsilon}$ and $h_T = \sqrt{s^2 + s^{2 \varepsilon}}$; i.e.,
\begin{align*}
\displaystyle
\frac{H_{T}}{h_{T}} =2.
\end{align*}
In this case, the estimate \eqref{int552=ex} becomes
\begin{align*}
\displaystyle
 |\varphi - {I}_{T} \varphi|_{W^{1,p}(T)}
&\leq 2 C_1^{I} \sum_{i=1}^2  {h}_i \left |  \frac{\partial \varphi}{\partial r_i} \right|_{W^{1,p}(T)},
\end{align*}
which is the anisotropic interpolation error estimate.
\end{ex}

\begin{ex}[Dagger] \label{ex14=4}
Let $T \subset \mathbb{R}^2$ be the simplex with vertices $p_1 := (0,0)^{\top}$, $p_2 := (s,0)^{\top}$ and $p_3 := (s^{\delta},s^{\varepsilon})^{\top}$ with $1 \< \varepsilon \< \delta $. We then have $h_1 = \sqrt{(s - s^{\delta})^2 + s^{2 \varepsilon}}$, $h_2 = \sqrt{s^{2 \delta} + s^{2 \varepsilon}}$ and $h_T=s$; i.e.,
\begin{align*}
\displaystyle
\frac{H_{T}}{h_{T}} &= \frac{ \sqrt{(s- s^{\delta})^2 + s^{2 \varepsilon}} \sqrt{s^{2 \delta} + s^{2 \varepsilon}}}{\frac{1}{2} s^{1 + \varepsilon}} \leq c.
\end{align*}
In this case, the estimate \eqref{int552=ex} becomes
\begin{align*}
\displaystyle
 |\varphi - {I}_{T} \varphi|_{W^{1,p}(T)}
&\leq c \sum_{i=1}^2  {h}_i \left |  \frac{\partial \varphi}{\partial r_i} \right|_{W^{1,p}(T)},
\end{align*}
which is the anisotropic interpolation error estimate.
\end{ex}

\begin{ex}[Blade] \label{ex14=5}
Let $T \subset \mathbb{R}^2$ be the simplex with vertices $p_1 := (0,0)^{\top}$, $p_2 := (2s,0)^{\top}$ and $p_3 := (s ,s^{\varepsilon})^{\top}$ with $1 \< \varepsilon $. We then have $h_1 = h_2 = \sqrt{s^{2} + s^{2 \varepsilon}}$ and $h_T = 2s$; i.e.,
\begin{align*}
\displaystyle
\frac{H_{T}}{h_{T}} = \frac{s^{2} + s^{2 \varepsilon}}{s^{1 + \varepsilon}} \to \infty \quad \text{as $s \to 0$}.
\end{align*}
In this case, the estimate \eqref{int552=ex} becomes
\begin{align*}
\displaystyle
 |\varphi - {I}_{T} \varphi|_{W^{1,p}(T)}
&\leq c s^{2 - \varepsilon} | \varphi |_{W^{2,p}(T)}. 
\end{align*}
When $\varepsilon > 2$, this implies that the estimate diverges as $s \to 0$. 
\end{ex}

\begin{ex}[Dagger] \label{ex14=6}
Let $T \subset \mathbb{R}^2$ be the simplex with vertices $p_1 := (0,0)^{\top}$, $p_2 := (s,0)^{\top}$ and $p_3 := (s^{\delta},s^{\varepsilon})^{\top}$ with $1 \< \delta \< \varepsilon $. We then have $h_1 = \sqrt{(s - s^{\delta})^2 + s^{2 \varepsilon}}$, $h_2 = \sqrt{s^{2 \delta} + s^{2 \varepsilon}}$ and $h_T=s$; i.e.,
\begin{align*}
\displaystyle
\frac{H_{T}}{h_{T}} &= \frac{ \sqrt{(s- s^{\delta})^2 + s^{2 \varepsilon}} \sqrt{s^{2 \delta} + s^{2 \varepsilon}}}{\frac{1}{2} s^{1 + \varepsilon}} \leq c s^{\delta - \varepsilon} \to \infty \quad \text{as $s \to 0$}.
\end{align*}
In this case, the estimate \eqref{int552=ex} becomes
\begin{align*}
\displaystyle
 |\varphi - {I}_{T} \varphi|_{W^{1,p}(T)}
&\leq c  s^{\delta - \varepsilon} \left( s \left |  \frac{\partial \varphi}{\partial r_1} \right|_{W^{1,p}(T)} + s^{\delta}\left |  \frac{\partial \varphi}{\partial r_2} \right|_{W^{1,p}(T)}  \right) \\
&\leq c s^{1 + \delta - \varepsilon}  | \varphi |_{W^{2,p}(T)}. 
\end{align*}
When $\varepsilon - \delta > 1$, this implies that the estimate diverges as $s \to 0$.
\end{ex}

\begin{ex}[Right-angled triangle] \label{ex14=7}
Let $T \subset \mathbb{R}^2$ be the simplex with vertices $p_1 := (0,0)^{\top}$, $p_2 := (s,0)^{\top}$ and $p_3 := (0, \delta s)^{\top}$ with $\delta \ll 1$. We then have $h_1 = s$, $h_2 = \delta s$ and $h_T = s \sqrt{1 + \delta^2}$; i.e.,
\begin{align*}
\displaystyle
\frac{H_{T}}{h_{T}} =2.
\end{align*}
In this case, the estimate \eqref{int552=ex} becomes
\begin{align*}
\displaystyle
 |\varphi - {I}_{T} \varphi|_{W^{1,p}(T)}
&\leq c \sum_{i=1}^2  {h}_i \left |  \frac{\partial \varphi}{\partial r_i} \right|_{W^{1,p}(T)},
\end{align*}
which is the anisotropic interpolation error estimate.
\end{ex}

\begin{ex}[Blade] \label{ex14=8}
Let $T \subset \mathbb{R}^2$ be the simplex with vertices $p_1 := (0,0)^{\top}$, $p_2 := (2s,0)^{\top}$ and $p_3 := (s , \delta s)^{\top}$ with $\delta \ll 1$. We then have $h_1 = h_2 = s \sqrt{1 + \delta^2}$ and $h_T = 2s$; i.e.,
\begin{align*}
\displaystyle
\frac{H_{T}}{h_{T}} = \frac{s^2 (1 + \delta^2)}{\delta s^2} \leq \frac{c}{\delta},
\end{align*}
In this case, the estimate \eqref{int552=ex} becomes
\begin{align*}
\displaystyle
 |\varphi - {I}_{T} \varphi|_{W^{1,p}(T)}
&\leq \frac{c}{\delta} s | \varphi |_{W^{2,p}(T)}. 
\end{align*}
This implies that the estimate converges as $s \to 0$ and the error may be large.   Thus, even if anisotropic mesh partitioning is used, it is unlikely to improve calculation efficiency.
\end{ex}

\subsection{Examples that do not satisfy conditions  \eqref{functionals} in Theorem \ref{thrApel}}
The following lemma (\cite[Lemma 4]{ApeDob92}, \cite[Lemma 2.3]{Ape99}) gives a criterion for the existence of linear functionals satisfying conditions \eqref{functionals=b} and \eqref{functionals=c}.

\begin{lem}
Let $\mathbb{P}$ be an arbitrary polynomial space and $\beta$ be a multi-index. We set $j := \dim (\partial^{\beta} \mathbb{P})$. Assume that $I:\mathcal{C}^{\mu} (\widehat{T}) \to \mathbb{P}$, $\mu \in \mathbb{N}$, is a linear operator with $I \hat{\eta} = \hat{\eta}$ $\forall \hat{\eta} \in \mathbb{P}$. Then, there exist linear functionals $\mathscr{F}_i: \mathcal{C}^{\infty}(\widehat{T}) \to \mathbb{R}$, $i = 1,\ldots,j$, such that
\begin{align}
\displaystyle
&\mathscr{F}_i ( \partial^{\beta} (\hat{\varphi} - I \hat{\varphi}) ) = 0 \quad  \forall i = 1,\ldots,j, \quad \forall \hat{\varphi} \in \mathcal{C}^{\infty}({\widehat{T}}), \label{CEX1} \\
&\hat{\eta} \in \mathbb{P}, \quad \mathscr{F}_i(\partial^{\beta} \hat{\eta}) = 0 \quad \forall i = 1,\ldots,j \quad \Rightarrow \quad \partial^{\beta} \hat{\eta} = 0 \label{CEX2}
\end{align}
if and only if the condition
\begin{align}
\displaystyle
\hat{\varphi} \in \mathcal{C}^{\infty}(\widehat{T}), \quad \partial^{\beta} \hat{\varphi} = 0 \quad \Rightarrow \quad \partial^{\beta} I \hat{\varphi} = 0 \label{CEX3}
\end{align}
holds.	
\end{lem}

\begin{pf*}
A proof can be found in \cite[Lemma 4]{ApeDob92}.
\qed
\end{pf*}

If Condition \eqref{CEX3} is violated, estimate \eqref{ineq512} does not hold. This means that one cannot obtain the estimate \eqref{int552}, which is sharper than \eqref{ineq523}. 

The following are examples that do not satisfy \eqref{CEX3}. Let $\widehat{T} \subset \mathbb{R}^2$ be the reference element with vertices $\widehat{p}_1 := (0,0)^{\top}$, $\widehat{p}_2 := (1,0)^{\top}$, $\widehat{p}_3 := (0, 1)^{\top}$.  We set $\widehat{p}_{4} := (1/3,1/3)^{\top}$. We define the barycentric coordinates $ \lambda_i: \mathbb{R}^2 \to \mathbb{R}$, $i=1,\ldots,3$, on the reference element as
\begin{align*}
\displaystyle
\lambda_1 := 1 - \hat{x}_1 - \hat{x}_2, \quad \lambda_2 := \hat{x}_1, \quad \lambda_3 := \hat{x}_2, \quad (\hat{x}_1,\hat{x}_2)^{\top} \in \widehat{T}.
\end{align*}

\begin{ex}[$\mathbb{P}^1$ + bubble Finite Element]
As mentioned in Example \ref{ex526}, we define the local basis functions as
\begin{align*}
\displaystyle
\theta_4(x) &:= 27 \lambda_1(x) \lambda_2(x) \lambda_3(x), \\
\theta_i(x) &:= \lambda_i(x) - \frac{1}{3} \theta_4(x), \quad i=1,2,3.
\end{align*}
The interpolation operator $I_T^b$ defined by
\begin{align*}
\displaystyle
I^{b}:\mathcal{C} (\widehat{T}) \ni \hat{\varphi}  \mapsto I^{b} \hat{\varphi} := \sum_{i=1}^{4} \hat{\varphi} (\widehat{P}_{i}) \theta_i \in \Span \{ \theta_1, \theta_2, \theta_3 , \theta_4 \}.
\end{align*}
Let $\beta = (1,0)$. Setting $\hat{\varphi}(\hat{x}_{1},\hat{x}_{2}) := \hat{x}_{2}^2$, we have $\frac{\partial \hat{\varphi}}{\partial \hat{x}_1} = 0$. By simple calculation, we obtain
\begin{align*}
\displaystyle
 \frac{\partial }{\partial \hat{x}_1} I^{b} \hat{\varphi}
 &= \hat{\varphi} (\widehat{p}_{1})  \frac{\partial \theta_1}{\partial \hat{x}_1} + \hat{\varphi} (\widehat{p}_{2})  \frac{\partial \theta_1}{\partial \hat{x}_1} + \hat{\varphi} (\widehat{p}_{3})  \frac{\partial \theta_3}{\partial \hat{x}_1} + \hat{\varphi} (\widehat{p}_{4})  \frac{\partial \theta_4}{\partial \hat{x}_1} \\
 &= \frac{\partial \theta_3}{\partial \hat{x}_1} + \frac{1}{3^2} \frac{\partial \theta_4}{\partial \hat{x}_1}
 = - \frac{1 }{3} \frac{\partial \theta_4}{\partial \hat{x}_1} + \frac{1}{3^2} \frac{\partial \theta_4}{\partial \hat{x}_1} \not\equiv 0.
\end{align*}
Therefore, the condition \eqref{CEX3} is not satisfied. This implies that the error estimate \eqref{ineq512} on the reference element does not hold for the $\mathcal{P}^1$ + bubble finite element. 	
\end{ex}

\begin{ex}[$\mathbb{P}^3$ Hermite Finite Element]
Following \cite[Theorem 2.2.8]{Cia02}, we define the Hermite interpolation operator $I^H: H^3(T) \to \mathbb{P}^3$ as
\begin{align*}
\displaystyle
 I^{H} \hat{\varphi} &:= \sum_{i=1}^3 \left( -2 \lambda_i^3 + 3 \lambda_i^2 - 7 \lambda_i \sum_{1 \leq j \< k \leq 3, j \neq i, k \neq i} \lambda_j \lambda_k \right) \hat{\varphi}(\widehat{p}_i) + 27 \lambda_1 \lambda_2 \lambda_3 \hat{\varphi}(\widehat{p}_4) \\
 &\quad + \sum_{i=1}^3 \left( \sum_{j=1}^3 \lambda_i \lambda_j (2 \lambda_i + \lambda_j -1) (\widehat{p}_j^{(1)} - \widehat{p}_i^{(1)}) \right) \frac{\partial \hat{\varphi}}{\partial \hat{x}_1}(\widehat{p}_i) \\
  &\quad + \sum_{i=1}^3 \left( \sum_{j=1}^3 \lambda_i \lambda_j (2 \lambda_i + \lambda_j -1) (\widehat{p}_j^{(2)} - \widehat{p}_i^{(2)}) \right) \frac{\partial \hat{\varphi}}{\partial \hat{x}_2}(\widehat{p}_i),
\end{align*}
where $\widehat{p}^{(k)}_{i}$, $1 \leq k \leq 2$, are the components of a point $\widehat{p}_i \in \mathbb{R}^2$. Let $\beta = (1,0)$. Setting $\hat{\varphi}(\hat{x}_{1},\hat{x}_{2}) := \hat{x}_{2}^4$, we have $\frac{\partial \hat{\varphi}}{\partial \hat{x}_1} = 0$. Furthermore, by a simple calculation, i.e.,
\begin{align*}
\displaystyle
\frac{\partial }{\partial \hat{x}_1 }  (\lambda_1 \lambda_2 \lambda_3 )
&= - \hat{x}_2^2 - 2 \hat{x}_1 \hat{x}_2 - \hat{x}_2,\\
\frac{\partial }{\partial \hat{x}_1 } \left\{ \lambda_3 \lambda_1 (2 \lambda_3 + \lambda_1 -1) \right\}
&= - \hat{x}_2 + 2 \hat{x}_1 \hat{x}_2,\\
\frac{\partial }{\partial \hat{x}_1 } \left\{  \lambda_3 \lambda_2 (2 \lambda_3 + \lambda_2 - 1) \right\}
&= - \hat{x}_2 + 2 \hat{x}_1 \hat{x}_2 + 2 \hat{x}_2^2,
\end{align*}
we obtain
\begin{align*}
\displaystyle
 \frac{\partial }{\partial \hat{x}_1} I^{H} \hat{\varphi}
 &=  \frac{\partial }{\partial \hat{x}_1 } \left( -2 \lambda_3^3 + 3 \lambda_3^2 - 7 \lambda_3 \sum_{1 \leq j \< k \leq 3, j \neq 3, k \neq 3} \lambda_j \lambda_k \right) \hat{\varphi}(\widehat{p}_3) \\
 &\quad + 27 \frac{\partial }{\partial \hat{x}_1 } (\lambda_1 \lambda_2 \lambda_3 ) \hat{\varphi}(\widehat{p}_4) \\
 &\quad +  \frac{\partial }{\partial \hat{x}_1 } \left( \sum_{j=1}^3 \lambda_3 \lambda_j (2 \lambda_3 + \lambda_j -1) (\widehat{p}_j^{(2)} - \widehat{p}_3^{(2)}) \right) \frac{\partial \hat{\varphi}}{\partial \hat{x}_2}(\widehat{p}_3) \\
 &= -7 \frac{\partial }{\partial \hat{x}_1 }  (\lambda_1 \lambda_2 \lambda_3 ) \hat{\varphi}(\widehat{p}_3)  + 27 \frac{\partial }{\partial \hat{x}_1 } (\lambda_1 \lambda_2 \lambda_3 ) \hat{\varphi}(\widehat{p}_4) \\
 &\quad + \frac{\partial }{\partial \hat{x}_1 } \left\{ \lambda_3 \lambda_1 (2 \lambda_3 + \lambda_1 -1) (\widehat{p}_1^{(2)} - \widehat{p}_3^{(2)}) \right \} \frac{\partial \hat{\varphi}}{\partial \hat{x}_2}(\widehat{p}_3) \\
  &\quad + \frac{\partial }{\partial \hat{x}_1 } \left\{ \lambda_3 \lambda_2 (2 \lambda_3 + \lambda_2 - 1) (\widehat{p}_2^{(2)} - \widehat{p}_3^{(2)})  \right \} \frac{\partial \hat{\varphi}}{\partial \hat{x}_2}(\widehat{p}_3) \\
 &= - 7 ( - \hat{x}_2^2 - 2 \hat{x}_1 \hat{x}_2 - \hat{x}_2) + \frac{1}{3} ( - \hat{x}_2^2 - 2 \hat{x}_1 \hat{x}_2 - \hat{x}_2) \\
 &\quad + 8 (\hat{x}_2 - 2 \hat{x}_1 \hat{x}_2 - \hat{x}_2^2) \not\equiv 0.
\end{align*}
Here, we used
\begin{align*}
\displaystyle
&\hat{\varphi}(\widehat{p}_i) = 0, \quad  \frac{\partial \hat{\varphi}}{\partial \hat{x}_2}(\widehat{p}_i) = 0, \quad i=1,2, \\
&\hat{\varphi}(\widehat{p}_3) = 1, \quad \hat{\varphi}(\widehat{p}_4) = \frac{1}{3^4}, \quad \frac{\partial \hat{\varphi}}{\partial \hat{x}_2}(\widehat{p}_3) = 4, \\
& \widehat{p}_1^{(2)} - \widehat{p}_3^{(2)} = -1, \quad \widehat{p}_2^{(2)} - \widehat{p}_3^{(2)} = -1.
\end{align*}
Therefore, Condition \eqref{CEX3} is not satisfied. This implies that error estimate \eqref{ineq512} on the reference element does not hold for Hermitian finite elements.
\end{ex}

\subsection{Effect of the quantity $|T|_d^{\frac{1}{q} - \frac{1}{p}}$ in the interpolation error estimates for $d = 2,3$}
We consider the effect of the factor $|T|_d^{\frac{1}{q} - \frac{1}{p}}$. 

\subsubsection{Case that $q \> p$} \label{sec=qp}
When $q \> p$, the factor may affect the convergence order. In particular, the interpolation error estimate may diverge on anisotropic mesh partitions. 

Let $T \subset \mathbb{R}^2$ be the triangle with vertices $p_1 := (0,0)^{\top}$, $p_2 := (s,0)^{\top}$, $p_3 := (0, s^{\varepsilon})^{\top}$ for $0 \< s \ll 1$, $\varepsilon \geq 1 $, $s \in \mathbb{R}$ and $\varepsilon \in \mathbb{R}$. Then,
\begin{align*}
\displaystyle
\frac{h_{\max}}{h_{\min}} = s^{1 - \varepsilon}, \quad \frac{H_{T}}{h_{T}} =2, \quad |T|_2 = \frac{1}{2} s^{1 + \varepsilon}.
\end{align*}
Let $k=1$, $\ell = 2$,  $m = 1$, $q = 2$, and $ p \in (1,2)$. Then, $W^{1,p}({T}) \hookrightarrow L^2({T})$ and Theorem \ref{thr551} lead to
\begin{align*}
\displaystyle
| {\varphi} - I_{{T}} {\varphi} |_{H^{1}({T})}
\leq  c s^{ - (1+\varepsilon)\frac{2 - p}{2p}} \left( s \left | \frac{\partial \varphi}{\partial r_1} \right |_{W^{1,p}({T})} + s^{\varepsilon} \left | \frac{\partial \varphi}{\partial r_2} \right |_{W^{1,p}({T})} \right).
\end{align*}
When $\varepsilon = 1$ (the case of the isotropic element), we get
\begin{align*}
\displaystyle
| {\varphi} - I_{{T}} {\varphi} |_{H^{1}({T})}
\leq  c h_T^{ \frac{2(p-1)}{p}} |\varphi|_{W^{2,p}(T)}, \quad  \frac{2(p-1)}{p} \> 0.
\end{align*}
However, when  $\varepsilon \> 1$ (the case of the anisotropic element),  the estimate may diverge as $s \to 0$. Therefore, if $q \> p$, the convergence order of the interpolation operator may deteriorate.

We next set $m=0$, $\ell= 2$, $q= \infty$, and $p=2$. Let 
\begin{align*}
\displaystyle
\varphi(x,y) := x^2 + y^2.
\end{align*}
Let $I_{T}^{L}: \mathcal{C}^0(T) \to \mathbb{P}^1$ be the local Lagrange interpolation operator. For any nodes $p_i$ of $T$, because $I_{T}^{L} \varphi (p_i) = \varphi (p_i)$, we have
\begin{align*}
\displaystyle
I_{T}^{L} \varphi (x,y) = s x + s^{\varepsilon} y.
\end{align*}
It thus holds that
\begin{align*}
\displaystyle
(\varphi - I_{T}^{L} \varphi)(x,y)
= \left(x - \frac{s}{2} \right)^2 + \left(y - \frac{s^{\varepsilon}}{2} \right)^2 - \frac{1}{4}(s^2 + s^{2 \varepsilon}).
\end{align*}
We therefore have, because $H^2(T) \hookrightarrow L^{\infty}(T)$,
\begin{align*}
\displaystyle
\| {\varphi} - I_{{T}}^L {\varphi} \|_{L^{\infty}({T})} &= \frac{1}{4}(s^2 + s^{2 \varepsilon}), \quad
\sum_{|\gamma| = 2} \widetilde{\mathscr{H}}^{\gamma}  \| \partial_x^{\gamma}  {\varphi} \|_{L^{2}({T})} 
= 2  |T|_2^{\frac{1}{2}} ( s^2  + s^{2 \varepsilon} ),
\end{align*}
and thus,
\begin{align*}
\displaystyle
\frac{\| {\varphi} - I_{{T}}^L {\varphi} \|_{L^{\infty}({T})} }{|T|_2^{- \frac{1}{2}} \sum_{|\gamma| = 2} \widetilde{\mathscr{H}}^{\gamma}  \| \partial_x^{\gamma}  {\varphi} \|_{L^{2}({T})} } = \frac{1}{8}.
\end{align*}
This example implies that the convergence order is not optimal, but the estimate converges on anisotropic meshes.

\subsubsection{Case that $q \< p$}
We consider Theorem \ref{thr551}. Let $I_{T}^{L}: \mathcal{C}(T) \to \mathbb{P}^k$ ($k \in \mathbb{N}$) be the local Lagrange interpolation operator. Let $  {\varphi}  \in W^{\ell,\infty}(T) $ be such that $\ell \in \mathbb{N}$, $2 \leq \ell \leq k+1$. It then holds that, for any $m \in \{ 0, \ldots,\ell -1\}$ and $q \in [1,\infty]$,
\begin{align}
\displaystyle
| {\varphi} - I_{{T}}^L {\varphi} |_{W^{m,q}({T})} \leq c |T|_d^{\frac{1}{q} } \left( \frac{H_{T}}{h_{T}} \right)^m \sum_{|\gamma| = \ell-m} h^{\gamma}| \partial_{r}^{\gamma} \varphi |_{W^{ m ,\infty}(T)}. \label{rem8=56}
\end{align}
The convergence order is therefore improved by $|T|_d^{\frac{1}{q} }$. We do numerical tests to confirm this. Let $k=1$ and 
\begin{align*}
\displaystyle
\varphi(x,y,z) := x^2 + \frac{1}{4} y^2 + z^2.
\end{align*}
Let $s := \frac{1}{N}$, $N \in \mathbb{N}$ and $\varepsilon \in \mathbb{R}$, $1 \< \varepsilon$. We compute the convergence order with respect to the $H^1$ norm defined by
\begin{align*}
\displaystyle
Err_s^{\varepsilon}(H^1) := | \varphi - I_{T}^{L} \varphi |_{H^1(T)}.
\end{align*}
The convergence indicator $r$ is defined by
\begin{align*}
\displaystyle
r = \frac{1}{\log(2)} \log \left( \frac{Err_s^{\varepsilon}(H^1)}{Err_{s/2}^{\varepsilon}(H^1)} \right).
\end{align*}

\begin{description}
  \item[(\Roman{lone})] Let $T \subset \mathbb{R}^3$ be the simplex with vertices $p_1 := (0,0,0)^{\top}$, $p_2 := (s,0,0)^{\top}$, $p_3 := (0,s^{\varepsilon},0)^{\top}$, and $p_4 := (0,0,s^{\delta})^{\top}$ ($1 \< \delta \leq  \varepsilon $), and $0 \< s \ll 1$, $s \in \mathbb{R}$. We then have $h_1 = \sqrt{s^2 + s^{2 \varepsilon}}$, $h_2 = s^{\varepsilon}$ and $h_3 := \sqrt{s^{2 \varepsilon} + s^{2 \delta}}$; i.e.,
\begin{align*}
\displaystyle
\frac{h_{\max}}{h_{\min}}  \leq c s^{1 - \varepsilon}, \quad \frac{H_{T}}{h_{T}} \leq c.
\end{align*}
From \eqref{rem8=56} with $m=1$, $\ell = 2$, and $q=2$, because $|T|_3 \approx s^{1+\varepsilon + \delta}$, we have the estimate
\begin{align*}
\displaystyle
&| {\varphi} - I_{{T}}^L {\varphi} |_{H^{1}({T})} \leq c  h_{T}^{\frac{3 + \varepsilon + \delta}{2}}.
\end{align*}
 Computational results are for the case that $\varepsilon = 3.0$ and $\delta = 2.0$ (Table \ref{modify=table1}).
  
\begin{table}[htbp]
\caption{Error of the local interpolation operator ($\varepsilon = 3.0, \delta = 2.0$)}
\centering
\begin{tabular}{l | l | l | l  } \hline
$N$ &  $s$  & $Err_s^{3.0}(H^1)$ & $r$   \\ \hline \hline
64 & 1.5625e-02 &  2.4336e-08  &           \\
128 & 7.8125e-03 & 1.5209e-09 &    4.00   \\
256 & 3.9062e-03  & 9.5053e-11  &    4.00      \\
\hline
\end{tabular}
\label{modify=table1}
\end{table}

\item[(\Roman{ltwo})] Let $T \subset \mathbb{R}^3$ be the simplex with vertices $p_1 := (0,0,0)^{\top}$, $p_2 := (s,0,0)^{\top}$, $p_3 := (\frac{s}{2},s^{\varepsilon},0)^{\top}$, and $p_4 := (0,0,s)^{\top}$ ($1 \< \varepsilon \leq 6$) and $0 \< s \ll 1$, $s \in \mathbb{R}$. We then have $h_1 = s$, $h_2 = \sqrt{ s^2/4 + s^{2 \varepsilon}}$ and $h_3 := s$; i.e.,
\begin{align*}
\displaystyle
\frac{h_{\max}}{h_{\min}} = \frac{t}{\sqrt{ s^2/4 + t^{2 \varepsilon}}} \leq c, \quad \frac{H_{T}}{h_{T}} \leq c s^{1 - \varepsilon}.
\end{align*}
From \eqref{rem8=56} with $m=1$, $\ell = 2$, and $q=2$, because $|T|_3 \approx s^{2+\varepsilon}$, we have the estimate
\begin{align*}
\displaystyle
&| {\varphi} - I_{{T}}^L {\varphi} |_{H^{1}({T})} \leq c  h_{T}^{3 - \frac{\varepsilon}{2}}.
\end{align*}
Computational results are for the cases that $\varepsilon = 3.0,6.0$ (Table \ref{modify=table2}).

\begin{table}[htbp]
\caption{Error of the local interpolation operator ($\varepsilon = 3.0,6.0$)}
\centering
\begin{tabular}{l | l | l | l| l| l  } \hline
$N$ &  $s$  & $Err_s^{3.0}(H^1)$ & $r$  & $Err_s^{6.0}(H^1)$ & $r$ \\ \hline \hline
64 & 1.5625e-02 &  1.9934e-04 &      &  1.0206e-01  &     \\
128 & 7.8125e-03 & 7.0477e-05   & 1.50  &  1.0206e-01 & 0\\
256 & 3.9062e-03  & 2.4917e-05  &   1.50 & 1.0206e-01  &  0    \\
\hline
\end{tabular}
\label{modify=table2}
\end{table}
\end{description}

\subsection{What happens if violating the maximum-angle condition?}
This subsection introduces two negative points by violating the maximum-angle condition. One is that it is practically disadvantageous. As an example, let $T \subset \mathbb{R}^2$ be the triangle with vertices $p_1 := (0,0)^{\top}$, $p_2 := (s,0)^{\top}$, $p_3 := (s/2, s^{\varepsilon})^{\top}$ for $0 \< s \ll 1$, $\varepsilon \geq 1 $, $s \in \mathbb{R}$ and $\varepsilon \in \mathbb{R}$. From Theorem \ref{thr551} with $k=1$, $\ell = 2$,  $m = 1$, $p = q = 2$, we have
\begin{align*}
\displaystyle
| {\varphi} - I_{{T}} {\varphi} |_{H^{1}({T})}
\leq  c s^{ 2 -  \varepsilon} \left | \frac{\partial \varphi}{\partial r_1} \right |_{H^{1}({T})} + s \left | \frac{\partial \varphi}{\partial r_2} \right |_{H^{1}({T})} .
\end{align*}
Even if one wants to reduce the step size in a specific direction ($y$-axis direction), the interpolation error may diverge as $s \to 0$ when $\varepsilon \> 2$. This loses the benefits of using anisotropic meshes.

Another is that violating the condition makes it challenging to show mathematical validity in the finite element method. One of the answers can be found in \cite{Apeet21}. That is,  the maximum-angle condition is sufficient to do numerical calculations safely.

\section{Lagrange Interpolation Error Estimates}

\subsection{One-dimensional Lagrange Interpolation}
Let $\Omega := (0,1) \subset \mathbb{R}$. For $N \in \mathbb{N}$, let $\mathbb{T}_h = \{  0 = x_0 \< x_1 \< \cdots \< x_N \< x_{N+1} = 1 \}$ be a mesh of $\overline{\Omega}$ such as
\begin{align*}
\displaystyle
\overline{\Omega} := \bigcup_{i=1}^N I_i, \quad \Int I_i \cap \Int I_j = \emptyset \quad \text{for $i \neq j$},
\end{align*}
where $I_i := [x_i , x_{i+1}]$ for $0 \leq i \leq N$. We denote $h_i := x_{i+1} - x_i $ for $0 \leq i \leq N$. For $\widehat{T} := [0,1] \subset \mathbb{R}$ and $\widehat{P} := \mathbb{P}^k$ with $k \in \mathbb{N}$, let $\{ \widehat{T} , \widehat{P} , \widehat{\Sigma} \}$ be the reference Lagrange finite element, e.g., see \cite{ErnGue04}. The corresponding interpolation operator is defined as
\begin{align*}
\displaystyle
I_{\widehat{T}}^k: \mathcal{C}(\widehat{T}) \ni \hat{v} \mapsto I_{\widehat{T}}^k (\hat{v}) := \sum_{m=0}^k \hat{v}(\hat{\xi}_m) \widehat{\mathcal{L}}_m^k,
\end{align*}
where $\hat{\xi}_m := \frac{m}{k}$ and $\{ \widehat{\mathcal{L}}_0^k , \ldots ,  \widehat{\mathcal{L}}_k^k  \}$ is the Lagrange polynomials associated with the nodes $\{ \hat{\xi}_0, \ldots, \hat{\xi}_k\}$. For $i \in \{ 0,\ldots,N\}$, we consider the affine transformations
\begin{align*}
\displaystyle
\Phi_i : \widehat{T} \ni t \mapsto x = x_i + t h_i \in I_i.
\end{align*}
For $\hat{v} \in \mathcal{C}(\widehat{T})$, we set $\hat{v} = v \circ \Phi_i$.

\begin{thr} \label{thr591}
Let $1 \leq p \leq \infty$ and assume that there exists a nonnegative integer $k$ such that
\begin{align*}
\displaystyle
\mathcal{P}^{k} = \widehat{{P}} \subset W^{k+1,p}(\widehat{T}) \subset \mathcal{C}(\widehat{T}).
\end{align*}
Let $\ell$ ($0 \leq \ell \leq k$) be such that $W^{\ell+1,p}(\widehat{T}) \subset \mathcal{C}(\widehat{T})$ with continuous embedding. Furthermore, assume that $\ell, m \in \mathbb{N} \cup \{ 0 \}$ and $p , q \in [1,\infty]$ such that $0 \leq m \leq \ell + 1$ and 
\begin{align*}
\displaystyle
W^{\ell +1,p}(\widehat{T}) \hookrightarrow W^{m,q} (\widehat{T}).
\end{align*}
It then holds that, for any $v \in W^{\ell+1,p}(I_i)$ with $\hat{v} = v \circ \Phi_i$,
\begin{align}
\displaystyle
 |v - {I}_{I_i}^k v |_{W^{m,q}(I_i)}
&\leq c h_{i}^{\frac{1}{q} - \frac{1}{p} + \ell+1-m} | v |_{W^{\ell+1,p}(I_i)}.\label{one59}
\end{align}	
\end{thr}

\begin{pf*}
We only show the outline of the proof. Scaling argument yields
\begin{align*}
\displaystyle
 |v - {I}_{I_i}^k v |_{W^{m,q}(I_i)}
 &=  h_i^{- m + \frac{1}{q}} | \hat{v} - {I}_{\hat{T}} \hat{v} |_{W^{m,q}(\widehat{T})}, \\
  |\hat{v}|_{W^{\ell+1,p}(\widehat{T})}
  &= h_i^{\ell+1 - \frac{1}{p}} |v|_{W^{\ell+1,p}(I_i)}.
\end{align*}
Using the Sobolev embedding theorem and the Bramble--Hilbert--type lemma, we have
\begin{align*}
\displaystyle
 | \hat{v} - {I}_{\hat{T}} \hat{v} |_{W^{m,q}(\widehat{T})}
 \leq c |\hat{v}|_{W^{\ell+1,p}(\widehat{T})}.
\end{align*}
Therefore, we obtain the estimate \eqref{one59}.
\qed
\end{pf*}

\begin{rem}
The assumptions of Theorem \ref{thr591} are standard; that is, there is no need to show the existence of functionals such as Theorem \ref{thr551}. Furthermore,  the quantity $h_{\max}/h_{\min}$ that deteriorates the convergent order does not appear in \eqref{one59}.
\end{rem}

\begin{rem}
If we set $x_j := \frac{j}{N+1}$, $j=0,1,\ldots,N,N+1$, the mesh  $\mathbb{T}_h$ is said to be the uniform mesh. If we set $x_j := g \left( \frac{j}{N+1} \right)$, $j=1,\ldots,N,N+1$ with a grading function $g$, the mesh $\mathbb{T}_h$ is said to be the graded mesh with respect to $x=0$, see \cite{BabSur94}. In particular, when one sets $g(y) := y^{\varepsilon}$ ($\varepsilon \> 0$), the mesh is called the radical mesh.
\end{rem}

\begin{rem}[Optimal order]
If $p = q$, it is possible to have the optimal error estimates even if the scale is different for each element. In the one-dimensional case, when $q \> p$, the convergence order of the interpolation operator may deteriorate, see Section \ref{sec=qp}.
\end{rem}

\subsection{Lagrange Finite Element}
Let $\widehat{T} \subset \mathbb{R}^d$ be the reference element defined in Sections \ref{reference2d} and \ref{reference3d}. Let $\alpha$ be a multi-index. For $k \in \mathbb{N}$, we define the set of Lagrange nodes as
\begin{align*}
\displaystyle
\mathscr{P} &:= \{ \widehat{p}_i \}_{i=1}^{N^{(2,k)}} := \left\{ \left(\frac{i_1}{k} ,\frac{i_2}{k} \right)^{\top} \in \mathbb{R}^2 \right\}_{0 \leq i_1 + i_2 \leq k} = \left\{ \frac{1}{k} \alpha \in \mathbb{R}^2 \right \}_{|\alpha| \leq k}, \quad \text{if $d=2$}, \\
\mathscr{P} &:=  \{ \widehat{p}_i \}_{i=1}^{N^{(3,k)}} := \widehat{T} \cap \left\{ \left(\frac{i_1}{k} ,\frac{i_2}{k} ,\frac{i_3}{k} \right)^{\top} \in \mathbb{R}^3 \right\}_{0 \leq i_1 , i_2, i_3 \leq k}, \quad \text{if $d=3$}.
\end{align*}
The Lagrange finite element on the reference element is defined by the triple $\{ \widehat{T}, \widehat{P}, \widehat{\Sigma\}}$ as follows.
\begin{enumerate}
 \item $\widehat{P} := \mathbb{P}^k(\widehat{T})$;
 \item $\widehat{\Sigma}$ is a set $\{ \hat{\chi}_{i} \}_{1 \leq i \leq N^{(d,k)}}$ of $N^{(d,k)}$ linear forms $\{ \hat{\chi}_{i} \}_{1 \leq i \leq N^{(d,k)}}$ with its components such that, for any $\hat{q} \in \widehat{P}$,
\begin{align}
\displaystyle
\hat{\chi}_{i}(\hat{q}) := \hat{q}( \widehat{p}_i) \quad \forall i \in \{ 1,\ldots, N^{(d,k)} \}. \label{Lag11}
\end{align}
\end{enumerate}
The nodal basis functions associated with the degrees of freedom by \eqref{Lag11} are defined as
\begin{align}
\displaystyle
\hat{\theta}_i (\widehat{p}_j) = \delta_{ij} \quad \forall i,j \in \{ 1,\ldots, N^{(d,k)} \}. \label{Lag12}
\end{align}
It then holds that $\hat{\chi}_{i} (\hat{\theta}_j) = \delta_{ij}$ for any $i,j \in \{ 1, \ldots ,d+1 \}$. Setting $V(\widehat{T}) := \mathcal{C}(\widehat{T})$ or $V(\widehat{T}) := W^{s,p}(\widehat{T})$ with $p \in [1,\infty]$ and $ps \> d$ ($s \geq d$ if $p=1$), the local operator ${I}_{\widehat{T}}^{L}$ is defined as
\begin{align}
\displaystyle
I_{\widehat{T}}^{L}: V(\widehat{T})  \ni \hat{\varphi}  \mapsto I_{\widehat{T}}^{L} \hat{\varphi} := \sum_{i=1}^{N^{(d,k)}} \hat{\varphi}(\widehat{p}_i) \hat{\theta}_i \in \widehat{P}. \label{Lag13}
\end{align}
By analogous argument in Section \ref{FEMG}, the Lagrange finite elements $\{ \widetilde{T} , \widetilde{P} , \widetilde{\Sigma}\}$ and $\{ T , {P} , \Sigma \}$ are constructed. The local shape functions are $\tilde{\theta}_{i} = \psi_{\widehat{T}}^{-1}(\hat{\theta}_i)$ and $\theta_{i} = \psi_{\widetilde{T}}^{-1}(\tilde{\theta}_i)$ for any  $i \in \{ 1, \ldots , N^{(d,k)}\}$, and the associated local interpolation operators are respectively defined as
\begin{align}
\displaystyle
{I}_{\widetilde{T}}^{L} : V(\widetilde{T}) \ni \tilde{\varphi} \mapsto {I}^{L}_{\widetilde{T}} \tilde{\varphi} &:=  \sum_{i=1}^{N^{(d,k)}} \tilde{\varphi}(\widetilde{p}_i) \tilde{\theta}_i  \in \widetilde{P}, \label{Lag14} \\
{I}_{T}^{L} : V(T) \ni \varphi \mapsto {I}^{L}_{T} \varphi &:= \sum_{i=1}^{N^{(d,k)}} \varphi(p_i) {\theta}_i  \in {P}, \label{Lag15}
\end{align}
where $\tilde{p}_i = {\Phi}_{\widetilde{T}}(\hat{p}_i)$, $p_i = {\Phi}_{T}(\tilde{p}_i)$ for $i \in \{ 1, \ldots, N^{(d,k)}\}$.

\subsection{Local Interpolation Error Estimates}
We first introduce the following lemmata.

\begin{lem}[$d=2$]
Let $\beta$ be a multi-index with $m := |\beta|$ and $\hat{\varphi} \in \mathcal{C}(\widehat{T})$ a function such that $\partial^{\beta}_{\hat{x}} \hat{\varphi} \in W^{\ell - m,p}(\widehat{T})$, where $\ell,m \in \mathbb{N}_0$, $p \in [1,\infty]$ are such that $0 \leq m \leq \ell \leq k+1$ and 
\begin{subequations} \label{Lag16}
\begin{align}
\displaystyle
&p=\infty \quad \text{if $m=0$ and $\ell=0$},  \label{Lag16a} \\
&p \> 2 \quad \text{if $m=0$ and $\ell=1$},  \label{Lag16b} \\
&m \< \ell \quad \text{if $\beta_1=0$ or $\beta_2 = 0$, and $m \> 0$}.  \label{Lag16c}
\end{align}
\end{subequations}
Fix $q \in [1,\infty]$ such that $W^{\ell - m,p}(\widehat{T}) \hookrightarrow L^q(\widehat{T})$. Let $I_{\widehat{T}} := I_{\widehat{T}}^L$.  It then holds that
\begin{align}
\displaystyle
\| \partial_{\hat{x}}^{\beta} (\hat{\varphi} - I_{\widehat{T}}^L \hat{\varphi}) \|_{L^q(\widehat{T})}
&\leq c | \partial_{\hat{x}}^{\beta}  \hat{\varphi} |_{W^{\ell  - m ,p}(\widehat{T})}.\label{Lag17}
\end{align}
\end{lem}

\begin{pf*}
We follow \cite[Lemma 2.4]{Ape99}. We first give proofs in some particular cases: $k=1,2$.

Let $k=1$. Let $m=0$, that is, $\beta = (0,0)$. We then have $j = \dim \mathbb{P}^1 = 3$.  From the Sobolev embedding theorem (Theorem \ref{thr=intro1}), we have $W^{\ell,p}(\widehat{T}) \subset \mathcal{C}^0(\widehat{T})$ with $1 \leq p \leq \infty$ and $2 \< \ell p$. Under this condition, we use
\begin{align*}
\displaystyle
\mathscr{F}_i(\hat{\varphi}) := \hat{\varphi} (\widehat{p}_i), \quad \hat{\varphi} \in W^{\ell,p}(\widehat{T}), \quad i=1,\ldots, 3.
\end{align*}
It then holds that
\begin{align*}
\displaystyle
|\mathscr{F}_i(\hat{\varphi})| \leq \| \hat{\varphi} \|_{\mathcal{C}^0(\widehat{T})} \leq c \| \hat{\varphi} \|_{W^{\ell,p}(\widehat{T})},
\end{align*}
which means $\mathscr{F}_i \in W^{\ell,p}(\widehat{T})^{\prime}$ for $i = 1,\ldots,3$, that is, \eqref{functionals=a} is satisfied. Furthermore, we have
\begin{align*}
\displaystyle
\mathscr{F}_i ( I_{\widehat{T}}^{L} \hat{\varphi} ) = ( I_{\widehat{T}}^{L} \hat{\varphi}) (\widehat{p}_i) =  \hat{\varphi} ({\widehat{p}_i}) = \mathscr{F}_i (\hat{\varphi}), \quad i = 1,\ldots,3,
\end{align*}
which satisfies \eqref{functionals=b}. For all $\hat{\eta} \in \mathbb{P}^{1}$, if $\mathscr{F}_i( \hat{\eta}) = 0$ for $i = 1,\ldots,3$, it obviously holds $\hat{\eta} = 0$. This means that \eqref{functionals=c} is satisfied.

Let $m=1$. We set $\beta = (1,0)$. We then have $j = \dim( \partial^{\beta} \mathbb{P}^1 ) = 1$. We consider a functional
\begin{align*}
\displaystyle
\mathscr{F}_1(\hat{\varphi}) := \int_{0}^{1} \hat{\varphi}(\hat{x}_1,0) d \hat{x}_1, \quad \hat{\varphi} \in W^{2,p}(\widehat{T}), \quad 1 \< p.
\end{align*}
We set $\widehat{I} := \{ \hat{x} \in \widehat{T}; \ \hat{x}_2 = 0 \}$. The continuity is then shown by the trace theorem (e.g., see Theorem \ref{thr=intro2}): if $1 = m \< \ell$,
\begin{align*}
\displaystyle
|\mathscr{F}_1(\hat{\varphi}) | \leq \| \hat{\varphi} \|_{L^1(\widehat{I})} \leq c \| \hat{\varphi} \|_{W^{\ell - 1,p}(\widehat{T})},
\end{align*}
which means $\mathscr{F}_1 \in W^{\ell -1,p}(\widehat{T})^{\prime}$, that is, \eqref{functionals=a} is satisfied. Furthermore, it holds that
\begin{align*}
\displaystyle
\mathscr{F}_1 (\partial^{(1,0)} (\hat{\varphi} - I_{\widehat{T}}^{L} \hat{\varphi}))
&= \int_{0}^{1} \frac{\partial}{\partial \hat{x}_1} (\hat{\varphi} - I_{\widehat{T}}^{L} \hat{\varphi})(\hat{x}_1,0) d \hat{x}_1 \\
&= \left[ \hat{\varphi} - I_{\widehat{T}}^{L} \hat{\varphi} \right]_{(0,0)}^{(1,0)} = 0,
\end{align*}
which satisfy \eqref{functionals=b}. Let $\hat{\eta} := a \hat{x}_1 + b \hat{x}_2 + c$. We then have
\begin{align*}
\displaystyle
\mathscr{F}_1 (\partial^{(1,0)} \hat{\eta}) = a.
\end{align*}
If $\mathscr{F}_1 (\partial^{(1,0)} \hat{\eta}) = 0$, $a = 0$. This implies that $\partial^{(1,0)} \hat{\eta} = 0$. This means that \eqref{functionals=c} is satisfied.

By analogous argument, the case $\beta = (0,1)$ holds.

Let $k=2$. Let $m=0$, that is, $\beta = (0,0)$. We then have $j = \dim \mathbb{P}^1 = 6$. Because $\dim \mathbb{P}^2 = 6$, we can show as in the case $k=1$ and $\beta = (0,0)$.

Let $\beta := (1,0)$. We define three functionals as
\begin{align*}
\displaystyle
\mathscr{F}_1(\hat{\varphi}) 
&:= \int_{0}^{\frac{1}{2}} \hat{\varphi}(\hat{x}_1,0) d \hat{x}_1,\\
\mathscr{F}_2(\hat{\varphi}) 
&:= \int_{\frac{1}{2}}^{1} \hat{\varphi}(\hat{x}_1,0) d \hat{x}_1,\\
\mathscr{F}_3(\hat{\varphi}) 
&:= \int_{0}^{\frac{1}{2}} \hat{\varphi}(\hat{x}_1, 1/2) d \hat{x}_1.
\end{align*}
We then show \eqref{functionals=a} and \eqref{functionals=b} as above. Let $\hat{\eta} \in \mathbb{P}^2$ be such that
\begin{align}
\displaystyle
\mathscr{F}_i(\partial_{\hat{x}}^{\beta} \hat{\eta}) = 0, \quad i=1,2,3. \label{Lag18}
\end{align}
We set the polynomial:
\begin{align}
\displaystyle
\hat{q} &:= \hat{\eta} - \hat{\eta}(1,0) \cdot 2 \left( \hat{x}_2 - \frac{1}{2}\right)(\hat{x}_2 - 1) - \hat{\eta}\left( \frac{1}{2},\frac{1}{2} \right) \cdot [ -4 \hat{x}_2 (\hat{x}_2 - 1)] \notag \\
&\quad - \hat{\eta}(0,1) \cdot 2 \hat{x}_2 \left( \hat{x}_2 - \frac{1}{2} \right) \in \mathcal{P}^2. \label{Lag19}
\end{align}
This has the following properties:
\begin{align}
\displaystyle
\partial_{\hat{x}} \hat{\eta} = \partial_{\hat{x}} \hat{q}, \quad \hat{q}(1,0) = \hat{q} \left( \frac{1}{2}, \frac{1}{2} \right) = \hat{q}(0,1) = 0. \label{Lag110}
\end{align}
We thus have
\begin{align*}
\displaystyle
0 = \mathscr{F}_3(\partial_{\hat{x}}^{\beta} \hat{\eta}) = \mathscr{F}_3(\partial_{\hat{x}}^{\beta} \hat{q}) = \hat{q} \left( \frac{1}{2}, \frac{1}{2} \right) - \hat{q} \left( 0, \frac{1}{2} \right),
\end{align*}
hence, $\hat{q} \left( 0, \frac{1}{2} \right) = 0$. By similar way, 
\begin{align*}
\displaystyle
0 &= \mathscr{F}_2(\partial_{\hat{x}}^{\beta} \hat{\eta}) = \mathscr{F}_2(\partial_{\hat{x}}^{\beta} \hat{q}) = \hat{q}(1,0) - \hat{q}\left( \frac{1}{2} , 0 \right),\\
0 &= \mathscr{F}_1(\partial_{\hat{x}}^{\beta} \hat{\eta}) = \mathscr{F}_1(\partial_{\hat{x}}^{\beta} \hat{q}) = \hat{q}\left( \frac{1}{2} , 0 \right) - \hat{q}(0,0),
\end{align*}
thus, $\hat{q} \left( \frac{1}{2} , 0 \right) = 0$ and $\hat{q}(0,0)  = 0$. Therefore, $\hat{q} \equiv 0$. Together with \eqref{Lag19}, we have $\hat{q} = \hat{q}(\hat{x}_2)$, $\partial_{\hat{x}}^{\beta} \hat{\eta} = 0$.
\qed
\end{pf*}

\begin{lem}[$d=3$]
Let $\beta$ be a multi-index with $m := |\beta|$ and $\hat{\varphi} \in \mathcal{C}(\widehat{T})$ a function such that $\partial^{\beta}_{\hat{x}} \hat{\varphi} \in W^{\ell - m,p}(\widehat{T})$, where $\ell,m \in \mathbb{N}_0$, $p \in [1,\infty]$ are such that $0 \leq m \leq \ell \leq k+1$ and 
\begin{subequations} \label{Lag18}
\begin{align}
\displaystyle
&p=\infty \quad \text{if $m=0$ and $\ell=0$},  \label{Lag18a} \\
&p \> \frac{3}{\ell} \quad \text{if $m=0$ and $\ell=1,2$},  \label{Lag18b} \\
&m \< \ell \quad \text{if $\beta_1=0$, $\beta_2 = 0$, or $\beta_3 = 0$},  \label{Lag18c} \\
&p \> 2 \quad \text{if $\beta \in \{ (\ell-1,0,0); (0,\ell-1,0); (0,0,\ell -1) \}$}. \label{Lag18d}
\end{align}
\end{subequations}
Fix $q \in [1,\infty]$ such that $W^{\ell - m,p}(\widehat{T}) \hookrightarrow L^q(\widehat{T})$. Let $I_{\widehat{T}} := I_{\widehat{T}}^L$. It then holds that
\begin{align}
\displaystyle
\| \partial_{\hat{x}}^{\beta} (\hat{\varphi} - I_{\widehat{T}} \hat{\varphi}) \|_{L^q(\widehat{T})}
&\leq c | \partial_{\hat{x}}^{\beta}  \hat{\varphi} |_{W^{\ell  - m ,p}(\widehat{T})}.\label{Lag19}
\end{align}

\end{lem}

\begin{pf*}
A proof can be found in \cite[Lemma 2.6]{Ape99}.
\qed
\end{pf*}

We have the following new Lagrange interpolation error estimates.

\begin{thr} \label{thr573}
Let $\{ \widehat{T} , \widehat{{P}} , \widehat{\Sigma} \}$ be the Lagrange finite element with  $V(\widehat{T}) := \mathcal{C}({\widehat{T}})$ and $\widehat{{P}} := \mathbb{P}^{k}(\widehat{T})$ with $k \geq 1$. Let $I_{\widehat{T}} := I_{\widehat{T}}^L$.  Let  $m \in \mathbb{N}_0$, $\ell \in \mathbb{N}$, and $p \in \mathbb{R}$ be such that $0 \leq m \leq \ell \leq k+1$ and
\begin{align*}
\displaystyle
&d=2: \ 
 \begin{cases}
p \in (2,\infty] \quad \text{if $m=0$, $\ell = 1$},\\
p \in [1,\infty] \quad \text{if $m=0$, $\ell \geq 2$ or $m \geq 1$, $\ell - m \geq 1$},
\end{cases} \\
&d=3: \ \begin{cases}
p \in \left(\frac{3}{\ell}, \infty \right] \quad \text{if $m=0$, $\ell=1,2$},\\
p \in (2,\infty] \quad \text{if $m \geq 1$,  $\ell-m = 1$}, \\
p \in [1,\infty] \quad \text{if $m=0$, $\ell \geq 3$ or $m \geq 1$,  $\ell-m \geq 2$}.
\end{cases}
\end{align*}
Setting $q \in [1,\infty)$ be such that
\begin{align}
\displaystyle
W^{\ell-m,p}(\widehat{T}) \hookrightarrow L^q(\widehat{T}), \label{Sobolev511}
\end{align}
that is $(\ell-m) - \frac{d}{p} \geq - \frac{d}{q}$. Then, for all $\hat{\varphi} \in W^{\ell ,p}(\widehat{T})$ with ${\varphi} := \hat{\varphi} \circ {\Phi}^{-1}$, we have
\begin{align}
\displaystyle
| {\varphi} - I_{{T}}^L {\varphi}|_{W^{m,q}({T})}
\leq  c |T|_d^{\frac{1}{q} - \frac{1}{p}} \left( \frac{H_{T}}{h_{T}} \right)^m \sum_{|\varepsilon| =  \ell-m} {h}^{\varepsilon} | \partial_r^{\varepsilon} \varphi |_{W^{m,p}(T)}. \label{Lag110}
\end{align}
In particular, if Condition \ref{Cond333} is imposed, it holds that,for all $\hat{\varphi} \in W^{\ell ,p}(\widehat{T})$ with ${\varphi} := \hat{\varphi} \circ {\Phi}^{-1}$,
\begin{align}
\displaystyle
| {\varphi} - I_{{T}}^L {\varphi}|_{W^{m,q}({T})}
\leq  C_2^{I} |T|_d^{\frac{1}{q} - \frac{1}{p}} \left( \frac{H_{T}}{h_{T}} \right)^m \sum_{|\varepsilon| = \ell-m} \widetilde{\mathscr{H}}^{\varepsilon}   | \partial_{\tilde{x}}^{\varepsilon}  ( \varphi \circ \Phi_{{T}}) |_{W^{m,p}(\Phi_{{T}}^{-1}(T))}. \label{Lag111}
\end{align}
Furthermore, for any $\hat{\varphi} \in \mathcal{C}(\widehat{T})$ with ${\varphi} := \hat{\varphi} \circ {\Phi}^{-1}$, it holds that
\begin{align*}
\displaystyle
\| {\varphi} - I_{{T}}^L {\varphi} \|_{L^{\infty}({T})}
\leq  c  \|  {\varphi} \|_{L^{\infty}(T)}.
\end{align*}
\end{thr}

\begin{pf*}
Proved in a similar way to Theorem \ref{thr551}.
\qed
\end{pf*}

\subsection{Global Interpolation Error Estimates}
Recall the space $V_h^n$ with $n=1$ (see \eqref{space=Vh}). We consider the space
\begin{align}
\displaystyle
V_h^L := \left\{ \varphi_h \in V_h^1: \ [\![ \varphi_h ]\!]_F = 0 \ \forall F \in \mathcal{F}_h^i  \right\} \subset H^1(\Omega). \label{space=Lag}
\end{align}
We also define the global interpolation $I_h^L$ to space $V_h^L$ as
\begin{align*}
\displaystyle
(I_h^L \varphi)|_{T} := I^L_{{T}} (\varphi|_{T}) \quad \forall T \in \mathbb{T}_h, \quad \forall \varphi \in \mathcal{C}(\overline{\Omega}).
\end{align*}

\begin{cor}
Suppose that the assumptions of Theorem \ref{thr573} are satisfied. We impose  Condition \ref{Cond1=sec6}. Let $I_h^L$ be the corresponding global Lagrange interpolation operator. It then holds that, for any $ \varphi \in W^{\ell ,p}(\Omega)$;
\begin{description}	
  \item[(\Roman{lone})] if Condition \ref{Cond333} is not imposed, 
\begin{align}
\displaystyle
| {\varphi} - I_{{h}}^L {\varphi}|_{W^{m,q}({\Omega})}
&\leq c  \sum_{T \in \mathbb{T}_h} |T|_d^{\frac{1}{q} - \frac{1}{p}} \sum_{|\varepsilon| =  \ell-m} {h}^{\varepsilon} | \partial_r^{\varepsilon} \varphi |_{W^{m,p}(T)}. \label{gint566}
\end{align}
   \item[(\Roman{ltwo})] if Condition \ref{Cond333} is imposed, 
\begin{align}
\displaystyle
| {\varphi} - I_{{h}}^L {\varphi}|_{W^{m,q}({\Omega})}
&\leq c  \sum_{T \in \mathbb{T}_h} |T|_d^{\frac{1}{q} - \frac{1}{p}} \sum_{|\varepsilon| =  \ell-m} \widetilde{\mathscr{H}}^{\varepsilon}  | \partial_{\tilde{x}}^{\varepsilon}   ( \varphi \circ \Phi_{T}) |_{W^{ m ,p}(\Phi_{T}^{-1}(T))}. \label{gint567}
\end{align}
\end{description}
\end{cor}

\begin{pf*}
This corollary is proved in the same argument as Corollary \ref{newglobal=cor}.
\qed
\end{pf*}

\section{$L^2$-orthogonal Projection}
This section considers error estimates of the $L^2$-orthogonal projection, e.g., for standard argument, see \cite[Section 1.4.3]{ErnGue04} and \cite[Section 11.5.3]{ErnGue21a}. 

\subsection{Finite Element}
Let $k \in \mathbb{N}_0$. Let $\widehat{T} \subset \mathbb{R}^d$ be the reference element defined in Section \ref{reference}. Let $\widehat{P}$ be a finite-dimensional space such that $\mathbb{P}^k \subset \widehat{P} \subset W^{k+1,\infty}(\widehat{T})$. The $L^2$-orthogonal projection onto $\widehat{P}$ is the linear operator ${\Pi}_{\widehat{T}}^{k}: L^1(\widehat{T}) \to \widehat{P}$ defined as
\begin{align}
\displaystyle
\int_{\widehat{T}} ({\Pi}_{\widehat{T}}^{k} \hat{\varphi} - \hat{\varphi}) \hat{q} d \hat{x} = 0 \quad \forall \hat{q} \in \widehat{P}, \quad \forall \hat{\varphi} \in L^1(\widehat{T}). \label{L2ortho=ref}
\end{align}
Because ${\Pi}_{\widehat{T}}^{k} \hat{\varphi} - \hat{\varphi}$ and ${\Pi}_{\widehat{T}}^{k} \hat{\varphi} - \hat{q}$ are $L^2$-orthogonal for any $\hat{q} \in \widehat{P}$, the Pythagorean identity yields 
\begin{align*}
\displaystyle
\| \hat{\varphi} - \hat{q} \|^2_{L^2(\widehat{T})}
= \| \hat{\varphi} - {\Pi}_{\widehat{T}}^{k} \hat{\varphi} \|^2_{L^2(\widehat{T})} + \| {\Pi}_{\widehat{T}}^{k} \hat{\varphi} - \hat{q} \|^2_{L^2(\widehat{T})}.
\end{align*}
This implies that
\begin{align*}
\displaystyle
\Pi_{\widehat{T}}^{k} \hat{\varphi} = \argmin_{\hat{q} \in \widehat{P}} \| \hat{\varphi} - \hat{q} \|_{L^2(\widehat{T})}.
\end{align*}
Therefore, $\widehat{P}$ is pointwise invariant under $\Pi_{\widehat{T}}^{k}$.

Let $\Phi_{\widetilde{T}}: \widehat{T} \to \widetilde{T}$ and $\Phi_{T}: \widetilde{T} \to T$ be the two affine mappings defined in Section \ref{two=step}. For any ${T} \in \mathbb{T}_h$ with $\widetilde{T} = {\Phi}_{\widetilde{T}}(\widehat{T})$ and $T ={\Phi}_{T} (\widetilde{T})$, let $\hat{\varphi} :=  \tilde{\varphi} \circ {\Phi_{\widetilde{T}}}$ and $ \tilde{\varphi} := \varphi \circ {\Phi_T}$. Furthermore, we set
\begin{align*}
\displaystyle
 \widetilde{P} &:= \{ \psi_{\widehat{T}}^{-1}(\hat{q}) ; \ \hat{q} \in \widehat{{P}}\}, \\
 {P} &:= \{ \psi_{\widetilde{T}}^{-1}(\tilde{q}) ; \ \tilde{q} \in \widetilde{{P}}\}.
\end{align*}
The $L^2$-orthognal projections onto $\widehat{P}$ and $P$ are respectively the linear operators $\Pi^k_{\widetilde{T}}: L^1(\widetilde{T}) \to \widetilde{P}$ and $\Pi^k_{{T}}: L^1({T}) \to {P}$ defined as
\begin{align*}
\displaystyle
&\int_{\widetilde{T}} ({\Pi}_{\widetilde{T}}^{k} \tilde{\varphi} - \tilde{\varphi}) \tilde{q} d \tilde{x} = 0 \quad \forall \tilde{q} \in \widetilde{P}, \quad \forall \tilde{\varphi} \in L^1(\widetilde{T}), \\
&\int_{{T}} ({\Pi}_{{T}}^{k} {\varphi} - {\varphi}) {q} d {x} = 0 \quad \forall {q} \in {P}, \quad \forall {\varphi} \in L^1({T}).
\end{align*}
Then,  $\widetilde{P}$ and $P$ are respectively pointwise invariant under $\Pi_{\widetilde{T}}^{k}$ and $\Pi_T^k$.

\subsection{Local Interpolation Error Estimates}
We have the following stability estimate of the projection $\Pi_{\widehat{T}}^k$.
\begin{lem}
Let $q \in [1,\infty)$. It holds that
\begin{align}
\displaystyle
\| \Pi_{\widehat{T}}^k \hat{\varphi} \|_{L^q(\widehat{T})} 
\leq c \| \hat{\varphi} \|_{L^q(\widehat{T})} \quad \forall \hat{\varphi} \in L^{q}(\widehat{T}). \label{chapp821}
\end{align}
\end{lem}

\begin{pf*}
Because all the norms in the finite-dimensional space $\widehat{P}$ are equivalent, there exist $\hat{c}_1$ and $\hat{c}_2$, depending on $\widehat{T}$, such that
\begin{align}
\displaystyle
\| \Pi_{\widehat{T}}^k \hat{\varphi} \|_{L^q(\widehat{T})} 
&\leq \hat{c}_1 \| \Pi_{\widehat{T}}^k \hat{\varphi} \|_{L^2(\widehat{T})}, \label{L2ortho=eqnorm1} \\
\| \Pi_{\widehat{T}}^k \hat{\varphi} \|_{L^{q^{\prime}}(\widehat{T})} 
&\leq  \hat{c}_2  \| \Pi_{\widehat{T}}^k \hat{\varphi} \|_{L^q(\widehat{T})},\label{L2ortho=eqnorm2}
\end{align}
where $\frac{1}{q} + \frac{1}{q^{\prime}} = 1$. Then,
\begin{align*}
\displaystyle
\| \Pi_{\widehat{T}}^k \hat{\varphi} \|_{L^q(\widehat{T})}^2
&\leq c \| \Pi_{\widehat{T}}^k \hat{\varphi} \|_{L^2(\widehat{T})}^2 = c \int_{\widehat{T}}  \hat{\varphi}  \Pi_{\widehat{T}}^k \hat{\varphi} d\hat{x} \\
&\leq c  \| \hat{\varphi} \|_{L^q(\widehat{T})}  \| \Pi_{\widehat{T}}^k \hat{\varphi} \|_{L^{q^{\prime}}(\widehat{T})} \\
&\leq c  \| \hat{\varphi} \|_{L^q(\widehat{T})}  \| \Pi_{\widehat{T}}^k \hat{\varphi} \|_{L^{q}(\widehat{T})},
\end{align*}
where we used \eqref{L2ortho=eqnorm1}, \eqref{L2ortho=eqnorm2}, \eqref{L2ortho=ref} with $\hat{q} := \Pi_{\widehat{T}}^k \hat{\varphi}$, and the H\"older's inequality with $\frac{1}{q} + \frac{1}{q^{\prime}} = 1$. This proves the target inequality.
\qed
\end{pf*}

The following theorem gives an anisotropic error estimate of the projection $\Pi_T^k$.

\begin{thr} \label{thr822}
For $k \in \mathbb{N}_0$, let $\ell \in \mathbb{N}_0$ be such that $0 \leq \ell \leq k$. Let $p \in [1,\infty)$ and $q \in [1,\infty)$ be such that
\begin{align}
\displaystyle
W^{1,p}({T}) \hookrightarrow L^q({T}), \label{Sobolev511}
\end{align}
that is $1 - \frac{d}{p} \geq - \frac{d}{q}$. It then holds that, for any $\hat{\varphi} \in W^{\ell+1,p}(\widehat{T})$ with ${\varphi} := \hat{\varphi} \circ {\Phi}^{-1}$,
\begin{align}
\displaystyle
\| \Pi_{T}^k \varphi - \varphi \|_{L^q(T)} \leq c |T|_d^{\frac{1}{q} - \frac{1}{p}} \sum_{|\epsilon| =  \ell+ 1} {h}^{\varepsilon} \| \partial_r^{\varepsilon} \varphi \|_{L^{p}(T)}. \label{chap822}
\end{align}
In particular, if Condition \ref{Cond333} is imposed, it holds that, for any $\hat{\varphi} \in W^{\ell+1,p}(\widehat{T})$ with ${\varphi} := \hat{\varphi} \circ {\Phi}^{-1}$,
\begin{align}
\displaystyle
\| \Pi_{T}^k \varphi - \varphi \|_{L^q(T)} \leq c |T|_d^{\frac{1}{q} - \frac{1}{p}} \sum_{|\epsilon| = \ell+ 1} \widetilde{\mathscr{H}}^{\varepsilon} \| \partial_{\tilde{x}}^{\varepsilon} (\varphi \circ \Phi_{{T}} ) \|_{L^{p}(\Phi_{{T}}^{-1}(T))}. \label{chap823}
\end{align}
\end{thr}

\begin{pf*}
Using the scaling argument, we have
\begin{align}
\displaystyle
\| \Pi_{T}^k \varphi - \varphi \|_{L^q(T)}
&= c |\det({A}_{\widetilde{T}})|^{\frac{1}{q}} \| \Pi_{\widehat{T}}^k \hat{\varphi} - \hat{\varphi} \|_{L^q(\widehat{T})}. \label{chap824}
\end{align}
where we used $|\det (A_T)| = 1$. For any $\hat{\eta} \in \mathbb{P}^{\ell} \subset \widehat{P}$ with $0 \leq \ell \leq k$, from the triangle inequality and $\Pi_{\widehat{T}}^k \hat{\eta} = \hat{\eta}$, we have
\begin{align}
\displaystyle
\| \Pi_{\widehat{T}}^k \hat{\varphi} - \hat{\varphi} \|_{L^q(\widehat{T})}
&\leq \| \Pi_{\widehat{T}}^k  (\hat{\varphi} - \hat{\eta}) \|_{L^q(\widehat{T})}  + \| \hat{\eta} - \hat{\varphi} \|_{L^q(\widehat{T})}. \label{chap825}
\end{align}
Using \eqref{chapp821} for the first term on the right-hand side of \eqref{chap825}, we have
\begin{align}
\displaystyle
\| \Pi_{\widehat{T}}^k  (\hat{\varphi} - \hat{\eta}) \|_{L^q(\widehat{T})} 
\leq c \| \hat{\varphi} - \hat{\eta} \|_{L^q(\widehat{T})}. \label{chap826}
\end{align}
Using the Sobolev embedding theorem for the second term on the right-hand side of \eqref{chap825} and \eqref{chap826}, we obtain
\begin{align}
\displaystyle
\| \hat{\varphi} - \hat{\eta} \|_{L^q(\widehat{T})} \leq c \| \hat{\varphi} - \hat{\eta} \|_{W^{1,p}(\widehat{T})}. \label{chap827}
\end{align}
Combining \eqref{chap824}, \eqref{chap825}, \eqref{chap826}, and \eqref{chap827}, we have
\begin{align}
\displaystyle
\| \Pi_{T}^k \varphi - \varphi \|_{L^q(T)}
&\leq c |\det({A}_{\widetilde{T}})|^{\frac{1}{q}} \inf_{\hat{\eta} \in \mathcal{P}^{\ell}} \| \hat{\varphi} - \hat{\eta} \|_{W^{1,p}(\widehat{T})}. \label{chap828}
\end{align}
From the Bramble--Hilbert-type lemma (e.g., see Subsection \ref{sec=BHL}), there exists a constant $\hat{\eta}_{\beta} \in \mathbb{P}^{\ell}$ such that, for any $\hat{\varphi} \in W^{\ell+1,p}(\widehat{T})$,
\begin{align}
\displaystyle
| \hat{\varphi} - \hat{\eta}_{\beta} |_{W^{t,p}(\widehat{T})} \leq C^{BH}(\widehat{T}) |\hat{\varphi}|_{W^{\ell+1,p}(\widehat{T})}, \quad t=0,1. \label{chap829}
\end{align}
If Condition \ref{Cond333} is not imposed, using \eqref{scaling2} ($m=0$) and \eqref{chap829}, we then have
\begin{align}
\displaystyle
\| \hat{\varphi} - \hat{\eta}_{\beta} \|_{W^{1,p}(\widehat{T})}
&\leq c |\hat{\varphi}|_{W^{\ell+1,p}(\widehat{T})} \notag \\
&\leq c |\det({A}_{\widetilde{T}})|^{-\frac{1}{p}} \sum_{|\epsilon| =  \ell+ 1} {h}^{\varepsilon} \| \partial_r^{\varepsilon} \varphi \|_{L^{p}(T)}.  \label{chap8210}
\end{align}
If Condition \ref{Cond333} is imposed, using \eqref{scaling3} ($m=0$) and \eqref{chap829}, we then have
\begin{align}
\displaystyle
\| \hat{\varphi} - \hat{\eta}_{\beta} \|_{W^{1,p}(\widehat{T})}
&\leq c |\hat{\varphi}|_{W^{\ell+1,p}(\widehat{T})} \notag \\
&\leq c |\det({A}_{\widetilde{T}})|^{-\frac{1}{p}} \sum_{|\epsilon| = \ell+ 1} \widetilde{\mathscr{H}}^{\varepsilon} \| \partial_{\tilde{x}}^{\varepsilon} \tilde{\varphi} \|_{L^{p}(\widetilde{T})}.  \label{chap8211}
\end{align}
Therefore, combining \eqref{chap828}, \eqref{chap8210}, and \eqref{chap8211} with \eqref{CN332}, we have \eqref{chap822} and \eqref{chap823} using $\widetilde{T} = \Phi^{-1}_{T}(T)$ and $\tilde{\varphi} = \varphi \circ \Phi_{T}$.
\qed
\end{pf*}

\subsection{Global Interpolation Error Estimates}
Recall the space $V_h^n$ with $n=1$ (see \eqref{space=Vh}). We define the standard discontinuous space as
\begin{align*}
\displaystyle
P_{dc,h}^{k} := V_h^1 = \left\{ p_h \in L^1(\Omega); \ p_h|_{T} \in P \quad \forall T \in \mathbb{T}_h  \right\}.
\end{align*}
We also define the global interpolation $\Pi_h^k$ to space $P_{dc,h}^{k}$ as
\begin{align*}
\displaystyle
(\Pi_h^k \varphi)|_{T} := \Pi_{{T}}^k (\varphi|_{T}) \quad \forall T \in \mathbb{T}_h, \quad \forall \varphi \in L^1(\Omega).
\end{align*}

\begin{cor}
Suppose that the assumptions of Theorem \ref{thr822} are satisfied. We impose Condition \ref{Cond1=sec6}. Let $\Pi_h^k$ be the corresponding global $L^2$-orthogonal projection. It then holds that, for any $ \varphi \in W^{\ell +1,p}(\Omega)$;
\begin{description}	
  \item[(\Roman{lone})] if Condition \ref{Cond333} is not imposed, 
\begin{align}
\displaystyle
\| \Pi_{h}^k \varphi - \varphi \|_{L^q(\Omega)}
&\leq c  \sum_{T \in \mathbb{T}_h} |T|_d^{\frac{1}{q} - \frac{1}{p}} \sum_{|\epsilon| =  \ell+ 1} {h}^{\varepsilon} \| \partial_r^{\varepsilon} \varphi \|_{L^{p}(T)}. \label{gint566}
\end{align}
   \item[(\Roman{ltwo})] if Condition \ref{Cond333} is imposed, 
\begin{align}
\displaystyle
\| \Pi_{h}^k \varphi - \varphi \|_{L^q(\Omega)}
&\leq c  \sum_{T \in \mathbb{T}_h} |T|_d^{\frac{1}{q} - \frac{1}{p}}  \sum_{|\epsilon| = \ell+ 1} \widetilde{\mathscr{H}}^{\varepsilon} \| \partial_{\tilde{x}}^{\varepsilon} (\varphi \circ \Phi_{{T}} ) \|_{L^{p}(\Phi_{{T}}^{-1}(T))}.  \label{gint567}
\end{align}
\end{description}
\end{cor}

\begin{pf*}
This corollary is proved in the same argument as Corollary \ref{newglobal=cor}.
\qed
\end{pf*}

\subsection{Another Estimate}

\begin{thr}  \label{thr943}
Let $T \subset \mathbb{R}^d$ be a simplex. Let $\Pi_T^0:L^2(T) \to \mathbb{P}^0(T)$ be the local $L^2$-projection defined by
\begin{align*}
\displaystyle
\Pi_{{T}}^{0} {\varphi} := \frac{1}{|{T}|_d} \int_{{T}} {\varphi} d {x} \quad \forall \varphi \in L^2(T).
\end{align*}
It then holds that 
\begin{align}
\displaystyle
\| \Pi^0_T \varphi - \varphi \|_{L^2(T)} &\leq \frac{h_T}{\pi} |\varphi|_{H^1(T)} \quad \forall \varphi \in H^1(T). \label{chap841}
\end{align}
\end{thr}

\begin{pf*}
For any $\varphi \in H^1(T)$, we set $w := \Pi^0_T \varphi - \varphi$. It then holds that
\begin{align*}
\displaystyle
\int_T w dx = \int_T ( \Pi^0_T \varphi - \varphi ) dx = \frac{1}{|T|_d} \int_T \varphi dx |T|_d - \int_T \varphi dx = 0.
\end{align*}
Therefore, using the Poincar\'e inequality \eqref{poincare}, we conclude \eqref{chap841}.
\qed	
\end{pf*}

\section{New Nonconforming FE Interpolation Error Estimates}

\subsection{Local Interpolation Error Estimates}
We introduce the following theorem using the error estimates of the $L^2$-orthogonal projection.

\begin{thr} \label{Thr2a}
Let $\alpha := (\alpha_1,\ldots,\alpha_d){^{\top}} \in \mathbb{N}_0^d$ be a multi-index and $k \in \mathbb{N}$. {Let $p \in [1,\infty)$ and $q \in [1,\infty)$ be such that \eqref{Sobolev511} holds.}
We define an interpolation operator $I_{T}: W^{k , p}(T) \to \mathbb{P}^{k}(T)$ that satisfies:
\begin{align}
\displaystyle
\partial_x^{\alpha} (I_{{T}} {\varphi}) = \Pi_T^0 (\partial_x^{\alpha} \varphi) \quad \forall \varphi \in W^{k, p}(T) \quad {\forall \alpha: \ |\alpha| \leq k}. \label{sp=cond} 
\end{align}
Then, for any $\hat{\varphi} \in W^{k+1,p}(\widehat{T})$ with ${\varphi} := \hat{\varphi} \circ {\Phi}^{-1}$ and any $\alpha$ with $|\alpha| \leq k$,
\begin{align}
\displaystyle
| I_{T} \varphi - \varphi |_{W^{k,q}({T})} &\leq  c |T|_d^{ \frac{1}{q} - \frac{1}{p} }  \sum_{i=1}^d  h_i \left | \frac{\partial\varphi}{\partial r_i} \right |_{W^{k,p}(T)}. \label{CR21=43} 
\end{align}
If Condition \ref{Cond333} is imposed, then:
\begin{align}
\displaystyle
| I_{T} \varphi - \varphi |_{W^{k,q}({T})} &\leq  c |T|_d^{ \frac{1}{q} - \frac{1}{p} }  \sum_{i=1}^d  \widetilde{\mathscr{H}}_i \left | \frac{\partial (\varphi \circ \Phi_T)}{\partial \tilde{x}_i} \right |_{W^{k,p}(\Phi_T^{-1}(T))}. \label{CR21=ass0} 
\end{align}
\end{thr}

\begin{pf*}
The error estimate of the $L^2$-orthogonal projection \eqref{chap822} with $\ell = 0$ yields
\begin{align*}
\displaystyle
| I_{T} \varphi - \varphi |_{W^{k,q}({T})}^q
&= \sum_{|\alpha| = k} \left \| \partial_x^{\alpha}  ( I_{T} \varphi - \varphi)  \right\|^q_{L^q(T)} \\
&=  \sum_{|\alpha| = k} \left \| \Pi_T^0 ( \partial_x^{\alpha} \varphi ) - \partial_x^{\alpha} \varphi \right\|^q_{L^q(T)} \\
&\leq c |T|_d^{\left( \frac{1}{q} - \frac{1}{p} \right) q} \sum_{|\alpha| = k} \sum_{i=1}^d  h_i^q \left \| \partial_x^{\alpha} \frac{\partial \varphi}{\partial r_i} \right \|_{L^p(T)}^q.
\end{align*}
Using the Jensen-type inequality \eqref{jensen}, we obtain
\begin{align*}
\displaystyle
| I_{T} \varphi - \varphi |_{W^{k,q}({T})}
&\leq c  |T|_d^{ \frac{1}{q} - \frac{1}{p}} \left( \sum_{i=1}^d  h_i^q \sum_{|\alpha| = k} \left \|  \partial_x^{\alpha} \frac{\partial\varphi}{\partial r_i} \right \|_{L^p(T)}^q \right)^{\frac{1}{q}} \\
&\leq c  |T|_d^{ \frac{1}{q} - \frac{1}{p}} \sum_{i=1}^d  h_i \left | \frac{\partial\varphi}{\partial r_i} \right |_{W^{k,p}(T)},
\end{align*}
which is the target inequality in Eq. \eqref{CR21=43}.

If Condition \ref{Cond333} is imposed, then \eqref{chap823} with $\ell = 0$ yields
\begin{align*}
\displaystyle
| I_{T} \varphi - \varphi |_{W^{k,q}({T})}^q
&\leq c |T|_d^{\left( \frac{1}{q} - \frac{1}{p} \right) q} \sum_{|\alpha| = k} \sum_{i=1}^d  \widetilde{\mathscr{H}}_i^q \left \| \partial_{\tilde{x}}^{\alpha} \frac{\partial (\varphi \circ \Phi_T)}{\partial \tilde{x}_i} \right \|_{L^p(\Phi_T^{-1}(T))}^q.
\end{align*}
Using the Jensen-type inequality in \eqref{jensen}, we obtained the target inequality in \eqref{CR21=ass0}.
\qed
\end{pf*}

\begin{note}
The operators that satisfy \eqref{sp=cond} exist; see Theorems \ref{thr=CR} and \ref{thr2=Morley}.
\end{note}

\subsection{CR Finite Element}
Let $\widehat{T} \subset \mathbb{R}^d$ be the reference element defined in Sections \ref{reference2d} and \ref{reference3d}. Let $\widehat{F}_i$ be the face of $\widehat{T}$ opposite to $\widehat{p}_i$. The CR finite element on the reference element is defined by the triple $\{ \widehat{T} , \widehat{P} , \widehat{\Sigma\}}$ as follows.
\begin{enumerate}
 \item $\widehat{P} := \mathbb{P}^1(\widehat{T})$;
 \item $\widehat{\Sigma}$ is a set $\{ \hat{\chi}_{i} \}_{1 \leq i \leq N^{(d,1)}}$ of $N^{(d,1)}$ linear forms $\{ \hat{\chi}_{i} \}_{1 \leq i \leq N^{(d,1)}}$ with its components such that, for any $\hat{q} \in \widehat{P}$,
\begin{align}
\displaystyle
\hat{\chi}_{i}(\hat{q}) := \frac{1}{| \widehat{F}_i |_{d-1}} \int_{\widehat{F}_i} \hat{q} d \hat{s} \quad \forall i \in \{ 1, \ldots ,d+1 \}. \label{CR911}
\end{align}
\end{enumerate}
The nodal basis functions associated with the degrees of freedom by \eqref{CR911} are defined as
\begin{align}
\displaystyle
\hat{\theta}_i(\hat{x}) := d \left( \frac{1}{d} - \hat{\lambda}_i(\hat{x}) \right) \quad \forall i \in \{ 1, \ldots,d+1 \}. \label{CR912}
\end{align}
It then holds that $\hat{\chi}_{i} (\hat{\theta}_j) = \delta_{ij}$ for any $i,j \in \{ 1, \ldots,d+1 \}$. Setting $V(\widehat{T}) := W^{1,1}(\widehat{T})$, the local operator ${I}_{\widehat{T}}^{CR}$ is defined as
\begin{align}
\displaystyle
I_{\widehat{T}}^{CR}: V(\widehat{T})  \ni \hat{\varphi}  \mapsto I_{\widehat{T}}^{CR} \hat{\varphi} := \sum_{i=1}^{d+1} \left( \frac{1}{| \widehat{F}_i|_{d-1}} \int_{\widehat{F}_i} \hat{\varphi} d \hat{s} \right)  \hat{\theta}_i \in \widehat{P}. \label{CR913}
\end{align}
By analogous argument in Section \ref{FEMG}, the CR finite elements $\{ \widetilde{T} , \widetilde{P} , \widetilde{\Sigma}\}$ and $\{ T , {P} , \Sigma \}$ are constructed. The local shape functions are $\tilde{\theta}_{i} = \psi_{\widehat{T}}^{-1}(\hat{\theta}_i)$ and $\theta_{i} = \psi_{\widetilde{T}}^{-1}(\tilde{\theta}_i)$ for any  $i \in \{ 1, \ldots , d+1\}$, and the associated local interpolation operators are respectively defined as
\begin{align}
\displaystyle
{I}_{\widetilde{T}}^{CR} : V(\widetilde{T}) \ni \tilde{\varphi} \mapsto {I}^{CR}_{\widetilde{T}} \tilde{\varphi} &:=  \sum_{i=1}^{d+1} \left( \frac{1}{| \widetilde{F}_i|_{d-1}} \int_{\widetilde{F}_i} \tilde{\varphi} d \tilde{s} \right)  \tilde{\theta}_i  \in \widetilde{P}, \label{CR914} \\
{I}^{CR}_{T} : V(T) \ni \varphi \mapsto {I}^{CR}_{T} \varphi &:=  \sum_{i=1}^{d+1} \left( \frac{1}{| {F}_{i}|_{d-1}} \int_{{F}_{i}} {\varphi} d{s} \right)  {\theta}_{i} \in {P}, \label{CR915T}
\end{align}
where $\{ \widetilde{F}_i :=  {\Phi}_{\widetilde{T}} (\widehat{F}_i)  \}_{i \in \{ 1, \ldots ,d+1\}}$ and $\{ F_{i} := {\Phi}_{T}(\widetilde{F}_i ) \}_{i \in \{ 1, \ldots,d+1\}}$.

\subsection{Local CR Interpolation Error Estimates}
We present anisotropic CR interpolation error estimates.

\begin{thr} \label{thr=CR}
Let $p \in [1,\infty)$ and $q \in [1,\infty)$ be such that \eqref{Sobolev511} holds. Then, 
\begin{align}
\displaystyle
| I_{T}^{CR} \varphi - \varphi |_{W^{1,q}({T})} &\leq  c |T|_d^{ \frac{1}{q} - \frac{1}{p} }  \sum_{i=1}^d  h_i \left | \frac{\partial\varphi}{\partial r_i} \right |_{W^{1,p}(T)}  \quad \forall {\varphi} \in W^{2,p}({T}), \label{CR921} \\
\| I_{T}^{CR} \varphi - \varphi \|_{L^{q}({T})} &\leq c |T|_d^{\frac{1}{q} - \frac{1}{p}} \sum_{i=1}^d {h}_i \left \| \frac{\partial \varphi}{\partial r_i} \right \|_{L^p(T)} \quad \forall {\varphi} \in W^{1,p}({T}), \label{CR922} \\
\| I_{T}^{CR} \varphi - \varphi \|_{L^{q}({T})} &\leq c |T|_d^{\frac{1}{q} - \frac{1}{p}} \sum_{|\varepsilon| = 2} h^{\varepsilon} \left\| \partial_{r}^{\varepsilon} \varphi  \right\|_{L^{p}(T)} \quad \forall {\varphi} \in W^{2,p}({T}). \label{CR922=b}
\end{align}
If Condition \ref{Cond333} is imposed, then:
\begin{align}
\displaystyle
| I_{T}^{CR} \varphi - \varphi |_{W^{1,q}({T})} &\leq  c |T|_d^{ \frac{1}{q} - \frac{1}{p} }  \sum_{i=1}^d  \widetilde{\mathscr{H}}_i \left | \frac{\partial (\varphi \circ \Phi_T)}{\partial \tilde{x}_i} \right |_{W^{1,p}(\Phi_T^{-1}(T))} \quad \forall \varphi \in W^{2,p}(T), \label{CR923} \\
\| I_{T}^{CR} \varphi - \varphi \|_{L^{q}(T)}
&\leq c |T|_d^{\frac{1}{q} - \frac{1}{p}} \sum_{i=1}^d \widetilde{\mathscr{H}}_i \left \| \frac{\partial }{\partial \tilde{x}_i} (\varphi \circ \Phi_{T}) \right \|_{L^p(\Phi_{T}^{-1}(T))} \quad \forall \varphi \in W^{1,p}(T), \label{CR924} \\
\| I_{T}^{CR} \varphi - \varphi \|_{L^{q}(T)}
&\leq c |T|_d^{\frac{1}{q} - \frac{1}{p}}  \sum_{|\varepsilon| = 2} \widetilde{\mathscr{H}}^{\varepsilon}  \| \partial_{\tilde{x}}^{\varepsilon} (  {\varphi} \circ \Phi_{T}) \|_{L^{p}({\Phi_{T}^{-1}(T)})} \quad \forall \varphi \in W^{2,p}(T), \label{CR924=b}
\end{align}
\end{thr}

\begin{pf*}
Only CR interpolation satisfies the condition \eqref{sp=cond} to prove the estimate \eqref{CR921} and \eqref{CR923}.

For ${\varphi} \in W^{2,p}({T})$, Green's formula and the definition of the CR interpolation imply that because $I_{T}^{CR} \varphi  \in \mathbb{P}^1$, 
\begin{align*}
\displaystyle
\frac{\partial }{\partial {x}_j} (I_{T}^{CR} \varphi)
&= \frac{1}{|T|} \int_{T}  \frac{\partial }{\partial {x}_j} (I_{T}^{CR} \varphi) dx = \frac{1}{|T|} \sum_{i=1}^{d+1} n_{{T}}^{(j)} \int_{F_i} I_{T}^{CR} \varphi ds \\
&= \frac{1}{|T|} \sum_{i=1}^{d+1} n_{{T}}^{(j)} \int_{F_i} \varphi ds 
= \frac{1}{|T|} \int_{T}  \frac{\partial \varphi}{\partial {x}_j} dx = \Pi_{T}^0 \left(  \frac{\partial \varphi}{\partial {x}_j} \right)
\end{align*}
for  $j=1,\ldots,d$, where $n_{{T}}$ denotes the outer unit normal vector to $T$ and $ n_{{T}}^{(j)}$ denotes the $j$th component of $n_T$. Therefore, from Theorem \ref{Thr2a}, the target inequalities \eqref{CR921} and \eqref{CR923} hold.

The (standard) scaling argument with $|\det(A_T)| = 1$ yields
\begin{align}
\displaystyle
\| I_{T}^{CR} \varphi - \varphi \|_{L^q(T)} 
&= |\det({A}_{\widetilde{T}})|^{\frac{1}{q}} \| I_{\widehat{T}}^{CR} \hat{\varphi} - \hat{\varphi} \|_{L^q(\widehat{T})}. \label{CR925}
\end{align}
For any $\hat{\eta} \in \mathbb{P}^{\ell}$ with $\ell \in \{ 0,1\}$, from the triangle inequality and $I_{\widehat{T}}^{CR} \hat{\eta} = \hat{\eta}$, we have
\begin{align}
\displaystyle
\| I_{\widehat{T}}^{CR} \hat{\varphi} - \hat{\varphi} \|_{L^q(\widehat{T})}
&\leq \| I_{\widehat{T}}^{CR}  (\hat{\varphi} - \hat{\eta}) \|_{L^q(\widehat{T})} + \| \hat{\eta} - \hat{\varphi} \|_{L^q(\widehat{T})}. \label{CR926}
\end{align}
Using the definition of the CR interpolation \eqref{CR913} and the trace theorem, we have
\begin{align}
\displaystyle
\| I_{\widehat{T}}^{CR}  (\hat{\varphi} - \hat{\eta}) \|_{L^q(\widehat{T})}
&\leq \sum_{i=1}^{d+1} \frac{1}{|\widehat{F}_i|_{d-1}} \int_{\widehat{F}_i} |\hat{\varphi} - \hat{\eta}| d\hat{s} \| \hat{\theta}_i \|_{L^q(\widehat{T})} \leq c \| \hat{\varphi} - \hat{\eta} \|_{W^{1,p}(\widehat{T})}. \label{CR927}
\end{align}
Using the Sobolev embedding theorem for the second term on the right-hand side of \eqref{CR926}, we obtain
\begin{align}
\displaystyle
\| \hat{\varphi} - \hat{\eta} \|_{L^q(\widehat{T})} \leq c \| \hat{\varphi} - \hat{\eta} \|_{W^{1,p}(\widehat{T})}. \label{cr=int=new1}
\end{align}
Combining \eqref{CR925}, \eqref{CR926}, \eqref{CR927} and \eqref{cr=int=new1}, we have
\begin{align}
\displaystyle
\| I_{T}^{CR} \varphi - \varphi \|_{L^q(T)}
&\leq C(\widehat{T}) |\det({A}_{\widetilde{T}})|^{\frac{1}{q}} \inf_{\hat{\eta} \in \mathbb{P}^{\ell}} \| \hat{\varphi} - \hat{\eta} \|_{W^{1,p}(\widehat{T})}. \label{CR928}
\end{align}
From the Bramble--Hilbert-type lemma (e.g., see Subsection \ref{sec=BHL}), there exists a constant $\hat{\eta}_{\beta} \in \mathbb{P}^{\ell}$ with $\ell \in \{ 0,1\}$ such that, for any $\hat{\varphi} \in W^{\ell+1,p}(\widehat{T})$,
\begin{align}
\displaystyle
| \hat{\varphi} - \hat{\eta}_{\beta} |_{W^{t,p}(\widehat{T})} \leq C^{BH}(\widehat{T}) |\hat{\varphi}|_{W^{\ell+1,p}(\widehat{T})}, \quad t=0,1. \label{chap829}
\end{align}
Thus, from \eqref{CR928} and \eqref{chap829},
\begin{align}
\displaystyle
\| I_{T}^{CR} \varphi - \varphi \|_{L^q(T)}
&\leq C(\widehat{T}) |\det({A}_{\widetilde{T}})|^{\frac{1}{q}} \inf_{\hat{\eta} \in \mathbb{P}^{\ell}} \| \hat{\varphi} - \hat{\eta} \|_{W^{1,p}(\widehat{T})} \notag \\
&\leq C(\widehat{T}) |\det({A}_{\widetilde{T}})|^{\frac{1}{q}}  \| \hat{\varphi} - \hat{\eta}_{\beta} \|_{W^{1,p}(\widehat{T})} \notag \\
&\leq C_1(\widehat{T})  |\det({A}_{\widetilde{T}})|^{\frac{1}{q}} |\hat{\varphi}|_{W^{\ell+1,p}(\widehat{T})}, \quad \ell = 0,1. \label{cr=int=new2}
\end{align}
Therefore, using \eqref{CN332}, \eqref{cr=int=new2} and \eqref{scaling2},
\begin{align*}
\displaystyle
\| I_{T}^{CR} \varphi - \varphi \|_{L^q(T)}
&\leq c |T|_{d}^{\frac{1}{q} - \frac{1}{p}}  \sum_{|\epsilon| =  \ell+1} {h}^{\varepsilon} \| \partial_r^{\epsilon} \varphi \|_{L^{p}(T)}, \quad \ell = 0,1,
\end{align*}
which leads to \eqref{CR922} and \eqref{CR922=b}. In particular, if Condition \ref{Cond333} is imposed, using \eqref{CN332}, \eqref{cr=int=new2} and \eqref{scaling3},
\begin{align*}
\displaystyle
\| I_{T}^{CR} \varphi - \varphi \|_{L^q(T)}
&\leq c |T|_{d}^{\frac{1}{q} - \frac{1}{p}}  \sum_{|\epsilon| =  \ell+1} \widetilde{\mathscr{H}}^{\varepsilon} \| \partial_{\tilde{x}}^{\varepsilon} \tilde{\varphi} \|_{L^{p}(\widetilde{T})}, \quad \ell = 0,1,
\end{align*}
which leads to \eqref{CR924} and \eqref{CR924=b}.
\qed
\end{pf*}

\subsection{Global CR Interpolation Error Estimates}
Recall the space $V_h^n$ with $n=1$ (see \eqref{space=Vh}). We define the CR finite element space as
\begin{align*}
\displaystyle
V_{h}^{CR} :=  \left\{ \varphi_h \in V_h^1; \ \int_F [\![ \varphi_h ]\!]_F ds = 0 \ \forall F \in \mathcal{F}_h^i  \right\}.
\end{align*}
We also define the global interpolation $I_h^{CR}: W^{1,1}(\Omega) \to V_{h}^{CR}$ as follows.
\begin{align*}
\displaystyle
(I_h^{CR} \varphi)|_{T} := I_{T}^{CR} (\varphi|_{T}) = \sum_{i=1}^{d+1} \left( \frac{1}{| {F}_{i}|_{d-1}} \int_{{F}_{i}} {\varphi} |_{T} d{s} \right)  {\theta}_{i} \quad \forall T \in \mathbb{T}_h, \ \forall \varphi \in W^{1,1}(\Omega).
\end{align*}

\begin{cor} \label{thr931}
Suppose that the assumptions of Theorem \ref{thr=CR} are satisfied. Let $I_h^{CR}$ be the corresponding global CR interpolation operator. Then,
\begin{align}
\displaystyle
| {\varphi} - I^{CR}_{{h}} {\varphi}|_{W^{1,q}({\mathbb{T}_h})}
&\leq c  \sum_{T \in \mathbb{T}_h} |T|_d^{ \frac{1}{q} - \frac{1}{p} }  \sum_{i=1}^d  h_i \left | \frac{\partial\varphi}{\partial r_i} \right |_{W^{1,p}(T)}  \quad \forall {\varphi} \in W^{2,p}({\Omega}), \label{CR931}\\
\| {\varphi} - I^{CR}_{{h}} {\varphi} \|_{L^{q}({\Omega})} 
&\leq c  \sum_{T \in \mathbb{T}_h}  |T|_d^{\frac{1}{q} - \frac{1}{p}} \sum_{i=1}^d {h}_i \left \| \frac{\partial \varphi}{\partial r_i} \right \|_{L^p(T)} \quad \forall {\varphi} \in W^{1,p}({\Omega}), \label{CR932} \\
\| {\varphi} - I^{CR}_{{h}} {\varphi} \|_{L^{q}({\Omega})} 
&\leq c  \sum_{T \in \mathbb{T}_h}  |T|_d^{\frac{1}{q} - \frac{1}{p}} \sum_{|\varepsilon| = 2} h^{\varepsilon} \left\| \partial_{r}^{\varepsilon} \varphi  \right\|_{L^{p}(T)} \quad \forall {\varphi} \in W^{2,p}({\Omega}). \label{CR932=b}
\end{align}
If Condition \ref{Cond333} is imposed, then:
\begin{align}
\displaystyle
| {\varphi} - I^{CR}_{{h}} {\varphi}|_{W^{1,q}({\mathbb{T}_h})}
&\leq c  \sum_{T \in \mathbb{T}_h} |T|_d^{ \frac{1}{q} - \frac{1}{p} }  \sum_{i=1}^d  \widetilde{\mathscr{H}}_i \left | \frac{\partial (\varphi \circ \Phi_T)}{\partial \tilde{x}_i} \right |_{W^{1,p}(\Phi_T^{-1}(T))} \quad \forall \varphi \in W^{2,p}(\Omega), \label{CR933}\\
\| {\varphi} - I^{CR}_{{h}} {\varphi} \|_{L^{q}({\Omega})} 
&\leq c  \sum_{T \in \mathbb{T}_h} |T|_d^{\frac{1}{q} - \frac{1}{p}} \sum_{i=1}^d \widetilde{\mathscr{H}}_i \left \| \frac{\partial }{\partial \tilde{x}_i} (\varphi \circ \Phi_{T}) \right \|_{L^p(\Phi_{T}^{-1}(T))} \quad \forall {\varphi} \in W^{1,p}({\Omega}), \label{CR934} \\
\| {\varphi} - I^{CR}_{{h}} {\varphi} \|_{L^{q}({\Omega})} 
&\leq c  \sum_{T \in \mathbb{T}_h} |T|_d^{\frac{1}{q} - \frac{1}{p}}  \sum_{|\varepsilon| = 2} \widetilde{\mathscr{H}}^{\varepsilon}  \| \partial_{\tilde{x}}^{\varepsilon} (  {\varphi} \circ \Phi_{T}) \|_{L^{p}({\Phi_{T}^{-1}(T)})} \quad \forall {\varphi} \in W^{2,p}({\Omega}). \label{CR934=b}
\end{align}
\end{cor}

\begin{pf*}
This corollary is proved in the same argument as Corollary \ref{newglobal=cor}.
\qed
\end{pf*}

\subsection{Another Estimate}

\begin{thr}
Let $T \subset \mathbb{R}^d$ be a simplex. Let $I_T^{CR}:H^1(T) \to \mathbb{P}^1(T)$ be the local CR interpolation operator defined as
\begin{align*}
\displaystyle
I_{T}^{CR}: H^{1}(T) \ni \varphi  \mapsto I_{T}^{CR} \varphi := \sum_{i=1}^{d+1} \left( \frac{1}{|F_i|_{d-1}} \int_{F_i} \varphi ds \right)  \theta_i \in \mathbb{P}^1. 
\end{align*}
It then holds that 
\begin{align}
\displaystyle
| I_{T}^{CR} \varphi - \varphi |_{H^1(T)}  &\leq \frac{h_{T}}{\pi} |\varphi|_{H^2(T)} \quad \forall \varphi \in H^2(T). \label{CR941}
\end{align}

\end{thr}

\begin{pf*}
Using \eqref{chap841}.
\begin{align}
\displaystyle
| I_{T}^{CR} \varphi - \varphi |_{H^{1}(T)}^2
&= \sum_{j=1}^d \left \|  \frac{\partial}{\partial x_j} ( I_{T}^{CR} \varphi - \varphi)  \right\|^2_{L^2(T)} \notag \\
&= \sum_{j=1}^d \left \| \Pi_T^0 \left(  \frac{\partial \varphi}{\partial {x}_j} \right) - \left(  \frac{\partial \varphi}{\partial {x}_j} \right) \right\|^2_{L^2(T)} \notag \\
&\leq  \left( \frac{h_{T}}{\pi} \right)^2  \sum_{i,j=1}^d  \left \| \frac{\partial^2 \varphi}{ \partial x_i \partial x_j} \right \|_{L^{2}(T)}^2 \notag \\
&= \left( \frac{h_{T}}{\pi} \right)^2 |\varphi|^2_{H^2(T)}, \notag
\end{align}
which conclude \eqref{CR941}.
\qed
\end{pf*}

\subsection{Nodal CR Interpolation Error Estimates}
Let $\widehat{T} \subset \mathbb{R}^d$ be the reference element defined in Sections \ref{reference2d} and \ref{reference3d}. Let $\widehat{F}_i$ be the face of $\widehat{T}$ opposite to $\widehat{p}_i$ and let $\hat{x}_{\widehat{F}_i}$ the barycentre of the face $\widehat{F}_i$. The (nodal) CR finite element on the reference element is defined by the triple $\{ \widehat{T} , \widehat{P} , \widehat{\Sigma\}}$ as follows.
\begin{enumerate}
 \item $\widehat{P} := \mathbb{P}^1(\widehat{T})$;
 \item $\widehat{\Sigma}$ is a set $\{ \hat{\chi}_{i} \}_{1 \leq i \leq N^{(d,1)}}$ of $N^{(d,1)}$ linear forms $\{ \hat{\chi}_{i} \}_{1 \leq i \leq N^{(d,1)}}$ with its components such that, for any $\hat{q} \in \widehat{P}$,
\begin{align}
\displaystyle
\hat{\chi}_{i}(\hat{p}) := \hat{q}(\hat{x}_{\widehat{F}_i}) \quad \forall i \in \{ 1, \ldots ,d+1 \}. \label{nCR1}
\end{align}
\end{enumerate}
The nodal basis functions associated with the degrees of freedom by \eqref{nCR1} are defined as
\begin{align}
\displaystyle
\hat{\theta}_i(\hat{x}) := d \left( \frac{1}{d} - \hat{\lambda}_i(\hat{x}) \right) \quad \forall i \in \{ 1, \ldots ,d+1 \}. \label{nCR2}
\end{align}
It then holds that $\hat{\chi}_{i} (\hat{\theta}_j) = \delta_{ij}$ for any $i,j \in \{ 1, \ldots ,d+1 \}$. Setting $V(\widehat{T}) := \mathcal{C}(\widehat{T})$ or $V(\widehat{T}) := W^{s,p}(\widehat{T})$ with $p \in [1,\infty]$ and $ps \> d$ ($s \geq d$ if $p=1$), the local operator ${I}_{\widehat{T}}^{nCR}$ is defined as
\begin{align}
\displaystyle
I_{\widehat{T}}^{nCR}: V(\widehat{T})  \ni \hat{\varphi}  \mapsto I_{\widehat{T}}^{nCR} \hat{\varphi} := \sum_{i=1}^{d+1} \hat{\varphi}(\hat{x}_{\widehat{F}_i}) \hat{\theta}_i \in \widehat{P}. \label{nCR3}
\end{align}
By analogous argument in Section \ref{FEMG}, the nodal CR finite elements $\{ \widetilde{T} , \widetilde{P} , \widetilde{\Sigma}\}$ and $\{ T , {P} , \Sigma \}$ are constructed. The local shape functions are $\tilde{\theta}_{i} = \psi_{\widehat{T}}^{-1}(\hat{\theta}_i)$ and $\theta_{i} = \psi_{\widetilde{T}}^{-1}(\tilde{\theta}_i)$ for any  $i \in \{ 1, \ldots , d+1 \}$, and the associated local interpolation operators are respectively defined as
\begin{align}
\displaystyle
{I}_{\widetilde{T}}^{nCR} : V(\widetilde{T}) \ni \tilde{\varphi} \mapsto {I}^{nCR}_{\widetilde{T}} \tilde{\varphi} &:=  \sum_{i=1}^{d+1} \tilde{\varphi}(\tilde{x}_{\widetilde{F}_i}) \tilde{\theta}_i  \in \widetilde{P}, \label{nCR4} \\
{I}_{T}^{nCR} : V(T) \ni \varphi \mapsto {I}^{nCR}_{T} \varphi &:= \sum_{i=1}^{d+1} \varphi(x_{F_i}) {\theta}_i  \in {P}, \label{nCR5T}
\end{align}
where $\{ \widetilde{F}_i :=  {\Phi}_{\widetilde{T}} (\widehat{F}_i)  \}_{i \in \{ 1, \ldots ,d+1\}}$,  $\{ F_{i} := {\Phi}_{T}(\widetilde{F}_i ) \}_{i \in \{ 1, \ldots,d+1\}}$,  $\tilde{x}_{\widetilde{F}_i} = {\Phi}_{\widetilde{T}}(\hat{x}_{\widehat{F}_i})$, $x_{F_i} = {\Phi}_{T}(\tilde{x}_{\widetilde{F}_i})$ for $i \in \{ 1, \ldots, d+1\}$.

\begin{cor}
Let $\{ \widehat{T} , \widehat{{P}} , \widehat{\Sigma} \}$ be the Crouzeix--Raviart finite element with  $V(\widehat{T}) := \mathcal{C}({\widehat{T}})$ and $\widehat{{P}} := \mathbb{P}^{1}(\widehat{T})$. Set $I_{\widehat{T}} ;= I_{\widehat{T}}^{nCR}$.  Let  $m \in \mathbb{N}_0$, $\ell \in \mathbb{N}$, and $p \in \mathbb{R}$ be such that
\begin{align*}
\displaystyle
&d=2: \ 
 \begin{cases}
p \in (2,\infty]  \quad \text{if $m=0$, $\ell = 1$},\\
p \in [1,\infty] \quad \text{if $m=0$, $\ell = 2$ or $m=1$, $\ell = 2$},
\end{cases} \\
&d=3: \ \begin{cases}
p \in \left(\frac{3}{\ell} , \infty \right] \quad \text{if $m=0$, $\ell=1,2$},\\
\displaystyle
p \in (2,\infty] \quad \text{if $m = 1$,  $\ell = 2$}.
\end{cases}
\end{align*}
Setting $q \in [1,\infty]$ such that $W^{\ell-m,p}(\widehat{T}) \hookrightarrow L^q(\widehat{T})$. Then, for all $\hat{\varphi} \in W^{\ell ,p}(\widehat{T})$ with ${\varphi} := \hat{\varphi} \circ {\Phi}^{-1}$, we have
\begin{align}
\displaystyle
| {\varphi} - I_{{T}}^{nCR} {\varphi}|_{W^{m,q}({T})}
\leq  c |T|_d^{\frac{1}{q} - \frac{1}{p}} \left( \frac{H_{T}}{h_{T}} \right)^m \sum_{|\varepsilon| =  \ell-m} {h}^{\varepsilon}  | \partial_r^{\varepsilon} \varphi |_{W^{m,p}(T)}, \label{nCR6}
\end{align}
In particular, if Condition \ref{Cond333} is imposed, it holds that, for all $\hat{\varphi} \in W^{\ell ,p}(\widehat{T})$ with ${\varphi} := \hat{\varphi} \circ {\Phi}^{-1}$, 
\begin{align}
\displaystyle
| {\varphi} - I_{{T}} {\varphi}|_{W^{m,q}({T})}
\leq c |T|_d^{\frac{1}{q} - \frac{1}{p}} \left( \frac{H_{T}}{h_{T}} \right)^m \sum_{|\varepsilon| = \ell-m} \widetilde{\mathscr{H}}^{\varepsilon}  | \partial_{\tilde{x}}^{\varepsilon}  ( \varphi \circ \Phi_{\widetilde{T}}) |_{W^{m,p}(\Phi_{\widetilde{T}}^{-1}(T))}. \label{nCR7}
\end{align}
Furthermore, for any $\hat{\varphi} \in \mathcal{C}(\widehat{T})$ with ${\varphi} := \hat{\varphi} \circ {\Phi}^{-1}$, it holds that
\begin{align*}
\displaystyle
\| {\varphi} - I_{{T}} {\varphi} \|_{L^{\infty}({T})}
\leq  c  \|  {\varphi} \|_{L^{\infty}(T)}.
\end{align*}
\end{cor}

\begin{pf*}
For $k=1$, we only introduce functionals $\mathscr{F}_i$ satisfying \eqref{functionals} in Theorem \ref{thr551} (or Theorem \ref{thrApel}) for each $\ell$ and $m$.

Let $m=0$, that is, $\beta = (0,\cdots,0) \in \mathbb{N}_0^d$. We then have $j = \dim \mathbb{P}^1 = d+1$.  From the Sobolev embedding theorem (Theorem \ref{thr=intro1}), we have $W^{\ell,p}(\widehat{T}) \subset \mathcal{C}^0(\widehat{T})$ with $1 \< p \leq \infty$, $d \< \ell p$ or $p=1$, $d \leq \ell$. Under this condition, we use
\begin{align*}
\displaystyle
\mathscr{F}_i(\hat{\varphi}) := \hat{\varphi} (\hat{x}_{\widehat{F}_i}), \quad \hat{\varphi} \in W^{\ell,p}(\widehat{T}), \quad i=1,\ldots, d+1.
\end{align*}
It then holds that
\begin{align*}
\displaystyle
|\mathscr{F}_i(\hat{\varphi})| \leq \| \hat{\varphi} \|_{\mathcal{C}^0(\widehat{T})} \leq c \| \hat{\varphi} \|_{W^{\ell,p}(\widehat{T})},
\end{align*}
which means $\mathscr{F}_i \in W^{\ell,p}(\widehat{T})^{\prime}$ for $i = 1,\ldots,d+1$, that is, \eqref{functionals=a} is satisfied. Furthermore, we have
\begin{align*}
\displaystyle
\mathscr{F}_i ( I_{\widehat{T}}^{nCR} \hat{\varphi} ) = ( I_{\widehat{T}}^{nCR} \hat{\varphi}) (\hat{x}_{\widehat{F}_i}) =  \hat{\varphi} (\hat{x}_{\widehat{F}_i}) = \mathscr{F}_i (\hat{\varphi}), \quad i = 1,\ldots,d+1,
\end{align*}
which satisfies \eqref{functionals=b}. For all $\hat{\eta} \in \mathbb{P}^{1}$, if $\mathscr{F}_i( \hat{\eta}) = 0$ for $i = 1,\ldots,d+1$, it obviously holds $\hat{\eta} = 0$. This means that \eqref{functionals=c} is satisfied.

Let $d=2$ and $m=1$ ($\ell = 2$). We set $\beta = (1,0)$. We then have $j = \dim( \partial^{\beta} \mathbb{P}^1 ) = 1$. We consider a functional
\begin{align*}
\displaystyle
\mathscr{F}_1(\hat{\varphi}) := \int_{0}^{\frac{1}{2}} \hat{\varphi}(\hat{x}_1,1/2) d \hat{x}_1, \quad \hat{\varphi} \in W^{2,p}(\widehat{T}), \quad 1 \< p.
\end{align*}
We set $\widehat{I} := \{ \hat{x} \in \widehat{T}; \ \hat{x}_2 = \frac{1}{2} \}$. The continuity is then shown by the trace theorem (e.g., see Theorem \ref{thr=intro2}):
\begin{align*}
\displaystyle
|\mathscr{F}_1(\hat{\varphi}) | \leq \| \hat{\varphi} \|_{L^1(\widehat{I})} \leq c \| \hat{\varphi} \|_{W^{1,p}(\widehat{T})},
\end{align*}
which means $\mathscr{F}_1 \in W^{2,p}(\widehat{T})^{\prime}$, that is, \eqref{functionals=a} is satisfied. Furthermore, it holds that
\begin{align*}
\displaystyle
\mathscr{F}_1 (\partial^{(1,0)} (\hat{\varphi} - I_{\widehat{T}}^{nCR} \hat{\varphi}))
&= \int_{0}^{\frac{1}{2}} \frac{\partial}{\partial \hat{x}_1} (\hat{\varphi} - I_{\widehat{T}}^{nCR} \hat{\varphi})(\hat{x}_1,1/2) d \hat{x}_1 \\
&= \left[ \hat{\varphi} - I_{\widehat{T}}^{nCR} \hat{\varphi} \right]_{(0,1/2)}^{(1/2,1/2)} = 0,
\end{align*}
which satisfy \eqref{functionals=b}. Let $\hat{\eta} := a \hat{x}_1 + b \hat{x}_2 + c$. We then have
\begin{align*}
\displaystyle
\mathscr{F}_1 (\partial^{(1,0)} \hat{\eta}) = \frac{1}{2} a.
\end{align*}
If $\mathscr{F}_1 (\partial^{(1,0)} \hat{\eta}) = 0$, $a = 0$. This implies that $\partial^{(1,0)} \hat{\eta} = 0$. This means that \eqref{functionals=c} is satisfied.

By analogous argument, the case $\beta = (0,1)$ holds.

Let $d=3$ and $m=1$ ($\ell = 2$). We consider Type (\roman{sone}) in Section \ref{reference3d} in detail. That is, the reference element is $\widehat{T} = \conv \{ 0,e_1, e_2,e_3 \}$. Here, $e_1, \ldots, e_3 \in \mathbb{R}^3$ are the canonical basis. We set $\beta = (1,0,0)$. We then have $j = \dim( \partial^{\beta} \mathcal{P}^1 ) = 1$. We consider a functional
\begin{align*}
\displaystyle
\mathscr{F}_1(\hat{\varphi}) := \int_{0}^{\frac{1}{3}} \hat{\varphi}(\hat{x}_1,1/3,1/3) d \hat{x}_1, \quad \hat{\varphi} \in W^{2,p}(\widehat{T}), \quad \frac{3}{2} \< p.
\end{align*}
We set $\widehat{I} := \{ \hat{x} \in \widehat{T}; \ \hat{x}_2 = \frac{1}{3}, \ \ \hat{x}_3 = \frac{1}{3} \}$. The continuity is then shown by the trace theorem:
\begin{align*}
\displaystyle
|\mathscr{F}_1(\hat{\varphi}) | \leq \| \hat{\varphi} \|_{L^1(\widehat{I})} \leq c \| \hat{\varphi} \|_{W^{2 ,p}(\widehat{T})} \quad \text{if $p \> 2$},
\end{align*}
which means $\mathscr{F}_1 \in W^{2,p}(\widehat{T})^{\prime}$, that is, \eqref{functionals=a} is satisfied. Furthermore, it holds that
\begin{align*}
\displaystyle
\mathscr{F}_1 (\partial^{(1,0,0)} (\hat{\varphi} - I_{\widehat{T}}^{nCR} \hat{\varphi}))
&= \left[ \hat{\varphi} - I_{\widehat{T}}^{nCR} \hat{\varphi} \right]_{(0,1/3,1/3)}^{(1/3,1/3,1/3)} = 0,
\end{align*}
which satisfy \eqref{functionals=b}. Let $\hat{\eta} := a \hat{x}_1 + b \hat{x}_2 + c \hat{x}_3 + d$. We then have
\begin{align*}
\displaystyle
\mathscr{F}_1 (\partial^{(1,0,0)} \hat{\eta}) = \frac{1}{3} a.
\end{align*}
If $\mathscr{F}_1 (\partial^{(1,0,0)} \hat{\eta}) = 0$, $a = 0$. This implies that $\partial^{(1,0,0)} \hat{\eta} = 0$. This means that \eqref{functionals=c} is satisfied.

By analogous argument, it holds the cases $\beta = (0,1,0), (0,0,1)$.

We consider Type (\roman{stwo}) in Section \ref{reference3d}. That is, the reference element is $\widehat{T} = \conv \{ 0,e_1, e_1 + e_2 , e_3 \}$. We set $\beta = (1,0,0)$. We then have $j = \dim( \partial^{\beta} \mathbb{P}^1 ) = 1$. We consider a functional
\begin{align*}
\displaystyle
\mathscr{F}_1(\hat{\varphi}) := \int_{\frac{1}{3}}^{\frac{2}{3}} \hat{\varphi}(\hat{x}_1,1/3,1/3) d \hat{x}_1, \quad \hat{\varphi} \in W^{2,p}(\widehat{T}).
\end{align*}
We can deduce the result by the similar argument with Type (\roman{sone}).

When $m = \ell = 0$, $p = \infty$ and $q \in [1,\infty]$, it holds that
\begin{align*}
\displaystyle
\| \hat{\varphi} - I_{\widehat{T}}^{nCR} \hat{\varphi} \|_{L^{q}(\widehat{T})}
&\leq c \| \hat{\varphi} \|_{L^{\infty}(\widehat{T})},
\end{align*}
because we have
\begin{align*}
\displaystyle
| ( I_{\widehat{T}}^{nCR}  \hat{\varphi} ) (\hat{x}) |
&\leq \sum_{i=1}^{d+1} | \hat{\varphi} (\hat{x}_{\widehat{F}_i})| | \hat{\theta}_i (\hat{x})|
\leq (d+1) \left( \max_{1 \leq i \leq d+1} \| \hat{\theta}_i \|_{L^{\infty}(\widehat{T})} \right) \| \hat{\varphi} \|_{L^{\infty}(\widehat{T})}.
\end{align*}
\qed
\end{pf*}

\subsection{Morley Finite Element}
Any dimensional Morley finite element is introduced in  \cite{WanJin06}.

Let $\widehat{T} \subset \mathbb{R}^d$, $d \in \{ 2,3\}$, be the reference element defined in Sections \ref{reference2d} and \ref{reference3d}. Let $\widehat{F}_i$, $1 \leq i \leq d+1$, be the $(d-1)$-dimensional subsimplex of $\widehat{T}$ without $\widehat{P}_i$ and $\widehat{S}_{i,j}$, $1 \leq i \< j \leq d+1$, the $(d-2)$-dimensional subsimplex of $\widehat{T}$ without $\widehat{P}_i$ and $\widehat{P}_j$. The $d$-dimensional Morley finite element on the reference element is defined by the triple $\{ \widehat{T} , \widehat{P} , \widehat{\Sigma\}}$ as 
\begin{enumerate}
 \item $\widehat{P} := \mathbb{P}^2(\widehat{T})$;
 \item $\widehat{\Sigma}$ is a set $\{ \hat{\chi}_{i} \}_{1 \leq i \leq N^{(d,2)}}$ of $N^{(d,2)} $ linear forms $\{ \hat{\chi}^{(1)}_{i,j} \}_{1 \leq i \< j \leq d+1} \cup \{ \hat{\chi}^{(2)}_{i} \}_{1 \leq i \leq d+1}$ with its components such that, for any $\hat{q} \in \widehat{P}$,
\begin{subequations} \label{Mor1}
\begin{align}
\displaystyle
\hat{\chi}^{(1)}_{i,j}(\hat{q}) &:= \frac{1}{|\widehat{S}_{i,j}|} \int_{\widehat{S}_{i,j}} \hat{q} d \hat{s}, \quad 1 \leq i \< j \leq d+1,  \label{Mor1=a}\\
\hat{\chi}^{(2)}_{i}(\hat{q}) &:= \frac{1}{|\widehat{F}_{i}|} \int_{\widehat{F}_{i}} \frac{\partial \hat{q}}{\partial \hat{n}_i} d \hat{s}, \quad 1 \leq i \leq d+1, \label{Mor1=b}
\end{align}
\end{subequations}
where $\frac{\partial }{\partial \hat{n}_i} = n_{\widehat{T},i} \cdot \nabla $, and $n_{\widehat{T},i}$ is the unit outer normal to $ \widehat{F}_i \subset \partial \widehat{T}$. For $d=2$, $\hat{\chi}^{(1)}_{i,j}(\hat{p})$ is interpreted as
\begin{align*}
\displaystyle
\hat{\chi}^{(1)}_{i,j}(\hat{q}) = \hat{q}(\widehat{p}_k), \quad k=1,2,3, \quad k \neq i,j.
\end{align*} 
\end{enumerate}
For a Morley finite {element}, $\widehat{\Sigma}$ is unisolvent (see \cite[Lemma 2]{WanJin06}). The nodal basis functions associated with the degrees of freedom provided by \eqref{Mor1} are defined as follows:
\begin{subequations} \label{Mor2}
\begin{align}
\displaystyle
\hat{\theta}_{i,j}^{(1)} &:= 1 - (d-1)(\hat{\lambda}_i + \hat{\lambda}_j) + d(d-1)\hat{\lambda}_i \hat{\lambda}_j \notag \\
&\quad - (d-1) (\nabla \hat{\lambda}_i)^T \nabla \hat{\lambda}_j \sum_{k = i,j} \frac{\hat{\lambda}_k (d \hat{\lambda}_k - 2)}{2 | \nabla \hat{\lambda}_k |_E^2}, \quad 1 \leq i \< j \leq d+1, \label{Mor2=a} \\
\hat{\theta}_{i}^{(2)} &:= \frac{\hat{\lambda}_i (d \hat{\lambda}_i - 2)}{2 | \nabla \hat{\lambda}_i |_E}, \quad 1 \leq i \leq d+1, \label{Mor2=b}
\end{align}
\end{subequations}
where $|\nabla \hat{\lambda}_i|_E$ denotes the Euclidean norm in $\mathbb{R}^d$. Subsequently, \cite[Theorem 1]{WanJin06} proved that, for $1 \leq i \< j \leq d+1$,
\begin{align}
\displaystyle
\hat{\chi}^{(1)}_{k,\ell}(\hat{\theta}_{i,j}^{(1)}) = \delta_{ik} \delta_{j \ell}, \quad 1 \leq k \< \ell \leq d+1,  \quad \hat{\chi}^{(2)}_{k}(\hat{\theta}_{i,j}^{(1)}) = 0, \quad 1 \leq k \leq d+1, \label{Mor3}
\end{align}
and, for $1 \leq i \leq d+1$,
\begin{align}
\displaystyle
\hat{\chi}^{(1)}_{k,\ell}(\hat{\theta}^{(2)}_{i}) = 0, \quad 1 \leq k \< \ell \leq d+1, \quad \hat{\chi}^{(2)}_{k}(\hat{\theta}^{(2)}_{i}) = \delta_{ik}, \quad 1 \leq k \leq d+1. \label{Mor4}
\end{align}
The local interpolation operator ${I}_{\widehat{T}}^{M}$ is defined by
\begin{align}
\displaystyle
{I}_{\widehat{T}}^{M}: W^{2,1}(\widehat{T}) \ni \hat{\varphi} \mapsto {I}_{\widehat{T}}^{M} \hat{\varphi} \in \widehat{P}, \label{Mor5}
\end{align}
with
\begin{align}
\displaystyle
{I}_{\widehat{T}}^{M} \hat{\varphi} := \sum_{1 \leq i \< j \leq d+1} \hat{\chi}^{(1)}_{i,j}(\hat{\varphi}) \hat{\theta}_{i,j}^{(1)} + \sum_{1 \leq i \leq d+1} \hat{\chi}^{(2)}_{i}(\hat{\varphi})  \hat{\theta}_{i}^{(2)}.  \label{Mor6}
\end{align}
Then, it holds that ${I}_{\widehat{T}}^{M} \hat{q} = \hat{q}$ for any $\hat{q} \in \widehat{P}$ and, for any $\hat{\varphi} \in W^{2,1}(\widehat{T})$,
\begin{subequations} \label{Mor7}
\begin{align}
\displaystyle
\hat{\chi}^{(1)}_{i,j}({I}_{\widehat{T}}^{M} \hat{\varphi}) &= \hat{\chi}^{(1)}_{i,j}( \hat{\varphi}), \quad 1 \leq i \< j \leq d+1, \label{Mor7=a}\\
\hat{\chi}^{(2)}_{i}({I}_{\widehat{T}}^{M} \hat{\varphi}) &= \hat{\chi}^{(2)}_{i}( \hat{\varphi}), \quad 1 \leq i \leq d+1.  \label{Mor7=b}
\end{align}
\end{subequations}
By analogous argument in Section \ref{FEMG}, the Morley finite elements $\{ \widetilde{T} , \widetilde{P} , \widetilde{\Sigma}\}$ and $\{ T , {P} , \Sigma \}$ are constructed. The local shape functions are
\begin{align*}
\displaystyle
&\tilde{\theta}_{i,j}^{(1)} = \psi_{\widehat{T}}^{-1}(\hat{\theta}_{i,j}^{(1)}), \quad 1 \leq i \< j \leq d+1, \quad \tilde{\theta}_{i}^{(2)} = \psi_{\widehat{T}}^{-1}(\hat{\theta}_{i}^{(2)}), \quad 1 \leq i \leq d+1, \\
&{\theta}_{i,j}^{(1)} = \psi_{\widetilde{T}}^{-1}(\tilde{\theta}_{i,j}^{(1)}), \quad 1 \leq i \< j \leq d+1, \quad {\theta}_{i}^{(2)} = \psi_{\widetilde{T}}^{-1}(\tilde{\theta}_{i}^{(2)}), \quad 1 \leq i \leq d+1.
\end{align*}
The associated local Morley interpolation operators are defined as
\begin{align}
\displaystyle
{I}_{\widetilde{T}}^{M}: W^{2,1}(\widetilde{T}) \ni \tilde{\varphi} \mapsto {I}_{\widetilde{T}}^{M} \tilde{\varphi} \in \widetilde{P}, \label{Mor8}
\end{align}
with, for any $\tilde{\varphi} \in W^{2,1}(\widetilde{T})$, 
\begin{subequations} \label{Mor9}
\begin{align}
\displaystyle
\tilde{\chi}^{(1)}_{i,j}({I}_{\widetilde{T}}^{M} \tilde{\varphi}) &= \tilde{\chi}^{(1)}_{i,j}( \tilde{\varphi}), \quad 1 \leq i \< j \leq d+1, \label{Mor9=a}\\
\tilde{\chi}^{(2)}_{i}({I}_{\widetilde{T}}^{M} \tilde{\varphi}) &= \tilde{\chi}^{(2)}_{i}( \tilde{\varphi}), \quad 1 \leq i \leq d+1,  \label{Mor9=b}
\end{align}
\end{subequations}
and
\begin{align}
\displaystyle
{I}_{{T}}^{M}: W^{2,1}({T}) \ni {\varphi} \mapsto {I}_{{T}}^{M} {\varphi} \in {P}, \label{Mor31}
\end{align}
with, for any ${\varphi} \in W^{2,1}({T})$, 
\begin{subequations} \label{Mor32}
\begin{align}
\displaystyle
{\chi}^{(1)}_{i,j}({I}_{{T}}^{M} {\varphi}) &= {\chi}^{(1)}_{i,j}( {\varphi}), \quad 1 \leq i \< j \leq d+1, \label{Mor32=a}\\
{\chi}^{(2)}_{i}({I}_{{T}}^{M} {\varphi}) &= {\chi}^{(2)}_{i}( {\varphi}), \quad 1 \leq i \leq d+1.  \label{Mor32=b}
\end{align}
\end{subequations}

\begin{rem}
The Morley FEM has not been defined uniquely. There are two versions: one defined in \cite{Mol68}, which is the original paper, and the other in \cite{ArnBre85,LasLea75,WanJin06}. {In original  Morley FEM, by normal derivatives on faces, the spans of the nodes are not preserved under push-forward. To overcome this difficulty, the mean value of the first normal derivative is used \cite{ArnBre85,LasLea75,WanJin06}. The original Morley interpolation error estimates are obtained using the modified Morley interpolation error estimates (see \cite{LasLea75}). }In this study, we used the Morley FEM introduced in  \cite{WanJin06}. 	
\end{rem}

\subsection{Local Morley Interpolation Error Estimates}
Using the idea of \cite[Lemma 1]{WanJin06}, the following lemma holds.

\begin{lem} \label{lem=new1}
Let $T \subset \mathbb{R}^d$ be a simplex. $n_{T,k}$ denotes the unit outer normal to the face $F_k$, $k  =1,\ldots,d+1$ of $T$, $S_1,\ldots,S_d$ are all $(d-2)$-dimensional subsimplexes of $F_k$. Let $v \in \mathcal{C}^1(T)$ be such that
\begin{align}
\displaystyle
\int_{S_{\ell}} v = 0, \quad \int_{F_k} \frac{\partial v}{\partial n_k} = 0, \label{Mor22}
\end{align}
for any $\ell=1,\ldots,d$ and  $k=1,\ldots,d+1$. It then holds that
\begin{align}
\displaystyle
 \int_{F_k} \frac{\partial v}{\partial x_i} = 0, \quad i=1,\ldots,d, \quad k=1,\ldots,d+1.  \label{Mor23}
\end{align}
\end{lem}

\begin{pf*}
 Let $v \in \mathcal{C}^1(T)$. Let $\xi \in \mathbb{R}^d$ be a constant vector, and let $\tau := \xi - (\xi \cdot n_{T,k}) n_{T,k}$. We have
\begin{align*}
\displaystyle
\tau \cdot n_{T,k}
&= \xi \cdot n_{T,k} - (\xi \cdot n_{T,k}) n_{T,k} \cdot n_{T,k} = 0,
\end{align*}
that is, $\tau$ is the tangent vector of $F_k$. Subsequently, from \eqref{Mor22} we obtain
\begin{align*}
\displaystyle
\int_{F_k} (\xi \cdot \nabla ) v 
= \int_{F_k}  \frac{\partial v}{\partial \tau} +  (\xi \cdot n_{T,k}) \int_{F_k} \frac{\partial v}{\partial n_k}
=  \int_{F_k}  \frac{\partial v}{\partial \tau}.
\end{align*}
Let $d=2$. Let $p_{k1}$ and $p_{k2}$ be the endpoints of the edge $F_k$, that is, $F_k = \overline{p_{k1} p_{k2}}$. Subsequently,  from \eqref{Mor22} we obtain
\begin{align}
\displaystyle
 \int_{F_k}  \frac{\partial v}{\partial \tau}
 &= \int_{s=0}^{s=|F_k|} \frac{d v}{ds} \left( \frac{|F_k| - s}{|F_k|}p_{k1} + \frac{s}{|F_k|} p_{k2} \right) = v(p_{k2}) - v(p_{k1}) = 0. \label{Mor21=2d}
\end{align}
Let $d=3$.  Let $\zeta^{(\ell)}$ be the unit outer normal of $S_{\ell}$ for $\ell=1,2,3$, From \eqref{Mor22}, the Gauss--Green formula yields
\begin{align}
\displaystyle
 \int_{F_k}  \frac{\partial v}{\partial \tau}
&= \sum_{\ell=1}^3 \tau \cdot \zeta^{(\ell)} \int_{S_{\ell}} v = 0. \label{Mor21}
\end{align}
From \eqref{Mor21=2d} and \eqref{Mor21}, it holds that for $d=2,3$
\begin{align}
\displaystyle
\int_{F_k} (\xi \cdot \nabla ) v = 0. \label{Mor21=23d}
\end{align}
Let $e_1, \ldots, e_d \in \mathbb{R}^d$ be a canonical basis. By setting $\xi := e_i$ in \eqref{Mor21=23d}, we obtain the desired result in \eqref{Mor23} under Assumption \eqref{Mor22}.
\qed
\end{pf*}

The anisotropic Morley interpolation error estimate is expressed as 

\begin{thr} \label{thr2=Morley}
Let $p \in [1,\infty)$ and $q \in [1,\infty)$ be such that \eqref{Sobolev511} holds true. Subsequently, for any ${\varphi} \in W^{3,p}({T}) \cap \mathcal{C}^1({T})$,we have
\begin{align}
\displaystyle
| I_{T}^{M} \varphi - \varphi |_{W^{2,q}({T})} &\leq  c |T|_d^{ \frac{1}{q} - \frac{1}{p} }  \sum_{i=1}^d  h_i \left | \frac{\partial\varphi}{\partial r_i} \right |_{W^{2,p}(T)}. \label{Mor24} 
\end{align}
If Condition \ref{Cond333} is imposed, then:
\begin{align}
\displaystyle
| I_{T}^{M} \varphi - \varphi |_{W^{2,q}({T})} &\leq  c |T|_d^{ \frac{1}{q} - \frac{1}{p} }  \sum_{i=1}^d  \widetilde{\mathscr{H}}_i \left | \frac{\partial (\varphi \circ \Phi_T)}{\partial \tilde{x}_i} \right |_{W^{2,p}(\Phi_T^{-1}(T) )}. \label{Mor24=ass0} 
\end{align}

\end{thr}

\begin{pf*}
Only Morley interpolation satisfies the condition \eqref{sp=cond}.

Let ${\varphi} \in W^{3,p}({T})  \cap \mathcal{C}^1(T)$ and set $v :=  I_{T}^{M} \varphi - \varphi$. Using the definition of the Morley interpolation operator \eqref{Mor32}, we obtain
\begin{align*}
\displaystyle
 \int_{{S}_{i,j}} {v} d {s} = 0, \quad 1 \leq i \< j \leq d+1, \quad
 \int_{{F}_{i}} \frac{\partial v}{ \partial n_i} ds = 0, \quad 1 \leq i \leq d+1.
\end{align*}
Therefore, from Lemma \ref{lem=new1}, we have
\begin{align}
\displaystyle
\int_{F_i} \frac{\partial v}{\partial x_k} = 0, \quad i=1,\ldots,d+1, \quad k = 1,\ldots,d. \label{Mor26}
\end{align}
From Green's formula and \eqref{Mor26}, it follows that, for $1 \leq j,k \leq d$,
\begin{align*}
\displaystyle
\int_{{T}} \frac{\partial^2 v}{\partial {x}_j \partial {x}_k} d {x}
&= \sum_{i=1}^{d+1}  n_{{T},j} \int_{{F}_i}  \frac{\partial v}{\partial {x}_k} d {s} = 0,
\end{align*}
which leads to
\begin{align*}
\displaystyle
 \frac{\partial^2}{\partial {x}_j \partial {x}_k} (I_{{T}}^{M} {\varphi} )
  = \frac{1}{|T|_d} \int_{T}  \frac{\partial^2 }{\partial {x}_j \partial {x}_k}  (I_{{T}}^{M} {\varphi} ) dx
 = \frac{1}{|T|_d} \int_{T}  \frac{\partial^2 \varphi}{\partial {x}_j \partial {x}_k} dx = \Pi_{T}^0 \left( \frac{\partial^2 \varphi}{\partial {x}_j \partial {x}_k} \right),
\end{align*}
because $ \frac{\partial^2}{\partial {x}_j \partial {x}_k} (I_{{T}}^{M} {\varphi}) \in \mathbb{P}^0(T)$, {Therefore, by Theorem \ref{Thr2a}, the target inequalities \eqref{Mor24} and \eqref{Mor24=ass0} hold.}
\qed
\end{pf*}

\subsection{Global Morley Interpolation Error Estimates}
Recall the space $V_h^n$ with $n=2$ (see \eqref{space=Vh}). the Morley finite element space is as follows:
\begin{align*}
\displaystyle
V_{h}^M &:=  \biggl \{ \varphi_h \in V_h^2: \ \int_F  \biggl {[} \! \biggl {[} \frac{ \partial \varphi_h}{  \partial n}  \biggr {]} \! \biggr {]} ds = 0 \ \forall F \in \mathcal{F}_h^i, \\
&\hspace{0.7cm}  \text{the integral average  of $\varphi_h$ over each $(d-2)$-dimensional} \\
&\hspace{0.7cm} \text{subsimplex of $T \in \mathbb{T}_h$ is continuous} \biggr \}, \\
V_{h0}^M &:=  \left\{ \varphi_h \in V_{h}^M; \ \text{degrees of freedom of $\varphi_h$ in \eqref{Mor1} vanish on $\partial \Omega$} \right \}.
\end{align*}
In particular, for $d=2$, the space $V_{h0}^M$ is described as
\begin{align*}
\displaystyle
V_{h0}^M &:=  \biggl \{ \varphi_h \in V_{h}^2: \ \int_F  \biggl {[} \! \biggl {[} \frac{ \partial \varphi_h}{  \partial n}  \biggr {]} \! \biggr {]} ds = 0 \ \forall F \in \mathcal{F}_h, \\
&\hspace{0.7cm} \text{$\varphi_h$ is continuous at each vertex in $\Omega$, $\varphi_h(p) = 0$, \ $p \in \partial \Omega$} \biggr \}.
\end{align*}
We also define the global interpolation $I_h^{M}: W^{2,1}(\Omega) \to V_{h}^{M}$ (or $I_h^{M}: W^{2,1}_{0}(\Omega) \to V_{h0}^{M}$) as follows.
\begin{align*}
\displaystyle
(I_h^{M} \varphi)|_{T} := I_{T}^{M} (\varphi|_{T})  \quad \forall T \in \mathbb{T}_h, \ \forall \varphi \in W^{2,1}(\Omega).
\end{align*}

\begin{cor} \label{thr1031}
Suppose that the assumptions of Theorem \ref{thr2=Morley} are satisfied. Let $I_h^{M}$ be the corresponding global Morley interpolation operator. It then holds that, for any ${\varphi} \in W^{3,p}({\Omega}) \cap \mathcal{C}^1(\overline{\Omega})$,
\begin{description}	
  \item[(\Roman{lone})] if Condition \ref{Cond333} is not imposed,
\begin{align}
\displaystyle
| I_{T}^{M} \varphi - \varphi |_{W^{2,q}({\mathbb{T}_h})} &\leq  c  \sum_{T \in \mathbb{T}_h}  |T|_d^{ \frac{1}{q} - \frac{1}{p} }  \sum_{i=1}^d  h_i \left | \frac{\partial\varphi}{\partial r_i} \right |_{W^{2,p}(T)}, \label{Mor33} 
\end{align}
  \item[(\Roman{ltwo})] if Condition \ref{Cond333} is imposed, 
\begin{align}
\displaystyle
| I_{T}^{M} \varphi - \varphi |_{W^{2,q}({\mathbb{T}_h})} &\leq  c  \sum_{T \in \mathbb{T}_h} |T|_d^{ \frac{1}{q} - \frac{1}{p} }  \sum_{i=1}^d  \widetilde{\mathscr{H}}_i \left | \frac{\partial (\varphi \circ \Phi_T)}{\partial \tilde{x}_i} \right |_{W^{2,p}(\Phi_T^{-1}(T) )}. \label{Mor34} 
\end{align}
\end{description}
\end{cor}

\begin{pf*}
This corollary is proved in the same argument as Corollary \ref{newglobal=cor}.
\qed
\end{pf*}

\subsection{Another Estimate}

\begin{thr}
Let $T \subset \mathbb{R}^d$ be a simplex. Let $I_T^{M}:H^2(T) \to \mathbb{P}^2(T)$ be the local Morley interpolation operator defined as
\begin{align*}
\displaystyle
{I}_{{T}}^{M}: H^{2}({T}) \ni {\varphi} \mapsto {I}_{{T}}^{M} {\varphi} \in \mathbb{P}^2(T), 
\end{align*}
with
\begin{subequations}
\begin{align*}
\displaystyle
{\chi}^{(1)}_{i,j}({I}_{{T}}^{M} {\varphi}) &= {\chi}^{(1)}_{i,j}( {\varphi}), \quad 1 \leq i \< j \leq d+1, \\
{\chi}^{(2)}_{i}({I}_{{T}}^{M} {\varphi}) &= {\chi}^{(2)}_{i}( {\varphi}), \quad 1 \leq i \leq d+1.  
\end{align*}
\end{subequations}
for any ${\varphi} \in H^{2}({T})$. It then holds that 
\begin{align}
\displaystyle
| I_{T}^{M} \varphi - \varphi |_{H^2(T)}  &\leq \frac{h_{T}}{\pi} |\varphi|_{H^3(T)} \quad \forall \varphi \in H^3(T). \label{Mor42}
\end{align}
\end{thr}

\begin{pf*}
Using \eqref{chap841}.
\begin{align*}
\displaystyle
| I_{T}^{M} \varphi - \varphi |_{H^2(T)}^2
&= \sum_{j,k=1}^d \left \|  \frac{\partial^2}{\partial x_j \partial x_k} ( I_{T}^{M} \varphi - \varphi)  \right\|^2_{L^2(T)} \notag \\
&= \sum_{j,k=1}^d \left \|  \Pi_T^0 \left( \frac{\partial^2 \varphi}{\partial {x}_j \partial {x}_k} \right) - \left( \frac{\partial^2 \varphi}{\partial {x}_j \partial {x}_k} \right)  \right\|^2_{L^2(T)} \\
&\leq \left( \frac{h_{T}}{\pi} \right)^2 \sum_{j,k=1}^d \left|   \frac{\partial^2 \varphi}{\partial {x}_j \partial {x}_k} \right|^2_{H^1(T)} = \left( \frac{h_{T}}{\pi} \right)^2 |\varphi|_{H^3(T)}^2,
\end{align*}
which conclude \eqref{Mor42}.
\qed
\end{pf*}

\section{New Scaling Argument: Part 2}

\subsection{Two-step Piola Transforms}
We adopt the following two-step Piola transformations.

\begin{defi}[Two-step Piola transforms] \label{piola=defi}
Let  $V(\widehat{T}) := \mathcal{C}(\widehat{T})^d$. The Piola transformation $\Psi := {\Psi}_{\widetilde{T}} \circ {\Psi}_{\widehat{T}} : V(\widehat{T}) \to V(T)$ is defined as
\begin{align}
\displaystyle
\Psi : V(\widehat{T})  &\to V(T) \label{piola} \\
\hat{v} &\mapsto v(x) :=  {\Psi}(\hat{v})(x) = \frac{1}{\det({A})} {A} \hat{v}(\hat{x}), \notag
\end{align}
with	two Piola transformations:
\begin{align*}
\displaystyle
{\Psi}_{\widehat{T}}: V(\widehat{T}) &\to V(\widetilde{T})\\
\hat{v} &\mapsto \tilde{v}(\tilde{x}) := {\Psi}_{\widehat{T}}(\hat{v})(\tilde{x}) := \frac{1}{\det ({A}_{\widetilde{T}})} {A}_{\widetilde{T}}\hat{v}(\hat{x}), \\
\Psi_{\widetilde{T}} :  V(\widetilde{T}) &\to V({T})  \\
\tilde{v} &\mapsto v(x) :=  {\Psi}_{\widetilde{T}} (\tilde{v})(x) :=  \frac{1}{\det ({{A}}_{T}) } {A}_{T} \tilde{v} (\tilde{x}).
\end{align*}
\end{defi}

\subsection{Property of the Piola Transformations}

\begin{lem}
If $\hat{v} \in \mathcal{C}^1(\widehat{T})^d$, then $v := \Psi \hat{v} \in \mathcal{C}^1(T)^d$ and it holds
\begin{align*}
\displaystyle
J_x v &=  \frac{1}{\det({{A}})} {{A}} \widehat{J}_{\hat{x}} \hat{v} A^{-1}, \\
\div v &= \frac{1}{\det({{A}})} \widehat{\div} \hat{v}, 
\end{align*}
where $J_x v$ and $ \widehat{J}_{\hat{x}} \hat{v}$ denote the Jacobian matrixes of $v$ and $\hat{v}$, respectively.
\end{lem}

\begin{pf*}
From the definition of the Piola transformation \eqref{piola}, we have
\begin{align*}
\displaystyle
J_x v (x) &= \frac{1}{\det({{A}})} {{A}} J_x (\hat{v} \circ \Phi^{-1})(x) = \frac{1}{\det({{A}})} {{A}} \widehat{J}_{\hat{x}} \hat{v}(\hat{x}) J_x \Phi^{-1} (x) \\
&= \frac{1}{\det({{A}})} {{A}} \widehat{J}_{\hat{x}} \hat{v}(\hat{x}) {{A}}^{-1}.
\end{align*}
Due to the property of the trace, we get
\begin{align*}
\displaystyle
\div v = \Tr (J_x v) &= \frac{1}{\det({{A}})} \Tr({{A}} \widehat{J}_{\hat{x}} \hat{v} {{A}}^{-1}) \\
&= \frac{1}{\det({{A}})} \Tr(\widehat{J}_{\hat{x}} \hat{v}) = \frac{1}{\det({{A}})} \widehat{\div} \hat{v}.
\end{align*}
\qed	
\end{pf*}

\begin{lem}
For $\hat{\varphi} \in \mathcal{C}^{1}(\widehat{T})$, $\hat{v} \in \mathcal{C}^{1}(\widehat{T})^d$ with $\varphi := \hat{\varphi} \circ \Phi^{-1}$ and $v := \Psi (\hat{v})$, it holds that
\begin{align}
\displaystyle
\int_{T} \div v \varphi dx &= \int_{\widehat{T}} \widehat{\div} \hat{v} \hat{\varphi} d \hat{x}, \label{piola=new2}\\
\int_{T} (v \cdot \nabla_x) \varphi dx &= \int_{\widehat{T}} ( \hat{v} \cdot \widehat{\nabla}_{\hat{x}} ) \hat{\varphi} d \hat{x},  \label{piola=new3}\\
\int_{\partial T} (v \cdot n ) \varphi ds &= \int_{\partial \widehat{T}} ( \hat{v} \cdot \hat{n}) \hat{\varphi} d \hat{s}.  \label{piola=new4}
\end{align}
\end{lem}

\begin{pf*}
Because $\det ({{A}})$ is positive, by a change a variable, 
\begin{align*}
\displaystyle
\int_{T} \div v \varphi dx
= \frac{1}{\det({{A}})}  \int_{\widehat{T}} \widehat{\div} \hat{v} \hat{\varphi} | \det ({{A}}) | d \hat{x},
\end{align*}
which leads to \eqref{piola=new2}.  Because 
\begin{align*}
\displaystyle
\nabla_x \varphi = {{A}}^{- T} \widehat{\nabla}_{\hat{x}} \hat{\varphi},
\end{align*}
we have
\begin{align*}
\displaystyle
\int_{T} (v \cdot \nabla_x) \varphi dx &= \frac{1}{\det({{A}})} \int_{\widehat{T}} ( {{A}} \hat{v} \cdot {{A}}^{- T} \widehat{\nabla}_{\hat{x}} ) \hat{\varphi} |\det({{A}})| d \hat{x} \\
&= \int_{\widehat{T}} [ ({{A}} \hat{v})^{T}  {{A}}^{- T} \widehat{\nabla}_{\hat{x}} ] \hat{\varphi}  d \hat{x} = \int_{\widehat{T}} ( \hat{v} \cdot \widehat{\nabla}_{\hat{x}} ) \hat{\varphi} d \hat{x},
\end{align*}
which is \eqref{piola=new3}.

From \eqref{piola=new2} and \eqref{piola=new3}, applying the Gauss--Green formula yields
\begin{align*}
\displaystyle
\int_{\partial T} (v \cdot n) \varphi ds
&= \int_{T} (\div v ) \varphi dx + \int_{T} (v \cdot \nabla_x ) \varphi dx \\
&= \int_{\widehat{T}} ( \widehat{\div} \hat{v} ) \hat{\varphi} d \hat{x} + \int_{\widehat{T}} ( \hat{v} \cdot \widehat{\nabla}_{\hat{x}} ) \hat{\varphi} d \hat{x} = \int_{\partial \widehat{T}} ( \hat{v} \cdot \hat{n} ) \hat{\varphi} d \hat{s},
\end{align*}
which is \eqref{piola=new4}.
\qed	
\end{pf*}

\subsection{Preliminalies}
\subsubsection{Calculations 1} \label{cal941}
We use the following calculations in \eqref{RT42}. Let $\hat{v} \in \mathcal{C}^{2}(\widehat{T})^d$ with $\tilde{v} = {\Psi}_{\widehat{T}} \hat{v}$ and ${v} = {\Psi}_{\widetilde{T}} \tilde{v}$. Using the definition of Piola transformations (Definition \ref{piola=defi}) yields, for $1 \leq i,k \leq d$,
\begin{align*}
\displaystyle
  \frac{\partial \hat{v}_k}{\partial \hat{x}_{i} }
 &= \det ({A}_{\widetilde{T}}) \sum_{\eta=1}^d [ {A}_{\widetilde{T}}^{-1}]_{k \eta} \sum_{{i_1^{(1)}} = 1}^d \frac{\partial \tilde{v}_{\eta}}{\partial \tilde{x}_{{i_1^{(1)}}}} \frac{\partial \tilde{x}_{{i_1^{(1)}}}}{\partial \hat{x}_{i} } \\
 &= \det ({A}_{\widetilde{T}}) \sum_{\eta=1}^d [ {A}_{\widetilde{T}}^{-1}]_{k \eta} \sum_{{i_1^{(1)}} = 1}^d \frac{\partial \tilde{v}_{\eta}}{\partial \tilde{x}_{{i_1^{(1)}}}} [A_{\widetilde{T}}]_{{{i_1^{(1)}}} {i}} \\
 &=  \det ({A}_{\widetilde{T}}) h_i h_k^{-1} \sum_{\eta=1}^d [ \widetilde{A}^{-1}]_{k \eta} \sum_{{i_1^{(1)}} = 1}^d \frac{\partial \tilde{v}_{\eta}}{\partial \tilde{x}_{{i_1^{(1)}}}} [\widetilde{A}]_{{{i_1^{(1)}}} {i}}, \\
\frac{\partial \tilde{v}_{\eta}}{\partial \tilde{x}_{{i_1^{(1)}}}}
 &= \det ({A}_T) \sum_{\nu = 1}^d[ {A}_T^{-1}]_{\eta \nu} \sum_{{i_1^{(0,1)}} = 1}^d \frac{\partial v_{\nu}}{\partial x_{{i_1^{(0,1)}}}} [{A}_T]_{{i_1^{(0,1)}} {i_1^{(1)}}},
\end{align*}
which leads to
\begin{align*}
\displaystyle
 \frac{\partial \hat{v}_k}{\partial \hat{x}_{i}}
 &= \det (A_{\widetilde{T}}) \det ({A}_T) h_k^{-1} \sum_{\eta,\nu=1}^d [ \widetilde{{A}}^{-1}]_{k \eta} [ {A}_T^{-1}]_{\eta \nu} \sum_{i_1^{(1)}, i_1^{(0,1)} = 1}^d h_i [\widetilde{{A}}]_{{{i_1^{(1)}}} {i}} [{A}_T]_{{i_1^{(0,1)}} {i_1^{(1)}}} \frac{\partial v_{\nu}}{\partial x_{{i_1^{(0,1)}}}}.
 \end{align*}
By an analogous calculation, for $1 \leq i, j,k \leq d$,
\begin{align*}
\displaystyle
 \frac{\partial^2 \hat{v}_k}{\partial \hat{x}_{i} \partial \hat{x}_{j}}
 &= \det ({A}_{\widetilde{T}}) h_i h_j h_k^{-1} \sum_{\eta=1}^d [\widetilde{{A}}^{-1}]_{k \eta} \sum_{i_1^{(1)} , j_1^{(1)} =1}^d \frac{\partial^2 \tilde{v}_{\eta}}{\partial \tilde{x}_{i_1^{(1)}} \partial \tilde{x}_{j_1^{(1)}}} [\widetilde{{A}}]_{i_1^{(1)} i} [ \widetilde{{A}} ]_{j_1^{(1)} j}, \\
  \frac{\partial^2 \tilde{v}_{\eta}}{\partial \tilde{x}_{i_1^{(1)}} \partial \tilde{x}_{j_1^{(1)}}}
 &=  \det ({A}_T) \sum_{\nu = 1}^d [ {A}_T^{-1}]_{\eta \nu} \sum_{i_1^{(0,1)}, j_1^{(0,1)} = 1}^d \frac{\partial^2 v_{\nu}}{\partial x_{i_1^{(0,1)}} \partial x_{j_1^{(0,1)}}}  [{A}_T]_{{i_1^{(0,1)}} {i_1^{(1)}}}  [{A}_T]_{{j_1^{(0,1)}} {j_1^{(1)}}},
 \end{align*}
which leads to
\begin{align*}
\displaystyle
\frac{\partial^2 \hat{v}_k}{\partial \hat{x}_{i} \partial \hat{x}_{j}} 
&= \det (A_{\widetilde{T}}) \det ({A}_T) h_k^{-1} \sum_{\eta,\nu=1}^d [ \widetilde{{A}}^{-1}]_{k \eta} [ {A}_T^{-1}]_{\eta \nu} \\
&\quad \sum_{i_1^{(1)}, i_1^{(0,1)} = 1}^d h_i [\widetilde{{A}}]_{{{i_1^{(1)}}} {i}} [{A}_T]_{{i_1^{(0,1)}} {i_1^{(1)}}} \sum_{j_1^{(1)}, j_1^{(0,1)} = 1}^d h_j [\widetilde{{A}}]_{{{j_1^{(1)}}} {j}} [{A}_T]_{{j_1^{(0,1)}} {j_1^{(1)}}} \frac{\partial^2 v_{\nu}}{\partial x_{i_1^{(0,1)}} \partial x_{j_1^{(0,1)}}}.
\end{align*}
For any multi-indices $\beta$ and $\gamma$, for $1 \leq k \leq d$, Let $\hat{v} \in \mathcal{C}^{|\beta| + |\gamma|}(\widehat{T})^d$ with $\tilde{v} = {\Psi}_{\widehat{T}} \hat{v}$ and ${v} = {\Psi}_{\widetilde{T}} \tilde{v}$. Then,
\begin{align*}
\displaystyle
\partial^{\beta + \gamma}_{\hat{x}}  \hat{v}_k &= \frac{\partial^{|\beta| + |\gamma|}}{\partial \hat{x}_1^{\beta_1} \cdots \partial \hat{x}_d^{\beta_d} \partial \hat{x}_1^{\gamma_1} \cdots \partial \hat{x}_d^{\gamma_d}}  \hat{v}_k \notag\\
&\hspace{-1.5cm}  =  \det (A_{\widetilde{T}}) \det ({A}_T) h_k^{-1} \sum_{\eta,\nu=1}^d [ \widetilde{{A}}^{-1}]_{k \eta} [ {A}_T^{-1}]_{\eta \nu} \\
&\hspace{-1.2cm}  \underbrace{\sum_{i_1^{(1)},i_1^{(0,1)} = 1}^d h_1 [\widetilde{{A}}]_{i_1^{(1)} 1} [{A}_{T}]_{i_1^{(0,1)} i_1^{(1)}} \cdots   \sum_{i_{\beta_1}^{(1)},i_{\beta_1}^{(0,1)} = 1}^d h_1 [\widetilde{{A}}]_{i_{\beta_1}^{(1)} 1} [{A}_{T}]_{i_{\beta_1}^{(0,1)} i_{\beta_1}^{(1)}}   }_{\beta_1 \text{times}} \cdots \notag \\
&\hspace{-1.2cm} \underbrace{ \sum_{i_1^{(d)} , i_1^{(0,d)} = 1}^d h_d [\widetilde{{A}}]_{i_1^{(d)} d} [{A}_{T}]_{i_1^{(0,d)} i_1^{(d)}}  \cdots \sum_{i_{\beta_d}^{(d)} , i_{\beta_d}^{(0,d)}= 1}^d h_d [\widetilde{{A}}]_{i_{\beta_d}^{(d)} d} [{A}_{T}]_{i_{\beta_d}^{(0,d)} i_{\beta_d}^{(d)}}  }_{\beta_d \text{times}}  \notag \\
&\hspace{-1.2cm}  \underbrace{\sum_{j_1^{(1)} , j_1^{(0,1)} = 1}^d h_1 [\widetilde{{A}}]_{j_1^{(1)} 1} [{A}_{T}]_{j_1^{(0,1)} j_1^{(1)}}  \cdots  \sum_{j_{\gamma_1}^{(1)} , j_{\gamma_1}^{(0,1)}= 1}^d h_1 [\widetilde{{A}}]_{j_{\gamma_1}^{(1)} 1} [{A}_{T}]_{j_{\gamma_1}^{(0,1)} j_{\gamma_1}^{(1)}} }_{\gamma_1 \text{times}}  \cdots   \notag \\
&\hspace{-1.2cm}  \underbrace{\sum_{j_1^{(d)} , j_1^{(0,d)} = 1}^d h_d [\widetilde{{A}}]_{j_1^{(d)} d} [{A}_{T}]_{j_1^{(0,d)} j_1^{(d)}}  \cdots \sum_{j_{\gamma_d}^{(d)} , j_{\gamma_d}^{(0,d)}= 1}^d h_d [\widetilde{{A}}]_{j_{\gamma_d}^{(d)} d} [{A}_{T}]_{j_{\gamma_d}^{(0,d)} j_{\gamma_d}^{(d)}}  }_{\gamma_d \text{times}} \notag \\
&\hspace{-1.2cm}  \underbrace{\frac{\partial^{\beta_1}}{\partial {x}_{i_1^{(0,1)}}^{} \cdots \partial {x}_{i_{\beta_1}^{(0,1)}}^{}}}_{\beta_1 \text{times}} \cdots \underbrace{\frac{\partial^{\beta_d}}{\partial {x}_{i_1^{(0,d)}}^{} \cdots \partial {x}_{i_{\beta_d}^{(0,d)}}^{}}}_{\beta_d \text{times}}  \underbrace{ \frac{\partial^{\gamma_1}}{ \partial {x}_{j_1^{(0,1)}}^{} \cdots \partial {x}_{j_{\gamma_1}^{(0,1)}}^{} } }_{\gamma_1 \text{times}} \cdots \underbrace{\frac{\partial^{\gamma_d}}{ \partial {x}_{j_1^{(0,d)}}^{} \cdots \partial {x}_{j_{\gamma_d}^{(0,d)}}^{} }}_{\gamma_d \text{times}} v_{\nu}.
\end{align*}
Using \eqref{Anorm}, \eqref{vec=ATvec}, \eqref{CN331c} and \eqref{CN332}, we have, for $1 \leq i,k \leq d$,
\begin{align*}
\displaystyle
\left| \frac{\partial \hat{v}_k}{\partial \hat{x}_{i}} \right|
&\leq |\det (A_{\widetilde{T}}) |  |\det ({A}_T) | h_k^{-1} \sum_{\eta,\nu=1}^d |  [ \widetilde{{A}}^{-1}]_{k \eta} |  | [ {A}_T^{-1}]_{\eta \nu} | \left| h_i  \sum_{i_1^{(1)}, i_1^{(0,1)} = 1}^d [{A}_T]_{{i_1^{(0,1)}} {i_1^{(1)}}} (\tilde{r}_i)_{i_1^{(1)}} \frac{\partial v_{\nu}}{\partial x_{{i_1^{(0,1)}}}} \right| \\
 &\leq  c |\det (A_{\widetilde{T}}) | h_k^{-1} \| \widetilde{A}^{-1} \|_2 \sum_{\nu=1}^d h_i \left| \frac{\partial v_{\nu}}{\partial r_i}  \right|,
\end{align*}
and, for $1 \leq i,j,k \leq d$,
\begin{align*}
\displaystyle
\left| \frac{\partial^2 \hat{v}_k}{\partial \hat{x}_{i} \partial \hat{x}_{j}} \right|
&\leq |\det (A_{\widetilde{T}})| | \det ({A}_T) | h_k^{-1} \sum_{\eta,\nu=1}^d | [ \widetilde{{A}}^{-1}]_{k \eta} | | [ {A}_T^{-1}]_{\eta \nu} | \\
&\quad \left | \sum_{i_1^{(1)}, i_1^{(0,1)} = 1}^d h_i [\widetilde{{A}}]_{{{i_1^{(1)}}} {i}} [{A}_T]_{{i_1^{(0,1)}} {i_1^{(1)}}} \sum_{j_1^{(1)}, j_1^{(0,1)} = 1}^d h_j [\widetilde{{A}}]_{{{j_1^{(1)}}} {j}} [{A}_T]_{{j_1^{(0,1)}} {j_1^{(1)}}} \frac{\partial^2 v_{\nu}}{\partial x_{i_1^{(0,1)}} \partial x_{j_1^{(0,1)}}} \right| \\
&\leq c  |\det (A_{\widetilde{T}}) | h_k^{-1} \| \widetilde{A}^{-1} \|_2 \sum_{\nu=1}^d h_i h_j \left| \frac{\partial^2 v_{\nu}}{\partial r_i {\partial r_j} } \right|.
\end{align*}

\subsubsection{Calculations 2} \label{cal941b}
We use the following calculations in \eqref{RT43}. Let $\hat{v} \in \mathcal{C}^{2}(\widehat{T})^d$ with $\tilde{v} = {\Psi}_{\widehat{T}} \hat{v}$ and ${v} = {\Psi}_{\widetilde{T}} \tilde{v}$. Using the definition of Piola transformations (Definition \ref{piola=defi}) yields, for $1 \leq i,k \leq d$,
\begin{align*}
\displaystyle
  \frac{\partial \hat{v}_k}{\partial \hat{x}_{i} }
 &=  \det ({A}_{\widetilde{T}}) h_k^{-1} \sum_{\eta=1}^d [ \widetilde{A}^{-1}]_{k \eta} \sum_{{i_1^{(1)}} = 1}^d h_i  [\widetilde{A}]_{{{i_1^{(1)}}} {i}} \frac{\partial \tilde{v}_{\eta}}{\partial \tilde{x}_{{i_1^{(1)}}}},
\end{align*}
and, for $1 \leq i,j,k \leq d$,
\begin{align*}
\displaystyle
 \frac{\partial^2 \hat{v}_k}{\partial \hat{x}_{i} \partial \hat{x}_{j}}
 &= \det ({A}_{\widetilde{T}}) h_k^{-1} \sum_{\eta=1}^d [\widetilde{{A}}^{-1}]_{k \eta} \sum_{i_1^{(1)}  =1}^d  h_i[\widetilde{{A}}]_{i_1^{(1)} i}  \sum_{j_1^{(1)} =1}^d h_j[ \widetilde{{A}} ]_{j_1^{(1)} j} \frac{\partial^2 \tilde{v}_{\eta}}{\partial \tilde{x}_{i_1^{(1)}} \partial \tilde{x}_{j_1^{(1)}}}.
\end{align*}
For any multi-indices $\beta$ and $\gamma$, for $1 \leq k \leq d$, Let $\hat{v} \in \mathcal{C}^{|\beta| + |\gamma|}(\widehat{T})^d$ with $\tilde{v} = {\Psi}_{\widehat{T}} \hat{v}$ and ${v} = {\Psi}_{\widetilde{T}} \tilde{v}$. Then,
\begin{align*}
\displaystyle
\partial^{\beta + \gamma}_{\hat{x}}  \hat{v}_k &= \frac{\partial^{|\beta| + |\gamma|}}{\partial \hat{x}_1^{\beta_1} \cdots \partial \hat{x}_d^{\beta_d} \partial \hat{x}_1^{\gamma_1} \cdots \partial \hat{x}_d^{\gamma_d}}  \hat{v}_k \notag\\
&\hspace{-1.5cm}  =  \det ({A}_{\widetilde{T}}) h_k^{-1} \sum_{\eta=1}^d [ \widetilde{A}^{-1}]_{k \eta} \\
&\hspace{-1.2cm}  \underbrace{\sum_{i_1^{(1)} = 1}^d h_1 [\widetilde{{A}}]_{i_1^{(1)} 1} \cdots   \sum_{i_{\beta_1}^{(1)} = 1}^d h_1 [\widetilde{{A}}]_{i_{\beta_1}^{(1)} 1}   }_{\beta_1 \text{times}} \cdots \underbrace{ \sum_{i_1^{(d)}  = 1}^d h_d [\widetilde{{A}}]_{i_1^{(d)} d}  \cdots \sum_{i_{\beta_d}^{(d)} = 1}^d h_d [\widetilde{{A}}]_{i_{\beta_d}^{(d)} d}  }_{\beta_d \text{times}}  \notag \\
&\hspace{-1.2cm}  \underbrace{\sum_{j_1^{(1)}  = 1}^d h_1 [\widetilde{{A}}]_{j_1^{(1)} 1} \cdots  \sum_{j_{\gamma_1}^{(1)}= 1}^d h_1 [\widetilde{{A}}]_{j_{\gamma_1}^{(1)} 1} }_{\gamma_1 \text{times}}  \cdots  \underbrace{\sum_{j_1^{(d)} = 1}^d h_d [\widetilde{{A}}]_{j_1^{(d)} d}  \cdots \sum_{j_{\gamma_d}^{(d)} = 1}^d h_d [\widetilde{{A}}]_{j_{\gamma_d}^{(d)} d}  }_{\gamma_d \text{times}} \notag \\
&\hspace{-1.2cm}  \underbrace{\frac{\partial^{\beta_1}}{\partial \tilde{x}_{i_1^{(1)}} \cdots \partial \tilde{x}_{i_{\beta_1}^{(1)}}}}_{\beta_1 \text{times}} \cdots \underbrace{\frac{\partial^{\beta_d}}{\partial \tilde{x}_{i_1^{(d)}} \cdots \partial \tilde{x}_{i_{\beta_d}^{(d)}}}}_{\beta_d \text{times}} \underbrace{ \frac{\partial^{\gamma_1}}{ \partial \tilde{x}_{j_1^{(1)}} \cdots \partial \tilde{x}_{j_{\gamma_1}^{(1)}} } }_{\gamma_1 \text{times}} \cdots \underbrace{\frac{\partial^{\gamma_d}}{ \partial \tilde{x}_{j_1^{(d)}} \cdots \partial \tilde{x}_{j_{\gamma_d}^{(d)}} }}_{\gamma_d \text{times}} \tilde{v}_{\eta}.
\end{align*}
Using Section \ref{add=new=cond=10} and \eqref{Anorm}, we have, for $1 \leq i,k \leq d$,
\begin{align*}
\displaystyle
 \left| \frac{\partial \hat{v}_k}{\partial \hat{x}_{i} } \right|
 &\leq | \det ({A}_{\widetilde{T}}) | h_k^{-1} \sum_{\eta=1}^d | [ \widetilde{A}^{-1}]_{k \eta} | \sum_{{i_1^{(1)}} = 1}^d h_i  | [\widetilde{A}]_{{{i_1^{(1)}}} {i}} | \left| \frac{\partial \tilde{v}_{\eta}}{\partial \tilde{x}_{{i_1^{(1)}}}} \right| \\
 &\leq c  | \det ({A}_{\widetilde{T}}) | h_k^{-1} \| \widetilde{{A}}^{-1} \|_{2}  \sum_{\eta=1}^d \sum_{{i_1^{(1)}} = 1}^d \widetilde{\mathscr{H}}_{i_1^{(1)}} \left| \frac{\partial \tilde{v}_{\eta}}{\partial \tilde{x}_{{i_1^{(1)}}}} \right|,
\end{align*}
and, for $1 \leq i,j,k \leq d$,
\begin{align*}
\displaystyle
\left| \frac{\partial^2 \hat{v}_k}{\partial \hat{x}_{i} \partial \hat{x}_{j}} \right|
&\leq | \det ({A}_{\widetilde{T}}) | h_k^{-1} \sum_{\eta=1}^d | [\widetilde{{A}}^{-1}]_{k \eta} | \sum_{i_1^{(1)}  =1}^d  h_i | [\widetilde{{A}}]_{i_1^{(1)} i} |  \sum_{j_1^{(1)} =1}^d h_j | [ \widetilde{{A}} ]_{j_1^{(1)} j} | \left|  \frac{\partial^2 \tilde{v}_{\eta}}{\partial \tilde{x}_{i_1^{(1)}} \partial \tilde{x}_{j_1^{(1)}}} \right| \\
&\leq c  | \det ({A}_{\widetilde{T}}) | h_k^{-1} \| \widetilde{{A}}^{-1} \|_{2}  \sum_{\eta=1}^d \sum_{i_1^{(1)},j_1^{(1)} = 1}^d \widetilde{\mathscr{H}}_{i_1^{(1)}} \widetilde{\mathscr{H}}_{j_1^{(1)}} \left|  \frac{\partial^2 \tilde{v}_{\eta}}{\partial \tilde{x}_{i_1^{(1)}} \partial \tilde{x}_{j_1^{(1)}}} \right|.
\end{align*}


\subsection{Main Results}

\begin{lem}
Let $p \in [1,\infty)$. Let $T\in \mathbb{T}_h$ satisfy Condition \ref{cond1} or Condition \ref{cond2} with $T ={\Phi}_{T} (\widetilde{T})$ and $\widetilde{T} = {\Phi}_{\widetilde{T}}(\widehat{T})$. Then,  for any $\hat{v} = (\hat{v}_1,\ldots,\hat{v}_d)^{\top} \in L^{p}(\widehat{T})^d$ with $\tilde{v} = (\tilde{v}_1,\ldots,\tilde{v}_d)^{\top} := \Psi_{\widehat{T}} \hat{v}$ and ${v} = ({v}_1,\ldots,{v}_d)^{\top} := \Psi_{\widetilde{T}} \tilde{v}$, 
\begin{align}
\displaystyle
\| v \|_{L^p({T})^d} \leq c |\det(A_{\widetilde{T}})|^{\frac{1-p }{p}} \| \widetilde{{A}} \|_{2} \left(  \sum_{j=1}^d h_j^p \| \hat{v}_j \|_{L^p(\widehat{T})}^p \right)^{1/p}. \label{RT41}
\end{align}
\end{lem}

\begin{pf*}
Because the space $\mathcal{C}(\widehat{T})^d$ is dense in the space $L^{p}(\widehat{T})^d$, we show \eqref{RT41} for $\hat{v} \in \mathcal{C}(\widehat{T})^d$ with $\tilde{v} = {\Psi}_{\widehat{T}} \hat{v}$ and  ${v} = {\Psi}_{\widetilde{T}} \tilde{v}$. From the definition of the Piola transformation, for $i=1,\ldots,d$,
\begin{align*}
\displaystyle
 v_i(x) &=  \frac{1}{\det ({{A}}_{T}) } \sum_{k=1}^d [A_T]_{ik} \tilde{v}_k (\tilde{x}), \quad
 \tilde{v}_i(\tilde{x}) = \frac{1}{\det ({A}_{\widetilde{T}})} \sum_{j=1}^d [\widetilde{A}]_{ij} h_j \hat{v}_j(\hat{x}),
\end{align*}
which leads to
\begin{align*}
\displaystyle
 v_i(x) &=  \frac{1}{\det ({{A}}_{T}) \det ({A}_{\widetilde{T}})} \sum_{j,k=1}^d  [A_T]_{ik} [\widetilde{A}]_{kj} h_j \hat{v}_j(\hat{x}).
\end{align*}
If $1 \leq p \< \infty$, using \eqref{Anorm}, \eqref{CN331c} and \eqref{CN332}, 
\begin{align*}
\displaystyle
\| v \|^p_{L^p(T)^d} &= \sum_{i=1}^d \| v_i \|_{L^p(T)}^p \leq c |\det(A_{\widetilde{T}})|^{1-p} \| \widetilde{{A}} \|_2^p \sum_{j=1}^d  h_{j}^p \| \hat{v}_j \|^p_{L^p(\widehat{T})},
\end{align*}
which leads to \eqref{RT41}.
\qed
\end{pf*}

\begin{lem}
Let $p \in [1,\infty)$. Let $T\in \mathbb{T}_h$ satisfy Condition \ref{cond1} or Condition \ref{cond2} with $T ={\Phi}_{T} (\widetilde{T})$ and $\widetilde{T} = {\Phi}_{\widetilde{T}}(\widehat{T})$. Let $\ell , m \in \mathbb{N}_0$ and  $k \in \mathbb{N}$ with $1 \leq k \leq d$. Let $\beta := (\beta_1,\ldots,\beta_d) \in \mathbb{N}_0^d$ and $\gamma := (\gamma_1,\ldots,\gamma_d) \in \mathbb{N}_0^d$ be multi-indices with $|\beta| = \ell$ and $|\gamma| = m$, respectively.  Then, for any $\hat{v} = (\hat{v}_1,\ldots,\hat{v}_d)^{\top} \in W^{ |\beta|+|\gamma|,p}(\widehat{T})^d$ with $\tilde{v} = (\tilde{v}_1,\ldots,\tilde{v}_d)^{\top} := \Psi_{\widehat{T}} \hat{v}$ and ${v} = (v_1, \ldots,v_d)^{\top} := \Psi_{\widetilde{T}} \tilde{v}$, 
\begin{align}
\displaystyle
\left \|  \partial_{\hat{x}}^{\beta}  \partial_{\hat{x}}^{\gamma} \hat{v}_k \right \|_{L^{p}(\widehat{T})} &\leq c  |\det ({A}_{\widetilde{T}})|^{\frac{p-1}{p}} h_k^{-1} \| \widetilde{{A}}^{-1} \|_{2}  \sum_{|\varepsilon| = \ell + m} h^{\varepsilon} \left \| \partial_{r}^{ \varepsilon} v  \right \|_{L^p(T)^d}. \label{RT42}
\end{align}
If Condition \ref{Cond333} is imposed, then
\begin{align}
\displaystyle
\left \|  \partial_{\hat{x}}^{\beta} \partial_{\hat{x}}^{\gamma} \hat{v}_k \right \|_{L^{p}(\widehat{T})} \leq c  |\det ({A}_{\widetilde{T}})|^{\frac{p-1}{p}} h_k^{-1} \| \widetilde{{A}}^{-1} \|_{2} \sum_{|\varepsilon| = \ell + m} \widetilde{\mathscr{H}}^{\varepsilon} \| \partial_{\tilde{x}}^{\varepsilon} (\Psi_{\widetilde{T}}^{-1}  v ) \|_{L^p( \Phi_{T}^{-1} (T))^d}.  \label{RT43}
\end{align}
\end{lem}

\begin{pf*}
Because the space $\mathcal{C}^{\ell+m}(\widehat{T})^d$ is dense
in the space $W^{\ell+m,p}(\widehat{T})^d$,  we show \eqref{RT42} and \eqref{RT43} for $\hat{v} \in \mathcal{C}^{\ell+m}(\widehat{T})^d$ with $\tilde{v} = {\Psi}_{\widehat{T}} \hat{v}$ and ${v} = {\Psi}_{\widetilde{T}} \tilde{v}$. 

\textbf{Proof of \eqref{RT42}.} Using \eqref{Anorm}, \eqref{CN331c} and \eqref{CN332}, through a simple calculation, we have, for $1 \leq k \leq d$,
\begin{align*}
\displaystyle
| \partial^{\beta + \gamma}_{\hat{x}}  \hat{v}_k| &= \left| \frac{\partial^{|\beta| + |\gamma|}}{\partial \hat{x}_1^{\beta_1} \cdots \partial \hat{x}_d^{\beta_d} \partial \hat{x}_1^{\gamma_1} \cdots \partial \hat{x}_d^{\gamma_d}}  \hat{v}_k \right| \notag\\
&\hspace{-1.5cm}  \leq c  |\det (A_{\widetilde{T}}) |  \| \widetilde{A}^{-1} \|_2 h_k^{-1} \sum_{\nu=1}^d h^{\beta} h^{\gamma} \notag\\
&\hspace{-1.2cm} \Biggl | \underbrace{\sum_{i_1^{(1)},i_1^{(0,1)} = 1}^d  [{A}_{T}]_{i_1^{(0,1)} i_1^{(1)}} (\tilde{r}_1)_{i_1^{(1)}} \cdots   \sum_{i_{\beta_1}^{(1)},i_{\beta_1}^{(0,1)} = 1}^d  [{A}_{T}]_{i_{\beta_1}^{(0,1)} i_{\beta_1}^{(1)}} (\tilde{r}_1)_{i_{\beta_1}^{(1)}}  }_{\beta_1 \text{times}} \cdots \notag \\
&\hspace{-1.2cm} \underbrace{ \sum_{i_1^{(d)} , i_1^{(0,d)} = 1}^d   [{A}_{T}]_{i_1^{(0,d)} i_1^{(d)}} (\tilde{r}_d)_{i_{1}^{(d)}}  \cdots \sum_{i_{\beta_d}^{(d)} , i_{\beta_d}^{(0,d)}= 1}^d [{A}_{T}]_{i_{\beta_d}^{(0,d)} i_{\beta_d}^{(d)}} (\tilde{r}_d)_{i_{\beta_d}^{(d)}}  }_{\beta_d \text{times}}  \notag \\
&\hspace{-1.2cm}  \underbrace{\sum_{j_1^{(1)} , j_1^{(0,1)} = 1}^d   [{A}_{T}]_{j_1^{(0,1)} j_1^{(1)}} (\tilde{r}_1)_{j_{1}^{(1)}}   \cdots  \sum_{j_{\gamma_1}^{(1)} , j_{\gamma_1}^{(0,1)}= 1}^d  [{A}_{T}]_{j_{\gamma_1}^{(0,1)} j_{\gamma_1}^{(1)}} (\tilde{r}_1)_{ j_{\gamma_1}^{(1)}} }_{\gamma_1 \text{times}}  \cdots   \notag \\
&\hspace{-1.2cm}  \underbrace{\sum_{j_1^{(d)} , j_1^{(0,d)} = 1}^d  [{A}_{T}]_{j_1^{(0,d)} j_1^{(d)}} (\tilde{r}_d)_{j_1^{(d)}}  \cdots \sum_{j_{\gamma_d}^{(d)} , j_{\gamma_d}^{(0,d)}= 1}^d [{A}_{T}]_{j_{\gamma_d}^{(0,d)} j_{\gamma_d}^{(d)}} (\tilde{r}_d)_{j_{\gamma_d}^{(d)}} }_{\gamma_d \text{times}} \notag \\
&\hspace{-1.2cm}  \underbrace{\frac{\partial^{\beta_1}}{\partial {x}_{i_1^{(0,1)}}^{} \cdots \partial {x}_{i_{\beta_1}^{(0,1)}}^{}}}_{\beta_1 \text{times}} \cdots \underbrace{\frac{\partial^{\beta_d}}{\partial {x}_{i_1^{(0,d)}}^{} \cdots \partial {x}_{i_{\beta_d}^{(0,d)}}^{}}}_{\beta_d \text{times}}  \underbrace{ \frac{\partial^{\gamma_1}}{ \partial {x}_{j_1^{(0,1)}}^{} \cdots \partial {x}_{j_{\gamma_1}^{(0,1)}}^{} } }_{\gamma_1 \text{times}} \cdots \underbrace{\frac{\partial^{\gamma_d}}{ \partial {x}_{j_1^{(0,d)}}^{} \cdots \partial {x}_{j_{\gamma_d}^{(0,d)}}^{} }}_{\gamma_d \text{times}} v_{\nu} \Biggr | \\
&\hspace{-1.5cm}  \leq c  |\det (A_{\widetilde{T}}) |  \| \widetilde{A}^{-1} \|_2 h_k^{-1}  \sum_{\nu=1}^d \sum_{|\varepsilon| = |\beta|+|\gamma|} h^{\varepsilon} | \partial_{r}^{\varepsilon} v_{\nu} |.
\end{align*}
Because $1 \leq p \< \infty$, it holds that, for $1 \leq k \leq d$,
\begin{align*}
\displaystyle
\left \|  \partial^{\beta}_{\hat{x}} \partial^{\gamma}_{\hat{x}}  \hat{v}_k  \right \|_{L^{p}(\widehat{T})}^p 
&\leq c  | \det ({A}_{\widetilde{T}}) |^{p-1} \| \widetilde{{A}}^{-1} \|_2^p h_k^{-p} \sum_{|\varepsilon| = |\beta|+|\gamma|} {h}^{\varepsilon p} \int_{T} |\partial^{\varepsilon}_{r} v|^p dx,
\end{align*}
which leads to \eqref{RT42} together with \eqref{jensen}.

\textbf{Proof of \eqref{RT43}.}
Using  Section \ref{add=new=cond=10} and \eqref{Anorm}, through a simple calculation, we have, for $1 \leq k \leq d$,
\begin{align*}
\displaystyle
| \partial^{\beta + \gamma}_{\hat{x}} \hat{v}_k |
&= \left|  \frac{\partial^{|\beta| + |\gamma|}}{\partial \hat{x}_1^{\beta_1} \cdots \partial \hat{x}_d^{\beta_d} \partial \hat{x}_1^{\gamma_1} \cdots \partial \hat{x}_d^{\gamma_d}}  \hat{v}_k \right| \\
&\hspace{-1.5cm}  \leq c  |\det (A_{\widetilde{T}}) |  \| \widetilde{A}^{-1} \|_2 h_k^{-1} \\
&\hspace{-1.2cm} \sum_{\eta=1}^d   \underbrace{\sum_{i_1^{(1)} = 1}^d  \cdots   \sum_{i_{\beta_1}^{(1)} = 1}^d  }_{\beta_1 \text{times}} \cdots \underbrace{ \sum_{i_1^{(d)}  = 1}^d  \cdots \sum_{i_{\beta_d}^{(d)} = 1}^d  }_{\beta_d \text{times}}   \underbrace{\sum_{j_1^{(1)}  = 1}^d  \cdots  \sum_{j_{\gamma_1}^{(1)}= 1}^d  }_{\gamma_1 \text{times}}  \cdots  \underbrace{\sum_{j_1^{(d)} = 1}^d   \cdots \sum_{j_{\gamma_d}^{(d)} = 1}^d  }_{\gamma_d \text{times}} \notag \\
&\hspace{-1.2cm}  \underbrace{\widetilde{\mathscr{H}}_{i_1^{(1)}} \cdots \widetilde{\mathscr{H}}_{i_{\varepsilon_1}^{(1)}}}_{\beta_1 \text{times}} \cdots  \underbrace{ \widetilde{\mathscr{H}}_{i_1^{(d)}} \cdots \widetilde{\mathscr{H}}_{i_{\varepsilon_d}^{(d)}}}_{\beta_d \text{times}}  \underbrace{\widetilde{\mathscr{H}}_{j_1^{(1)}} \cdots \widetilde{\mathscr{H}}_{j_{\varepsilon_1}^{(1)}}}_{\gamma_1 \text{times}} \cdots  \underbrace{ \widetilde{\mathscr{H}}_{j_1^{(d)}} \cdots \widetilde{\mathscr{H}}_{j_{\varepsilon_d}^{(d)}}}_{\gamma_d \text{times}} \notag \\
&\hspace{-1.2cm} \Biggl |  \underbrace{\frac{\partial^{\beta_1}}{\partial \tilde{x}_{i_1^{(1)}} \cdots \partial \tilde{x}_{i_{\beta_1}^{(1)}}}}_{\beta_1 \text{times}} \cdots \underbrace{\frac{\partial^{\beta_d}}{\partial \tilde{x}_{i_1^{(d)}} \cdots \partial \tilde{x}_{i_{\beta_d}^{(d)}}}}_{\beta_d \text{times}} \underbrace{ \frac{\partial^{\gamma_1}}{ \partial \tilde{x}_{j_1^{(1)}} \cdots \partial \tilde{x}_{j_{\gamma_1}^{(1)}} } }_{\gamma_1 \text{times}} \cdots \underbrace{\frac{\partial^{\gamma_d}}{ \partial \tilde{x}_{j_1^{(d)}} \cdots \partial \tilde{x}_{j_{\gamma_d}^{(d)}} }}_{\gamma_d \text{times}} \tilde{v}_{\eta} \Biggr| \\
&\hspace{-1.5cm}  \leq c |\det (A_{\widetilde{T}}) |  \| \widetilde{A}^{-1} \|_2 h_k^{-1} \sum_{\eta=1}^d \sum_{|\varepsilon| = |\beta|+|\gamma|} \widetilde{\mathscr{H}}^{\varepsilon} |\partial^{\varepsilon}_{\tilde{x}} \tilde{v}_{\eta}|.
\end{align*}
Because $1 \leq p \< \infty$, it holds that, for $1 \leq k \leq d$,
\begin{align*}
\displaystyle
\left \|  \partial^{\beta}_{\hat{x}} \partial^{\gamma}_{\hat{x}}  \hat{v}_k  \right \|_{L^{p}(\widehat{T})}^p 
&\leq c  | \det ({A}_{\widetilde{T}}) |^{p-1} \| \widetilde{{A}}^{-1} \|_2^p h_k^{-p} \sum_{|\varepsilon| = |\beta|+|\gamma|} \widetilde{\mathscr{H}}^{\varepsilon p} \int_{\widetilde{T}} |\partial^{\varepsilon}_{\tilde{x}} \tilde{v}|^p d\tilde{x},
\end{align*}
which leads to \eqref{RT43} together with \eqref{jensen}.
\qed
\end{pf*}

\begin{rem} \label{ex=01}
In inequality \eqref{RT43}, it is possible to obtain the estimates in $T$ by specifically determining the matrix ${A}_{T}$.

Let $\hat{v} \in \mathcal{C}^{1}(\widehat{T})^d$ with $\tilde{v} = {\Psi}_{\widehat{T}} \hat{v}$ and ${v} = {\Psi}_{\widetilde{T}} \tilde{v}$. Using \eqref{Anorm}, \eqref{CN331c}, \eqref{CN332} and the definition of Piola transformations (Definition \ref{piola=defi}), we have, for $1 \leq i,k \leq d$,
\begin{align*}
\displaystyle
\left| \frac{\partial \hat{v}_k}{\partial \hat{x}_{i}}\right|
 &\leq c  |\det ({A}_{\widetilde{T}}) | \| \widetilde{{A}}^{-1} \|_2 h_k^{-1} \sum_{\nu=1}^d \sum_{i_1^{(1)}, i_1^{(0,1)} = 1}^d \widetilde{\mathscr{H}}_{i_1^{(1)}} | {A}_T]_{{i_1^{(0,1)}} {i_1^{(1)}}} |  \left| \frac{\partial v_{\nu}}{\partial x_{{i_1^{(0,1)}}}} \right |.
 \end{align*}
Let $d=3$. We define the matrix ${A}_{T}$ as 
\begin{align*}
\displaystyle
\mathcal{A}_{T} := 
\begin{pmatrix}
\cos \frac{\pi}{2}  & - \sin \frac{\pi}{2} & 0\\
 \sin \frac{\pi}{2} & \cos \frac{\pi}{2} & 0 \\
 0 & 0 & 1
\end{pmatrix}.
\end{align*}
We then have
\begin{align*}
\displaystyle
\left| \frac{\partial \hat{v}_k}{\partial \hat{x}_{i}}\right|
 &\leq c  |\det ({A}_{\widetilde{T}}) | \| \widetilde{{A}}^{-1} \|_2 h_k^{-1} \sum_{\nu=1}^3 \left( \widetilde{\mathscr{H}}_{1}  \left| \frac{\partial v_{\nu}}{\partial x_{2}} \right | +  \widetilde{\mathscr{H}}_{2}  \left| \frac{\partial v_{\nu}}{\partial x_{1}} \right | +  \widetilde{\mathscr{H}}_{3}  \left| \frac{\partial v_{\nu}}{\partial x_{3}} \right | \right).
 \end{align*}
Because $1 \leq p \< \infty$, it holds that, for $1 \leq i, k \leq 3$,
\begin{align*}
\displaystyle
\left \| \frac{\partial \hat{v}_k}{\partial \hat{x}_{i}}\right \|^p_{L^p(\widehat{T})}
&\leq c  |\det ({A}_{\widetilde{T}}) |^{p-1} \| \widetilde{{A}}^{-1} \|^p_2 h_k^{-p} \\
&\quad \times \left( \widetilde{\mathscr{H}}_1^p \left \| \frac{\partial v}{\partial x_2} \right \|^p_{L^p(T)} + \widetilde{\mathscr{H}}_2^p \left \| \frac{\partial v}{\partial x_1} \right \|^p_{L^p(T)} + \widetilde{\mathscr{H}}_3^p \left \| \frac{\partial v}{\partial x_3} \right \|^p_{L^p(T)} \right).
 \end{align*}
\end{rem}

The following two lemmata are divided into the element on $\mathfrak{T}^{(2)}$ or $\mathfrak{T}_1^{(3)}$ and the element on $\mathfrak{T}_2^{(3)}$. 

\begin{lem} \label{lem1142}
Let $p \in [1,\infty)$. Let $T\in \mathbb{T}_h$ satisfy Condition \ref{cond1} or Condition \ref{cond2} with $T ={\Phi}_{T} (\widetilde{T})$ and $\widetilde{T} = {\Phi}_{\widetilde{T}}(\widehat{T})$, where $\widetilde{T} \in \mathfrak{T}^{(2)}$ or $\widetilde{T} \in \mathfrak{T}_1^{(3)}$. Let $\beta := (\beta_1,\ldots,\beta_d) \in \mathbb{N}_0^d$ be a multi-index with $|\beta| = \ell$. Then, for any $\hat{v} \in W^{\ell+1,p}(\widehat{T})^d$ with $\tilde{v}= {\Psi}_{\widehat{T}} \hat{v}$ and ${v} = {\Psi}_{\widetilde{T}} \tilde{v}$,  
\begin{align}
\displaystyle
\left \| \partial^{\beta}_{\hat{x}} \nabla_{\hat{x}} \cdot \hat{v} \right \|_{L^p(\widehat{T}_1)}
&\leq c |\det ({A}_{\widetilde{T}})|^{\frac{p-1}{p}} \sum_{|\varepsilon| = \ell} h^{\varepsilon} \left \|  \partial_r^{\varepsilon} \nabla \cdot v \right \|_{L^p(T)}. \label{RT12}
\end{align}
If Condition \ref{Cond333} is imposed, it holds that
\begin{align}
\displaystyle
\left \| \partial^{\beta}_{\hat{x}}  \nabla_{\hat{x}} \cdot \hat{v}  \right \|_{L^p(\widehat{T}_1)}
&\leq c |\det ({A}_{\widetilde{T}})|^{\frac{p-1}{p}} \sum_{|\varepsilon| = \ell} \widetilde{\mathscr{H}}^{\varepsilon} \left \|  \partial^{\varepsilon}_{\tilde{x}} \nabla_{\tilde{x}} \cdot (\Psi_{\widetilde{T}}^{-1}  v) \right \|_{L^p(\Phi_{T}^{-1}(T))}. \label{RT13}
\end{align}
\end{lem}

\begin{pf*}
Because the space $\mathcal{C}^{\ell+1}(\widehat{T})^d$ is dense
in the space $W^{\ell+1,p}(\widehat{T})^d$,  we show \eqref{RT12} and \eqref{RT13} for $\hat{v} \in \mathcal{C}^{\ell+1}(\widehat{T})^d$ with $\tilde{v} = {\Psi}_{\widehat{T}} \hat{v}$ and ${v} = {\Psi}_{\widetilde{T}} \tilde{v}$. 

By a simple calculation from Sections \ref{cal941} and \ref{cal941b},
\begin{align*}
\displaystyle
\nabla_{\hat{x}} \cdot \hat{v} &= \sum_{k=1}^d \frac{\partial \hat{v}_k}{\partial \hat{x}_{k}} \\
 &= \det ({A}_{\widetilde{T}}) \det ({A}_T) \sum_{k,\eta,\nu,i_1^{(1)}, i_1^{(0,1)}=1}^d  [ \widetilde{{A}}^{-1}]_{k \eta }[\widetilde{{A}}]_{{{i_1^{(1)}}} {k}}  [ {A}_T^{-1}]_{\eta \nu} [{A}_T]_{{i_1^{(0,1)}} {i_1^{(1)}}} \frac{\partial v_{\nu}}{\partial x_{{i_1^{(0,1)}}}} \\
  &= \det ({A}_{\widetilde{T}}) \det ({A}_T) \nabla \cdot v, \\
\frac{\partial}{\partial \hat{x}_i} \nabla_{\hat{x}} \cdot \hat{v} &= \sum_{k=1}^d \frac{\partial ^2 \hat{v}_k}{\partial \hat{x}_i \partial \hat{x}_{k}} \\
&=  \det ({A}_{\widetilde{T}}) \det ({A}_T)  h_i  \sum_{i_1^{(1)}, i_1^{(0,1)} = 1}^d [\widetilde{{A}}]_{{{i_1^{(1)}}} {i}} [{A}_T]_{{i_1^{(0,1)}} {i_1^{(1)}}}  \\
&\quad \sum_{k,\eta,\nu,j_1^{(1)}, j_1^{(0,1)}=1}^d  [ \widetilde{{A}}^{-1}]_{k \eta}[\widetilde{{A}}]_{{{j_1^{(1)}}} {k}}  [{A}_T^{-1}]_{\eta \nu}  [{A}_T]_{{j_1^{(0,1)}} {j_1^{(1)}}} \frac{\partial^2 v_{\nu}}{\partial x_{i_1^{(0,1)}} \partial x_{j_1^{(0,1)}}} \\
&=  \det ({A}_{\widetilde{T}}) \det ({A}_T)  h_i  \sum_{i_1^{(1)}, i_1^{(0,1)} = 1}^d [\widetilde{{A}}]_{{{i_1^{(1)}}} {i}} [{A}_T]_{{i_1^{(0,1)}} {i_1^{(1)}}} \frac{\partial (\nabla \cdot v)}{\partial x_{i_1^{(0,1)}}}.
 \end{align*}
For a general derivative ${\partial}_{\hat{x}}^{\beta} \nabla_{\hat{x}} \cdot \hat{v}$ with order $|\beta| = \ell$, we obtain
 \begin{align*}
\displaystyle
\partial^{\beta}_{\hat{x}} \nabla_{\hat{x}} \cdot \hat{v} &= \frac{\partial^{|\beta|}}{\partial \hat{x}_1^{\beta_1} \cdots \partial \hat{x}_d^{\beta_d}} \nabla_{\hat{x}} \cdot \hat{v} \notag\\
&\hspace{-1.5cm}  = \det ({A}_{\widetilde{T}}) \det ({A}_T)   \notag\\
&\hspace{-1.2cm}  \underbrace{\sum_{i_1^{(1)},i_1^{(0,1)} = 1}^d h_1 [\widetilde{{A}}]_{i_1^{(1)} 1} [{A}_{T}]_{i_1^{(0,1)} i_1^{(1)}} \cdots   \sum_{i_{\beta_1}^{(1)},i_{\beta_1}^{(0,1)} = 1}^d h_1 [\widetilde{{A}}]_{i_{\beta_1}^{(1)} 1} [{A}_{T}]_{i_{\beta_1}^{(0,1)} i_{\beta_1}^{(1)}}   }_{\beta_1 \text{times}} \cdots \notag \\
&\hspace{-1.2cm} \underbrace{ \sum_{i_1^{(d)} , i_1^{(0,d)} = 1}^d h_d [\widetilde{{A}}]_{i_1^{(d)} d} [{A}_{T}]_{i_1^{(0,d)} i_1^{(d)}}  \cdots \sum_{i_{\beta_d}^{(d)} , i_{\beta_d}^{(0,d)}= 1}^d h_d [\widetilde{{A}}]_{i_{\beta_d}^{(d)} d} [{A}_{T}]_{i_{\beta_d}^{(0,d)} i_{\beta_d}^{(d)}}  }_{\beta_d \text{times}}  \notag \\
&\hspace{-1.2cm}  \underbrace{\frac{\partial^{\beta_1}}{\partial {x}_{i_1^{(0,1)}}^{} \cdots \partial {x}_{i_{\beta_1}^{(0,1)}}^{}}}_{\beta_1 \text{times}} \cdots \underbrace{\frac{\partial^{\beta_d}}{\partial {x}_{i_1^{(0,d)}}^{} \cdots \partial {x}_{i_{\beta_d}^{(0,d)}}^{}}}_{\beta_d \text{times}} \nabla \cdot v \\
&\hspace{-1.5cm}  =  \det ({A}_{\widetilde{T}}) \det ({A}_T)   \notag\\
&\hspace{-1.2cm}  \underbrace{\sum_{i_1^{(1)},i_1^{(0,1)} = 1}^d h_1 [{A}_{T}]_{i_1^{(0,1)} i_1^{(1)}} (\tilde{r}_1)_{i_1^{(1)}} \cdots   \sum_{i_{\beta_1}^{(1)},i_{\beta_1}^{(0,1)} = 1}^d h_1 [{A}_{T}]_{i_{\beta_1}^{(0,1)} i_{\beta_1}^{(1)}} (\tilde{r}_1)_{i_{\beta_1}^{(1)}}  }_{\beta_1 \text{times}} \cdots \notag \\
&\hspace{-1.2cm} \underbrace{ \sum_{i_1^{(d)} , i_1^{(0,d)} = 1}^d h_d [{A}_{T}]_{i_1^{(0,d)} i_1^{(d)}} (\tilde{r}_d)_{i_1^{(d)}}  \cdots \sum_{i_{\beta_d}^{(d)} , i_{\beta_d}^{(0,d)}= 1}^d h_d [{A}_{T}]_{i_{\beta_d}^{(0,d)} i_{\beta_d}^{(d)}} (\tilde{r}_d)_{i_{\beta_d}^{(d)}} }_{\beta_d \text{times}}  \notag \\
&\hspace{-1.2cm}  \underbrace{\frac{\partial^{\beta_1}}{\partial {x}_{i_1^{(0,1)}}^{} \cdots \partial {x}_{i_{\beta_1}^{(0,1)}}^{}}}_{\beta_1 \text{times}} \cdots \underbrace{\frac{\partial^{\beta_d}}{\partial {x}_{i_1^{(0,d)}}^{} \cdots \partial {x}_{i_{\beta_d}^{(0,d)}}^{}}}_{\beta_d \text{times}} \nabla \cdot v.
\end{align*}
It then holds that, using \eqref{CN332} and \eqref{jensen},
\begin{align*}
\displaystyle
| \partial^{\beta}_{\hat{x}} \nabla_{\hat{x}} \cdot \hat{v} |
&\leq c |\det ({A}_{\widetilde{T}})| \sum_{|\varepsilon| = \ell} h^{\varepsilon} |\partial_r^{\varepsilon} \nabla \cdot v|,
\end{align*}
which leads to
\begin{align*}
\displaystyle
\| \partial^{\beta}_{\hat{x}} \nabla_{\hat{x}} \cdot \hat{v} \|_{L^p(\widehat{T})}
&\leq c |\det ({A}_{\widetilde{T}})|^{\frac{p-1}{p}}\sum_{|\varepsilon| = \ell} h^{\varepsilon} \| \partial_r^{\varepsilon} \nabla \cdot v \|_{L^p(T)}.
\end{align*}

Using an analogous argument, if Condition \ref{Cond333} is imposed, for a general derivative ${\partial}_{\hat{x}}^{\beta} \nabla_{\hat{x}} \cdot \hat{v}$ with order $|\beta| = \ell$, we obtain
\begin{align*}
\displaystyle
\partial^{\beta}_{\hat{x}} \nabla_{\hat{x}} \cdot \hat{v} &= \frac{\partial^{|\beta|}}{\partial \hat{x}_1^{\beta_1} \cdots \partial \hat{x}_d^{\beta_d}} \nabla_{\hat{x}} \cdot \hat{v} \notag\\
&\hspace{-1.5cm}  = \det ({A}_{\widetilde{T}})  \notag \\
&\hspace{-1.2cm}  \underbrace{\sum_{i_1^{(1)} = 1}^d h_1 [\widetilde{{A}}]_{i_1^{(1)} 1} \cdots   \sum_{i_{\beta_1}^{(1)} = 1}^d h_1 [\widetilde{{A}}]_{i_{\beta_1}^{(1)} 1}   }_{\beta_1 \text{times}} \cdots \underbrace{ \sum_{i_1^{(d)}  = 1}^d h_d [\widetilde{{A}}]_{i_1^{(d)} d}  \cdots \sum_{i_{\beta_d}^{(d)} = 1}^d h_d [\widetilde{{A}}]_{i_{\beta_d}^{(d)} d}  }_{\beta_d \text{times}}  \notag \\
&\hspace{-1.2cm}  \underbrace{\frac{\partial^{\beta_1}}{\partial \tilde{x}_{i_1^{(1)}} \cdots \partial \tilde{x}_{i_{\beta_1}^{(1)}}}}_{\beta_1 \text{times}} \cdots \underbrace{\frac{\partial^{\beta_d}}{\partial \tilde{x}_{i_1^{(d)}} \cdots \partial \tilde{x}_{i_{\beta_d}^{(d)}}}}_{\beta_d \text{times}} \nabla_{\tilde{x}} \cdot \tilde{v}.
\end{align*}
It then holds that
\begin{align*}
\displaystyle
&\left| \partial^{\beta}_{\hat{x}} \nabla_{\hat{x}} \cdot \hat{v} \right| \\
&\quad \leq c |\det ({A}_{\widetilde{T}})| \\
&\quad \quad  \underbrace{\sum_{i_1^{(1)} = 1}^d h_1 | [\widetilde{{A}}]_{i_1^{(1)} 1} | \cdots   \sum_{i_{\beta_1}^{(1)} = 1}^d h_1 | [\widetilde{{A}}]_{i_{\beta_1}^{(1)} 1} |  }_{\beta_1 \text{times}} \cdots \underbrace{ \sum_{i_1^{(d)}  = 1}^d h_d | [\widetilde{{A}}]_{i_1^{(d)} d} |  \cdots \sum_{i_{\beta_d}^{(d)} = 1}^d h_d | [\widetilde{{A}}]_{i_{\beta_d}^{(d)} d} |  }_{\beta_d \text{times}} \\
&\quad \quad \Biggl | \underbrace{\frac{\partial^{\beta_1}}{\partial \tilde{x}_{i_1^{(1)}} \cdots \partial \tilde{x}_{i_{\beta_1}^{(1)}}}}_{\beta_1 \text{times}} \cdots \underbrace{\frac{\partial^{\beta_d}}{\partial \tilde{x}_{i_1^{(d)}} \cdots \partial \tilde{x}_{i_{\beta_d}^{(d)}}}}_{\beta_d \text{times}} \nabla_{\tilde{x}} \cdot \tilde{v} \Biggr | \\
&\quad \leq c |\det ({A}_{\widetilde{T}})| \sum_{|\varepsilon| = \ell} \widetilde{\mathscr{H}}^{\varepsilon} | \partial_{\tilde{x}}^{\varepsilon}  \nabla_{\tilde{x}} \cdot \tilde{v} |,
\end{align*}
which leads to
\begin{align*}
\displaystyle
\| \partial^{\beta}_{\hat{x}} \nabla_{\hat{x}} \cdot \hat{v} \|_{L^p(\widehat{T}_1)}
&\leq c |\det ({A}_{\widetilde{T}})|^{\frac{p-1}{p}}\sum_{|\varepsilon| = \ell} \widetilde{\mathscr{H}}^{\varepsilon} \| \partial_{\tilde{x}}^{\varepsilon} \nabla_{\tilde{x}} \cdot \tilde{v} \|_{L^p(\widetilde{T})}.
\end{align*}
\qed
\end{pf*}

\begin{lem} \label{lem1143}
Let $p \in [1,\infty)$ and $d=3$. Let $T\in \mathbb{T}_h$ satisfy Condition \ref{cond2} with $T ={\Phi}_{T} (\widetilde{T})$ and $\widetilde{T} = {\Phi}_{\widetilde{T}}(\widehat{T})$, where $\widetilde{T} \in \mathfrak{T}_2^{(3)}$. Let $\ell \in \mathbb{N}_0$ and  $k \in \mathbb{N}$ with $1 \leq k \leq 3$. Let $\beta := (\beta_1,\beta_2,\beta_3) \in \mathbb{N}_0^3$ be a multi-index with $|\beta| = \ell$. Then,  for any $\hat{v} \in W^{\ell+1,p}(\widehat{T})^3$ with $\tilde{v}= {\Psi}_{\widehat{T}} \hat{v}$ and ${v} = {\Psi}_{\widetilde{T}} \tilde{v}$,  
\begin{align}
\displaystyle
\left \| \partial_{\hat{x}}^{\beta} \frac{\partial \hat{v}_k}{\partial \hat{x}_k} \right\|_{L^p(\widehat{T})} 
&\leq c  |\det ({A}_{\widetilde{T}})|^{\frac{p-1}{p}} \| \widetilde{{A}}^{-1} \|_2 \sum_{|\varepsilon| = \ell} h^{\varepsilon} \left \| \partial^{\varepsilon}_r \frac{\partial v}{\partial r_k} \right \|_{L^p(T)^3}.  \label{RT14}
\end{align}
If Condition \ref{Cond333} is imposed, it holds that
\begin{align}
\displaystyle
\left \| \partial^{\beta}_{\hat{x}}  \frac{\partial \hat{v}_k}{\partial \hat{x}_{k}} \right \|_{L^p(\widehat{T})}
&\leq c  |\det ({A}_{\widetilde{T}})|^{\frac{p-1}{p}} \| \widetilde{{A}}^{-1} \|_2 \sum_{|\varepsilon| = \ell} \widetilde{\mathscr{H}}^{\varepsilon} \left \| \partial^{\varepsilon}_{\tilde{x}} \frac{\partial (\Psi_{\widetilde{T}}^{-1} v)}{\partial \tilde{r}_k} \right \|_{L^p(\Phi_{T}^{-1}(T))^3}. \label{RT14b}
\end{align}

\end{lem}

\begin{pf*}
Because the space $\mathcal{C}^{\ell+1}(\widehat{T})^3$ is dense
in the space $W^{\ell+1,p}(\widehat{T})^3$,  we show \eqref{RT14} and \eqref{RT14b} for $\hat{v} \in \mathcal{C}^{\ell+1}(\widehat{T})^3$ with $\tilde{v} = {\Psi}_{\widehat{T}} \hat{v}$ and ${v} = {\Psi}_{\widetilde{T}} \tilde{v}$. 

By a simple calculation from Section \ref{cal941}, for $1 \leq i,k \leq 3$,
\begin{align*}
\displaystyle
 \frac{\partial \hat{v}_k}{\partial \hat{x}_{k}}
 &= \det ({A}_{\widetilde{T}}) \det ({A}_T) \sum_{\eta,\nu,i_1^{(1)}, i_1^{(0,1)}=1}^3 [ \widetilde{{A}}^{-1}]_{k \eta} [\widetilde{{A}}]_{{{i_1^{(1)}}} {k}} [ {A}_T^{-1}]_{\eta \nu} [{A}_T]_{{i_1^{(0,1)}} {i_1^{(1)}}} \frac{\partial v_{\nu}}{\partial x_{{i_1^{(0,1)}}}} \\
 &= \det ({A}_{\widetilde{T}}) \det ({A}_T)   \sum_{\eta,\nu,i_1^{(1)}, i_1^{(0,1)}=1}^3 [ \widetilde{{A}}^{-1}]_{k \eta}  [ {A}_T^{-1}]_{\eta \nu} [{A}_T]_{{i_1^{(0,1)}} {i_1^{(1)}}} (\tilde{r}_k)_{i_1^{(1)}} \frac{\partial v_{\nu}}{\partial x_{{i_1^{(0,1)}}}} \\
  &= \det ({A}_{\widetilde{T}}) \det ({A}_T)   \sum_{\eta,\nu=1}^3 [ \widetilde{{A}}^{-1}]_{k \eta}  [ {A}_T^{-1}]_{\eta \nu} \frac{\partial v_{\nu}}{\partial r_{k}}, \\
 \frac{\partial^2 \hat{v}_k}{\partial \hat{x}_{i} \partial \hat{x}_{k}} 
&= \det ({A}_{\widetilde{T}}) \det ({A}_T)  \sum_{\eta,\nu,i_1^{(1)}, i_1^{(0,1)},j_1^{(1)}, j_1^{(0,1)}=1}^3 [ \widetilde{{A}}^{-1}]_{k \eta} [ {A}_T^{-1}]_{\eta \nu} \\
&\quad  h_i [\widetilde{{A}}]_{{{i_1^{(1)}}} {i}} [\widetilde{{A}}]_{{{j_1^{(1)}}} {k}} [{A}_T]_{{i_1^{(0,1)}} {i_1^{(1)}}} [{A}_T]_{{j_1^{(0,1)}} {j_1^{(1)}}} \frac{\partial^2 v_{\nu}}{\partial x_{i_1^{(0,1)}} \partial x_{j_1^{(0,1)}}} \\
&= \det ({A}_{\widetilde{T}}) \det ({A}_T)  \sum_{\eta,\nu,i_1^{(1)}, i_1^{(0,1)},j_1^{(1)}, j_1^{(0,1)}=1}^3 [ \widetilde{{A}}^{-1}]_{k \eta} [ {A}_T^{-1}]_{\eta \nu} \\
&\quad  h_i [{A}_T]_{{i_1^{(0,1)}} {i_1^{(1)}}} (\tilde{r}_i)_{i_1^{(1)}} [{A}_T]_{{j_1^{(0,1)}} {j_1^{(1)}}} (\tilde{r}_k)_{j_1^{(1)}} \frac{\partial^2 v_{\nu}}{\partial x_{i_1^{(0,1)}} \partial x_{j_1^{(0,1)}}} \\
&= \det ({A}_{\widetilde{T}}) \det ({A}_T)  \sum_{\eta,\nu=1}^3 [ \widetilde{{A}}^{-1}]_{k \eta} [ {A}_T^{-1}]_{\eta \nu} h_i \frac{\partial^2 v_{\nu}}{\partial r_{i} \partial r_{k}}.
 \end{align*}
For a general derivative ${\partial}_{\hat{x}}^{\beta} \frac{\partial \hat{v}_k}{\partial \hat{x}_{k}}$ ($1 \leq k \leq 3$) with order $|\beta| = \ell$, we obtain
 \begin{align*}
\displaystyle
\partial^{\beta}_{\hat{x}}  \frac{\partial \hat{v}_k}{\partial \hat{x}_{k}} &= \frac{\partial^{|\beta| }}{\partial \hat{x}_1^{\beta_1} \partial \hat{x}_2^{\beta_2} \partial \hat{x}_3^{\beta_3} }   \frac{\partial \hat{v}_k}{\partial \hat{x}_{k}} \notag\\
&\hspace{-1.5cm}  = \det ({A}_{\widetilde{T}}) \det ({A}_T)  \sum_{\eta,\nu=1}^3 [ \widetilde{{A}}^{-1}]_{k \eta} [ {A}_T^{-1}]_{\eta \nu}  \notag\\
&\hspace{-1.2cm}  \underbrace{\sum_{i_1^{(1)},i_1^{(0,1)} = 1}^3 h_1  [{A}_{T}]_{i_1^{(0,1)} i_1^{(1)}} (\tilde{r}_1)_{i_1^{(1)}} \cdots   \sum_{i_{\beta_1}^{(1)},i_{\beta_1}^{(0,1)} = 1}^3 h_1 [{A}_{T}]_{i_{\beta_1}^{(0,1)} i_{\beta_1}^{(1)}} (\tilde{r}_1)_{i_{\beta_1}^{(1)}}   }_{\beta_1 \text{times}} \cdots \notag \\
&\hspace{-1.2cm} \underbrace{ \sum_{i_1^{(3)} , i_1^{(0,3)} = 1}^3 h_3  [{A}_{T}]_{i_1^{(0,3)} i_1^{(3)}} (\tilde{r}_3)_{i_1^{(3)}}  \cdots \sum_{i_{\beta_3}^{(3)} , i_{\beta_3}^{(0,3)}= 1}^3 h_3 [{A}_{T}]_{i_{\beta_3}^{(0,3)} i_{\beta_3}^{(3)}} (\tilde{r}_3)_{i_{\beta_3}^{(3)}} }_{\beta_d \text{times}}  \notag \\
&\hspace{-1.2cm}  \underbrace{\frac{\partial^{\beta_1}}{\partial {x}_{i_1^{(0,1)}}^{} \cdots \partial {x}_{i_{\beta_1}^{(0,1)}}^{}}}_{\beta_1 \text{times}} \cdots \underbrace{\frac{\partial^{\beta_3}}{\partial {x}_{i_1^{(0,3)}}^{} \cdots \partial {x}_{i_{\beta_d}^{(0,d)}}^{}}}_{\beta_d \text{times}}   \frac{\partial v_{\nu}}{\partial r_k} \\
&\hspace{-1.5cm}  = \det ({A}_{\widetilde{T}}) \det ({A}_T) \sum_{\eta,\nu=1}^3 [ \widetilde{{A}}^{-1}]_{k \eta} [ {A}_T^{-1}]_{\eta \nu}   h^{\beta}  \underbrace{\frac{\partial^{\beta_1}}{\partial {r}_{1}^{} \cdots \partial {r}_{1}^{}}}_{\beta_1 \text{times}} \cdots \underbrace{\frac{\partial^{\beta_3}}{\partial {r}_{3}^{} \cdots \partial {r}_{3}^{}}}_{\beta_3 \text{times}}   \frac{\partial v_{\nu}}{\partial r_k}.
\end{align*}
It then holds that, using \eqref{Anorm}, \eqref{CN331c} and \eqref{CN332},
\begin{align*}
\displaystyle
\left | \partial^{\beta}_{\hat{x}}  \frac{\partial \hat{v}_k}{\partial \hat{x}_{k}} \right|
&\leq c |\det ({A}_{\widetilde{T}})| \| \widetilde{{A}}^{-1} \|_2 \sum_{\nu=1}^3 \sum_{|\varepsilon| = |\beta|} h^{\varepsilon} \left| \partial^{\varepsilon}_r \frac{\partial v_{\nu}}{\partial r_k} \right|,
\end{align*}
which leads to, using \eqref{jensen},
\begin{align*}
\displaystyle
\left \| \partial^{\beta}_{\hat{x}}  \frac{\partial \hat{v}_k}{\partial \hat{x}_{k}} \right \|_{L^p(\widehat{T}_2)}
&\leq c  |\det ({A}_{\widetilde{T}})|^{\frac{p-1}{p}} \| \widetilde{{A}}^{-1} \|_2 \sum_{|\varepsilon| = |\beta|} h^{\varepsilon} \left \| \partial^{\varepsilon}_r \frac{\partial v}{\partial r_k} \right \|_{L^p(T)^3}.
\end{align*}

If Condition \ref{Cond333} is imposed, by a simple calculation from Section \ref{cal941b}, for $1 \leq i,k \leq 3$,
\begin{align*}
\displaystyle
 \frac{\partial \hat{v}_k}{\partial \hat{x}_{k}}
&= \det ({A}_{\widetilde{T}})  \sum_{\eta=1}^3 [ \widetilde{{A}}^{-1}]_{k \eta} \sum_{i_1^{(1)} = 1}^3 [\widetilde{{A}}]_{{{i_1^{(1)}}} {k}} \frac{\partial \tilde{v}_{\eta}}{\partial \tilde{x}_{{i_1^{(1)}}}} \\
&= \det ({A}_{\widetilde{T}})  \sum_{\eta=1}^3 [ \widetilde{{A}}^{-1}]_{k \eta} \sum_{i_1^{(1)} = 1}^3 (\tilde{r}_k)_{i_1^{(1)}} \frac{\partial \tilde{v}_{\eta}}{\partial \tilde{x}_{{i_1^{(1)}}}} \\
&= \det ({A}_{\widetilde{T}})  \sum_{\eta=1}^3 [ \widetilde{{A}}^{-1}]_{k \eta} \frac{\partial \tilde{v}_{\eta}}{\partial \tilde{r}_{k}}, \\
\frac{\partial^2 \hat{v}_k}{\partial \hat{x}_{i} \partial \hat{x}_{k}} 
&= \det ({A}_{\widetilde{T}})  \sum_{\eta,i_1^{(1)}=1}^3 [ \widetilde{{A}}^{-1}]_{k \eta} h_i  [\widetilde{{A}}]_{{{i_1^{(1)}}} {i}} \sum_{j_1^{(1)}=1}^3 [\widetilde{{A}}]_{{{j_1^{(1)}}} {k}} \frac{\partial^2 \tilde{v}_{\eta}}{\partial \tilde{x}_{i_1^{(1)}} \partial \tilde{x}_{j_1^{(1)}}} \\
&= \det ({A}_{\widetilde{T}})  \sum_{\eta,i_1^{(1)}=1}^3 [ \widetilde{{A}}^{-1}]_{k \eta} h_i  [\widetilde{{A}}]_{{{i_1^{(1)}}} {i}} \sum_{j_1^{(1)}=1}^3  (\tilde{r}_k)_{j_1^{(1)}} \frac{\partial^2 \tilde{v}_{\eta}}{\partial \tilde{x}_{i_1^{(1)}} \partial \tilde{x}_{j_1^{(1)}}} \\
&= \det ({A}_{\widetilde{T}})  \sum_{\eta,i_1^{(1)}=1}^3 [ \widetilde{{A}}^{-1}]_{k \eta} h_i  [\widetilde{{A}}]_{{{i_1^{(1)}}} {i}} \frac{\partial^2 \tilde{v}_{\eta}}{\partial \tilde{x}_{i_1^{(1)}} \partial \tilde{r}_{k}^s}.
 \end{align*}
For a general derivative ${\partial}_{\hat{x}}^{\beta} \frac{\partial \hat{v}_k}{\partial \hat{x}_{k}}$ ($1 \leq k \leq 3$) with order $|\beta| = \ell$, we obtain
\begin{align*}
\displaystyle
\partial^{\beta}_{\hat{x}}  \frac{\partial \hat{v}_k}{\partial \hat{x}_{k}} &= \frac{\partial^{|\beta| }}{\partial \hat{x}_1^{\beta_1} \partial \hat{x}_2^{\beta_2} \partial \hat{x}_3^{\beta_3} }   \frac{\partial \hat{v}_k}{\partial \hat{x}_{k}} \notag\\
&\hspace{-1.5cm}  = \det ({A}_{\widetilde{T}}) \sum_{\eta=1}^3 [ \widetilde{{A}}^{-1}]_{k \eta}   \notag \\
&\hspace{-1.2cm}  \underbrace{\sum_{i_1^{(1)} = 1}^3 h_1 [\widetilde{{A}}]_{i_1^{(1)} 1} \cdots   \sum_{i_{\beta_1}^{(1)} = 1}^3 h_1 [\widetilde{{A}}]_{i_{\beta_1}^{(1)} 1}   }_{\beta_1 \text{times}} \cdots \underbrace{ \sum_{i_1^{(3)}  = 1}^3 h_3 [\widetilde{{A}}]_{i_1^{(3)} 3}  \cdots \sum_{i_{\beta_3}^{(3)} = 1}^3 h_3 [\widetilde{{A}}]_{i_{\beta_3}^{(3)} 3}  }_{\beta_3 \text{times}}  \notag \\
&\hspace{-1.2cm}  \underbrace{\frac{\partial^{\beta_1}}{\partial \tilde{x}_{i_1^{(1)}} \cdots \partial \tilde{x}_{i_{\beta_1}^{(1)}}}}_{\beta_1 \text{times}} \cdots \underbrace{\frac{\partial^{\beta_3}}{\partial \tilde{x}_{i_1^{(3)}} \cdots \partial \tilde{x}_{i_{\beta_3}^{(3)}}}}_{\beta_3 \text{times}} \frac{\partial \tilde{v}_{\eta}}{ \partial \tilde{r}_{k}}.
\end{align*}
It then holds that, using \eqref{Anorm},
\begin{align*}
\displaystyle
&\left| \partial^{\beta}_{\hat{x}}  \frac{\partial \hat{v}_k}{\partial \hat{x}_{k}} \right| \\
&\quad \leq | \det ({A}_{\widetilde{T}}) | \| \widetilde{{A}}^{-1} \|_2 \sum_{\eta=1}^3 \\
&\quad \quad  \underbrace{\sum_{i_1^{(1)} = 1}^3 h_1 | [\widetilde{{A}}]_{i_1^{(1)} 1} | \cdots   \sum_{i_{\beta_1}^{(1)} = 1}^3 h_1 | [\widetilde{{A}}]_{i_{\beta_1}^{(1)} 1} |  }_{\beta_1 \text{times}} \cdots \underbrace{ \sum_{i_1^{(3)}  = 1}^3 h_3 | [\widetilde{{A}}]_{i_1^{(3)} 3} | \cdots \sum_{i_{\beta_3}^{(3)} = 1}^3 h_3 | [\widetilde{{A}}]_{i_{\beta_d}^{(3)} 3} |  }_{\beta_3 \text{times}} \\
&\quad \quad  \Biggl | \underbrace{\frac{\partial^{\beta_1}}{\partial \tilde{x}_{i_1^{(1)}} \cdots \partial \tilde{x}_{i_{\beta_1}^{(1)}}}}_{\beta_1 \text{times}} \cdots \underbrace{\frac{\partial^{\beta_3}}{\partial \tilde{x}_{i_1^{(3)}} \cdots \partial \tilde{x}_{i_{\beta_3}^{(3)}}}}_{\beta_3 \text{times}} \frac{\partial \tilde{v}_{\eta}}{ \partial \tilde{r}_{k}} \Biggr | \\
&\quad \leq c | \det ({A}_{\widetilde{T}}) | \| \widetilde{{A}}^{-1} \|_2 \sum_{\eta=1}^3 \sum_{|\varepsilon| = |\beta|} \widetilde{\mathscr{H}}^{\varepsilon} \left | \partial_{\tilde{x}}^{\varepsilon} \frac{\partial \tilde{v}_{\eta}}{ \partial \tilde{r}_{k}} \right|,
\end{align*}
which leads to, using \eqref{jensen},
\begin{align*}
\displaystyle
\left \| \partial^{\beta}_{\hat{x}}  \frac{\partial \hat{v}_k}{\partial \hat{x}_{k}} \right \|_{L^p(\widehat{T})}
&\leq c  |\det ({A}_{\widetilde{T}})|^{\frac{p-1}{p}} \| \widetilde{{A}}^{-1} \|_2 \sum_{|\varepsilon| = |\beta|} \widetilde{\mathscr{H}}^{\varepsilon} \left \| \partial^{\varepsilon}_{\tilde{x}} \frac{\partial \tilde{v}}{\partial \tilde{r}_k} \right \|_{L^p(\widetilde{T})^3}.
\end{align*}
\qed
\end{pf*}

\section{New RT Interpolation Error Estimates}

\subsection{RT Finite Element}
Let $T \subset \mathbb{R}^d$ be a simplex. We define a space as
\begin{align*}
\displaystyle
\mathbb{R}^k(\partial T) := \{ \varphi_h \in L^2(\partial T): \ \varphi_h |_F \in \mathbb{P}^k(F) \ \forall F \in \mathcal{F}_{T} \}.
\end{align*}

Let $\widehat{T} \subset \mathbb{R}^d$ be the reference element defined in Sections \ref{reference2d} and \ref{reference3d}. Let $\widehat{F}_i$ be the face of $\widehat{T}$ opposite to $\widehat{p}_i$.  The RT finite element on the reference element is defined by the triple $\{ \widehat{T} , \widehat{P} , \widehat{\Sigma\}}$ as follows:
 \begin{enumerate}
 \item $\widehat{P} := \mathbb{RT}^k(\widehat{T})$;
 \item $\widehat{\Sigma}$ is a set $\{ \hat{\chi}_{i} \}_{1 \leq i \leq N^{(RT)}}$ of $N^{(RT)}$ linear forms with its components such that, for any $\hat{r} \in \widehat{P}$, 
\begin{align}
\displaystyle
& \int_{\widehat{F}} \hat{r} \cdot \hat{n}_{\widehat{F}} \hat{q}_k d\hat{s}, \quad \forall \hat{q}_k \in \mathbb{R}^k(\partial \widehat{T}), \label{RT1112} \\
 &\int_{\widehat{T}} \hat{r} \cdot \hat{q}_{k-1} d\hat{x}, \quad \forall \hat{q}_{k-1} \in \mathbb{P}^{k-1}(\widehat{T})^d, \label{RT1113}
\end{align}
shere $\hat{n}_{\widehat{F}}$ denotes the outer unit normal vector of $\widehat{T}$ on the face $\widehat{F}$. Note that for $k=0$, the local degrees of freedom of type \eqref{RT1113} are violated.
 \end{enumerate}
For the simplicial RT element in $\mathbb{R}^d$, it holds that
\begin{align}
\displaystyle
\dim \mathbb{RT}^k(\widehat{T}) =
\begin{cases}
(k+1)(k+3) \quad \text{if $d=2$}, \\
\frac{1}{2}(k+1)(k+2)(k+4) \quad \text{if $d=3$}.
\end{cases} \label{RT1114}
\end{align}
The RT finite element with the local degrees of freedom with respect to \eqref{RT1112} and \eqref{RT1113} is unisolvent; for example, see \cite[Proposition 2.3.4]{BofBreFor13}.

 We set the domain of the local RT interpolation to $V(\widehat{T}) := W^{1,1}(\widehat{T})^d$; for example, see Theorem \ref{low=trace}. The local RT interpolation $I_{\widehat{T}}^{RT^k}: V(\widehat{T}) \to \widehat{P}$ is then defined as follows: For any $\hat{v} \in V(\widehat{T})$,
\begin{align}
\displaystyle
\int_{\widehat{F}} ( I_{\widehat{T}}^{RT^k} \hat{v} - \hat{v} ) \cdot \hat{n}_{\widehat{F}} \hat{q}_k d\hat{s} &= 0 \quad \forall \hat{q}_k \in\mathbb{R}^k(\partial \widehat{T}), \label{RT1115}
\end{align}
and if $k \geq 1$,
\begin{align}
\displaystyle
\int_{\widehat{T}} ( I_{\widehat{T}}^{RT^k} \hat{v} - \hat{v} ) \cdot \hat{q}_{k-1} d\hat{x} &= 0 \quad \forall \hat{q}_{k-1} \in \mathbb{P}^{k-1}(\widehat{T})^d. \label{RT1116}
\end{align}
In particular, when $k=0$, the degrees of freedom by \eqref{RT1112} are describe as
\begin{align}
\displaystyle
\hat{\chi}_i (\hat{r}) := \int_{\widehat{F}_i} \hat{r} \cdot \hat{n}_{\widehat{F}_i} d \hat{s} \quad \forall \hat{r} \in \mathbb{RT}^0(\widehat{T}), \quad \forall i \in \{ 1 , \ldots , d+1\}. \label{RT1117}
\end{align}
The local shape functions are as follows. 
\begin{align*}
\displaystyle
\hat{\theta}_{i}^{RT^0}(x) := \frac{\iota_{\widehat{F}_{i},\widehat{T}}}{d |\widehat{T}|_d} (\hat{x} - \hat{p}_{i}) \quad \forall i \in \{ 1, \ldots , d+1 \},
\end{align*}
where $\iota_{\widehat{F}_{i},\widehat{T}} := 1$ if $\hat{n}_{\widehat{F}_i}$ points outwards, and $ - 1$ otherwise \cite[Chapter 14]{ErnGue21a}. Indeed, $\hat{\theta}_{i}^{RT^0} \in \mathbb{RT}^0(\widehat{T})$ and $\hat{\chi}_i(\hat{\theta}_{j}^{RT^0}) = \delta_{ij}$ for any $i , j\in \{ 1 , \ldots , d+1\}$. The local RT interpolation $I_{\widehat{T}}^{RT^0}: V(\widehat{T}) \to \mathbb{RT}^0(\widehat{T})$ is then described as
\begin{align}
\displaystyle
 I_{\widehat{T}}^{RT^0}:V(\widehat{T}) \ni \hat{v} \mapsto I_{\widehat{T}}^{RT^0} \hat{v} := \sum_{i=1}^{d+1} \left( \int_{\widehat{F}_i} \hat{v} \cdot \hat{n}_{\widehat{F}_i} d \hat{s} \right) \hat{\theta}_{i}^{RT^0}  \in \mathbb{RT}^0(\widehat{T}).  \label{RT1119}
\end{align}

Let $\Phi_{\widetilde{T}}: \widehat{T} \to \widetilde{T}$ and $\Phi_{T}: \widetilde{T} \to T$ be the affine mappings defined in Section \ref{two=step}. Let ${\Psi}_{\widehat{T}}: V(\widehat{T}) \to V(\widetilde{T})$ and $\Psi_{\widetilde{T}} :  V(\widetilde{T}) \to V({T})$ be the Piola transformations defined in Definition \ref{piola=defi}. The triples $\{ \widetilde{T} , \widetilde{P} , \widetilde{\Sigma} \}$ and $\{ {T} , P , {\Sigma} \}$  are defined as
\begin{align*}
\displaystyle
\begin{cases}
\displaystyle
\widetilde{T} = {\Phi}_{\widetilde{T}}(\widehat{T}); \\
\displaystyle
\widetilde{P} = \{ {\Psi}_{\widehat{T}}(\hat{q}) ; \ \hat{q} \in \widehat{P} \}; \\
\displaystyle
\widetilde{\Sigma} = \{ \{ \tilde{\chi}_{i} \}_{1 \leq i \leq N^{(RT)}}; \ \tilde{\chi}_{i} = \hat{\chi}_i( \widehat{\Psi}^{-1}(\tilde{q})), \forall \tilde{q} \in \widetilde{P}, \hat{\chi}_i \in \widehat{\Sigma}  \};
\end{cases}
\end{align*}
and 
\begin{align*}
\displaystyle
\begin{cases}
\displaystyle
T = {\Phi}_{T}(\widetilde{T}); \\
\displaystyle
P = \{ {\Psi}_{\widetilde{T}}(\tilde{q}) ; \ \tilde{q} \in \widetilde{P}  \}; \\
\displaystyle
\Sigma = \{ \{ \chi_{i} \}_{1 \leq i \leq N^{(RT)}}; \ \chi_{i} = \tilde{\chi}_i( {\Psi}_{\widetilde{T}}^{-1}(q)), \forall q \in P, \tilde{\chi}_i \in \widetilde{\Sigma} \}.
\end{cases}
\end{align*}
The triples $\{ \widetilde{T} , \widetilde{P} , \widetilde{\Sigma} \}$ and $\{ {T} ,P , {\Sigma} \}$ are then the RT finite elements. Furthermore, let 
\begin{align}
\displaystyle
I_{\widetilde{T}}^{RT^k}: V(\widetilde{T}) \to \mathbb{RT}^k(\widetilde{T}) \label{RT11110}
\end{align}
and
\begin{align}
\displaystyle
I_{{T}}^{RT^k}: V(T) \to \mathbb{RT}^k(T) \label{RT11111}
\end{align}
be the associated local RT interpolation defined in \eqref{RT1115} and \eqref{RT1116}, respectively.

\begin{rem}
Let $T \subset \mathbb{R}^d$ be a simplex. Let $v \in H(\div;T)$ and $q \in \mathbb{R}^k(\partial T)$. Let $\gamma^d: H(\div;T) \to H^{- \frac{1}{2}} (\partial T)$ be the trace operator. Then, $\gamma^d (v) \cdot n_T \in H^{- \frac{1}{2}} (\partial T) $, see Theorem \ref{thrB22}. We consider 
\begin{align*}
\displaystyle
\int_{\partial T }( \gamma^d (v) \cdot n_T ) q ds.
\end{align*}
Then, we cannot take the integral over an edge $F$ of $\partial T$. Because functions $q \in \mathbb{R}^k(\partial T)$ do not belong in $H^{ \frac{1}{2}} (\partial T)$.

 As an example, we introduce the following remark (\cite[Remark 2.5.1]{BofBreFor13}).
Given a function $\chi \in H^{-\frac{1}{2}}(\partial T)$, even if we are allowed to take
\begin{align*}
\displaystyle
\int_{\partial T} \chi ds := \langle \chi , \psi \rangle \quad \text{with $\psi \equiv 1$},
\end{align*}
we cannot take the integral over an edge $F$ of $\partial T$. Because the function identically equal to $1$ on the whole boundary $\partial T$ belongs to  $H^{\frac{1}{2}}(\partial T)$, while the function that is equal to $1$ on the edge $F$ and $0$ on the rest of $\partial T$ does not belong to  $H^{\frac{1}{2}}(\partial T)$.
\end{rem}

\begin{prop} \label{prop1111}
For any $\hat{v} \in H^{1}(\widehat{T})^d$ with $v := \Psi (\hat{v})$, it holds that
\begin{align*}
\displaystyle
 {\Psi} ( I_{\widehat{T}}^{RT^k} \hat{v} ) =  I_{T}^{RT^k} ({\Psi} \hat{v}),
\end{align*}
that is, the diagrams
\begin{align*}
\displaystyle
\xymatrix{
V(T)\ar[r]^-{{\Psi}_{\widetilde{T}}^{-1}}\ar[d]_-{I_{T}^{RT^k}}\ar@{}|{}&V(\widetilde{T})\ar[r]^-{{\Psi}_{\widehat{T}}^{-1}}\ar[d]_-{{I}_{\widetilde{T}}^{RT^k}}\ar@{}|{}&V(\widehat{T})\ar[d]^-{{I}_{\widehat{T}}^{RT^k}} \\
P\ar[r]_-{{\Psi}_{\widetilde{T}}^{-1}}&\widetilde{{P}}\ar[r]_-{{\Psi}_{\widehat{T}}^{-1}}&\widehat{{P}}
}
\end{align*}
commute.
\end{prop}

\begin{pf*}
A proof can be found in \cite[Lemma 3.4]{BofBre08}. 
\qed
\end{pf*}

\begin{lem} \label{rtl2=com}
Let $T \in \mathbb{T}_h$. Let $V^{RT}({T}) := W^{1,1}({T})^d$ and $V^{L^2}({T}) := L^1(T)$. 
For $k \in \mathbb{N}_0$, let $I_{{T}}^{RT^k}: V^{RT}(T) \to \mathbb{RT}^k(T)$ and $\Pi^k_{{T}}: V^{L^2}(T) \to {\mathbb{P}^k(T)}$ be the RT interpolation operator and the $L^2$-orthogonal projection, respectively. Then, the following diagram commutes:
\begin{align*}
\displaystyle
\xymatrix{
V^{RT}(T) \ar[r]^-{\nabla \cdot}\ar[d]_-{{I}_{{T}}^{RT^k}}\ar@{}|{}&V^{L^2}(T) \ar[d]^-{{\Pi}_{{T}}^{k}} \\
 \mathbb{RT}^k(T) \ar[r]_-{\nabla \cdot}& \mathbb{P}^k(T)
}
\end{align*}
In other words,
\begin{align}
\displaystyle
\nabla \cdot (I_{{T}}^{RT^k} v) = \Pi_T^k (\nabla \cdot v) \quad \forall v \in V^{RT}(T). \label{rt=div=com}
\end{align}
\end{lem}

\begin{pf*}
A proof can be found in \cite[Lemma 16.2]{ErnGue21a}.
\qed
\end{pf*}

\begin{lem}
Let $T \in \mathbb{T}_h$ and $q \in \mathbb{RT}^k(T)$. Then,
\begin{align}
\displaystyle
\div q &\in \mathbb{P}^k(T), \label{B33} \\
q \cdot n |_{\partial T} &\in \mathbb{R}^k(\partial T).  \label{B34}
\end{align}
\end{lem}

\begin{pf*}
A proof can be found in \cite[Proposition 2.3.3]{BofBreFor13}.
\qed
\end{pf*}

\begin{lem}
The RT finite element with the nodal values in \eqref{RT1112} and \eqref{RT1113} is unisolvent.
\end{lem}

\begin{pf*}
A proof can be found in \cite[Proposition 2.3.4]{BofBreFor13}.
\qed
\end{pf*}

\subsection{Remarks on the Anisotropic RT Interpolation Error Estimate}
We consider the simplex $\widehat{T} \subset \mathbb{R}^2$ with vertices $\widehat{p}_1 := (0,0)^{\top}$, $\widehat{p}_2 := (1,0)^{\top}$ and $\widehat{p}_3 := (0,1)^{\top}$. For $1 \leq i \leq 3$, let $\widehat{F}_i$ be the face of $\widehat{T}$ opposite to $\widehat{p}_i$. The RT interpolation of $\hat{v}$ is defined as
\begin{align*}
\displaystyle
 I_{\widehat{T}}^{RT^0} \hat{v}
  = \sum_{i=1}^3 \left( \int_{\widehat{F}_i} \hat{v} \cdot \hat{n}_i d \hat{s} \right) \hat{\theta}_i \in \mathbb{RT}^0,
\end{align*}
where 
\begin{align*}
\displaystyle
\hat{\theta}_i := \frac{1}{2 |\widehat{T}|} (\hat{x} - \widehat{p}_i), \quad \hat{x} = (\hat{x}_1,\hat{x}_2)^{\top}.
\end{align*}
Setting $\hat{v} := (0, \hat{x}_2^2)^{\top}$  yields 
\begin{align*}
\displaystyle
I_{\widehat{T}}^{RT^0} \hat{v}
&= \frac{1}{\sqrt{2}} \left( \int_{\widehat{F}_1} \hat{x}_2^2 d \hat{s} \right) (\hat{x}_1,\hat{x}_2)^{\top} - \left( \int_{\widehat{F}_3} \hat{x}_2^2 d \hat{s} \right)  (\hat{x}_1,\hat{x}_2 - 1)^{\top} \\
&= \frac{1}{3} (\hat{x}_1,\hat{x}_2)^{\top}.
\end{align*}
This implies that $( I_{\widehat{T}}^{RT^0} \hat{v} )_1 - \hat{v}_1 \not\equiv 0$ for any $\hat{x} \in \mathbb{R}^2$ and the following component-wise stability does not hold:
\begin{align*}
\displaystyle
 \| ( I_{\widehat{T}}^{RT^0}  \hat{v} )_1  \|_{L^2(\widehat{T})} \leq c | \hat{v}_1|_{H^{1}(\widehat{T})}.
\end{align*}
In other words, $( I_{\widehat{T}}^{RT^0}  \hat{v} )_1$ depends on both $\hat{v}_1$ and $\hat{v}_2$. Meanwhile, setting $\hat{v} := (0, \hat{x}_1^2)^{\top}$  yields $ I_{\widehat{T}}^{RT^0} \hat{v} = \frac{1}{3} (0,1)^{\top}$. A key observation is that if $\hat{r} := (0,g(\hat{x}_1))^{\top}$, then $( I_{\widehat{T}}^{RT^0}  \hat{r} )_1 = 0$. In the next section, we introduce component-wise stabilities of the RT interpolation on the reference element by \cite{AcoApe10}.

\subsection{Component-wise Stability of the RT interpolation on the Reference Element}
\textcolor{red}{We will use symbols only used in this subsection.}

\subsubsection{Two-dimensional case}
Let $\widehat{T} \subset \mathbb{R}^2$ be the reference triangle with vertices $\widehat{A}_1 := (1,0)^{\top}$, $\widehat{A}_2 := (0,1)^{\top}$, and $\widehat{A}_3 := (0,0)^{\top}$ with $\hat{N}_1 := (-1,0)^{\top}$, $\hat{N}_2:= (0,-1)^{\top}$, and $\hat{N}_3 := \frac{1}{\sqrt{2}}(1,1)^{\top}$. For $1 \leq i \leq 3$, let $\widehat{E}_i$ be the edge of $\widehat{T}$ opposite to $\widehat{A}_i$. 

We use the same notation for a function of some variable as for its extension to $\widehat{T}$ as a function independent of the other variable. For example, $f(\hat{x}_2)$ denotes a function defined on $\widehat{E}_1$ as well as one is defined in $\widehat{T}$. Furthermore, the same notation is used to denote a polynomial $\hat{p}_k$ on an edge and a polynomial in two variables such that its restriction to that edge agrees with $\hat{q}_k$. For example, for $\hat{q}_k \in \mathbb{P}^k(\widehat{E}_3)$, we write $\hat{q}_k(1-\hat{x}_2,\hat{x}_2)$.

\begin{lem}
Let $\hat{f}_i \in L^p(\widehat{E}_i)$, $i \in \{ 1,2\}$. If
\begin{align*}
\displaystyle
\hat{u}(\hat{x}) = (\hat{f}_1(\hat{x}_2),0)^{\top}, \quad \hat{v}(\hat{x}) = (0 , \hat{f}_2(\hat{x}_1))^{\top},
\end{align*}
then there exist polynomials $\hat{q}_i \in \mathcal{P}^k(\widehat{E}_i)$, $i \in \{ 1,2\}$, such that
\begin{align*}
\displaystyle
I^{RT^k}_{\widehat{T}} \hat{u} = (\hat{q}_1(\hat{x}_2) , 0)^{\top}, \quad I^{RT^k}_{\widehat{T}} \hat{v} = (0 , \hat{q}_2(\hat{x}_1) )^{\top}.
\end{align*}

\end{lem}

\begin{pf*}
A proof is provided in \cite[Lemma 3.2]{AcoApe10} (also see Lemma \ref{lem1133}) for the case $d=3$. The estimate in the case $d=2$ can be proved analogously.
\qed
\end{pf*}

\begin{lem}
For $k \in \mathbb{N}_0$, there exists a constant $c$ such that, for all $\hat{u} = (\hat{u}_1, \hat{u}_2)^{\top} \in W^{1,p}(\widehat{T})^2$, 
\begin{align}
\displaystyle
\| (I_{\widehat{T}}^{RT^k} \hat{u})_i \|_{L^p(\widehat{T})} \leq c \left( \| \hat{u}_i \|_{W^{1,p}(\widehat{T})} + \| \widehat{\div} \hat{u} \|_{L^p(\widehat{T})} \right), \quad i=1, 2. \label{RT1131}
\end{align}
\end{lem}

\begin{pf*}
The proof is provided in \cite[Lemma 3.3]{AcoApe10} (also see Lemma \ref{lem1134}) for the case $d=3$. The estimate in the case $d=2$ can be proved analogously.
\qed
\end{pf*}

\subsubsection{Three-dimensional case: Type \roman{sone}}
Let $\widehat{T} \subset \mathbb{R}^3$ be the reference triangle with vertices $\widehat{A}_1 := (1,0,0)^{\top}$, $\widehat{A}_2 := (0,1,0)^{\top}$, $\widehat{A}_3 := (0,0,1)^{\top}$, and $\widehat{A}_4 := (0,0,0)^{\top}$ with $\widehat{N}_1 := (-1,0,0)^{\top}$, $\widehat{N}_2:= (0,-1,0)^{\top}$, $\widehat{N}_3:= (0,0,-1)^{\top}$, and $\widehat{N}_4 := \frac{1}{\sqrt{3}}(1,1,1)^{\top}$. For $1 \leq i \leq 4$, let $\widehat{E}_i$ be the edge of $\widehat{T}$ opposite to $\widehat{A}_i$. 

In the two-dimensional case, we use the same notation for a function of some variable as for its extension to $\widehat{T}$ as a function independent of the other variable. For example, $f(\hat{x}_2,\hat{x}_3)$ denotes a function define on $\widehat{E}_1$ as well as one is defined in $\widehat{T}$. Furthermore, the same notation is used to denote a polynomial $\hat{p}_k$ on an edge and a polynomial in two variables such that its restriction to that edge agrees with $\hat{p}_k$. For example, for $\hat{p}_k \in \mathbb{P}^k(\widehat{E}_4)$, we write $\hat{p}_k(1-\hat{x}_2 - \hat{x}_3,\hat{x}_2,\hat{x}_3)$.

\begin{lem} \label{lem1133}
Let $k \in \mathbb{N}_0$. Let $\hat{f}_i \in L^p(\widehat{E}_i)$, $i \in \{ 1,2,3\}$. If
\begin{align*}
\displaystyle
\hat{u}(\hat{x}) &= (\hat{f}_1(\hat{x}_2,\hat{x}_3),0,0)^{\top}, \quad \hat{v}(\hat{x}) = (0 , \hat{f}_2(\hat{x}_1,\hat{x}_3),0)^{\top}, \\
\hat{w}(\hat{x}) &= (0, 0 , \hat{f}_3(\hat{x}_1,\hat{x}_2))^{\top},
\end{align*}
then there exist polynomials $\hat{q}_i \in \mathbb{P}^k(\widehat{E}_i)$, $i \in \{ 1,2,3\}$, such that
\begin{align*}
\displaystyle
I^{RT^k}_{\widehat{T}} \hat{u} &= (\hat{q}_1(\hat{x}_2,\hat{x}_3), 0 , 0)^{\top}, \quad I^{RT^k}_{\widehat{T}} \hat{v} = (0 , \hat{q}_2(\hat{x}_1,\hat{x}_3),0 )^{\top}, \\
 I^{RT^k}_{\widehat{T}} \hat{w} &= (0 , 0, \hat{q}_3(\hat{x}_1,\hat{x}_2))^{\top}.
\end{align*}
\end{lem}

\begin{pf*}
We follow \cite[Lemma 3.2]{AcoApe10}.  Because $\widehat{\div} \hat{u} = 0$, from the definition of the RT interpolation and the Green's formula, we have, for any $\hat{r}_k \in \mathbb{P}^k(\widehat{T})$,
\begin{align*}
\displaystyle
0 &= \int_{\widehat{T}} \hat{r}_k \widehat{\div} \hat{u} d \hat{x} \\
&= \sum_{i=1}^4 \int_{\widehat{E}_i} (\hat{r}_k \widehat{N}_{i}) \cdot \hat{u} d \hat{s} - \int_{\widehat{T}} (\hat{u} \cdot \widehat{\nabla} ) \hat{r}_k d \hat{x} \\
&= \sum_{i=1}^4 \int_{\widehat{E}_i} (\hat{r}_k \widehat{N}_{i}) \cdot (I_{\widehat{T}}^{RT^k} \hat{u}) d \hat{s} - \int_{\widehat{T}} ((I_{\widehat{T}}^{RT^k} \hat{u}) \cdot \widehat{\nabla} ) \hat{r}_k d \hat{x} \\
&= \int_{\widehat{T}} \hat{r}_k \widehat{\div} (I_{\widehat{T}}^{RT^k} \hat{u}) d \hat{x},
\end{align*}
which leads to $\widehat{\div} (I_{\widehat{T}}^{RT^k} \hat{u}) = 0$. Therefore, form the property of the RT interpolation, $I_{\widehat{T}}^{RT^k} \hat{u} \in \mathbb{P}^k(\widehat{T})^3$, e.g. see \cite[Lemma 3.1]{BofBre08}. 

Using \eqref{RT1115} for $i=2,3$, and $\hat{u}_2 = \hat{u}_3 = 0$, we have
\begin{align*}
\displaystyle
\int_{\widehat{E}_i} (I_{\widehat{T}}^{RT^k} \hat{u})_{i} \hat{r}_k d \hat{s} = 0 \quad \forall \hat{r}_k \in \mathbb{P}^k(\widehat{E}_i), \quad i=2,3.
\end{align*}
Setting $\hat{r}_k :=  (I_{\widehat{T}}^{RT^k} \hat{u})_{i}$, we obtain that $ (I_{\widehat{T}}^{RT^k} \hat{u})_{i} |_{\widehat{E}_i} = 0$ for $i=2,3$. 

For $k=0$, because  $I_{\widehat{T}}^{RT^0} \hat{u} \in \mathbb{P}^0(\widehat{T})^3$ and $ (I_{\widehat{T}}^{RT^0} \hat{u})_{i} |_{\widehat{E}_i} = 0$ for $i=2,3$, it holds that $(I_{\widehat{T}}^{RT^0} \hat{u})_i = 0$ in $\widehat{T}$ for $i=2,3$. This implies that the first result holds.

For $k \geq 1$, there exists a polynomial $\hat{r}_{i} \in \mathbb{P}^{k-1}(\widehat{T})$, $i=2,3$, such that $(I_{\widehat{T}}^{RT^k} \hat{u})_{i} = \hat{x}_{i} \hat{r}_{i}$. Using \eqref{RT1116} for $i=2,3$, and $\hat{u}_2 = \hat{u}_3 = 0$, we have, for $i=2,3$,
\begin{align*}
\displaystyle
\int_{\widehat{T}} (I_{\widehat{T}}^{RT^k} \hat{u})_i \hat{r}_i d \hat{x} = 0. \quad \text{as $\hat{q}_{k-1} := (0,\hat{r}_i,0)^{\top}$ in \eqref{RT1116}},
\end{align*}
which leads to
\begin{align*}
\displaystyle
\int_{\widehat{T}} (I_{\widehat{T}}^{RT^k} \hat{u})_i^2 d \hat{x} 
&= \int_{\widehat{T}} \hat{x}_{i} \hat{x}_{i} \hat{r}_i^2 d \hat{x}
\leq \|  \hat{x}_{i} \|_{L^{\infty}(\widehat{T})}  \int_{\widehat{T}}  \hat{x}_{i} \hat{r}_i^2 d \hat{x} = 0.
\end{align*}
Note that $\hat{x}_i \geq 0$ in $\widehat{T}$ for $i=2,3$. We hence conclude that $(I_{\widehat{T}}^{RT^k} \hat{u})_i = 0$ in $\widehat{T}$ for $i=2,3$. 

Because $\widehat{\div} (I_{\widehat{T}}^{RT} \hat{u}) = 0$, it follows that 
\begin{align*}
\displaystyle
\frac{\partial (I_{\widehat{T}}^{RT^k} \hat{u})_1}{\partial \hat{x}_1} = 0.
\end{align*}
This means that $(I_{\widehat{T}}^{RT^k} \hat{u})_1$ is independent of $\hat{x}_1$.

The other two results are analogous.
\qed
\end{pf*}

\begin{lem}  \label{lem1134}
For $k \in \mathbb{N}_0$, there exists a constant $c$ such that, for all $\hat{u} = (\hat{u}_1, \hat{u}_2, \hat{u}_3)^{\top} \in W^{1,p}(\widehat{T})^3$, 
\begin{align}
\displaystyle
\| (I_{\widehat{T}}^{RT^k} \hat{u})_i \|_{L^p(\widehat{T})} \leq c \left( \| \hat{u}_i \|_{W^{1,p}(\widehat{T})} + \| \widehat{\div} \hat{u} \|_{L^p(\widehat{T})} \right), \quad i=1, 2, 3. \label{RT1132}
\end{align}

\end{lem}

\begin{pf*}
Only shown when $k=0$. When $k \geq 1$. see \cite[Lemma 3.3]{AcoApe10}.

From Lemma \ref{lem1133}, if
\begin{align*}
\displaystyle
\hat{v} := (\hat{u}_1 , \hat{u}_2 - \hat{u}_2(\hat{x}_1,0,\hat{x}_3) ,   \hat{u}_3 - \hat{u}_3(\hat{x}_1,\hat{x}_2,0))^{\top},
\end{align*}
it holds that
\begin{align*}
\displaystyle
I_{\widehat{T}}^{RT^k} \hat{v}
&= I_{\widehat{T}}^{RT^k} \hat{u} - I_{\widehat{T}}^{RT^k}  (0,  \hat{u}_2(\hat{x}_1,0,\hat{x}_3) ,  0 )^{\top} - I_{\widehat{T}}^{RT^k}  (0,0 ,  \hat{u}_3(\hat{x}_1,\hat{x}_2,0))^{\top},
\end{align*}
and thus, $(I_{\widehat{T}}^{RT^k} \hat{v})_1 = (I_{\widehat{T}}^{RT^k} \hat{u})_1$.

Let $k=0$. Because $\hat{v}_2|_{\widehat{E}_2} = 0$ and $\hat{v}_3|_{\widehat{E}_3} = 0$, $I_{\widehat{T}}^{RT^0} \hat{v}$ is determined by the equations
\begin{subequations} \label{RT1133}
\begin{align}
\displaystyle
\int_{\widehat{E}_1} (I_{\widehat{T}}^{RT^0} \hat{v})_1 d \hat{s} &= \int_{\widehat{E}_1} \hat{v}_1 d \hat{s}, \label{RT1133a}\\
\int_{\widehat{E}_2} (I_{\widehat{T}}^{RT^0} \hat{v})_2 d \hat{s} &=0, \label{RT1133b}\\
\int_{\widehat{E}_3} (I_{\widehat{T}}^{RT^0} \hat{v})_3 d \hat{s} &=0, \label{RT1133c}\\
\int_{\widehat{E}_4} \{ (I_{\widehat{T}}^{RT^0} \hat{v})_1 + (I_{\widehat{T}}^{RT^0} \hat{v})_2 + (I_{\widehat{T}}^{RT^0} \hat{v})_3 \} d \hat{s} &= \int_{\widehat{E}_4} ( \hat{v}_1 + \hat{v}_2 + \hat{v}_3 ) d \hat{s}. \label{RT1133d}
\end{align}
\end{subequations}
From the divergence formula and the definition of $\hat{v}$, we have
\begin{align}
\displaystyle
\int_{\widehat{T}} \widehat{\div} \hat{v} d \hat{x}
&= \int_{\partial \widehat{T}} \hat{v} \cdot \hat{n} d \hat{s}
= \frac{1}{\sqrt{3}} \int_{\widehat{E}_4} (\hat{v}_1 + \hat{v}_2 + \hat{v}_3) d \hat{s} +  \int_{\partial \widehat{T} \backslash \widehat{E}_4} \hat{v} \cdot \hat{n} d \hat{s} \notag \\
&=  \frac{1}{\sqrt{3}} \int_{\widehat{E}_4} (\hat{v}_1 + \hat{v}_2 + \hat{v}_3) d \hat{s} + \int_{\widehat{E}_1} \hat{v}_1 d \hat{s}.  \label{RT1134}
\end{align}
Because $\hat{u}_1 = \hat{v}_1$, $\widehat{\div} \hat{u} = \widehat{\div} \hat{v}$, $(I_{\widehat{T}}^{RT^0} \hat{u})_1 = (I_{\widehat{T}}^{RT^0} \hat{v})_1$, \eqref{RT1133}, \eqref{RT1134}, the definition of the Raviart--Thomas interpolation, and the trace theorem, we have
\begin{align*}
\displaystyle
\| (I_{\widehat{T}}^{RT^0} \hat{u})_1 \|_{L^p(\widehat{T})}
&= \| (I_{\widehat{T}}^{RT^0} \hat{v})_1 \|_{L^p(\widehat{T})} \\
&\leq \sum_{i=1}^4 \left| \int_{\widehat{E}_i}  \hat{v} \cdot \widehat{N}_i  d \hat{s} \right| \| (\hat{\theta}_i)_1 \|_{L^p(\widehat{T})}  \\
&\leq c \left| \int_{\widehat{E}_1} \hat{v}_1 d \hat{s}  +  \frac{1}{\sqrt{3}} \int_{\widehat{E}_4} ( \hat{v}_1 + \hat{v}_2 + \hat{v}_3 ) d \hat{s} \right| \\
&\leq c \left( \| \hat{u}_1 \|_{W^{1,p}(\widehat{T})} + \| \widehat{\div} \hat{u} \|_{L^p(\widehat{T})} \right),
\end{align*}
which is the desired result for $k=0$. By analogous argument, the estimates for $ (I_{\widehat{T}}^{RT^0} \hat{u})_i$, $i=2,3$, can be proved.

\qed
\end{pf*}

\subsubsection{Three-dimensional case: Type \roman{stwo}}
Let $\widehat{T} \subset \mathbb{R}^3$ be the reference triangle with vertices $\widehat{A}_1 := (1,0,0)^{\top}$, $\widehat{A}_2 := (1,1,0)^{\top}$, $\widehat{A}_3 := (0,0,1)^{\top}$, and $\widehat{A}_4 := (0,0,0)^{\top}$ with $\widehat{N}_1 := \frac{1}{\sqrt{2}} (-1,1,0)^{\top}$, $\widehat{N}_2:= (0,-1,0)^{\top}$, $\widehat{N}_3:= (0,0,-1)^{\top}$, and $\widehat{N}_4 := \frac{1}{\sqrt{2}}(1,0,1)^{\top}$. For $1 \leq i \leq 4$, let $\widehat{E}_i$ be the edge of $\widehat{T}$ opposite to $\widehat{A}_i$ and with $\overline{E}_1$ the projection of $\widehat{E}_1$ onto the plane given by $\hat{x}_1 = 0$. 

\begin{lem} \label{lem1135}
Let $k \in \mathbb{N}_0$. Let $\hat{f}_1 \in L^p(\overline{E}_1)$, and $\hat{f}_i \in L^p(\widehat{E}_i)$, $i \in \{ 2,3\}$. If
\begin{align*}
\displaystyle
\hat{u}(\hat{x}) &= (\hat{f}_1(\hat{x}_2,\hat{x}_3),0,0)^{\top}, \quad \hat{v}(\hat{x}) = (0 , \hat{f}_2(\hat{x}_1,\hat{x}_3),0)^{\top}, \\
\hat{w}(\hat{x}) &= (0, 0 , \hat{f}_3(\hat{x}_1,\hat{x}_2))^{\top},
\end{align*}
then there exist polynomials  $\hat{q}_1 \in \mathbb{P}^k(\overline{E}_1)$, and $\hat{q}_i \in \mathbb{P}^k(\widehat{E}_i)$, $i \in \{ 2,3\}$, such that
\begin{align*}
\displaystyle
I^{RT^k}_{\widehat{T}} \hat{u} &= (\hat{q}_1(\hat{x}_2,\hat{x}_3), 0 , 0)^{\top}, \quad I^{RT^k}_{\widehat{T}} \hat{v} = (0 , \hat{q}_2(\hat{x}_1,\hat{x}_3),0 )^{\top}, \\
 I^{RT^k}_{\widehat{T}} \hat{w} &= (0 , 0, \hat{q}_3(\hat{x}_1,\hat{x}_2))^{\top}.
\end{align*}
\end{lem}

\begin{pf*}
We follow \cite[Lemma 4.2]{AcoApe10}. The proof is similar to that of Lemma \ref{lem1133}. We prove the second equality. The other two follow in an analogous argument.

Because $\widehat{\div} \hat{v} = 0$, from the definition of the RT interpolation and the Green's formula, we have $\widehat{\div} (I_{\widehat{T}}^{RT^k} \hat{v}) = 0$. Therefore, form the property of the RT interpolation, $I_{\widehat{T}}^{RT^k} \hat{v} \in \mathbb{P}^k(\widehat{T})^3$. Using \eqref{RT1115} for $i=3$, and $\hat{v}_3 = 0$, we have
\begin{align*}
\displaystyle
\int_{\widehat{E}_3} (I_{\widehat{T}}^{RT^k} \hat{v})_{3} \hat{p}_k d \hat{s} = 0 \quad \forall \hat{p}_k \in \mathbb{P}^k(\widehat{E}_3).
\end{align*}
Setting $\hat{p}_k :=  (I_{\widehat{T}}^{RT^k} \hat{v})_{3}$, we obtain that $ (I_{\widehat{T}}^{RT^k} \hat{v})_{3} |_{\widehat{E}_3} = 0$.

Let $k=0$. Because  $I_{\widehat{T}}^{RT^0} \hat{v} \in \mathbb{P}^0(\widehat{T})^3$ and $ (I_{\widehat{T}}^{RT^0} \hat{v})_{3} |_{\widehat{E}_3} = 0$, it holds that $(I_{\widehat{T}}^{RT^0} \hat{v})_3 = 0$ in $\widehat{T}$. Using \eqref{RT1115} for $i=4$, and $\hat{v}_1 = \hat{v}_3 = 0$, we have
\begin{align*}
\displaystyle
\int_{\widehat{E}_4} \{  (I_{\widehat{T}}^{RT^0} \hat{v})_{1} +  (I_{\widehat{T}}^{RT^0} \hat{v})_{3} \}d \hat{s} = 0,
\end{align*}
which leads to $ (I_{\widehat{T}}^{RT^0} \hat{v})_{1}|_{\widehat{E}_4} = 0$. It then holds that that $(I_{\widehat{T}}^{RT^0} \hat{v})_1 = 0$ in $\widehat{T}$. This implies that the second result holds.

Let $k \geq 1$. As in the proof of Lemma \ref{lem1133}, we obtain that $(I_{\widehat{T}}^{RT^k} \hat{v})_3 = 0$ in $\widehat{T}$.  Using \eqref{RT1115} for $i=4$, and $\hat{v}_1 = \hat{v}_3 = 0$, we have
\begin{align*}
\displaystyle
\int_{\widehat{E}_4} \{  (I_{\widehat{T}}^{RT^k} \hat{v})_{1} +  (I_{\widehat{T}}^{RT^k} \hat{v})_{3} \} \hat{p}_k d \hat{s} = 0 \quad \forall \hat{p}_k \in \mathbb{P}^k(\widehat{E}_4),
\end{align*}
which implies that $\{  (I_{\widehat{T}}^{RT^k} \hat{v})_{1} +  (I_{\widehat{T}}^{RT^k} \hat{v})_{3} \} |_{\widehat{E}_4} = 0$, and hence
\begin{align*}
\displaystyle
(I_{\widehat{T}}^{RT^k} \hat{v})_{1} +  (I_{\widehat{T}}^{RT^k} \hat{v})_{3} = (1 - \hat{x}_1 - \hat{x}_3) \hat{r}
\end{align*}
for some $\hat{r} \in \mathbb{P}^{k-1}(\widehat{T})$. Using \eqref{RT1116} and $\hat{v}_1 = \hat{v}_3 = 0$, we have
\begin{align*}
\displaystyle
\int_{\widehat{T}}  \{  (I_{\widehat{T}}^{RT^k} \hat{v})_{1} +  (I_{\widehat{T}}^{RT^k} \hat{v})_{3} \} \hat{r} d \hat{x} = 0. \quad \text{as $\hat{q}_{k-1} := (\hat{r},0,\hat{r})^{\top}$ in \eqref{RT1116}},
\end{align*}
which leads to
\begin{align*}
\displaystyle
\int_{\widehat{T}}  \{ (I_{\widehat{T}}^{RT^k} \hat{v})_{1} +  (I_{\widehat{T}}^{RT^k} \hat{v})_{3} \}^2  d \hat{x}
&= \int_{\widehat{T}} (1 - \hat{x}_1 - \hat{x}_3)^2 \hat{r}^2 d \hat{x} \\
&\leq \| 1 - \hat{x}_1 - \hat{x}_3 \|_{L^{\infty}(\widehat{T})} \int_{\widehat{T}} (1 - \hat{x}_1 - \hat{x}_3) \hat{r}^2 d \hat{x} = 0.
\end{align*}
Note that $1 - \hat{x}_1 - \hat{x}_3 \geq 0$ in $\widehat{T}$. We hence have $(I_{\widehat{T}}^{RT^k} \hat{v})_{1} +  (I_{\widehat{T}}^{RT^k} \hat{v})_{3} = 0$ in $\widehat{T}$. Because we know $(I_{\widehat{T}}^{RT^k} \hat{v})_3 = 0$ in $\widehat{T}$, we conclude that $(I_{\widehat{T}}^{RT^k} \hat{v})_1 = 0$ in $\widehat{T}$.

Because $\widehat{\div} (I_{\widehat{T}}^{RT^k} \hat{v}) = 0$, it follows that 
\begin{align*}
\displaystyle
\frac{\partial (I_{\widehat{T}}^{RT^k} \hat{v})_2}{\partial \hat{x}_2} = 0.
\end{align*}
This means that $(I_{\widehat{T}}^{RT^k} \hat{v})_2$ is independent of $\hat{x}_2$.
\qed
\end{pf*}

\begin{lem} \label{lem1136}
For $k \in \mathbb{N}_0$, there exists a constant $c$ such that, for all $\hat{u} = (\hat{u}_1,\hat{u}_2, \hat{u}_3)^{\top} \in W^{1,p}(\widehat{T})^3$, 
\begin{subequations} \label{RT11311}
\begin{align}
\displaystyle
\| (I_{\widehat{T}}^{R^kT} \hat{u})_i \|_{L^p(\widehat{T})} 
&\leq c \left( \| \hat{u}_1 \|_{W^{1,p}(\widehat{T})} + \left\| \frac{\partial \hat{u}_2}{\partial \hat{x}_2} \right\|_{L^p(\widehat{T})} + \left\| \frac{\partial \hat{u}_3}{\partial \hat{x}_3} \right\|_{L^p(\widehat{T})}  \right), \label{RT11311a} \\
\| (I_{\widehat{T}}^{RT^k} \hat{u})_i \|_{L^p(\widehat{T})} 
&\leq c \left( \| \hat{u}_i \|_{W^{1,p}(\widehat{T})} + \left\| \div \hat{u} \right\|_{L^p(\widehat{T})} \right), \quad i=2,3. \label{RT11311b}
\end{align}
\end{subequations}
In particular,
\begin{align}
\displaystyle
\| (I_{\widehat{T}}^{RT^k} \hat{u})_i \|_{L^p(\widehat{T})} 
&\leq c \left( \| \hat{u}_i \|_{W^{1,p}(\widehat{T})} + \sum_{j=1, j\neq i}^3 \left\| \frac{\partial \hat{u}_j}{\partial \hat{x}_j} \right\|_{L^p(\widehat{T})} \right), \quad i=1,2,3. \label{RT11312}
\end{align}
\end{lem}

\begin{pf*}
Only shown when $k=0$. When $k \geq 1$. see \cite[Lemma 4.3]{AcoApe10}. We prove the estimates  \eqref{RT11311a} and  \eqref{RT11311b} with $i=2$. The other one follows in an analogous argument.

\textbf{Case for \eqref{RT11311a}.} 
From Lemma \ref{lem1135}, if
\begin{align*}
\displaystyle
\hat{v} := (\hat{u}_1 , \hat{u}_2 - \hat{u}_2(\hat{x}_1,0,\hat{x}_3) ,   \hat{u}_3 - \hat{u}_3(\hat{x}_1,\hat{x}_2,0))^{\top},
\end{align*}
it holds that
\begin{align*}
\displaystyle
I_{\widehat{T}}^{RT^k} \hat{v}
&= I_{\widehat{T}}^{RT^k} \hat{u} - I_{\widehat{T}}^{RT^k}  (0,  \hat{u}_2(\hat{x}_1,0,\hat{x}_3) ,  0 )^{\top} - I_{\widehat{T}}^{RT^k}  (0,0 ,  \hat{u}_3(\hat{x}_1,\hat{x}_2,0))^{\top},
\end{align*}
and thus, $(I_{\widehat{T}}^{RT^k} \hat{v})_1 = (I_{\widehat{T}}^{RT^k} \hat{u})_1$.

Let $k=0$. Because $\hat{v}_2|_{\widehat{E}_2} = 0$ and $\hat{v}_3|_{\widehat{E}_3} = 0$, $I_{\widehat{T}}^{RT^0} \hat{v}$ is determined by the equations
\begin{subequations} \label{RT11313}
\begin{align}
\displaystyle
\int_{\widehat{E}_1} \{ - (I_{\widehat{T}}^{RT^0} \hat{v})_1 + (I_{\widehat{T}}^{RT^0} \hat{v})_2 \} d \hat{s} &= \int_{\widehat{E}_1} ( - \hat{v}_1 + \hat{v}_2 ) d \hat{s}, \label{RT11313a}\\
\int_{\widehat{E}_2} (I_{\widehat{T}}^{RT^0} \hat{v})_2 d \hat{s} &=0, \label{RT11313b}\\
\int_{\widehat{E}_3} (I_{\widehat{T}}^{RT^0} \hat{v})_3 d \hat{s} &=0, \label{RT11313c}\\
\int_{\widehat{E}_4} \{ (I_{\widehat{T}}^{RT^0} \hat{v})_1  + (I_{\widehat{T}}^{RT^0} \hat{v})_3 \} d \hat{s} &= \int_{\widehat{E}_4} ( \hat{v}_1 + \hat{v}_3 ) d \hat{s}. \label{RT11313d}
\end{align}
\end{subequations}
From the divergence formula and the definition of $\hat{v}$, we have
\begin{align}
\displaystyle
\int_{\widehat{T}} \widehat{\div} (\hat{v}_1,\hat{v}_2,0)^{\top} d \hat{x}
&= \int_{\partial \widehat{T}}  (\hat{v}_1,\hat{v}_2,0)^{\top}  \cdot \hat{n} d \hat{s} \notag\\
&= \frac{1}{\sqrt{2}} \int_{\widehat{E}_1} (- \hat{v}_1 + \hat{v}_2 ) d \hat{s} +  \frac{1}{\sqrt{2}} \int_{ \widehat{E}_4} \hat{v}_1 d \hat{s}, \label{RT11314}
\end{align}
and
\begin{align}
\displaystyle
\int_{\widehat{T}} \widehat{\div} (\hat{v}_1, 0 , \hat{v}_3)^{\top} d \hat{x}
&= \int_{\partial \widehat{T}}  (\hat{v}_1, 0 , \hat{v}_3)^{\top}  \cdot \hat{n} d \hat{s} \notag\\
&= \frac{1}{\sqrt{2}} \int_{\widehat{E}_1} (- \hat{v}_1 ) d \hat{s} +  \frac{1}{\sqrt{2}} \int_{ \widehat{E}_4} ( \hat{v}_1 + \hat{v}_3 ) d \hat{s}, \label{RT11315}
\end{align}
Because $\hat{u}_1 = \hat{v}_1$, $\frac{\partial \hat{u}_j}{\partial \hat{x}_j} = \frac{\partial \hat{v}_j}{\partial \hat{x}_j}$ for $j=2,3$, $(I_{\widehat{T}}^{RT^0} \hat{u})_1 = (I_{\widehat{T}}^{RT^0} \hat{v})_1$, \eqref{RT11313}, \eqref{RT11314}, \eqref{RT11315}, the definition of the RT interpolation, and the trace theorem, we have
\begin{align*}
\displaystyle
\| (I_{\widehat{T}}^{RT^0} \hat{u})_1 \|_{L^p(\widehat{T})}
&= \| (I_{\widehat{T}}^{RT^0} \hat{v})_1 \|_{L^p(\widehat{T})} \\
&\leq \sum_{j=1}^4 \left| \int_{\widehat{E}_j} \hat{v} \cdot \widehat{N}_j  d \hat{s} \right| \| (\hat{\theta}_j)_1 \|_{L^p(\widehat{T})}  \\
&\leq c \left| \frac{1}{\sqrt{2}} \int_{\widehat{E}_1} ( - \hat{v}_1 + \hat{v}_2 ) d \hat{s} +  \frac{1}{\sqrt{2}} \int_{\widehat{E}_4} ( \hat{v}_1 + \hat{v}_3 ) d \hat{s} \right| \\
&\leq c \left( \| \hat{u}_1 \|_{W^{1,p}(\widehat{T})} + \left\| \frac{\partial \hat{u}_2}{\partial \hat{x}_2} \right\|_{L^p(\widehat{T})} + \left\| \frac{\partial \hat{u}_3}{\partial \hat{x}_3} \right\|_{L^p(\widehat{T})} \right),
\end{align*}
which is the desired result for $k=0$.

\textbf{Case for \eqref{RT11311b} with $i=2$.} 
From Lemma \ref{lem1135}, if
\begin{align*}
\displaystyle
\hat{v} := (\hat{u}_1  - \hat{u}_1(\hat{x}_2,\hat{x}_2,\hat{x}_3)  , \hat{u}_2,   \hat{u}_3 - \hat{u}_3(\hat{x}_1,\hat{x}_2,0))^{\top},
\end{align*}
it holds that
\begin{align*}
\displaystyle
I_{\widehat{T}}^{RT^k} \hat{v}
&= I_{\widehat{T}}^{RT^k} \hat{u} - I_{\widehat{T}}^{RT^k}  ( \hat{u}_1(\hat{x}_2,\hat{x}_2,\hat{x}_3), 0  ,  0 )^{\top} - I_{\widehat{T}}^{RT^k}  (0,0 ,  \hat{u}_3(\hat{x}_1,\hat{x}_2,0))^{\top},
\end{align*}
and thus, $(I_{\widehat{T}}^{RT^k} \hat{v})_2 = (I_{\widehat{T}}^{RT^k} \hat{u})_2$.

Let $k=0$. Because $\hat{v}_1|_{\widehat{E}_1} = 0$ and $\hat{v}_3|_{\widehat{E}_3} = 0$, $I_{\widehat{T}}^{RT^0} \hat{v}$ is determined by the equations
\begin{subequations} \label{RT11323}
\begin{align}
\displaystyle
\int_{\widehat{E}_1} \{ - (I_{\widehat{T}}^{RT^0} \hat{v})_1 + (I_{\widehat{T}}^{RT^0} \hat{v})_2 \} d \hat{s} &= \int_{\widehat{E}_1} \hat{v}_2  d \hat{s}, \label{RT11323a}\\
\int_{\widehat{E}_2} (I_{\widehat{T}}^{RT^0} \hat{v})_2 d \hat{s} &= \int_{\widehat{E}_2} \hat{v}_2  d \hat{s}, \label{RT11323b}\\
\int_{\widehat{E}_3} (I_{\widehat{T}}^{RT^0} \hat{v})_3 d \hat{s} &=0, \label{RT11323c}\\
\int_{\widehat{E}_4} \{ (I_{\widehat{T}}^{RT^0} \hat{v})_1  + (I_{\widehat{T}}^{RT^0} \hat{v})_3 \} d \hat{s} &= \int_{\widehat{E}_4} ( \hat{v}_1 + \hat{v}_3 ) d \hat{s}. \label{RT11323d}
\end{align}
\end{subequations}
From the divergence formula and the definition of $\hat{v}$, we have
\begin{align}
\displaystyle
\int_{\widehat{T}} \widehat{\div} (\hat{v}_1, \hat{v}_2 , \hat{v}_3)^{\top} d \hat{x}
&= \int_{\partial \widehat{T}}  (\hat{v}_1, \hat{v}_2 , \hat{v}_3)^{\top}  \cdot \hat{n} d \hat{s} \notag \\
&= \frac{1}{\sqrt{2}} \int_{\widehat{E}_1} \hat{v}_2 d \hat{s} - \int_{\widehat{E}_2} \hat{v}_2 d \hat{s} + \frac{1}{\sqrt{2}} \int_{ \widehat{E}_4} ( \hat{v}_1 + \hat{v}_3 ) d \hat{s}, \label{RT11324}
\end{align}
Because $\hat{u}_2 = \hat{v}_2$, $\frac{\partial \hat{u}_j}{\partial \hat{x}_j} = \frac{\partial \hat{v}_j}{\partial \hat{x}_j}$ for $j=1,3$, $(I_{\widehat{T}}^{RT} \hat{u})_2 = (I_{\widehat{T}}^{RT} \hat{v})_2$, \eqref{RT11323}, \eqref{RT11324}, the definition of the RT interpolation, and the trace theorem, we have
\begin{align*}
\displaystyle
\| (I_{\widehat{T}}^{RT^0} \hat{u})_2 \|_{L^p(\widehat{T})}
&= \| (I_{\widehat{T}}^{RT^0} \hat{v})_2 \|_{L^p(\widehat{T})} \\
&\leq \sum_{j=1}^4 \left| \int_{\widehat{E}_j} \hat{v} \cdot \widehat{N}_j  d \hat{s} \right| \| (\hat{\theta}_j)_2 \|_{L^p(\widehat{T})}  \\
&\leq c \left| \frac{1}{\sqrt{2}} \int_{\widehat{E}_1} \hat{v}_2 d \hat{s} +  \int_{\widehat{E}_2} \hat{v}_2 d \hat{s} +  \frac{1}{\sqrt{2}} \int_{\widehat{E}_4} ( \hat{v}_1 + \hat{v}_3 ) d \hat{s} \right| \\
&\leq c \left( \| \hat{u}_2 \|_{W^{1,p}(\widehat{T})} + \left\|  \widehat{\div} \hat{u} \right\|_{L^p(\widehat{T})} \right),
\end{align*}
which is the desired result for $k=0$.

\qed
\end{pf*}

\subsection{Stability of the local RT interpolation}
The following two lemmata are divided into the element on $\mathfrak{T}^{(2)}$ or $\mathfrak{T}_1^{(3)}$ and the element on $\mathfrak{T}_2^{(3)}$. 

\begin{lem} \label{lem1151}
Let $p \in [1,\infty)$. Let $T\in \mathbb{T}_h$ satisfy Condition \ref{cond1} or Condition \ref{cond2} with $T ={\Phi}_{T} (\widetilde{T})$ and $\widetilde{T} = {\Phi}_{\widetilde{T}}(\widehat{T})$, where $\widetilde{T} \in \mathfrak{T}^{(2)}$ or $\widetilde{T} \in \mathfrak{T}_1^{(3)}$. Then, for any $\hat{v} \in W^{1,p}(\widehat{T})^d$ with $\tilde{v}= {\Psi}_{\widehat{T}} \hat{v}$ and ${v} = {\Psi}_{\widetilde{T}} \tilde{v}$,  
\begin{align}
\displaystyle
\| {I_{T}^{RT^k} v} \|_{L^p({T})^d} 
&\leq c \left[ \frac{H_{T}}{h_{T}} \left( \| v \|_{L^p(T)^d} + \sum_{|\varepsilon|=1} h^{\varepsilon} \left \| \partial_r^{\varepsilon} v \right \|_{L^p(T)^d} \right ) + h_{T} \| \nabla \cdot {v} \|_{L^p({T})} \right]. \label{RT51}
\end{align}
\end{lem}

\begin{pf*}
From \eqref{RT41},
\begin{align}
\displaystyle
\| {I_{T}^{RT^k} v} \|_{L^p({T})^d}
&\leq c |\det({A}_{\widetilde{T}})|^{- \frac{p - 1}{p}} \| \widetilde{{A}} \|_{2} \left(  \sum_{j=1}^d h_j^p \| (I_{\widehat{T}}^{RT^k} \hat{v})_j \|_{L^p(\widehat{T})}^p \right)^{1/p}. \label{RT53}
\end{align}
The component-wise stability \eqref{RT1131} for $2d$ or \eqref{RT1132} for $3d$ yields
\begin{align}
\displaystyle
\sum_{j=1}^d h_j^p \| (I_{\widehat{T}}^{RT^k} \hat{v})_j \|_{L^p(\widehat{T})}^p
&\leq c \sum_{j=1}^d h_j^p \left( \| \hat{v}_j \|_{W^{1,p}(\widehat{T})}^p + \| \nabla_{\hat{x}} \cdot \hat{v}  \|_{L^p(\widehat{T})}^p \right). \label{RT54}
\end{align}
From \eqref{RT42} with $\ell=0$ and $m \in \{0,1 \}$,
\begin{align}
\displaystyle
\| \hat{v}_j \|_{W^{1,p}(\widehat{T})}^p 
&= \| \hat{v}_j \|_{L^{p}(\widehat{T})}^p + \sum_{k=1}^d \left \| \frac{\partial \hat{v}_j}{\partial \hat{x}_k}  \right \|^p_{L^p(\widehat{T})} \notag \\
&\leq c |\det({A}_{\widetilde{T}})|^{p-1} \|  \widetilde{{A}}^{-1} \|_2^p h_j^{-p} \left[ \| v \|^p_{L^p(T)^d} + \left(  \sum_{ |\varepsilon| =1} h^{\varepsilon} \left \| \partial_r^{\varepsilon} v \right \|_{L^p(T)^d} \right)^p \right]. \label{RT55}
\end{align}
From \eqref{RT12} with $\ell=0$,
\begin{align}
\displaystyle
\left \|  \nabla_{\hat{x}} \cdot \hat{v}  \right \|_{L^p(\widehat{T}_1)}
\leq c  |\det ({A}_{\widetilde{T}})|^{\frac{p-1}{p}} \left \|  \nabla \cdot v \right \|_{L^p(T)}. \label{RT56}
\end{align}
Combining the above inequalities \eqref{RT53}, \eqref{RT54}, \eqref{RT55}, and \eqref{RT56} with \eqref{CN331b} and \eqref{jensen} yields
\begin{align*}
\displaystyle
\| {I_{T}^{RT^k} v} \|_{L^p({T})^d}
&\leq c \left[ \frac{H_{T}}{h_{T}} \left( \| v\|_{L^p(T)^d} + \sum_{|\varepsilon|=1} h^{\varepsilon} \left \| \partial_r^{\varepsilon} v \right \|_{L^p(T)^d} \right ) + h_{T} \| \nabla \cdot {v} \|_{L^p({T})} \right],
\end{align*}
which is the desired estimate.
\qed
\end{pf*}

\begin{lem} \label{lem1152}
Let $p \in [1,\infty)$ and $d=3$. Let $T\in \mathbb{T}_h$ satisfy Condition \ref{cond2} with $T ={\Phi}_{T} (\widetilde{T})$ and $\widetilde{T} = {\Phi}_{\widetilde{T}}(\widehat{T})$, where $\widetilde{T} \in \mathfrak{T}_2^{(3)}$. Then,  for any $\hat{v} \in W^{\ell+1,p}(\widehat{T})^3$ with $\tilde{v}= {\Psi}_{\widehat{T}} \hat{v}$ and ${v} = {\Psi}_{\widetilde{T}} \tilde{v}$,  
\begin{align}
\displaystyle
\| {I_{T}^{RT^k} v} \|_{L^p({T})^3} 
&\leq c \frac{H_{T}}{h_{T}} \left[  \| v \|_{L^p(T)^3}  + h_T \sum_{k=1}^3 \left \| \frac{\partial v}{\partial r_k} \right \|_{L^p(T)^3} \right]. \label{RT58}
\end{align}

\end{lem}

\begin{pf*}
The component-wise stability \eqref{RT11311} yields
\begin{align}
\displaystyle
\sum_{j=1}^3 h_j^p \| (I_{\widehat{T}}^{RT^k} \hat{v})_j \|_{L^p(\widehat{T})}^p
&\leq c \sum_{j=1}^3 h_j^p \left( \| \hat{v}_j \|_{W^{1,p}(\widehat{T})}^p +\sum_{k=1, k \neq j}^3 \left\| \frac{\partial \hat{v}_k}{\partial \hat{x}_k} \right\|_{L^p(\widehat{T})}^p  \right). \label{RT59}
\end{align}
From \eqref{RT14} with $\ell=0$,
\begin{align}
\displaystyle
\left \|  \frac{\partial \hat{v}_k}{\partial \hat{x}_k} \right\|_{L^p(\widehat{T})} 
&\leq c  |\det ({A}_{\widetilde{T}})|^{\frac{p-1}{p}} \| \widetilde{{A}}^{-1} \|_2 \left \| \frac{\partial v}{\partial r_k} \right \|_{L^p(T)^3}.  \label{RT510}
\end{align}
By analogous argument in Lemma \ref{lem1151},
\begin{align}
\displaystyle
\| \hat{v}_j \|_{W^{1,p}(\widehat{T})}^p 
&= \| \hat{v}_j \|_{L^{p}(\widehat{T})}^p + \sum_{k=1}^3 \left \| \frac{\partial \hat{v}_j}{\partial \hat{x}_k}  \right \|^p_{L^p(\widehat{T})} \notag \\
&\quad \leq c |\det({A}_{\widetilde{T}})|^{p-1} \|  \widetilde{{A}}^{-1} \|_2^p h_j^{-p} \left[ \| v \|^p_{L^p(T)^3} + \left(  \sum_{ |\varepsilon| =1} h^{\varepsilon} \left \| \partial_r^{\varepsilon} v \right \|_{L^p(T)^3} \right)^p \right]. \label{RT511}
\end{align}
Combining the above inequalities \eqref{RT53}, \eqref{RT59}, \eqref{RT510}, and \eqref{RT511} with \eqref{CN331b} and \eqref{jensen} yields
\begin{align*}
\displaystyle
\| {I_{T}^{RT^k} v} \|_{L^p({T})^3} 
&\leq c |\det({A}_{\widetilde{T}})|^{- \frac{p - 1}{p}} \| \widetilde{{A}} \|_{2} \left(  \sum_{j=1}^d h_j^p \| (I_{\widehat{T}}^{RT^k} \hat{v})_j \|_{L^p(\widehat{T})}^p \right)^{1/p} \\
&\leq c \frac{H_{T}}{h_{T}} \left[  \| v \|_{L^p(T)^3} + \sum_{ |\varepsilon| =1} h^{\varepsilon} \left \| \partial_r^{\varepsilon} v \right \|_{L^p(T)^3} + \sum_{j=1}^3 h_j \sum_{k=1, k \neq j}^3 \left \| \frac{\partial v}{\partial r_k} \right \|_{L^p(T)^3} \right],
\end{align*}
which is the desired result.
\qed
\end{pf*}

\subsection{Local RT Interpolation Error Estimates}
The following two theorems are divided into the element on $\mathfrak{T}^{(2)}$ or $\mathfrak{T}_1^{(3)}$ and the element on $\mathfrak{T}_2^{(3)}$.

\begin{thr} \label{thr1161}
Let $p \in [1,\infty)$. Let $T\in \mathbb{T}_h$ satisfy Condition \ref{cond1} or Condition \ref{cond2} with $T ={\Phi}_{T} (\widetilde{T})$ and $\widetilde{T} = {\Phi}_{\widetilde{T}}(\widehat{T})$, where $\widetilde{T} \in \mathfrak{T}^{(2)}$ or $\widetilde{T} \in \mathfrak{T}_1^{(3)}$. For $k \in \mathbb{N}_0$, let $\{ {T} , \mathbb{RT}^k({T}) , {\Sigma} \}$ be the RT finite element and $I_{T}^{RT^k}$ the local interpolation operator defined in \eqref{RT11111}. Let $\ell$ be such that $0 \leq \ell \leq k$. Then, for any $\hat{v} \in W^{\ell+1,p}(\widehat{T})^d$ with $\tilde{v}= {\Psi}_{\widehat{T}} \hat{v}$ and ${v} = {\Psi}_{\widetilde{T}} \tilde{v}$,  
\begin{align}
\displaystyle
\| I_{T}^{RT^k} v - v \|_{L^p(T)^d} 
&\leq  c \left( \frac{H_{T}}{h_{T}} \sum_{|\varepsilon| = \ell+ 1} h^{\varepsilon} \left \| \partial_r^{\varepsilon} v \right \|_{L^p(T)^d} +  h_{T} \sum_{|\beta| = \ell} h^{\beta} \| \partial_r^{\beta} \nabla \cdot {v} \|_{L^{p}({T})} \right). \label{RT61}
\end{align}
If Condition \ref{Cond333} is imposed, it holds that
\begin{align}
\displaystyle
\| I_{T}^{RT^k} v - v \|_{L^p(T)^d}
&\leq c  \Biggl( \frac{H_{T}}{h_{T}}  \sum_{|\varepsilon| = \ell + 1} \widetilde{\mathscr{H}}^{\varepsilon} \| \partial^{\varepsilon}_{\tilde{x}} (\Psi_{\widetilde{T}}^{-1} v) \|_{L^p(\Phi_{T}^{-1}(T))^d} \notag \\
&\quad +  h_{T} \sum_{|\beta| = \ell} \widetilde{\mathscr{H}}^{\beta} \| \partial^{\beta}_{\tilde{x}} \nabla_{\tilde{x}} \cdot (\Psi_{\widetilde{T}}^{-1} v) \|_{L^{p}(\Phi_{T}^{-1}(T))} \Biggr). \label{RT62}
\end{align}
\end{thr}

\begin{pf*}
Let $\hat{v} \in W^{\ell+1,p}(\widehat{T})^d$. Let ${I}_{\widehat{T}}^{RT^k}$ be the local interpolation operators on $\widehat{T}$ defined by \eqref{RT1115} and \eqref{RT1116}. If ${q} \in \mathbb{P}^{\ell} ({T})^d \subset \mathbb{RT}^k({T})$, then $I_{{T}}^{RT} {q} = {q}$. 

We set ${\mathfrak{Q}}^{(\ell+1)} {v} := ({Q}^{(\ell+1)} {v}_1, \ldots , {Q}^{(\ell+1)} {v}_d)^{\top} \in \mathbb{P}^{\ell}({T})^d$, where ${Q}^{(\ell+1)} {v}_j$ is defined by \eqref{BH=2} for any $j$. We then obtain
\begin{align}
\displaystyle
\| I_{T}^{RT^k} v - v \|_{L^p(T)^d} &\leq \| I_{T}^{RT^k} (v - \mathfrak{Q}^{(\ell+1)} {v}) \|_{L^p(T)^d} + \| \mathfrak{Q}^{(\ell+1)} {v} - v \|_{L^p(T)^d}. \label{RT63}
\end{align}
The inequality \eqref{RT41} for the first term on the right-hand side of \eqref{RT63} yield
\begin{align}
\displaystyle
\| I_{T}^{RT^k} (v - \mathfrak{Q}^{(\ell+1)} {v}) \|_{L^p(T)^d}
&\leq c  |\det({A}_{\widetilde{T}})|^{\frac{1-p }{p}} \| \widetilde{{A}} \|_{2} \left(  \sum_{j=1}^d h_j^p \|  \{  I_{\widehat{T}}^{RT^k} (\hat{v} -  \widehat{\mathfrak{Q}}^{(\ell+1)} \hat{v}) \}_j \|_{L^p(\widehat{T})}^p \right)^{1/p}.\label{RT64}
\end{align}
The component-wise stability \eqref{RT1131} for $2d$ or \eqref{RT1132} for $3d$ yields
\begin{align}
\displaystyle
&\sum_{j=1}^d h_j^p \|  \{  I_{\widehat{T}_1}^{RT^k} (\hat{v} -  \widehat{\mathfrak{Q}}^{(\ell+1)} \hat{v}) \}_j \|_{L^p(\widehat{T})}^p \notag \\
&\quad \leq c \sum_{j=1}^d h_j^p \left( \| \hat{v}_j - \widehat{Q}^{(\ell+1)} \hat{v}_j \|_{W^{1,p}(\widehat{T})}^p + \| \nabla_{\hat{x}} \cdot ( \hat{v} - \widehat{\mathfrak{Q}}^{(\ell+1)} \hat{v}) \|_{L^p(\widehat{T})}^p \right). \label{RT65}
\end{align}
The inequality \eqref{RT41} for the second term on the right-hand side of \eqref{RT63} yields
\begin{align}
\displaystyle
\| \mathfrak{Q}^{(\ell+1)} {v} - v \|_{L^p(T)^d} 
&\leq c  |\det({A}_{\widetilde{T}})|^{\frac{1-p }{p}} \| \widetilde{{A}} \|_{2} \left(  \sum_{j=1}^d h_j^p \|  \widehat{Q}^{(\ell+1)} \hat{v}_j - \hat{v}_j \|_{L^p(\widehat{T})}^p \right)^{1/p}.\label{RT66}
\end{align}
The Bramble--Hilbert-type lemma (Lemma \ref{BH=1}) and \eqref{RT42},
\begin{align}
\displaystyle
\| \hat{v}_j - \widehat{Q}^{(\ell+1)} \hat{v}_j \|_{W^{1,p}(\widehat{T})}^p 
&= \| \hat{v}_j - \widehat{Q}^{(\ell+1)} \hat{v}_j \|_{L^{p}(\widehat{T})}^p + \sum_{k=1}^d \left \| \frac{\partial}{\partial \hat{x}_k} ( \hat{v}_j - \widehat{Q}^{(\ell+1)} \hat{v}_j ) \right \|_{L^{p}(\widehat{T})}^p \notag\\
&\leq c \left( \sum_{|\gamma| = \ell+ 1} \left \|  \partial_{\hat{x}}^{\gamma} \hat{v}_j \right \|_{L^{p}(\widehat{T})}^p + \sum_{k=1}^d \sum_{|\beta| = \ell} \left\|  \partial_{\hat{x}}^{\beta} \frac{\partial \hat{v}_j}{\partial \hat{x}_k} \right\|_{L^{p}(\widehat{T})}^p \right) \notag\\
&\leq c  |\det ({A}_{\widetilde{T}})|^{p-1} h_j^{-p} \| \widetilde{{A}}^{-1} \|_{2}^p \left( \sum_{|\varepsilon| = \ell+ 1} h^{\varepsilon} \left \| \partial_r^{\varepsilon} v \right \|_{L^p(T)^d} \right)^p. \label{RT67}
\end{align}
Because from \cite[Proposition 4.1.17]{BreSco08} it holds that
\begin{align}
\displaystyle
\widehat{\div} ( \widehat{\mathfrak{Q}}^{(\ell+1)} \hat{v}) = \widehat{Q}^{\ell}(\widehat{\div} \hat{v}). \label{RT68}
\end{align}
From the Bramble--Hilbert-type lemma (Lemma \ref{BH=1}) and \eqref{RT12},
\begin{align}
\displaystyle
 \|  \nabla_{\hat{x}} \cdot ( \hat{v} - \widehat{\mathfrak{Q}}^{(\ell+1)} \hat{v})) \|_{L^p(\widehat{T})}^p
  &=  \| \nabla_{\hat{x}} \cdot \hat{v} -\widehat{Q}^{\ell}(\nabla_{\hat{x}} \cdot \hat{v}) \|_{L^p(\widehat{T})}^p \notag\\
 &\leq  \| \nabla_{\hat{x}} \cdot \hat{v} -\widehat{Q}^{\ell}( \nabla_{\hat{x}} \cdot \hat{v}) \|_{W^{\ell,p}(\widehat{T})}^p \notag\\
 &\leq c | \nabla_{\hat{x}} \cdot \hat{v} |_{W^{\ell,p}(\widehat{T})}^p 
 = c \sum_{|\beta| = \ell} \| \partial^{\beta} \nabla_{\hat{x}} \cdot \hat{v} \|_{L^{p}(\widehat{T})}^p \notag\\
& \leq  c  |\det ({A}_{\widetilde{T}})|^{{p-1}} \left( \sum_{|\varepsilon| = \ell} h^{\varepsilon} \left \|  \partial_r^{\varepsilon} \nabla \cdot v \right \|_{L^p(T)} \right)^p. \label{RT69}
\end{align}
Combining \eqref{RT64}, \eqref{RT65}, \eqref{RT67}, and \eqref{RT69} with \eqref{CN331b} yields
\begin{align}
\displaystyle
\| I_{T}^{RT^k} (v - \mathfrak{Q}^{(\ell+1)} {v}) \|_{L^p(T)^d}
&\leq c \left( \frac{H_{T}}{h_{T}} \sum_{|\varepsilon| = \ell+ 1} h^{\varepsilon} \left \|  \partial_r^{\varepsilon} v \right \|_{L^p(T)^d} +  h_{T} \sum_{|\beta| = \ell} h^{\beta} \| \partial_r^{\beta} \nabla \cdot {v} \|_{L^{p}({T})} \right). \label{RT610}
\end{align}
Furthermore, using a similar argument, from the Bramble--Hilbert-type lemma (Lemma \ref{BH=1}), \eqref{RT42}, and \eqref{RT66} together with \eqref{CN331b},
\begin{align}
\displaystyle
\| \mathfrak{Q}^{(\ell+1)} {v} - v \|_{L^p(T)^d}
&\leq c \frac{H_{T}}{h_{T}} \sum_{|\varepsilon| = \ell+ 1} h^{\varepsilon} \left \| \partial_r^{\varepsilon} v \right \|_{L^p(T)^d}. \label{RT611}
\end{align}
Therefore, from \eqref{RT63}, \eqref{RT610}, and \eqref{RT611}, we have \eqref{RT61}.

\textbf{Case in which Condition \ref{Cond333} is imposed.} 
From the Bramble--Hilbert-type lemma (Lemma \ref{BH=1}) and \eqref{RT43},
\begin{align}
\displaystyle
\| \hat{v}_j - \widehat{Q}^{(\ell+1)} \hat{v}_j \|_{W^{1,p}(\widehat{T})}^p 
&\leq c |\hat{v}_j|_{W^{\ell+1,p}(\widehat{T}_1)}^p + c \sum_{k=1}^d \left | \frac{\partial \hat{v}_j}{\partial \hat{x}_k} \right |_{W^{\ell,p}(\widehat{T})}^p \notag\\
&= c \left( \sum_{|\gamma| = \ell+ 1} \left \|  \partial_{\hat{x}}^{\gamma} \hat{v}_j \right \|_{L^{p}(\widehat{T})}^p + \sum_{k=1}^d \sum_{|\beta| = \ell} \left\|  \partial_{\hat{x}}^{\beta} \frac{\partial \hat{v}_j}{\partial \hat{x}_k} \right\|_{L^{p}(\widehat{T})}^p \right) \notag\\
&\leq c  |\det ({A}_{\widetilde{T}})|^{{p-1}} h_j^{-p} \| \widetilde{{A}}^{-1} \|_{2}^p \left( \sum_{|\varepsilon| = \ell + 1} \widetilde{\mathscr{H}}^{\varepsilon} \| \partial^{\varepsilon}_{\tilde{x}} \tilde{v} \|_{L^p(\widetilde{T})^d} \right)^p. \label{RT612}
\end{align}
Because \eqref{RT68}, from the Bramble--Hilbert-type lemma (Lemma \ref{BH=1}) and \eqref{RT13},
\begin{align}
\displaystyle
 \| \nabla_{\hat{x}} \cdot ( \hat{v} - \widehat{\mathfrak{Q}}^{(\ell+1)} \hat{v})) \|_{L^p(\widehat{T})}^p
&\leq c  |\det ({A}_{\widetilde{T}})|^{{p-1}} \left( \sum_{|\varepsilon| = \ell} \widetilde{\mathscr{H}}^{\varepsilon} \left \|  \partial^{\varepsilon}_{\tilde{x}} \nabla_{\tilde{x}} \cdot \tilde{v} \right \|_{L^p(\widetilde{T})} \right)^p. \label{RT613}
\end{align}
Combining \eqref{RT64}, \eqref{RT65}, \eqref{RT612}, and \eqref{RT613} with \eqref{CN331b} yields
\begin{align}
\displaystyle
\| I_{T}^{RT^k} (v - \mathfrak{Q}^{(\ell+1)} {v}) \|_{L^p(T)^d}
&\leq c \left( \frac{H_{T}}{h_{T}}  \sum_{|\varepsilon| = \ell + 1} \widetilde{\mathscr{H}}^{\varepsilon} \| \partial^{\varepsilon}_{\tilde{x}} \tilde{v} \|_{L^p(\widetilde{T})^d} +  h_{T} \sum_{|\beta| = \ell} \widetilde{\mathscr{H}}^{\beta} \| \partial^{\beta}_{\tilde{x}} \nabla_{\tilde{x}} \cdot \tilde{v} \|_{L^{p}(\widetilde{T})} \right). \label{RT614}
\end{align}
Furthermore, using a similar argument, from the Bramble--Hilbert-type lemma (Lemma \ref{BH=1}), \eqref{RT43}, and \eqref{RT66} together with \eqref{CN331b},
\begin{align}
\displaystyle
\| \mathfrak{Q}^{(\ell+1)} {v} - v \|_{L^p(T)^d}
&\leq c \frac{H_{T}}{h_{T}}  \sum_{|\varepsilon| = \ell + 1} \widetilde{\mathscr{H}}^{\varepsilon} \| \partial^{\varepsilon}_{\tilde{x}} \tilde{v} \|_{L^p(\widetilde{T})^d}, \label{RT615}
\end{align}
Therefore, from \eqref{RT63}, \eqref{RT614}, and \eqref{RT615}, we have \eqref{RT62}.
\qed
\end{pf*}

\begin{thr}
Let $p \in [1,\infty)$ and $d=3$. Let $T\in \mathbb{T}_h$ satisfy Condition \ref{cond2} with $T ={\Phi}_{T} (\widetilde{T})$ and $\widetilde{T} = {\Phi}_{\widetilde{T}}(\widehat{T})$, where $\widetilde{T} \in \mathfrak{T}_2^{(3)}$. For $k \in \mathbb{N}_0$, let $\{ {T} , \mathbb{RT}^k({T}) , {\Sigma} \}$ be the RT finite element and $I_{T}^{RT^k}$ the local interpolation operator defined in \eqref{RT11111}. Let $\ell$ be such that $0 \leq \ell \leq k$. Then, for any $\hat{v} \in W^{\ell+1,p}(\widehat{T})^d$ with $\tilde{v}= {\Psi}_{\widehat{T}} \hat{v}$ and ${v} = {\Psi}_{\widetilde{T}} \tilde{v}$, 
\begin{align}
\displaystyle
&\| I_{T}^{RT^k} v - v \|_{L^p(T)^3} 
\leq c \frac{H_{T}}{h_{T}} \Biggl(  h_T \sum_{k=1}^3 \sum_{|\varepsilon| = \ell} h^{\varepsilon} \left \| \partial_{r}^{\varepsilon} \frac{\partial v}{\partial r_k}  \right \|_{L^p(T)^3} \Biggr). \label{RT616}
\end{align}
If Condition \ref{Cond333} is imposed, it holds that
\begin{align}
\displaystyle
\| I_{T}^{RT^k} v - v \|_{L^p(T)^3}
&\leq c \frac{H_{T}}{h_{T}} \Biggl( \sum_{|\varepsilon| = \ell+ 1} \widetilde{\mathscr{H}}^{\varepsilon} \| \partial^{\varepsilon}_{\tilde{x}} (\Psi_{\widetilde{T}}^{-1}  v ) \|_{L^p( \Phi_{T}^{-1} (T))^3} \notag \\
&\quad \quad + h_T \sum_{k=1}^3 \sum_{|\varepsilon| = \ell} \widetilde{\mathscr{H}}^{\varepsilon} \left \| \partial^{\varepsilon}_{\tilde{x}} \frac{\partial (\Psi_{\widetilde{T}}^{-1} v)}{\partial \tilde{r}_k^s} \right \|_{L^p(\Phi_{T}^{-1}(T))^3} \Biggr). \label{RT616b}
\end{align}
\end{thr}

\begin{pf*}
An analogous proof of Theorem \ref{thr1161} yields the desired result \eqref{RT616}, where we use Lemma \ref{lem1136} instead of Lemma \ref{lem1134}, and Lemma \ref{lem1143} instead of Lemma \ref{lem1142}.

Let $\hat{v} \in W^{\ell+1,p}(\widehat{T}_2)^3$. Let ${I}_{\widehat{T}}^{RT^k}$ be the local interpolation operators on $\widehat{T}$ defined by \eqref{RT1115} and \eqref{RT1116}. If ${q} \in \mathbb{P}^{\ell} ({T})^3 \subset \mathbb{RT}^k({T})$, then $I_{{T}}^{RT^k} {q} = {q}$. 

We set ${\mathfrak{Q}}^{(\ell+1)} {v} := ({Q}^{(\ell+1)} {v}_1, {Q}^{(\ell+1)} {v}_2 , {Q}^{(\ell+1)} {v}_3)^{\top} \in \mathbb{P}^{\ell}({T})^3$, where ${Q}^{(\ell+1)} {v}_j^s$ is defined by \eqref{BH=2} for any $j$. We then obtain
\begin{align}
\displaystyle
\| I_{T}^{RT^k} v - v \|_{L^p(T)^3} &\leq \| I_{T}^{RT^k} (v - \mathfrak{Q}^{(\ell+1)} {v}) \|_{L^p(T)^3} + \| \mathfrak{Q}^{(\ell+1)} {v} - v \|_{L^p(T)^3}. \label{RT617}
\end{align}
The inequality \eqref{RT41} for the first term on the right-hand side of \eqref{RT617} yields
\begin{align}
\displaystyle
\| I_{T}^{RT^k} (v - \mathfrak{Q}^{(\ell+1)} {v}) \|_{L^p(T)^3}
&\leq c  |\det({A}_{\widetilde{T}})|^{\frac{1-p }{p}} \| \widetilde{{A}} \|_{2} \left(  \sum_{j=1}^3 h_j^p \|  \{  I_{\widehat{T}}^{RT^k} (\hat{v} -  \widehat{\mathfrak{Q}}^{(\ell+1)} \hat{v}) \}_j \|_{L^p(\widehat{T})}^p \right)^{1/p}.\label{RT618}
\end{align}
The component-wise stability \eqref{RT11312} for $3d$ yields
\begin{align}
\displaystyle
&\sum_{j=1}^3 h_j^p \|  \{  I_{\widehat{T}}^{RT^k} (\hat{v} -  \widehat{\mathfrak{Q}}^{(\ell+1)} \hat{v}) \}_j \|_{L^p(\widehat{T})}^p \notag \\
&\quad \leq c \sum_{j=1}^3 h_j^p \left( \| \hat{v}_j - \widehat{Q}^{(\ell+1)} \hat{v}_j \|_{W^{1,p}(\widehat{T})}^p +\sum_{k=1, k \neq j}^3 \left\| \frac{\partial}{\partial \hat{x}_k} ( \hat{v} - \widehat{\mathfrak{Q}}^{(\ell+1)} \hat{v})_k \right\|_{L^p(\widehat{T})}^p \right). \label{RT619}
\end{align}
The inequality \eqref{RT41} for the second term on the right-hand side of \eqref{RT617} yields
\begin{align}
\displaystyle
\| \mathfrak{Q}^{(\ell+1)} {v} - v \|_{L^p(T)^3}
&\leq c  |\det({A}_{\widetilde{T}})|^{\frac{1-p }{p}} \| \widetilde{{A}} \|_{2} \left(  \sum_{j=1}^d h_j^p \|  \widehat{Q}^{(\ell+1)} \hat{v}_j - \hat{v}_j \|_{L^p(\widehat{T})}^p \right)^{1/p}.\label{RT623}
\end{align}
From the Bramble--Hilbert-type lemma (Lemma \ref{BH=1}) and \eqref{RT42}, we have
\begin{align}
\displaystyle
\| \hat{v}_j - \widehat{Q}^{(\ell+1)} \hat{v}_j \|_{W^{1,p}(\widehat{T})}^p
&= \| \hat{v}_j - \widehat{Q}^{(\ell+1)} \hat{v}_j \|_{L^{p}(\widehat{T})}^p + \sum_{k=1}^d \left \| \frac{\partial}{\partial \hat{x}_k} ( \hat{v}_j - \widehat{Q}^{(\ell+1)} \hat{v}_j ) \right \|_{L^{p}(\widehat{T})}^p \notag\\
&\leq c \left( \sum_{|\gamma| = \ell+ 1} \left \|  \partial_{\hat{x}}^{\gamma} \hat{v}_j \right \|_{L^{p}(\widehat{T})}^p + \sum_{k=1}^3 \sum_{|\beta| = \ell} \left\|  \partial_{\hat{x}}^{\beta} \frac{\partial \hat{v}_j}{\partial \hat{x}_k} \right\|_{L^{p}(\widehat{T})}^p \right) \notag\\
&\leq c  |\det ({A}_{\widetilde{T}})|^{p-1} h_j^{-p} \| \widetilde{{A}}^{-1} \|_{2}^p \left( \sum_{|\varepsilon| = \ell+ 1} h^{\varepsilon} \left \|  \partial^{\varepsilon}_r v \right \|_{L^p(T)^3} \right)^p.  \label{RT620}
\end{align}
From the Bramble--Hilbert-type lemma (Lemma \ref{BH=1}) and  \eqref{RT14}, we have
\begin{align}
\displaystyle
\sum_{k=1, k \neq j}^3 \left \| \frac{\partial}{\partial \hat{x}_k} ( \hat{v}_k - \widehat{{Q}}^{(\ell+1)} \hat{v}_k)) \right \|_{L^p(\widehat{T})}^p \notag 
&= \sum_{k=1, k \neq j}^3 \left \| \frac{\partial \hat{v}_k }{\partial \hat{x}_k} - \widehat{{Q}}^{(\ell+1)} \left( \frac{\partial \hat{v}_k }{\partial \hat{x}_k} \right ) \right \|_{L^p(\widehat{T})}^p \\
&\hspace{-4.5cm} \leq c \sum_{k=1, k \neq j}^3  \sum_{|\beta| = \ell} \left\|  \partial_{\hat{x}}^{\beta} \frac{\partial \hat{v}_k}{\partial \hat{x}_k} \right\|_{L^{p}(\widehat{T})}^p \notag \\
&\hspace{-4.5cm} \leq c  |\det ({A}_{\widetilde{T}})|^{p-1} \| \widetilde{{A}}^{-1} \|^p_2 \sum_{k=1, k \neq j}^3 \sum_{|\varepsilon| = \ell} h^{\varepsilon p} \left \| \partial_{r}^{\varepsilon} \frac{\partial v}{\partial r_k}  \right \|_{L^p(T)^3}^p. \label{RT621}
\end{align}
Gathering \eqref{RT618}, \eqref{RT619}, \eqref{RT620} and \eqref{RT621} together with \eqref{jensen} and \eqref{CN331b} yields
\begin{align}
\displaystyle
&\| I_{T}^{RT^k} (v - \mathfrak{Q}^{(\ell+1)} {v}) \|_{L^p(T)^3} \notag\\
&\quad \leq c  \frac{H_{T}}{h_{T}} \Biggl ( \sum_{|\varepsilon| = \ell+ 1} h^{\varepsilon} \left \|  \partial^{\varepsilon}_r v \right \|_{L^p(T)^3} + \sum_{j=1}^3 h_j \sum_{k=1, k \neq j}^3 \sum_{|\varepsilon| = \ell} h^{\varepsilon} \left \| \partial_{r}^{\varepsilon} \frac{\partial v}{\partial r_k}  \right \|_{L^p(T)^3} \Biggr)  \notag \\
&\quad \leq c  \frac{H_{T}}{h_{T}} \Biggl ( \sum_{|\varepsilon| = \ell+ 1} h^{\varepsilon} \left \|  \partial^{\varepsilon}_r v \right \|_{L^p(T)^3} + h_T \sum_{k=1}^3 \sum_{|\varepsilon| = \ell} h^{\varepsilon} \left \| \partial_{r}^{\varepsilon} \frac{\partial v}{\partial r_k}  \right \|_{L^p(T)^3} \Biggr). \label{RT622}
\end{align}
From the Bramble--Hilbert-type lemma (Lemma \ref{BH=1}) and \eqref{RT42},
\begin{align}
\displaystyle
\| \hat{v}_j - \widehat{Q}^{(\ell+1)} \hat{v}_j \|_{L^{p}(\widehat{T})}^p
&\leq c \sum_{|\gamma| = \ell+ 1} \left \|  \partial_{\hat{x}}^{\gamma} \hat{v}_j \right \|_{L^{p}(\widehat{T})}^p \notag\\
&\leq c  |\det ({A}_{\widetilde{T}})|^{p-1} h_j^{-p} \| \widetilde{{A}}^{-1} \|_{2}^p \left( \sum_{|\varepsilon| = \ell+ 1} h^{\varepsilon} \left \| \partial_r^{\varepsilon} v \right \|_{L^p(T)^3} \right)^p. \label{RT624}
\end{align}
From \eqref{RT623} and \eqref{RT624} together with \eqref{CN331b}, we obtain
\begin{align}
\displaystyle
\| \mathfrak{Q}^{(\ell+1)} {v} - v \|_{L^p(T)^3}
&\leq c  \frac{H_{T}}{h_{T}} \sum_{|\varepsilon| = \ell+ 1} h^{\varepsilon} \left \|  \partial^{\varepsilon}_r v \right \|_{L^p(T)^3}. \label{RT625}
\end{align}
Therefore, from \eqref{RT617} and \eqref{RT622}, \eqref{RT625}, we have \eqref{RT616}.

\textbf{Case in which Condition \ref{Cond333} is imposed.} 
From the Bramble--Hilbert-type lemma (Lemma \ref{BH=1}) and \eqref{RT43}, we have
\begin{align}
\displaystyle
\| \hat{v}_j - \widehat{Q}^{(\ell+1)} \hat{v}_j \|_{W^{1,p}(\widehat{T})}^p
&= \| \hat{v}_j - \widehat{Q}^{(\ell+1)} \hat{v}_j \|_{L^{p}(\widehat{T})}^p + \sum_{k=1}^d \left \| \frac{\partial}{\partial \hat{x}_k} ( \hat{v}_j - \widehat{Q}^{(\ell+1)} \hat{v}_j ) \right \|_{L^{p}(\widehat{T})}^p \notag\\
&\leq c \left( \sum_{|\gamma| = \ell+ 1} \left \|  \partial_{\hat{x}}^{\gamma} \hat{v}_j \right \|_{L^{p}(\widehat{T})}^p + \sum_{k=1}^3 \sum_{|\beta| = \ell} \left\|  \partial_{\hat{x}}^{\beta} \frac{\partial \hat{v}_j}{\partial \hat{x}_k} \right\|_{L^{p}(\widehat{T})}^p \right) \notag\\
&\leq c  |\det ({A}_{\widetilde{T}})|^{p-1} h_j^{-p} \| \widetilde{{A}}^{-1} \|_{2}^p \left(\sum_{|\varepsilon| = \ell+1} \widetilde{\mathscr{H}}^{\varepsilon} \| \partial_{\tilde{x}}^{\varepsilon} (\Psi_{\widetilde{T}}^{-1}  v ) \|_{L^p( \Phi_{T}^{-1} (T))^3} \right)^p.  \label{RT9627}
\end{align}
From the Bramble--Hilbert-type lemma (Lemma \ref{BH=1}) and  \eqref{RT14}, we have
\begin{align}
\displaystyle
\sum_{k=1, k \neq j}^3 \left \| \frac{\partial}{\partial \hat{x}_k} ( \hat{v}_k - \widehat{{Q}}^{(\ell+1)} \hat{v}_k)) \right \|_{L^p(\widehat{T})}^p \notag 
&= \sum_{k=1, k \neq j}^3 \left \| \frac{\partial \hat{v}_k }{\partial \hat{x}_k} - \widehat{{Q}}^{(\ell+1)} \left( \frac{\partial \hat{v}_k }{\partial \hat{x}_k} \right ) \right \|_{L^p(\widehat{T})}^p \\
&\hspace{-4.5cm} \leq c \sum_{k=1, k \neq j}^3  \sum_{|\beta| = \ell} \left\|  \partial_{\hat{x}}^{\beta} \frac{\partial \hat{v}_k}{\partial \hat{x}_k} \right\|_{L^{p}(\widehat{T})}^p \notag \\
&\hspace{-4.5cm} \leq c  |\det ({A}_{\widetilde{T}})|^{p-1} \| \widetilde{{A}}^{-1} \|^p_2 \sum_{k=1, k \neq j}^3 \sum_{|\varepsilon| = \ell} \widetilde{\mathscr{H}}^{\varepsilon} \left \| \partial^{\varepsilon}_{\tilde{x}} \frac{\partial (\Psi_{\widetilde{T}}^{-1} v)}{\partial \tilde{r}_k^s} \right \|_{L^p(\Phi_{T}^{-1}(T))^3}. \label{RT9628}
\end{align}
Gathering \eqref{RT618}, \eqref{RT619}, \eqref{RT9627} and \eqref{RT9628} together with \eqref{jensen} and \eqref{CN331b} yields
\begin{align}
\displaystyle
\| I_{T}^{RT^k} (v - \mathfrak{Q}^{(\ell+1)} {v}) \|_{L^p(T)^3}
&\leq c  \frac{H_{T}}{h_{T}} \Biggl ( \sum_{|\varepsilon| = \ell+ 1} \widetilde{\mathscr{H}}^{\varepsilon} \| \partial_{\tilde{x}}^{\varepsilon} (\Psi_{\widetilde{T}}^{-1}  v ) \|_{L^p( \Phi_{T}^{-1} (T))^3} \notag \\
&\quad \quad + \sum_{j=1}^3 h_j \sum_{k=1, k \neq j}^3 \sum_{|\varepsilon| = \ell} \widetilde{\mathscr{H}}^{\varepsilon} \left \| \partial^{\varepsilon}_{\tilde{x}} \frac{\partial (\Psi_{\widetilde{T}}^{-1} v)}{\partial \tilde{r}_k} \right \|_{L^p(\Phi_{T}^{-1}(T))^3} \Biggr)  \notag \\
&\leq c  \frac{H_{T}}{h_{T}} \Biggl ( \sum_{|\varepsilon| = \ell+ 1} \widetilde{\mathscr{H}}^{\varepsilon} \| \partial_{\tilde{x}}^{\varepsilon} (\Psi_{\widetilde{T}}^{-1}  v ) \|_{L^p( \Phi_{T}^{-1} (T))^3} \notag \\
&\quad \quad+ h_T \sum_{k=1}^3 \sum_{|\varepsilon| = \ell} \widetilde{\mathscr{H}}^{\varepsilon} \left \| \partial^{\varepsilon}_{\tilde{x}} \frac{\partial (\Psi_{\widetilde{T}}^{-1} v)}{\partial \tilde{r}_k} \right \|_{L^p(\Phi_{T}^{-1}(T))^3} \Biggr). \label{RT9629}
\end{align}
From the Bramble--Hilbert-type lemma (Lemma \ref{BH=1}) and \eqref{RT43},
\begin{align}
\displaystyle
\| \hat{v}_j - \widehat{Q}^{(\ell+1)} \hat{v}_j \|_{L^{p}(\widehat{T})}^p 
&\leq c \sum_{|\gamma| = \ell+ 1} \left \|  \partial_{\hat{x}}^{\gamma} \hat{v}_j \right \|_{L^{p}(\widehat{T})}^p \notag\\
&\leq c  |\det ({A}_{\widetilde{T}})|^{p-1} h_j^{-p} \| \widetilde{{A}}^{-1} \|_{2}^p \left(\sum_{|\varepsilon| = \ell+1} \widetilde{\mathscr{H}}^{\varepsilon} \| \partial_{\tilde{x}}^{\varepsilon} (\Psi_{\widetilde{T}}^{-1}  v ) \|_{L^p( \Phi_{T}^{-1} (T))^3} \right)^p. \label{RT9630}
\end{align}
From \eqref{RT623} and \eqref{RT9630} together with \eqref{CN331b}, we obtain
\begin{align}
\displaystyle
\| \mathfrak{Q}^{(\ell+1)} {v} - v \|_{L^p(T)^3}
&\leq c  \frac{H_{T}}{h_{T}} \sum_{|\varepsilon| = \ell+ 1} \widetilde{\mathscr{H}}^{\varepsilon} \| \partial_{\tilde{x}}^{\varepsilon} (\Psi_{\widetilde{T}}^{-1}  v ) \|_{L^p( \Phi_{T}^{-1} (T))^3}. \label{RT6231}
\end{align}
Therefore, from \eqref{RT617} and \eqref{RT9629}, \eqref{RT6231}, we have \eqref{RT616b}.
\qed
\end{pf*}

\subsection{Global RT Interpolation Error Estimates}
We define a broken finite element space as
\begin{align*}
\displaystyle
RT^k(\mathbb{T}_h) := \left\{ v_h \in L^1(\Omega)^d; \ v_h|_T \in \mathbb{RT}^k({T}) \ \forall T \in \mathbb{T}_h \right \}.
\end{align*}
The corresponding (global) RT finite element space is defined as
\begin{align*}
\displaystyle
V^{RT^k}_{h} &:= \{ v_h \in RT^k(\mathbb{T}_h) ;  \  [\![ v_h \cdot n ]\!]_F = 0, \ \forall F \in \mathcal{F}_h^i \}.
\end{align*}

\begin{lem}
It holds that
\begin{align*}
\displaystyle
V^{RT^k}_h \subset H(\div; \Omega).
\end{align*}
\end{lem}

\begin{pf*} 
Let $v_h \in V^{RT^k}_h$. Because its restriction to every $T \in \mathbb{T}_h$ is a polynomial, it is differentiable in the classical sense. Let us consider the function $w_h \in L^2(\Omega)$ defined on $T$ by $w_h|_{T} = \div(v_h)|_{T}$. Let $\varphi \in C_{0}^{\infty}(\Omega)$. Then, using the Green formula yields
\begin{align*}
\displaystyle
\int_{\Omega} w_h \varphi dx &= \sum_{T \in \mathbb{T}_h} \int_{T} w_h \varphi dx \\
&= - \sum_{T \in \mathbb{T}_h} \int_{T} (v_h)|_{T} \cdot \nabla \varphi dx + \sum_{F \in \mathcal{F}_h^i} \int_F [\![ v_h \cdot n ]\!]_F \varphi ds.
\end{align*}
Because $ [\![ v_h \cdot n ]\!]_F = 0$,
\begin{align*}
\displaystyle
\int_{\Omega} w_h \varphi dx = - \int_{\Omega} (v_h \cdot  \nabla) \varphi dx.
\end{align*}
Therefore, the distributional divergence of $v_h$ is $w_h$. Because $w_h \in L^2(\Omega)$, $\div v_h \in L^2(\Omega)$.
\qed
\end{pf*}

We define the global RT  interpolation $I_h^{RT^k} : W^{1,1}(\Omega) \to V^{RT^k}_{h}$ as
\begin{align*}
\displaystyle
(I_h^{RT^k} v )|_{T} = I_{T}^{RT^k} (v|_{T}) \quad \forall T \in \mathbb{T}_h, \quad \forall v \in W^{1,1}(\Omega).
\end{align*}

\begin{cor}[de Rham complex] \label{cor1171}
The following diagram commutes:
\begin{align*}
\displaystyle
\xymatrix{
W^{1,1}(\Omega)^d \ar[r]^-{\nabla \cdot} \ar[d]_-{{I}_{h}^{RT^k}}&L^1(\Omega)\ar[d]^-{{\Pi}_{h}^{k}} \\
V_h^{RT^k}\ar[r]_-{\nabla \cdot}& P_{dc,h}^{k}
}
\end{align*}
In other words, it holds that
\begin{align}
\displaystyle
\div (I_{h}^{RT^k} v) = \Pi_{h}^k (\div v) \quad \forall v \in W^{1,1}(\Omega)^d. \label{RT171}
\end{align}
\end{cor}

\begin{pf*}
Combine Lemma \ref{rtl2=com}.
\qed	
\end{pf*}

\begin{cor}[Stability] \label{RTstability}
Let $p \in [1,\infty)$.  We impose Condition \ref{Cond1=sec6} with $h \leq 1$. Then,
\begin{align*}
\displaystyle
\| {I_{h}^{RT^k} v} \|_{L^p({\Omega})^d}
&\leq c \| v\|_{W^{1,p}(\Omega)^d} \quad \forall v \in W^{1,p}(\Omega)^d. 
\end{align*}
\end{cor}

\begin{pf*}
Lemmata \ref{lem1151} and \ref{lem1152} yield
\begin{align*}
\displaystyle
\| {I_{h}^{RT^k} v} \|_{L^p({\Omega})^d}^p
&= \sum_{T \in \mathbb{T}_h} \| {I_{T}^{RT^k} v} \|_{L^p({T})^d}^p
\leq c  \sum_{T \in \mathbb{T}_h} \| v\|_{W^{1,p}(T)^d}^p = c  \| v\|_{W^{1,p}(\Omega)^d}^p,
\end{align*}
which leads to the desired result. 
\qed
\end{pf*}

\begin{cor} \label{thr1173}
Let $p \in [1,\infty)$.  We impose Condition \ref{Cond1=sec6} with $h \leq 1$. Let $\ell$ be such that $0 \leq \ell \leq k$. Then, for any $ v \in W^{\ell+1,p}(\Omega)^d$,  if all mesh elements $T\in \mathbb{T}_h$ satisfy Condition \ref{cond1} or Condition \ref{cond2} with $T ={\Phi}_{T} (\widetilde{T})$ and $\widetilde{T} = {\Phi}_{\widetilde{T}}(\widehat{T})$, where $\widetilde{T} \in \mathfrak{T}^{(2)}$ or $\widetilde{T} \in \mathfrak{T}_1^{(3)}$,
\begin{align}
 \displaystyle
\| I_{h}^{RT^k} v - v \|_{L^p(\Omega)^d} 
&\leq c \sum_{T \in \mathbb{T}_h} \sum_{|\varepsilon| = \ell+ 1} h^{\varepsilon} \left \| \partial_r^{\varepsilon} v \right \|_{L^p(T)^d} + h \left( \sum_{T \in \mathbb{T}_h}  \sum_{|\beta| = \ell} h^{\beta p} \| \partial_r^{\beta} \nabla \cdot {v} \|_{L^{p}({T})}^p \right)^{1/p}.\label{RT72}
\end{align}
Furthermore, if Condition \ref{Cond333} is imposed, 
\begin{align}
\displaystyle
\| I_{h}^{RT^k} v - v \|_{L^p(\Omega)^d} 
&\leq  c \sum_{T \in \mathbb{T}_h}  \Biggl(  \sum_{|\varepsilon| = \ell + 1} \widetilde{\mathscr{H}}^{\varepsilon} \| \partial^{\varepsilon}_{\tilde{x}} (\Psi_{\widetilde{T}}^{-1} v) \|_{L^p(\Phi_{T}^{-1}(T))^d} \notag \\
&\quad \quad  + h_{T} \sum_{|\beta| = \ell} \widetilde{\mathscr{H}}^{\beta} \| \partial^{\beta}_{\tilde{x}} \nabla_{\tilde{x}} \cdot (\Psi_{\widetilde{T}}^{-1} v) \|_{L^{p}(\Phi_{T}^{-1}(T))} \Biggr). \label{RT73}
\end{align}
Let $d=3$. For any $ v \in W^{\ell+1,p}(\Omega)^d$, if all mesh elements $T\in \mathbb{T}_h$ satisfy Condition \ref{cond2} with $T ={\Phi}_{T} (\widetilde{T})$ and $\widetilde{T} = {\Phi}_{\widetilde{T}}(\widehat{T})$, where $\widetilde{T} \in \mathfrak{T}_2^{(3)}$,
\begin{align}
\displaystyle
&\| I_{h}^{RT^k} v - v \|_{L^p(\Omega)^3} 
\leq c h  \Biggl( \sum_{T \in \mathbb{T}_h} \sum_{k=1}^3 \sum_{|\varepsilon| = \ell} h^{\varepsilon p} \left \| \partial_{r}^{\varepsilon} \frac{\partial v}{\partial r_k}  \right \|_{L^p(T)^3}^p \Biggr)^{1/p}. \label{RT74}
\end{align}
Furthermore, if Condition \ref{Cond333} is imposed, 
\begin{align}
\displaystyle
\| I_{h}^{RT^k} v - v \|_{L^p(\Omega)^3}
&\leq c  \sum_{T \in \mathbb{T}_h} \Biggl( \sum_{|\varepsilon| = \ell+ 1} \widetilde{\mathscr{H}}^{\varepsilon} \| \partial_{\tilde{x}}^{\varepsilon} (\Psi_{\widetilde{T}}^{-1}  v ) \|_{L^p( \Phi_{T}^{-1} (T))^3} \notag \\
&\quad \quad + h_T \sum_{k=1}^3 \sum_{|\varepsilon| = \ell} \widetilde{\mathscr{H}}^{\varepsilon} \left \| \partial_{\tilde{x}}^{\varepsilon} \frac{\partial (\Psi_{\widetilde{T}}^{-1} v)}{\partial \tilde{r}_k} \right \|_{L^p(\Phi_{T}^{-1}(T))^3} \Biggr). \label{RT74b}
\end{align}
\end{cor}

\begin{pf*}
This corollary is proved in the same argument as Corollary \ref{newglobal=cor}.
\qed
\end{pf*}

\section{Inverse Inequalities on Anisotropic Meshes}
This section presents some limited results for the inverse inequalities.

\begin{lem}
Let $\widehat{P} := \mathbb{P}^k$ with $k \in \mathbb{N}$. Let $p,q \in [1,\infty]$. Then, if $d=2$, there exists a positive constant $C^{IV,2}$, independent of $h_{T}$ and ${T}$, such that, for all ${\varphi}_h \in{P} = \{  \hat{\varphi}_h \circ {\Phi}^{-1} ; \ \hat{\varphi}_h \in \widehat{P} \}$,
\begin{align}
\displaystyle
\left \| \frac{\partial \varphi_h}{\partial x_i} \right \|_{L^q(T)}
&\leq 
C^{IV,2} |T|_d^{\frac{1}{q} - \frac{1}{p}} \left( \frac{1}{\widetilde{\mathscr{H}}_1}  |(A_T^{-1})_{1i}|  + \frac{2}{\widetilde{\mathscr{H}}_2} |(A_T^{-1})_{2i}| \right) \| {\varphi}_h \|_{L^p({T})}, \quad \text{$i=1,2$}. \label{inv=2d}
\end{align}
Let $d=3$. If Condition \ref{Cond333} is imposed, there exist positive constants $C^{IV,3}$, independent of $h_{T}$ and ${T}$, such that, for all ${\varphi}_h \in{P} = \{  \hat{\varphi}_h \circ {\Phi}^{-1} ; \ \hat{\varphi}_h \in \widehat{P} \}$, 
\begin{align}
\displaystyle
&\left \| \frac{\partial \varphi_h}{\partial x_i} \right \|_{L^q(T)} \notag \\
&\leq 
C^{IV,2} |T|_d^{\frac{1}{q} - \frac{1}{p}} \left( \frac{1}{\widetilde{\mathscr{H}}_1}  |(A_T^{-1})_{1i}|  + \frac{2}{\widetilde{\mathscr{H}}_2} |(A_T^{-1})_{2i}| + \frac{2(M+1)}{\widetilde{\mathscr{H}}_3} |(A_T^{-1})_{3i}|  \right) \| {\varphi}_h \|_{L^p({T})}, \label{inv=3d}
\end{align}
for $i=1,2,3$.
\end{lem}

\begin{pf*}
Let ${\varphi}_h \in{P}$. Then, for $ i= 1,\ldots,d$,
\begin{align*}
\displaystyle
\frac{\partial  {\varphi}_h}{\partial {x}_i} 
= \sum_{j,k=1}^d \frac{\partial \hat{\varphi}_h}{\partial \hat{x}_j} h_j^{-1} \widetilde{A}^{-1}_{jk} (A_T^{-1})_{ki},
\end{align*}
which leads to
\begin{align*}
\displaystyle
\left \| \frac{\partial \varphi_h}{\partial x_i} \right \|_{L^q(T)}^q
&= \int_{T} \left| \frac{\partial \varphi_h}{\partial x_i} \right|^q dx = |\det({A})| \int_{\widehat{T}} \left|  \widehat{\frac{\partial \varphi_h}{\partial x_i}} \right|^q d\hat{x} \\
&\leq c |\det({A})|  \sum_{j,k=1}^d  h_j^{-q} |\widetilde{A}^{-1}_{jk} (A_T^{-1})_{ki}|^q \left \| \frac{\partial \hat{\varphi}_h}{\partial \hat{x}_j} \right \|_{L^q(\widehat{T})}^q.
\end{align*}
Using the Jensen-type inequality \eqref{jensen} yields
\begin{align*}
\displaystyle
\left \| \frac{\partial \varphi_h}{\partial x_i} \right \|_{L^q(T)}
&\leq c  |\det({A})|^{\frac{1}{q}}  \sum_{j,k=1}^d  h_j^{-1} |\widetilde{A}^{-1}_{jk} (A_T^{-1})_{ki}| \left \| \frac{\partial \hat{\varphi}_h}{\partial \hat{x}_j} \right \|_{L^q(\widehat{T})}.
\end{align*}
All the norms in $\widehat{P}$ are equivalent, that is, there exists a positive constant $C^E$ depending on $\widehat{T}$ and $s \in \mathbb{N}_0$ such that
\begin{align}
\displaystyle
\| \hat{\varphi}_h \|_{W^{s,\infty}(\widehat{T})} \leq C^E  \| \hat{\varphi}_h \|_{L^1(\widehat{T})} \quad \forall \hat{\varphi}_h \in \widehat{P}. \label{inv1}
\end{align}
Using the standard scaling argument and \eqref{inv1}, we have
\begin{align*}
\displaystyle
 \left \| \frac{\partial  \hat{\varphi}_h}{\partial \hat{x}_j} \right \|_{L^q(\widehat{T})} &\leq \|  \hat{\varphi}_h \|_{W^{1,q}(\widehat{T})} \leq c  \| \hat{\varphi}_h \|_{L^p(\widehat{T})} \leq c |\det({A})|^{- \frac{1}{p}} \| {\varphi}_h \|_{L^p({T})}, 
\end{align*}
which leads to
\begin{align}
\displaystyle
\left \| \frac{\partial \varphi_h}{\partial x_i} \right \|_{L^q(T)}
&\leq c  |T|_d^{\frac{1}{q} - \frac{1}{p}} \left( \sum_{j,k=1}^d  h_j^{-1} |\widetilde{A}^{-1}_{jk} (A_T^{-1})_{ki}| \right) \| {\varphi}_h \|_{L^p({T})}, \label{inv2}
\end{align}
where we used the fact that $|\det (A)| = \frac{|T|_d}{|\widehat{T}|_d}$.

Let $d=2$. By the simple calculation,
\begin{align*}
\displaystyle
\widetilde{{A}}^{-1} =
\begin{pmatrix}
1 & - \frac{s}{t} \\
0 & \frac{1}{t} \\
\end{pmatrix}.
\end{align*}
Because \eqref{sym=mthscrH}, $|s| \leq 1$ and $h_2 \leq h_1$, we have
\begin{align}
\displaystyle
 \sum_{j,k=1}^2  h_j^{-1} |\widetilde{A}^{-1}_{jk} (A_T^{-1})_{ki}| 
 \leq \frac{1}{\widetilde{\mathscr{H}}_1}  |(A_T^{-1})_{1i}|  + \frac{2}{\widetilde{\mathscr{H}}_2} |(A_T^{-1})_{2i}|, \quad \text{$i=1,2.$}  \label{inv3}
 \end{align}
Gatherring \eqref{inv2} and \eqref{inv3}, we have the target inequalities for $ i= 1,2$.

Let $d=3$.  By the simple calculation,
\begin{align*}
\displaystyle
\widetilde{{A}}_1^{-1} =
\begin{pmatrix}
1 & - \frac{s_1}{t_1} & \frac{s_1 s_{22} - t_1 s_{21}}{t_1 t_2} \\
0 & \frac{1}{t_1} & - \frac{s_{22}}{t_1 t_2}\\
0 & 0  &  \frac{1}{t_2} \\
\end{pmatrix}, \
\widetilde{{A}}_2^{-1} =
\begin{pmatrix}
1 & \frac{s_1}{t_1} & \frac{- s_1 s_{22} - t_1 s_{21}}{t_1 t_2} \\
0 & \frac{1}{t_1}  & - \frac{s_{22}}{t_1 t_2} \\
0 & 0  & \frac{1}{t_2}\\
\end{pmatrix}.
\end{align*}
Recall that  $\widetilde{A} \in \{ \widetilde{A}_1 , \widetilde{A}_2 \}$. 
If Condition \ref{Cond333} is imposed, there exists a positive constant $M$ independent of $h_T$ such that $|s_{22}| \leq M \frac{h_2 t_1}{h_3}$. Because \eqref{sym=mthscrH}, $|s_1| \leq 1$, $|s_{21}| \leq 1$, $|s_{21}| \leq 1$ and $h_2 \leq h_3 \leq h_1$, we have
\begin{align*}
\displaystyle
&\frac{1}{h_1} = \frac{1}{\widetilde{\mathscr{H}}_1}, \quad \frac{|s_1|}{h_1 t_1} \leq \frac{1}{h_2 t_1} \leq \frac{1}{\widetilde{\mathscr{H}}_2}, \quad \frac{|s_{22}|}{h_2 t_1 t_2} \leq \frac{M}{h_3 t_2} = \frac{M}{\widetilde{\mathscr{H}}_3},\\
&\frac{|s_1 s_{22}| + |t_1 s_{21}|}{h_1 t_1 t_2} \leq M \frac{h_2}{h_1 h_3 t_2} + \frac{1}{h_1 t_2} \leq \frac{M+1}{h_3 t_2} = \frac{M+1}{\widetilde{\mathscr{H}}_3}.
\end{align*}
Using these inequalities,
\begin{align}
\displaystyle
 &\sum_{j,k=1}^3  h_j^{-1} |\widetilde{A}^{-1}_{jk} (A_T^{-1})_{ki}| \notag\\
&\quad \leq
\frac{1}{\widetilde{\mathscr{H}}_1}  |(A_T^{-1})_{1i}|  + \frac{2}{\widetilde{\mathscr{H}}_2} |(A_T^{-1})_{2i}| + \frac{2(M+1)}{\widetilde{\mathscr{H}}_3} |(A_T^{-1})_{3i}|,  \quad \text{$i=1,2,3$}.  \label{inv4}
\end{align}
Gatherring \eqref{inv2} and \eqref{inv4}, we have the target inequalities for $ i= 1,2,3$.
\qed
\end{pf*}

\begin{ex}
Let $d=2$. We set
\begin{align*}
\displaystyle
A_{r} \coloneq 
\begin{pmatrix}
-1  & 0 \\
 0 & 1
\end{pmatrix}, \quad 
A_{\theta} \coloneq 
\begin{pmatrix}
\cos \theta  & - \sin \theta \\
 \sin \theta & \cos \theta
\end{pmatrix}.
\end{align*}
Here, $A_r$ is a reflection matrix with respect to the $y$-axis and $A_{\theta}$
 is a rotation matrix. We define ${A}_{T} \in O(2)$ as $A_T := A_r A_{\theta}$, where  $\theta$ denotes the angle. Then,
\begin{align*}
\displaystyle
A_{T} =
\begin{pmatrix}
-1  & 0 \\
 0 & 1
\end{pmatrix} \quad  \text{if $\theta = 0$}, \quad 
A_{T} =
\begin{pmatrix}
0  & 1 \\
 1 & 0
\end{pmatrix} \quad  \text{if $\theta = \frac{\pi}{2}$}, \quad 
A_{T} =
\begin{pmatrix}
1  & 0 \\
 0 & -1
\end{pmatrix} \quad  \text{if $\theta = \pi$}.
\end{align*}
The inverse inequality \eqref{inv=2d} can be written as
\begin{align*}
\displaystyle
\left \| \frac{\partial \varphi_h}{\partial x_i} \right \|_{L^q(T)}
&\leq c |T|_d^{\frac{1}{q} - \frac{1}{p}} \frac{1}{\widetilde{\mathscr{H}}_i} \| {\varphi}_h \|_{L^p({T})}, \quad \text{$i=1,2$, if $\theta = 0,\pi$},\\
\left \| \frac{\partial \varphi_h}{\partial x_i} \right \|_{L^q(T)}
&\leq c |T|_d^{\frac{1}{q} - \frac{1}{p}} \frac{1}{\widetilde{\mathscr{H}}_{i+1}} \| {\varphi}_h \|_{L^p({T})}, \quad \text{$i=1,2$, if $\theta = \frac{\pi}{2}$},
\end{align*}
where the indices $i$, $i+1$ have to be understood mod 2. 
\end{ex}

%
%

\end{document}